\documentclass{article}
\usepackage{bm, amssymb, amsmath, amsthm, accents, algorithmic, algorithm, appendix, breqn, bigints, booktabs, caption,  color, enumitem, graphicx, multirow, mathtools, mathrsfs, subcaption, parskip, tabularx, tikz, verbatim, algorithm,algorithmic,float}
\usepackage[authoryear]{natbib}
\usepackage{bm,bbm}
\usepackage{enumitem}
\usepackage{calligra}
\usepackage{calrsfs}

\newtheorem{theorem}{Theorem}
\newtheorem{lemma}[theorem]{Lemma}
\newtheorem{proposition}[theorem]{Proposition}
\newtheorem{remark}[theorem]{Remark}
\newtheorem{corollary}[theorem]{Corollary}
\newtheorem{definition}[theorem]{Definition}

\newtheorem{assumption}{\unskip} 
\renewcommand{\theassumption}{(A\arabic{assumption})}

\newtheorem{assumptionU}{\unskip} 
\renewcommand{\theassumptionU}{(U\arabic{assumptionU})}

\newcommand{\beq}{\begin{equation}}
\newcommand{\eeq}{\end{equation}}

\def\eq#1/{(\ref{e:#1})}
\newcommand{\si}{\sigma}

\newcommand{\al}{\alpha}

\newcommand{\de}{\delta}

\newcommand{\la}{\lambda}

\newcommand{\ff}{\infty}
\newcommand{\ra}{\rightarrow}

\newcommand{\ga}{\gamma}

\newcommand \QED{\hbox{\hskip6pt\vrule height 8pt width 6pt}}
\newcommand{\ep}{\epsilon}

\DeclarePairedDelimiterX{\inp}[2]{\langle}{\rangle}{#1, #2}

\newcommand{\R}{\mathbb{R}}

\newcommand{\beqq}{\begin{eqnarray}}
\newcommand{\eeqq}{\end{eqnarray}}
\newcommand{\ta}{\theta}
\newcommand{\bta}{\bm{\ta}}
\newcommand{\Ta}{\Theta}

\newcommand{\interior}[1]{{\kern0pt#1}^{\mathrm{o}}}

\newcommand{\argmin}{\underset{\bm{\theta} \in {\Theta}}{\mbox{argmin}}}

\newcommand{\triple}{(\Omega, \mathcal{F}, P)}

\def\Real{\mathbb R}

\def \Rm{\Real^m}

\usepackage{xcolor}
\usepackage{xurl}
\usepackage{etoolbox}
\usepackage[hidelinks]{hyperref}
\makeatletter
\apptocmd{\@floatboxreset}{\centering}{}{}
\makeatother
\newenvironment{keywords}{\par\noindent\textbf{Key words}: }{\par}

\renewcommand{\theassumptionU}{{\bf(U\arabic{assumptionU})}}
\def\theHassumptionU{\arabic{assumptionU}}
\def\theHassumption{\arabic{assumption}}
\newif\ifAVmarks
\AVmarksfalse
\newcommand{\AV}[1]{{\ifAVmarks\color{red}\fi#1}}

\begin{document}
\title{Private Minimum Hellinger Distance Estimation via Hellinger Distance Differential Privacy}
\author{Fengnan Deng \\ Department of Statistics \\ George Mason University \\ Fairfax, VA 22030 \and Anand N.\ Vidyashankar \\Department of Statistics \\George Mason University \\Fairfax, VA 22030 }
\date{ }
\maketitle
\begin{abstract}
Objective functions based on Hellinger distance yield robust and efficient estimators of model parameters. Motivated by privacy and regulatory requirements encountered in contemporary applications, we derive in this paper \emph{private minimum Hellinger distance estimators}. The estimators satisfy a new privacy constraint, namely, Hellinger differential privacy. Specifically, the proposed estimators are consistent and efficient while satisfying the stated Hellinger differential privacy guarantees. We demonstrate that Hellinger differential privacy shares several features of standard differential privacy while allowing reduced Gaussian perturbations relative to basic approximate-DP composition at a common final privacy target. Additionally, for computational purposes, we also develop Hellinger differentially private gradient descent and Newton-Raphson algorithms. We illustrate the behavior of our estimators in finite samples using numerical experiments and verify that they retain robustness properties under gross-error contamination.
\end{abstract}


\section{Introduction}

Recently, the adoption of AI (artificial intelligence)- whose success relies on data- by many scientific fields has brought an increased focus on data privacy. Many entities collect individually identifiable data to provide personalized services and share them with other organizations to improve and enhance the quality of the service. However, such data-sharing activities lead to increased data privacy concerns. Specifically, in healthcare, data are often shared with various research groups to facilitate treatment payment operations and improve the quality of care. The Health Insurance Portability and Accountability Act (HIPAA), the Health Information Technology for Economic and Clinical Health (HITECH), the California Consumer Privacy Act (CCPA), and the General Data Protection Regulation (GDPR) are some of the recent regulations that impact commerce.  Historically, anonymization techniques, like encrypting or removing personally identifiable information, have been widely used to ensure privacy protection. However, recent studies 
\citep{gymrek2013, homer2008,narayan2008,sweeney1997}
have shown that many of the existing anonymization methods are fragile and can lead to the leakage of private information. Specifically, an intruder might still be
able to identify individuals by cross-classifying categorical
variables in the dataset and matching them with some external
database. A need for a rigorous technical framework to measure and analyze the de-identification methods has long been noted (see \cite{DL86}). 

Differential privacy (DP) is a probabilistic framework that quantifies how individual privacy is preserved in a database when certain information is released by querying the database. The basic idea behind the DP framework is to measure the indistinguishability, using a parameter $\ep$, of two probability distributions of a dataset in the presence or absence of a record. Small values of $\ep$ correspond to hard distinguishability and high privacy. The distributions in question typically correspond to those of statistics obtained by querying the database.

In an interactive setting of the DP, a data warehouse provides a group of query functions that allow users to pose queries about the data and receive responses with noise added for privacy (referred to as a mechanism). In a non-interactive setting, the data warehouses offer a dataset with added noise, and users can apply any models and methods to this data. Controlling privacy breaches is challenging in the non-interactive setting (see \cite{Dwork2006b}, for instance). In contrast, privacy-preserving mechanisms that rely on specific query functions without direct access to the dataset enable assessment of privacy and utility. In this paper, we focus on the interactive approach and assume the existence of a trusted curator holding individuals' data in a database. The goal
of private inference is to protect individual data while simultaneously allowing statistical analysis of 
the entire database. An analyst can only access a model's perturbed summary statistics or outputs in such cases. While adding noise preserves privacy, the amount of noise needs to be small to ensure optimal statistical performance. 

This paper describes a new notion, namely Hellinger distance differential privacy (HDP)- a particular case of power divergence privacy (PDP)- and studies private estimation and inference for Minimum Hellinger Distance Estimators (MHDEs). The power divergence family parameterized by $\la$ encompasses R\'enyi divergence up to a logarithmic transformation. The classical $\ep-$DP can be obtained by taking the limit as $\la$ diverges to infinity in the PDP. Following the privacy literature, we describe privacy using the parameters $\la$ and $\ep$ and the terminology $(\la, \ep)-$PDP. When $\la=-\frac{1}{2}$, the power divergence reduces to twice the squared Hellinger distance between densities, and hence $(-\frac{1}{2},2\ep)-$PDP is referred to as $\ep-$HDP.

Some well-known differential privacy methods, such as $\epsilon-$DP, $(\epsilon,\delta)-$DP, and $(\alpha,\epsilon)-$R\'enyi differential privacy (RDP) (\cite{Mironov2017}) can be defined using statistical divergence measures and are subsumed in our $(\la, \ep)-$PDP. Specifically, for $\la>0$, $(\la, \ep)-$PDP is equivalent to $(\la+1,\frac{1}{\la}\log(\ep \la(\la+1)+1)-$RDP.
While the RDP focuses on the case $\la >0$, the case $ \la <0$ has several useful properties. Specifically, it turns out that the additive Gaussian mechanism has minimal variance when $\la=-\frac{1}{2}$. This observation motivates a detailed study of this case, namely, the HDP. We establish that HDP  has better composition and group privacy properties than PDP (see Theorem \ref{thm:HDP} and Theorem \ref{thm:HDP_group}). We also discuss the relationship between HDP and other privacy frameworks, which broadens the potential applications of HDP.

Another contribution of our paper is deriving private MHDE (PMHDE) estimators; these estimators are private, robust to model misspecification, and efficient. In contrast, M-estimators are not always efficient due to the boundedness of the score function. Additionally, assumptions such as convexity and Lipschitz property are typically used to derive the estimators and study the properties (\cite{avella2021privacy,chaudhuri2011differentially,chaudhuri2012convergence,chen2019renyi,slavkovic2012perturbed,wang2017selecting}).
This paper develops asymptotic properties of PMHDE obtained through private optimization algorithms without explicit assumptions of convexity and boundedness. While $\ep-$HDP is used for privacy guarantees, the methods also work for any $(\la,\ep)-$PDP for all real values of $\la$. Additionally, other approaches such as $\mu-$GDP (\cite{Dong2022}) and $\rho-$zCDP (\cite{Bun2016,Dwork2016}) can also be used for privacy guarantees. Some of these are described in Section \ref{sec:HDP}.


The sensitivity of the query function
plays a central role in DP investigations. In applications in point estimation, the query concerns the gradient of the loss function under appropriate model regularity. Since the second-order properties of estimators also rely on the first and second-order derivatives of the loss functions, it is reasonable to anticipate a link between the sensitivity and statistical efficiency. We make this precise in Sections \ref{sec:utility_PMHDE} and \ref{sec:efficiency_PMHDE}. Specifically, we obtain sharp estimates of the sensitivities of the gradient and Hessian of the Hellinger loss functions (see Theorem \ref{thm:sensitivity}) and use them to derive the limit distribution of the PMHDE. The arguments required for establishing this are subtle and involved. Finally, it is worth emphasizing here that we do not make convexity assumptions. Instead, we leverage the family regularity, the properties of the Hellinger objective function, and the $L_1$ properties of kernel density estimates to develop our results. 

Algorithms such as gradient descent (GD) and Newton-Raphson (NR) are typically used to obtain MHDE and other M-estimators (\cite{bassily2014, feldman2020,lee2018,loh2013regularized, song2013}). The private versions of these estimators are derived by optimizing the private (perturbed) objective functions obtained by adding an appropriate noise at every iteration. We study private versions of GD and NR, namely PGD and PNR, that output $\ep-$HDP counterparts of MHDE, namely PMHDE (see Section \ref{sec:private_opti}). While it is common to study the private versions for a fixed $\ep >0$ as $n \to \ff$, a unique aspect of our work is that we allow $\ep$ to decrease to zero as $n \to \ff$. By a judicious choice of this rate, we retain full privacy and efficiency. Alternatively, one could use similar perturbations to obtain $(\la, \ep)-$PDP counterparts of MHDE. Analysis of these algorithms is critical to study the properties of the PMHDE. These are investigated in Section \ref{sec:utility_PMHDE} and Section \ref{sec:efficiency_PMHDE}.
We begin in Section \ref{sec:HDP} with the background and notations of DP, while in Section \ref{sec:numerical}, we present several numerical experiments that evaluate the performance of our estimators under the true model and contamination. Section \ref{sec:ext and cr} contains some extensions and concluding remarks. The online supplementary material contains proofs in Section 6 and regularity conditions, additional calculations, and numerical results in Sections 7--11. References to Sections 6--11 below refer to this supplement.

\section{Background, notations, and Hellinger distance differential privacy} \label{sec:HDP}
Let $\{X_n: n \ge 1\}$ denote a collection of independent and identically distributed (i.i.d.) real-valued random variables defined on the probability space $(\Omega, \mathcal{F}, P)$ and set $\mathbf{X}=(X_1, X_2, \cdots X_n)$ so that 
$\mathbf{X}: \triple \mapsto (\mathbb{R}^n,\mathcal{B}(\mathbb{R}^n),P_X)$, where $P_X$ is the induced probability measure on the Borel subsets of $\R^n$. We denote by $\mathcal{D} \subset \R^n$, the database of the i.i.d. observations; that is, $\mathcal{D} = \{\mathbf{X}(\omega): \omega \in \Omega\}$. A query function $W(\cdot)$ is a statistic; namely a measurable mapping $W:(\mathcal{D},\mathcal{B}(\R^n))\mapsto (\Rm,\mathcal{B}(\mathbb{R}^m))$. We denote a typical element of $\mathcal{D}$, namely the dataset, by $D$, and the query applied to $D$ by $w \coloneqq W(D)$. In our applications, we will be interested in functions of the type, $f(w, D)$,  where $f(\cdot, \cdot)$ is a measurable mapping from $ (\R^m \times \R^n, \mathcal{B}(\R^m)\times \mathcal{B}(\R^n) ) \mapsto  (\Rm, \mathcal{B}(\Rm))$. A simple example of $f(w,D)$ is $f(w,D)=w$.

Next, we introduce one of the essential concepts of DP, namely the mechanism (or a randomized algorithm) denoted by $M$. In statistical terms, $M$ is a measurable mapping from $ (\R^m \times \R^n, \mathcal{B}(\R^m)\times \mathcal{B}(\R^n) ) \mapsto  (\Rm, \mathcal{B}(\Rm))$.
$M$ is said to be an additive mechanism, if $M(w, D) = f(w, D)+ \mathbf{Y}$, where $\mathbf{Y}=[Y_1,\cdots,Y_m] \in 
\Rm$, is a random vector (with i.i.d. components) representing the noise and independent of $(w,D)$. 
In here, $M(\cdot, \cdot)$ represents a private version of $f(\cdot, \cdot)$. Continuing with the above example, the additive mechanism will output $w+\mathbf{Y}$, a perturbed version of $w$. 

A critical component of the privacy measure is the sensitivity of the query function. Formally, in DP, it is based on two \emph{adjacent} datasets differing in one observation. Specifically, for $D,D'\in\mathcal{D}$, set $z_i=\bm{1}_{[x_i \ne x_i']}$, for $1 \le i \le n$. Then, the Hamming distance between $D$ and $D'$ is given by 
\begin{align*}
    ||D-D'||_H = \sum_{i=1}^n z_i.
\end{align*}
We say $D$ and $D'$ are adjacent if $||D-D'||_H = 1$. Define ${L}_1$ and ${L}_2$ sensitivity of a query function $W$ to be 
\begin{align*}
    \Delta_{L_1}W=\sup\limits_{||D-D'||_H = 1}||W(D)-W(D')||_1\\
    \Delta_{L_2}W=\sup\limits_{||D-D'||_H =1}||W(D)-W(D')||_2.
\end{align*}

\noindent We now turn to a precise description of some widely used DP measures. To this end, let $Q$ denote the conditional distribution of $M(w,D)$ given $w, D$. We start with $\epsilon-$differential privacy introduced in
\cite{Dwork2006b}.
\begin{definition}
A mechanism, ${M}$, satisfies $\epsilon-$differential privacy (DP) if for any $S \in \mathcal{B}(\Rm)$ and adjacent $D,D'$,
\begin{align}{\label{ep-DP}}
    Q(M(w,D)\in S) \leq e^{\epsilon}\cdot Q(M(w,D')\in S).
\end{align}
\end{definition}
\noindent When $M$ is an additive mechanism, that is,  $M(w, D)=W(D)+\mathbf{Y}$ and the random variables $Y_1, Y_2, \cdots, Y_m$ are i.i.d. with $Y_1 \sim Lap(0,\frac{\Delta_{L_1}W}{\epsilon})$, then $M$ satisfies $\ep-$DP. A natural case of an additive mechanism, namely the Gaussian mechanism, does not satisfy (\ref{ep-DP}). Hence, a relaxation of $\ep-$(DP) was studied in \cite{Dwork2006b} and is referred to as approximate DP or, more commonly, as $(\ep, \de)-$DP. 
\begin{definition}
    A mechanism $M$, satisfies $(\epsilon,\delta)-$differential privacy (DP) if for any possible output $S\subset \mathcal{B}(\Rm)$ and adjacent $D,D'$,
\begin{align*}
    Q(M(w,D)\in S)  \leq e^{\epsilon}\cdot Q(M(w,D')\in S)  + \delta.
\end{align*}
\end{definition}
In this case, the distribution of $\mathbf{Y}$ is $N_m(\mathbf{0}, \sigma^2 \mathbf{I}_{m\times m})$,
where, for  $0 < \epsilon, \delta <1$,
\begin{align*}
    \sigma^2=\frac{2\cdot (\Delta_{L_2}W)^2\cdot\log(1.25/\delta)}{\epsilon^2}.
\end{align*}

Other notions of relaxations of differential privacy have been studied in the literature. Some of the most commonly studied ones include concentrated differential privacy (CDP) (\cite{Dwork2016}), zero concentrated differential privacy (zCDP) (\cite{Bun2016}), Gaussian differential privacy (\cite{Dong2022}), and Renyi differential privacy (RDP) (\cite{Mironov2017}). Some of these approaches can be unified under the general notion of divergence. We first turn to the notion of Hellinger distance differential privacy (HDP) and a related generalization. In the following, for two random variables $X$ and $Y$ with distribution $P_1$ and $P_2$ respectively, we denote the divergence between them as $D(X,Y)$, which is equivalent to $D(P_1,P_2)$.
\begin{definition}{(Hellinger differential privacy)}
A mechanism $M$ is said to satisfy $\epsilon-$Hellinger differential privacy (HDP), for $\ep \in [0, 2]$, if for any adjacent $D,D' \in \mathcal{D}$,  
\begin{align*}
D_{HD}(M(w,D),M(w,D'))\leq \epsilon.
\end{align*}
where for two distributions $P_1, P_2$, with densities $p_1(\cdot)$ and $p_2(\cdot)$ with respect to (w.r.t.) the Lebesgue measure,
\begin{align*}
D_{HD}(P_1,P_2) = \int_{\mathbb{R}^m} \left(\sqrt{p_1(x)}-\sqrt{p_2(x)}\right)^2 \mathrm{d}x.
\end{align*}
\end{definition}

Hellinger distance is a member of the general class of divergences referred to as \emph{Power divergence}, introduced by \cite{Cressie1984}
and further analyzed in \cite{Read1988} for performing goodness of fit tests in multinomial and multivariate discrete data. The ideas were unified in the work of Lindsay (\cite{Lindsay1994}), who studied general divergences for robust and efficient estimation in parametric models (see also \cite{Basu2011}). Let $P_1$ and $P_2$ be two probability distributions possessing densities $p_1(\cdot)$ and $p_2(\cdot)$ on $\Rm$. The power divergence $D_{\lambda}(P_1,P_2)$ between $P_1$ and $P_2$, denoted by $D_{\la}(P_1,P_2)$ is defined as follows: for $\la \neq -1,0$
\begin{align*}
D_{\lambda}(P_1,P_2) = 
       \frac{1}{\lambda(\lambda+1)}\mathbf{E}_{X\sim p_2}\left[\left(\frac{p_1(X)}{p_2(X)}\right)^{\lambda+1}-1\right].  
\end{align*}
$D_0 (P_1,P_2)$ and $D_{-1}(P_1,P_2)$ are defined by taking the limits as $\la$ approaches $0$ or $-1$. A standard calculation shows that
\begin{align*}
   D_0(P_1,P_2)=D_{-1}(P_2,P_1)=KL(P_1,P_2), 
\end{align*}
where $KL(\cdot, \cdot)$ represents the Kullback-Leibler divergence. R\'enyi divergence is a particular case of power divergence when $\la >0$; specifically, setting $\alpha=(\la+1)$, the R\'enyi divergence of order $\alpha$ is given by
\begin{align*}
    D_{\alpha}^{(R)}(P_1,P_2)=\frac{1}{\lambda}\log\left[\lambda(\lambda+1)D_{\lambda}(P_1,P_2)+1\right].
\end{align*}
When $P_1$ and $P_2$ have the same support, the limit of $D_{\la}(P_1, P_2)$ exists and is referred to as the Max divergence and is given by
\begin{align*}
    D_{\infty}(P_1,P_2)  \coloneqq \lim\limits_{\lambda\to\infty}\frac{1}{\lambda}\log[\lambda(\lambda+1)D_{\lambda}(P_1,P_2)+1] = \max_{S\in Supp(P_1)}\log \frac{P_1(S)}{P_2(S)}.
\end{align*}
For $\delta\in(0,1)$, a $\delta-$relaxation  of the above max divergence is given by
\begin{align*}
    D_{\infty}^\delta(P_1,P_2) \coloneqq  \max_{S\in Supp(P_1):P_1(S)\geq\delta}\log \frac{P_1(S)-\delta}{P_2(S)}.
\end{align*}
Using the above notions, one can express the $\ep-$DP, $(\ep, \de)-$DP, and $(\ep, \al)-$RDP as follows: let $D,D'$ be adjacent. A mechanism, $M$ is $\epsilon-$DP if
$D_{\infty}(M(W(D)),M(W(D')))  \leq \epsilon$,  while it is $(\epsilon,\delta)-$DP if $D_{\infty}^{\delta}(M(W(D)),M(W(D'))) \leq \epsilon$.
It is said to be $(\alpha,\epsilon)-$RDP ($\alpha>1$)  if $D_{\alpha}^{(R)}(M(W(D)),M(W(D')))\leq \epsilon$. 
We now turn to describe a new privacy measure called \emph{Power Divergence Differential Privacy} (PDP).
\begin{definition}
Let $\la \in \Real$ and $ \ep > 0$ if $[\la(\la+1)] \ge 0 $ and $ 0 < \ep < -[\la(\la+1)]^{-1}$ otherwise. A mechanism $M$ is said to satisfy $(\lambda,\epsilon)-$Power differentially private (PDP) if for any fixed adjacent $D,D'$,  
\begin{align*}
D_{\lambda}(M(w,D)),M(w,D')))\leq \epsilon.
\end{align*}
\end{definition}

\begin{remark} \label{rem:PDP_relation}
It follows from the above definitions that $\epsilon-$HDP is equivalent to $(-\frac{1}{2},2\epsilon)-$PDP.
We emphasize here PDP is defined for any $\lambda\in\mathbb{R}$ and includes the RDP. Specifically, for $\la>0$, a standard calculation shows that $(\la,\ep)-$PDP is equivalent to $(\la+1,\frac{1}{\la}\log(\ep\la(\la+1)+1))-$RDP. Furthermore, if $M$ satisfies $(\la,\ep)-$PDP then it satisfies $(\alpha,\alpha\epsilon)-$RDP, where $\alpha=\la+1>1$. If $M$ is an additive Gaussian mechanism with common dataset-independent covariance $\sigma^2 I_m$, the same PDP guarantee for $\lambda>0$ also implies $\epsilon-$zCDP (details are in Section 9).
For any $\delta\in(0,1)$ and $\lambda>0$, $(\lambda,\epsilon)-$PDP implies
$(\frac{1}{\lambda}\log(\frac{\lambda(\lambda+1)\epsilon+1}{\delta}),\delta)-$DP.
The values $\lambda<-1$ fall outside the range used in the definition of RDP. Nevertheless, by reversing the ordered pair in the divergence identity, one obtains $(\frac{-1}{\lambda+1}\log(\frac{\lambda(\lambda+1)\epsilon+1}{\delta}),\delta)-$DP; see Section 9.
For additive Gaussian and scalar Laplace mechanisms with common dataset-independent noise scales, $M$ also satisfies $\mu-$GDP. Writing $\Delta_{L_i}W$ for the corresponding sensitivity, the smallest GDP parameters are, respectively,
\[
\mu_G=\frac{\Delta_{L_2}W}{\sigma},\qquad
\mu_L=2\Phi^{-1}\left(1-\frac12e^{-\Delta_{L_1}W/(2b)}\right),
\]
where $b$ is the scalar Laplace scale. For Gaussian noise and $\lambda(\lambda+1)\neq0$, the PDP bound gives
$\mu_G^2\leq\frac{2}{\lambda(\lambda+1)}\log\{1+\lambda(\lambda+1)\epsilon\}$, with the continuous bound $\mu_G^2\leq2\epsilon$ at $\lambda=0,-1$. The mechanism-specific GDP calculations are given in the supplementary material.
\end{remark}

\subsection{Alternative partition-based definitions}
\label{sec:pd-partitions}

The preceding definitions can also be expressed through finite-valued
versions of the released output. Let $\mathfrak P$ denote the collection
of finite Borel measurable partitions of $\mathbb R^m$. For
$\pi=\{A_1,\ldots,A_k\}\in\mathfrak P$, choose distinct representatives
$x_j\in A_j$ and put
\[
X_\pi=\sum_{j=1}^k x_j\mathbf1_{\{X\in A_j\}}.
\]
Let $P_i^\pi$ be the distribution of $X_\pi$ when $X$ has distribution
$P_i$, so that $p_i^\pi(x_j)=P_i(A_j)$. For positive cell probabilities,
the preceding power-divergence formula becomes, for $\lambda\ne-1,0$,
\begin{equation}
D_\lambda(P_1^\pi,P_2^\pi)
=\frac{\sum_{j=1}^k P_1(A_j)^{\lambda+1}P_2(A_j)^{-\lambda}-1}
{\lambda(\lambda+1)}.
\label{eq:pd-cell-formula}
\end{equation}
The above identity also applies when some cell probabilities are zero. When exactly one of the two cell probabilities is zero, the corresponding term in the sum is defined by letting that probability decrease to zero while holding the positive probability fixed. When both cell probabilities are zero, the term is assigned the value zero. These conventions allow the divergence to take the value positive infinity.
The same convention applies to the density-form divergence. The
logarithmic case is
\[
D_0(P_1^\pi,P_2^\pi)
=\sum_jP_1(A_j)\log\{P_1(A_j)/P_2(A_j)\},
\]
and $D_{-1}(P_1^\pi,P_2^\pi)=D_0(P_2^\pi,P_1^\pi)$.
We use $0\log(0/b)=0$, including $b=0$, and
$a\log(a/0)=+\infty$ for $a>0$.
Thus the value depends only on the cell probabilities, not on the
representatives.

Partition representations of information and $\chi^2$ divergences are
used in \citet{Bei01}. The following result is the power-divergence
case of the general finite-partition identity; see also
\citet{Gilardoni2009}.

\begin{proposition}[Partition equivalence]
\label{prop:pd-partition-equivalence}
For every $\lambda\in\mathbb R$ and probability measures $P_1,P_2$,
\begin{equation}
D_\lambda(P_1,P_2)
=\sup_{\pi\in\mathfrak P}D_\lambda(P_1^\pi,P_2^\pi).
\label{eq:pd-partition-equivalence}
\end{equation}
The equality holds in $[0,\infty]$. The partition divergence is
nondecreasing under refinement, and there is a refining sequence
$(\pi_r)$ for which
$D_\lambda(P_1^{\pi_r},P_2^{\pi_r})\uparrow D_\lambda(P_1,P_2)$.
Consequently, for the admissible parameters in the preceding definition,
$M$ satisfies $(\lambda,h)$-PDP if and only if
$D_\lambda(P_D^\pi,P_{D'}^\pi)\leq h$ for every finite measurable
partition and every ordered pair of adjacent datasets, where $P_D$
denotes the output law of $M$ at $D$.
\end{proposition}

\noindent\emph{Proof.} See Supplementary Section 6.1.

Although the divergences were introduced using Lebesgue densities,
the same statements hold for arbitrary probability measures when the
densities are replaced by Radon--Nikodym derivatives relative to a
common dominating measure, with the usual boundary conventions.

\begin{remark}[Partition formulation of HDP]
\label{rem:hdp-partitions}
Since $D_{-1/2}=2D_{HD}$, Proposition~\ref{prop:pd-partition-equivalence}
gives
\begin{equation}
D_{HD}(P_1,P_2)
=\sup_{\pi\in\mathfrak P}
\sum_{A\in\pi}
\left\{\sqrt{P_1(A)}-\sqrt{P_2(A)}\right\}^2.
\label{eq:hdp-partition-equivalence}
\end{equation}
Thus $h$-HDP is equivalent to requiring the sum in
\eqref{eq:hdp-partition-equivalence} to be at most $h$ for every finite
partition and every neighboring pair. In the nontrivial parameter
range $0<h<2$, this is also $(-1/2,2h)$-PDP.
Each partition represents a finite-valued observation of the release,
and the constraint applies jointly to all its cell probabilities.
Any finite collection of measurable output events generates such a
partition.
\end{remark}

\subsection{Properties and mechanisms of PDP}

In applications, one encounters multiple queries applied to datasets.
These are referred to as compositions.  Three commonly occurring compositions are: (i) Sequential composition, (ii) Adaptive composition, and (iii) Parallel composition. Parallel composition involves disjoint datasets, while sequential and adaptive compositions typically involve the same dataset.

Let $\{M_k: k \ge 1\}$ denote a sequence of mechanisms. The adaptive $n-$composition of mechanisms $M_1, \cdots, M_n$, denoted by $M^{(n)}$, represents the trajectory of the outputs from $n$ mechanisms and is defined recursively as follows: let $M^{(0)}(w,D)=w$, and
let $M^{(1)}(w, D)=[M_1(M^{(0)}(w),D)]$; for $n \ge 2$
\begin{align*}
    M^{(n)}(w,D)=[M^{(n-1)}(w,D), M_n(\langle\mathbf{e}_n, M^{(n-1)}(w)\rangle,D)],
\end{align*}
where $\mathbf{e}_n=(0,0, \cdots,1)_{1\times (n-1)}$ is the unit vector. We note here that the trajectory is useful for describing the composition property and for calculating the mechanisms. \emph{However, only the $n^{\text{th}}$ component of $M^{(n)}$ is released.} Next, turning to the sequential composition, it is given by
\begin{align*}
    M^{(n)}(w,D)=[M^{(n-1)}(w,D), M_n(w,D)]=[M_1(w,D),\cdots,M_n(w,D)].
\end{align*}
The parallel composition is defined by applying a sequence of queries on disjoint datasets, specifically,
\begin{align*}
    M^{(n)}(\mathbf{W}^{(n)},\mathbf{D}^{(n)})=[M_1(w_1,D_1),\cdots,M_n(w_n,D_n)],
\end{align*}
where $D_i\cap D_j=\emptyset$ for all $i\neq j$, and $\mathbf{D}^{(n)}=D_1\times\cdots\times D_n \in \mathcal{D}^{(n)} = \mathcal{D}_1\times\cdots\times\mathcal{D}_n$ and $\mathbf{W}^{(n)}=[w_1,\cdots,w_n]: \mathcal{D}^{(n)} \mapsto \Rm_1\times\cdots\times\Rm_n$. To extend the definition of adjacent datasets to $\mathbf{D}^{(n)}$, we extend the Hamming distance between  $\mathbf{D}^{(n)}$ and $\mathbf{D}^{(n)'}=D_1'\times\cdots\times D_n'$ as follows:
\begin{align*}
    ||\mathbf{D}^{(n)}-\mathbf{D}^{(n)'}||_H = \sum_{i=1}^n ||D_i-D_i'||_H.
\end{align*}
We say $\mathbf{D}^{(n)}$ and $\mathbf{D}^{(n)'}$ are adjacent if
\begin{align*}
    ||\mathbf{D}^{(n)}-\mathbf{D}^{(n)'}||_H=1.
\end{align*}
Below, by $M_i(W, D)|W=w$ we mean the mechanism $M_i$ acting on a given query $W=w$ and dataset $D$.
Our next result is concerned with the PDP properties of the compositions.
\begin{theorem} \label{thm:PDP}
~
\begin{enumerate}
\item (Adaptive composition) Let $M_1(w,D)$ satisfy $(\lambda,\epsilon_1)-$PDP and  $M_2(W,D)|W$ satisfy $(\lambda,\epsilon_2)-$PDP. Then the composition $M^{(2)}(w,D)=(M_1(w,D),M_2(M_1(w,D),D))$ satisfies $(\lambda,(\epsilon_1+\epsilon_2 + \lambda(\lambda+1)\epsilon_1\epsilon_2))-$PDP.
    \item (Sequential composition) Let $M_1(w,D)$ satisfy $(\lambda,\epsilon_1)-$PDP and $M_2(w,D)$ satisfy $(\lambda,\epsilon_2)-$PDP. Then the composition $M^{(2)}(w,D)=(M_1(w,D),M_2(w,D))$ satisfies $(\lambda,(\epsilon_1+\epsilon_2 + \lambda(\lambda+1)\epsilon_1\epsilon_2))-$PDP.
    \item (Parallel composition) 
    Let $M_1$ and $M_2$ satisfy $\ep_1$ and $\ep_2$ PDP  on two disjoint datasets $D_1$ and $D_2$ with distinct queries $w_1$ and $w_2$ respectively. Then, the parallel composition $M^{(2)}(\mathbf{W}^{(2)}, \mathbf{D}^{(2)})$ satisfies $(\lambda,\max\{\epsilon_1,\epsilon_2\})-$PDP.
\end{enumerate}
\end{theorem}
We now focus on the additive mechanism with Gaussian perturbation. Recalling that such a mechanism can be represented as
$M(w,D)=w+\mathbf{Y}$, $\mathbf{Y}\sim N(0,\sigma^2\cdot\mathbf{I})$ (or $m-$dimensional Laplace with independent components), our aim is to identify $\sigma^2$ (or $b$) to achieve $(\lambda,\epsilon)-$PDP. Our
next Theorem summarizes this result.
\begin{theorem} \label{thm:PDP_mech}
Let $\mathbf{Y} =(Y_1, Y_2, \cdots, Y_m)$. If $Y_i$'s are i.i.d. $N(0, \si^2 )$, then the choice 
\begin{align} \label{eq:additive_noise}
    \si^2 \coloneqq \si_{\lambda,\ep}^2=\begin{cases}
        \left(\Delta_{L_2} W\right)^2\cdot \frac{\lambda(\lambda+1)}{2\log (1+\lambda(\lambda+1)\epsilon)} & \text{if } \la(\la+1)  \neq 0 \\
        \left(\Delta_{L_2} W\right)^2 \cdot \frac{1}{2\epsilon} &\text{otherwise},
       \end{cases}
\end{align}
renders the mechanism $(\la, \ep)-$PDP. If $Y_i$'s are i.i.d. $Lap(0,b)$, then the choice
\begin{eqnarray}
    b \coloneqq b_{\la, \ep} = \begin{cases}
    \max\left\{ \frac{sign(\lambda)(\lambda+1)\Delta_{L_1} W}{\log(\lambda(\lambda+1)\epsilon+1)},  \frac{sign(\lambda+1)(\lambda)\Delta_{L_1} W}{\log(\lambda(\lambda+1)\epsilon+1)}\right\}, &\text{if }\la(\la+1)  \neq 0\\
    \frac{\Delta_{L_1}W}{\ep} &\text{otherwise},
    \end{cases}
\end{eqnarray}
renders the mechanism $(\la,\ep)-$PDP. Furthermore, choosing $\la=-\frac{1}{2}$ and replacing $\ep$ by $2\ep$, we obtain for Gaussian $\mathbf{Y}$ and Laplace $\mathbf{Y}$
\begin{align*}
    \si^2_{HDP, \ep}=\frac{\left(\Delta_{L_2} W\right)^2}{8\log(\frac{1}{1-0.5\epsilon})},\quad \text{and}~~ b_{HDP, \ep}=\frac{\Delta_{L_1}W}{2\log(\frac{1}{1-0.5\epsilon})} \quad \text{respectively}.
\end{align*}
\end{theorem}
The case $\lambda=-\frac{1}{2}$ is interesting for the following reasons:
\begin{enumerate}
\item When $\lambda=-\frac{1}{2}$, the Power divergence is twice the squared Hellinger distance, which is widely used in point estimation. Specifically, the minimum Hellinger distance estimator achieves robustness and efficiency when the model is correctly specified (\cite{Beran1977,Cheng2006}).
\item When using the additive mechanism with Gaussian perturbation, for a fixed numerical PDP budget $0<\epsilon<4$, $\la=-\frac{1}{2}$ minimizes $\si^2_{\la}$ for $\la \in \Real$. To see this, setting $t= \la(\la+1)\in(-\frac{1}{\ep},\infty)$, and $f(t)=\frac{t}{2\log(1+t\epsilon)}$, observe that $$f'(t)=\frac{\log(1+t\epsilon)-\frac{t\epsilon}{t\epsilon+1}}{2(\log(1+t\epsilon))^2}, ~~f'(0)=\frac14.$$ Next, setting $g(t)=\log(1+t\epsilon)-\frac{t\epsilon}{t\epsilon+1}$,  we observe that 
$g'(t)>0$ for  $t>0$, $g'(t)<0$ as $t<0$, $g'(0)=0$. Since $g(0)=0$, these derivative signs imply $g(t)\geq0$.
Therefore, $f'(t)\geq0$ for $t\neq0$. At zero, the continuous
extension $f(0)=1/(2\epsilon)$ gives $f'(0)=1/4$.
Hence $f(t)$ is non-decreasing in $t$. Since $t$ is quadratic in $\lambda$ and is minimized at $\lambda=-\frac{1}{2}$, $\sigma^2$ is minimized at $\lambda=-\frac{1}{2}$. 
\item For both adaptive and sequential composition, the privacy is maximized in the power divergence class at $\la=-\frac{1}{2}$.
\item When $\lambda=-\frac{1}{2}$, PDP is a symmetric divergence and has a simpler group privacy representation (see Theorem \ref{thm:HDP_group} below).
\end{enumerate}

We now turn to a more detailed analysis of the case $\la=-\frac{1}{2}$. As explained previously,  a sequence of analyses is performed on the same dataset, with each analysis using the information from the previous ones. If each analysis satisfies a certain privacy level, then the overall privacy guarantee for this sequence is given by the adaptive composition rule. The next result is a particular case of Theorem \ref{thm:PDP}. We state it separately to emphasize the choice of $\lambda$ and $\ep$.
\begin{theorem}{\label{thm:HDP}}
~
\begin{enumerate}  
\item (Adaptive composition) Let $M_1(w,D)$ satisfy $\epsilon_1-$HDP and $M_2(W,D)|W$ satisfy $\epsilon_2-$HDP. Then the composition $M^{(2)}(w,D)=(M_1(w,D),M_2(M_1(w,D),D))$ satisfies $(\epsilon_1+\epsilon_2-\frac{1}{2}\epsilon_1\epsilon_2)-$HDP.
\item (Sequential composition) Let $M_1(w,D)$ satisfy $\epsilon_1-$HDP and $M_2(w,D)$ satisfy $\epsilon_2-$HDP. Then the composition $M^{(2)}(w,D)=(M_1(w,D),M_2(w,D))$ satisfies $(\epsilon_1+\epsilon_2-\frac{1}{2}\epsilon_1\epsilon_2)-$HDP.
\item (Parallel composition) Let $M_1$ and $M_2$ satisfy $\ep_1$ and $\ep_2$ HDP on two disjoint datasets $D_1$ and $D_2$ with distinct queries $w_1$ and $w_2$ respectively. Then, the parallel composition 
    \begin{align*}
    M^{(2)}(w_1,w_2, D_1, D_2) \coloneqq (M_1(w_1,D_1),M_2(w_2,D_2))
    \end{align*}
satisfies $\max\{\epsilon_1,\epsilon_2\}-$HDP.
\end{enumerate}
\end{theorem}
It is known that the privacy levels degrade ($\ep$ increases) with the number of compositions. Our next Corollary provides a useful quantification of this degradation after $j-$compositions.
\begin{corollary} \label{cor:composition}
Let $h_1(x)=x$ and for all $j\geq1$, the mechanism $M_j$ satisfies $\ep-$HDP. Then, for both adaptive and sequential compositions and for all $j \ge 1$, $M^{(j)}$ satisfies $h_j(\ep)-$HDP, where 
\begin{eqnarray*}
    h_j(x)=x+h_{j-1}(x)-\frac{1}{2}x h_{j-1}(x).
\end{eqnarray*}
For the parallel compositions, $M^{(j)}$ satisfies $\ep-$HDP.
\end{corollary}

When applying compositions to HDP, a key ingredient is the post-processing property. Specifically, for any mechanism $M(\cdot, \cdot)$ satisfying $\ep-$HDP, and $g: \Rm \mapsto \Rm$, the mechanism $g\circ M$ also satisfies $\ep-$HDP. This follows immediately from the post-processing inequality of the Hellinger distance (see \cite{Wu2017}). A natural next question concerns the relationship between HDP and other privacy measures. This is described in the next proposition.

\begin{proposition}
\label{prop:HDP_relation}
If $M$ satisfies $\epsilon$-HDP with $0<\epsilon<1$,
then $M$ also satisfies $(\epsilon',\delta')$ differential privacy,
where $\epsilon'=0$ and $\delta'=\sqrt{\epsilon}$.
Furthermore, if $M$ is an additive Gaussian mechanism with
a common dataset-independent covariance $\sigma^2 I_m$, or a
scalar additive Laplace mechanism with a common dataset-independent
scale, then $M$ also satisfies $\mu$-GDP, where
\(
\mu=2\Phi^{-1}\left(\frac{1+\sqrt{\epsilon}}{2}\right)
\).
\end{proposition}


As illustrated above, the definition of differential privacy is based on the pairs of adjacent datasets. However, in practice, it is convenient to define adjacent datasets as those that differ by $k$ records. This is common in applications such as healthcare, where one must protect groups of individuals. To address this scenario, we define group privacy using $k-$neighbor datasets. We say $D,D'$ are $k-$neighbor datasets if there exists datasets $D=D_0,D_1,\cdots,D_k=D'$ such that $D_i$ and $D_{i+1}$ are adjacent or identical for $i=0,\cdots,k-1$. That is,
\begin{align}\label{eq:knbr}
    ||D-D'||_{H}  = k.
\end{align}
Our next Theorem shows that HDP has a simple characterization for evaluating group privacy. This is in sharp contrast to other values of $\la$ (see, for instance, \cite{Mironov2017}).
\begin{theorem}[Group privacy] \label{thm:HDP_group}
If a mechanism $M(\cdot, \cdot)$ is $\ep-$HDP, $D$ and $D'$ satisfy (\ref{eq:knbr}) then 
\begin{align*}
 D_{HD}(M(w, D), M(w, D')) \le k^2 \ep.
\end{align*}
That is, for any  $k-$neighbor datasets, the mechanism is $k^2\ep-$HDP.
\end{theorem}
We now turn to implementing the HDP via the additive mechanisms involving the Gaussian and Laplace perturbations. We separate the following proposition from Theorem \ref{thm:PDP_mech} to focus on the HDP case.

\begin{proposition} \label{prop:HDP_mech}
Let $M(w,D)=w+\mathbf{Y}$ be an additive mechanism and   for $i=1,2$, $\Delta_{L_i} W$ be the $L_i$ sensitivity of $W$. Then
\begin{enumerate}
\item If $Y\sim N(0,\sigma^2\cdot\mathbf{I})$, then, to achieve $\epsilon-$HDP, \[\sigma^2=\frac{\left(\Delta_{L_2} W\right)^2}{8\log(\frac{1}{1-0.5\epsilon})}.\]
\item If $\mathbf{Y}=(Y_1,\cdots,Y_m)$, $Y_i\sim Lap(0,b)$, then to achieve $\epsilon-$HDP, $b=\frac{\Delta_{L_1}W}{2\log(\frac{1}{1-0.5\epsilon})}$. 
\end{enumerate}
\end{proposition}
If $m=1$, then a sharper value of $b$ for the Laplace mechanism can be obtained by using Lemma 8.3
in Section 8
and solving
\[-2\left[e^{-\frac{\Delta_{L_1}w}{2b}}+\frac{\Delta_{L_1}w}{2b}e^{-\frac{\Delta_{L_1}w}{2b}}-1\right]=\epsilon.
\] 
The exact value of $b$ for multidimensional parameter space is not explicit, and the previous proposition provides an upper bound. An extension of the above additive mechanism for matrix-valued queries is referred to as symmetric matrix mechanism and is given below. We will use this mechanism in Section \ref{sec:private_opti} for obtaining PMHDE using the PNR algorithm and in Section \ref{sec:CI} for constructing private confidence intervals.

\begin{proposition}\label{prop:symmetric_HDP}
Let the query $w:\mathcal{D}\to \mathbb{R}^{m\times m}$ be a matrix-valued function and $w(D)$ is a symmetric matrix. Then the additive mechanism $M(w,D)=w(D)+E$ satisfies $\ep-$HDP, where $E$ is a random upper triangle matrix including diagonals, whose components are i.i.d. random variables with distribution $N(0,\sigma^2)$ or $Lap(0,b)$, and $\sigma^2$ and $b$ are chosen using Proposition \ref{prop:HDP_mech}.
\end{proposition}
\noindent The proof of this proposition is similar to that given on page 14 of \cite{Dwork2014} Algorithm 1 and hence is omitted.

One common use of the composition property, post-processing rules, and additive mechanisms is in the optimization algorithms. For instance, in parametric estimation problems, legal and regulatory requirements may need privacy-preserving parameter estimates. A common approach is to modify an existing optimization algorithm to obtain private parameter estimates. Widely used optimization algorithms, such as GD and NR algorithms, iteratively update the estimators. Using the additive mechanism, with Gaussian or Laplace perturbations, it is possible to achieve the required levels of privacy at each iteration and ensure that the final iteration produces a desired private estimator. These modified algorithms are called PGD and PNR algorithms. Several versions of these optimization algorithms have been explored in the context of M-estimators: \cite{avella2021privacy,chaudhuri2011differentially,chaudhuri2012convergence,chen2019renyi,dalenius1977privacy,slavkovic2012perturbed,wang2017selecting}. 

\subsection{Operational interpretations}
\label{sec:hdp-operational}
The partition formulation expresses HDP through probabilities of output
events. We now use it to examine disclosure risks, posterior beliefs, and
hypothesis-testing performance. The comparisons depend on the privacy
parameters and the disclosure quantity under consideration.

Throughout this subsection, let $P_1,P_2$ be the output laws for an
ordered neighboring pair, with densities $p_1,p_2$ relative to a common
dominating measure $\nu$. Let $\rho(P_1,P_2)=\int\sqrt{p_1p_2}\,d\nu$
and $\operatorname{TV}(P_1,P_2)=\sup_A|P_1(A)-P_2(A)|$.
For $0\leq h\leq2$, set $c_h=1-h/2$ and
$s_h=(h-h^2/4)^{1/2}$. Thus, $D_{HD}(P_1,P_2)\leq h$ if and only if
$\rho(P_1,P_2)\geq c_h$.

\paragraph{Event probabilities and disclosure risk.}
The binary partition $\{A,A^c\}$ gives constraints on the probability
of every output event, including rare disclosures.

\begin{proposition}[Event-probability bounds]
\label{prop:hdp-op-events}
Suppose $D_{HD}(P_1,P_2)\leq h$. For $b\in[0,1]$, put
\[
L_h(b)=\left(c_h\sqrt b-s_h\sqrt{1-b}\right)_+^2,
\qquad U_h(b)=1-L_h(1-b).
\]
Then, for every measurable $A$,
\begin{equation}
L_h\{P_2(A)\}\leq P_1(A)\leq U_h\{P_2(A)\}.
\label{eq:hdp-op-sharp-events}
\end{equation}
Both bounds are sharp over pairs satisfying the Hellinger constraint.
In particular,
\begin{equation}
P_2(A)=0\quad\Longrightarrow\quad P_1(A)\leq h-h^2/4,
\label{eq:hdp-op-null-event}
\end{equation}
and, when $P_2(A)>0$,
\begin{equation}
\left|\frac{P_1(A)}{P_2(A)}-1\right|
\leq2\sqrt{\frac{h}{P_2(A)}}+\frac{h}{P_2(A)}.
\label{eq:hdp-op-relative}
\end{equation}
\end{proposition}

\noindent\emph{Proof.} See Supplementary Section 6.

In particular, if $P_2(A)=0$, observing $A$ distinguishes the two
neighboring dataset hypotheses with certainty, and Proposition~\ref{prop:hdp-op-events}
bounds its probability under $P_1$ by $h-h^2/4$.
For $h>0$, this probability need not vanish, whereas pure DP requires
it to be zero. The distinction concerns the form of protection, not
the existence of protection. For events with $b=P_2(A)>0$,
the relative-probability bound in Proposition~\ref{prop:hdp-op-events}
shows that the ratio $P_1(A)/b$ tends to one when $h/b\to0$.
Uniform versions and rare-event consequences are given in Supplementary Section 12.

\paragraph{Posterior-odds interpretation.}
The event bounds also control how observing the released output can
change an adversary's prior odds between two neighboring dataset
hypotheses. Let  \(w\in(0,1)\) and \(1-w\) denote the prior probabilities for neighboring dataset hypotheses inducing \(P_1\) and \(P_2\), respectively. Bayes' rule then gives
$\text{posterior odds}/\text{prior odds}=p_1(Y)/p_2(Y)$.
For $t>0$, put
$B_t=\{p_1>e^tp_2\}\cup\{p_2>e^tp_1\}$. The HDP constraint implies
\begin{equation}
\max\{P_1(B_t),P_2(B_t)\}
\leq \min\left\{1,\frac{h}{(1-e^{-t/2})^2}\right\}.
\label{eq:hdp-op-likelihood}
\end{equation}
Outside $B_t$, the posterior-to-prior odds multiplier belongs to
$[e^{-t},e^t]$. The probability bound also holds under any mixture
$wP_1+(1-w)P_2$, independently of the prior odds; $B_t$ may depend
on the neighboring pair. Thus HDP controls the probability of a large
change in posterior odds rather than requiring the multiplier to be
bounded almost surely, as pure DP does. Approximate DP, $\mu$-GDP,
and zCDP also permit unbounded multipliers but provide their own
quantitative controls.

We next extend this two-hypothesis interpretation to arbitrary priors
on the database using the counterfactual posterior formulation of
\citet{KasiviswanathanSmith2014}. Let $D\sim\Pi$, where $\Pi$ is an
arbitrary prior on $\mathcal D$, and let $P_D$ be the corresponding
output law of $M$. The marginal output distribution induced by $\Pi$
is $\overline P_\Pi(A)=\mathbf E_{D\sim\Pi}[P_D(A)]
=\int_{\mathcal D}P_u(A)\,\Pi(du)$.
For individual $i$, let $D^{(-i)}$ replace that individual's record
by a fixed reference value. Given output $y$, the actual and
counterfactual posteriors $\Pi_y$ and $\Pi_y^{(-i)}$ both concern $D$,
but use the likelihoods $P_D$ and $P_{D^{(-i)}}$, respectively, with
the same prior $\Pi$. We assume standard Borel database and output
spaces and fresh mechanism randomization so that these posteriors
are well defined.

\begin{proposition}[Counterfactual posterior stability]
\label{prop:hdp-op-posterior}
Suppose $M$ is $h$-HDP, $D\sim\Pi$, and $Y$ has the marginal
output law $\overline P_\Pi$ defined above. For each individual $i$,
let $\Pi_Y$ and $\Pi_Y^{(-i)}$ be the corresponding actual and
counterfactual posteriors of $D$. Then
\begin{equation}
\mathbf E_{\overline P_\Pi}\left[
\operatorname{TV}(\Pi_Y,\Pi_Y^{(-i)})^2\right]\leq s_h^2,
\label{eq:hdp-op-posterior-second}
\end{equation}
Consequently, for every $\eta>0$,
\begin{equation}
\mathbf P_{\overline P_\Pi}\left\{
\sup_C|\Pi_Y(C)-\Pi_Y^{(-i)}(C)|>\eta
\right\}
\leq\min\{1,s_h^2/\eta^2\}.
\label{eq:hdp-op-posterior-tail}
\end{equation}
\end{proposition}

\noindent\emph{Proof.} See Supplementary Section 6.

Since the supremum in \eqref{eq:hdp-op-posterior-tail} is inside the
probability, the conclusion holds simultaneously for all measurable
assertions about the database. Side information $Z$ about $D$,
available to an adversary, can also be incorporated by taking
$\Pi=\mathcal L(D\mid Z=z)$ for almost every $z$, provided
$\mathcal L(Y\mid D,Z)=P_D$. The resulting bound quantifies the effect
of replacing an individual's record beyond the information already
available to the adversary through $Z$. For a database of size $n$,
a union bound gives $\min\{1,ns_h^2/\eta^2\}$ as an upper bound on
the probability that the posterior discrepancy exceeds $\eta$ for
at least one individual.

\paragraph{Hypothesis-testing interpretation.}
The preceding results quantify how a release affects an adversary's
posterior beliefs. A complementary question is how reliably it can
distinguish two neighboring datasets; see \citet{Dong2022}.
Let $\varphi:\mathbb R^m\to[0,1]$ be a test function for testing
$H:Y\sim P_1$ against $K:Y\sim P_2$, with size
$\alpha=\mathbf E_{P_1}\varphi$ and power
$1-\beta=\mathbf E_{P_2}\varphi$.

\begin{proposition}[Testing bounds]
\label{prop:hdp-op-testing}
If $D_{HD}(P_1,P_2)\leq h$, every such test satisfies
\begin{equation}
\beta\geq
\left(c_h\sqrt{1-\alpha}-s_h\sqrt\alpha\right)_+^2,
\qquad
\alpha+\beta\geq1-s_h.
\label{eq:hdp-op-testing}
\end{equation}
\end{proposition}

\begin{proof}
Applying Cauchy--Schwarz separately with the weights $\varphi$ and
$1-\varphi$ yields
\[
c_h\leq\rho(P_1,P_2)
\leq\sqrt{\alpha(1-\beta)}+\sqrt{(1-\alpha)\beta}.
\]
The same binary calculation as in
Proposition~\ref{prop:hdp-op-events} gives the first bound in
\eqref{eq:hdp-op-testing}; binary experiments attain it.
Also,
\[
\operatorname{TV}(P_1,P_2)
\leq\frac12\left\{
\int(\sqrt{p_1}-\sqrt{p_2})^2\,d\nu
\int(\sqrt{p_1}+\sqrt{p_2})^2\,d\nu
\right\}^{1/2}
\leq s_h.
\]
Since
$\alpha+\beta
=1-\{\mathbf E_{P_2}\varphi-\mathbf E_{P_1}\varphi\}
\geq1-\operatorname{TV}(P_1,P_2)$,
the second inequality follows.
\end{proof}

\noindent\emph{Sharpness and Bayesian interpretation.}
The first inequality in \eqref{eq:hdp-op-testing} is sharp: for every
prescribed size $\alpha$, suitable binary output distributions
attain equality under the Hellinger constraint. Suppose $H$ and $K$
each have prior probability $1/2$. Then the unconditional
probabilities of rejecting $H$ when it is true and accepting $H$
when $K$ is true are $\alpha/2$ and $\beta/2$, respectively.
Consequently, the marginal probability of correctly identifying the
true hypothesis is $1-(\alpha+\beta)/2$, which by
Proposition~\ref{prop:hdp-op-testing} is at most $(1+s_h)/2$.

Sharp conversions between HDP and approximate DP (Propositions 19
and 20), together with examples distinguishing their eventwise
guarantees, appear in Supplementary Section 6.

\addtocounter{theorem}{2}

\subsection{Exact calibration across privacy notions}
\label{sec:privacy-calibration}
The preceding subsection establishes operational interpretations of
HDP through event-probability bounds, posterior stability, and
hypothesis-testing guarantees. Sharp conversions between HDP and
approximate DP are developed in the supplementary material. For
neighboring output laws $P_1,P_2$, define the directional privacy
profiles $\delta_{12}(t)=\sup_A\{P_1(A)-e^tP_2(A)\}$ and
$\delta_{21}(t)=\sup_A\{P_2(A)-e^tP_1(A)\}$, $t\geq0$.
An $(\epsilon,\delta)$-DP guarantee requires both profiles to be
bounded by $\delta$ at $t=\epsilon$, uniformly over neighboring
datasets. In particular,
$\delta_{12}(0)=\delta_{21}(0)=\operatorname{TV}(P_1,P_2)$.

We next turn to exact calibration within the Gaussian family,
relating the noise variance to the attained Hellinger level and
approximate-DP profile. These calculations also allow quantitative
comparisons of noise requirements under specified composition rules.

\paragraph{Gaussian calibration.}
Consider the additive Gaussian mechanism
$M_\sigma(D)=W(D)+\sigma Z$, where $Z\sim N_m(0,I_m)$ is independent
of the dataset and $\sigma>0$ is dataset-independent. Assume
$0<\Delta:=\sup_{D\sim D'}\|W(D)-W(D')\|_2<\infty$.
For neighboring datasets $D_1,D_2$, write $w_i=W(D_i)$,
$P_i=N_m(w_i,\sigma^2I_m)$, $d=\|w_1-w_2\|_2$, and
$\mu_d=d/\sigma$. Let $\Phi$ and $\phi$ denote the standard normal
distribution function and density.

\begin{proposition}[Exact Gaussian privacy quantities]
\label{prop:cal-gaussian}
For $d>0$, the Hellinger distance and directional privacy profiles
satisfy
\begin{align}
D_{HD}(P_1,P_2)&=2(1-e^{-\mu_d^2/8}),
\label{eq:cal-gaussian-pair-h}\\
\delta_{12}(\epsilon)=\delta_{21}(\epsilon)
&=\delta_G(\epsilon;\mu_d),\qquad \epsilon\geq0.
\label{eq:cal-gaussian-pair-profile}
\end{align}
Here, for $\mu>0$,
$\delta_G(\epsilon;\mu)
=\Phi(\mu/2-\epsilon/\mu)
-e^\epsilon\Phi(-\mu/2-\epsilon/\mu)$.
\end{proposition}

\noindent\emph{Proof.} See Supplementary Section 6.

In particular,
$\operatorname{TV}(P_1,P_2)=\delta_G(0;\mu_d)
=2\Phi(\mu_d/2)-1$; when $d=0$, the laws coincide.
Both discrepancies increase with $\mu_d$, so their worst-case values
over neighboring datasets result from replacing $\mu_d$ by
$\mu=\Delta/\sigma$, even when the supremum is not attained.

For these Gaussian output laws, the power divergence is
\begin{equation}
D_\lambda(P_1,P_2)=
\begin{cases}
\displaystyle
\frac{\exp\{\lambda(\lambda+1)\mu_d^2/2\}-1}
{\lambda(\lambda+1)},&\lambda\notin\{-1,0\},\\[4pt]
\mu_d^2/2,&\lambda\in\{-1,0\}.
\end{cases}
\label{eq:cal-gaussian-power}
\end{equation}
In particular, $D_{-1/2}=2D_{HD}$, consistent with our
normalization. For $\alpha>1$, the corresponding R\'enyi divergence
is $D_\alpha^{(R)}(P_1,P_2)=\alpha\mu_d^2/2$ \citep{Mironov2017}.
These identities follow from the Gaussian likelihood-ratio moments
in the supplementary material.

\paragraph{Forward and reverse Gaussian calibration.}
Let $\sigma_{\mathrm{BW}}^2(\epsilon,\delta)$ denote the minimum
Gaussian variance at the specified approximate-DP parameters.

\begin{corollary}[Gaussian parameter correspondence]
\label{cor:cal-gaussian-correspondence}
Fix $\epsilon\geq0$ and $0<\delta<1$, and let $\mu_*>0$ be the
unique solution of $\delta_G(\epsilon;\mu_*)=\delta$.
The minimum Gaussian variance satisfying $(\epsilon,\delta)$-DP
and the corresponding attained HDP level are
\begin{equation}
\sigma_*^2:=\sigma_{\mathrm{BW}}^2(\epsilon,\delta)
=\frac{\Delta^2}{\mu_*^2},\qquad
h_{\mathrm{achieved}}=2(1-e^{-\mu_*^2/8}).
\label{eq:cal-gaussian-forward-variance}
\end{equation}
Conversely, for $0<h<2$, put
$\mu_h=\sqrt{-8\log(1-h/2)}$.
The minimum Gaussian variance satisfying $h$-HDP and its exact
approximate-DP profile are
\begin{equation}
\sigma_H^2(h)=\frac{\Delta^2}{\mu_h^2},\qquad
\delta_h^G(\epsilon)=\delta_G(\epsilon;\mu_h),\quad\epsilon\geq0.
\label{eq:cal-gaussian-reverse}
\end{equation}
At this variance, the mechanism satisfies $(\epsilon,\delta)$-DP
if and only if $\delta\geq\delta_h^G(\epsilon)$.
Every boundary pair $(\epsilon,\delta_h^G(\epsilon))$
has minimum Gaussian variance $\sigma_H^2(h)$.
\end{corollary}

\begin{proof}
For fixed $\epsilon\geq0$, the function
$\delta_G(\epsilon;\mu)$ is strictly increasing for $\mu>0$, since
$\partial_\mu\delta_G(\epsilon;\mu)
=\phi(\mu/2-\epsilon/\mu)>0$.
Its limits as $\mu\downarrow0$ and $\mu\to\infty$ are $0$ and $1$.
Thus $\delta_G(\epsilon;\Delta/\sigma)\leq\delta$ if and only if
$\sigma^2\geq\Delta^2/\mu_*^2$.
Proposition~\ref{prop:cal-gaussian} gives the corresponding
Hellinger level. Conversely, solving
$2(1-e^{-\Delta^2/(8\sigma^2)})\leq h$ gives
$\sigma_H^2(h)$. The remaining claims follow from the exact Gaussian
privacy profile and its strict monotonicity.
\end{proof}

At $\epsilon=0$, the Gaussian profile gives
$\delta_h^G(0)=2\Phi(\mu_h/2)-1$, the total variation distance.
For every finite $\epsilon$, a nonconstant Gaussian mechanism with
finite positive variance satisfies $\delta_h^G(\epsilon)>0$ and hence
does not satisfy pure $\epsilon$-DP. The reverse correspondence gives
an entire approximate-DP boundary; pairs above it are weaker
guarantees for the same mechanism. An upper bound on $\Delta$
gives certified rather than necessarily attained privacy levels.

\begin{figure}[H]
\centering
\begin{minipage}[t]{0.485\textwidth}
\centering
\includegraphics[width=\linewidth,trim=0 0 0 76bp,clip]{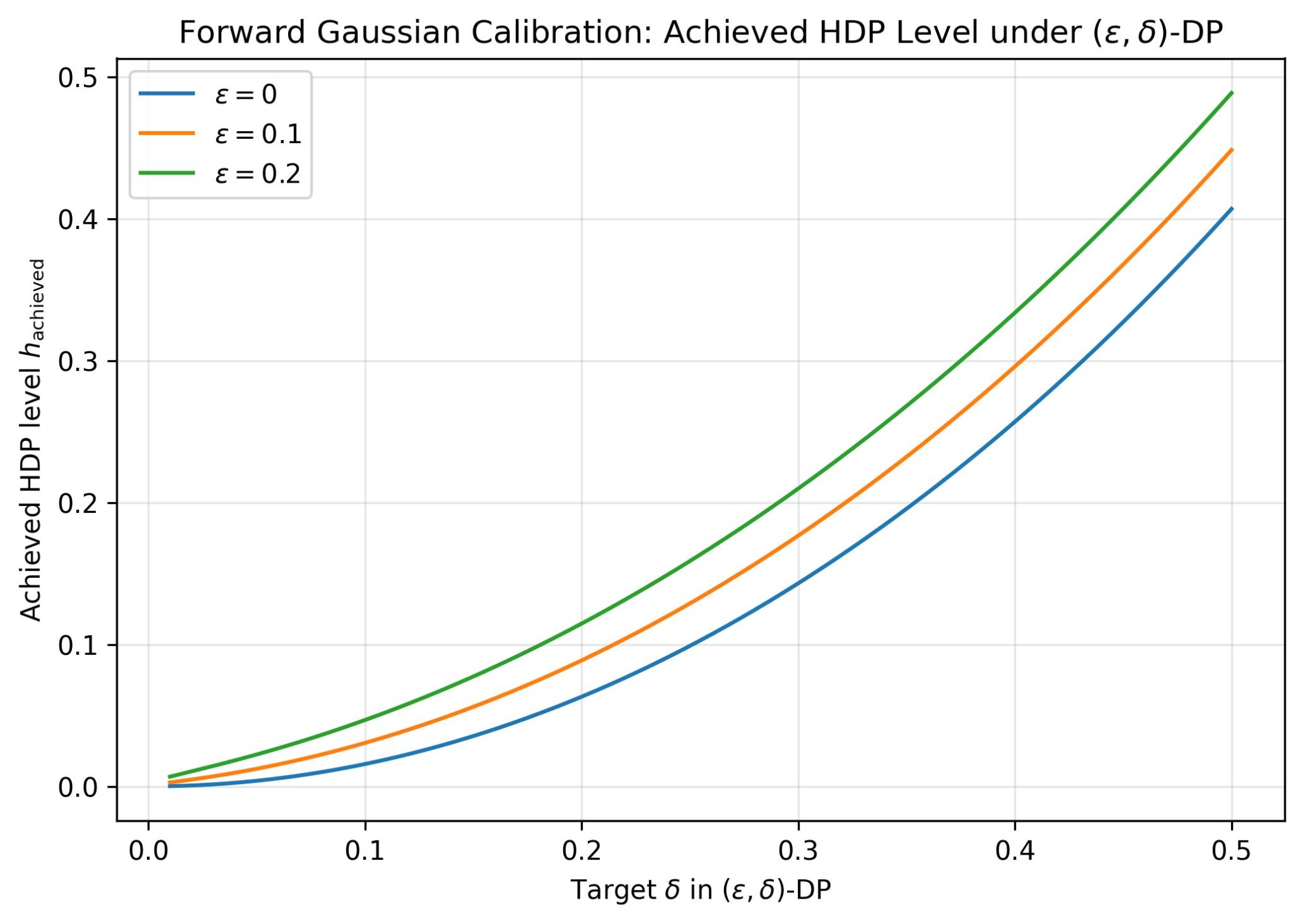}
\par\smallskip{\small (a) Gaussian: $(\epsilon,\delta)$-DP to HDP}
\end{minipage}\hfill
\begin{minipage}[t]{0.485\textwidth}
\centering
\includegraphics[width=\linewidth,trim=0 0 0 76bp,clip]{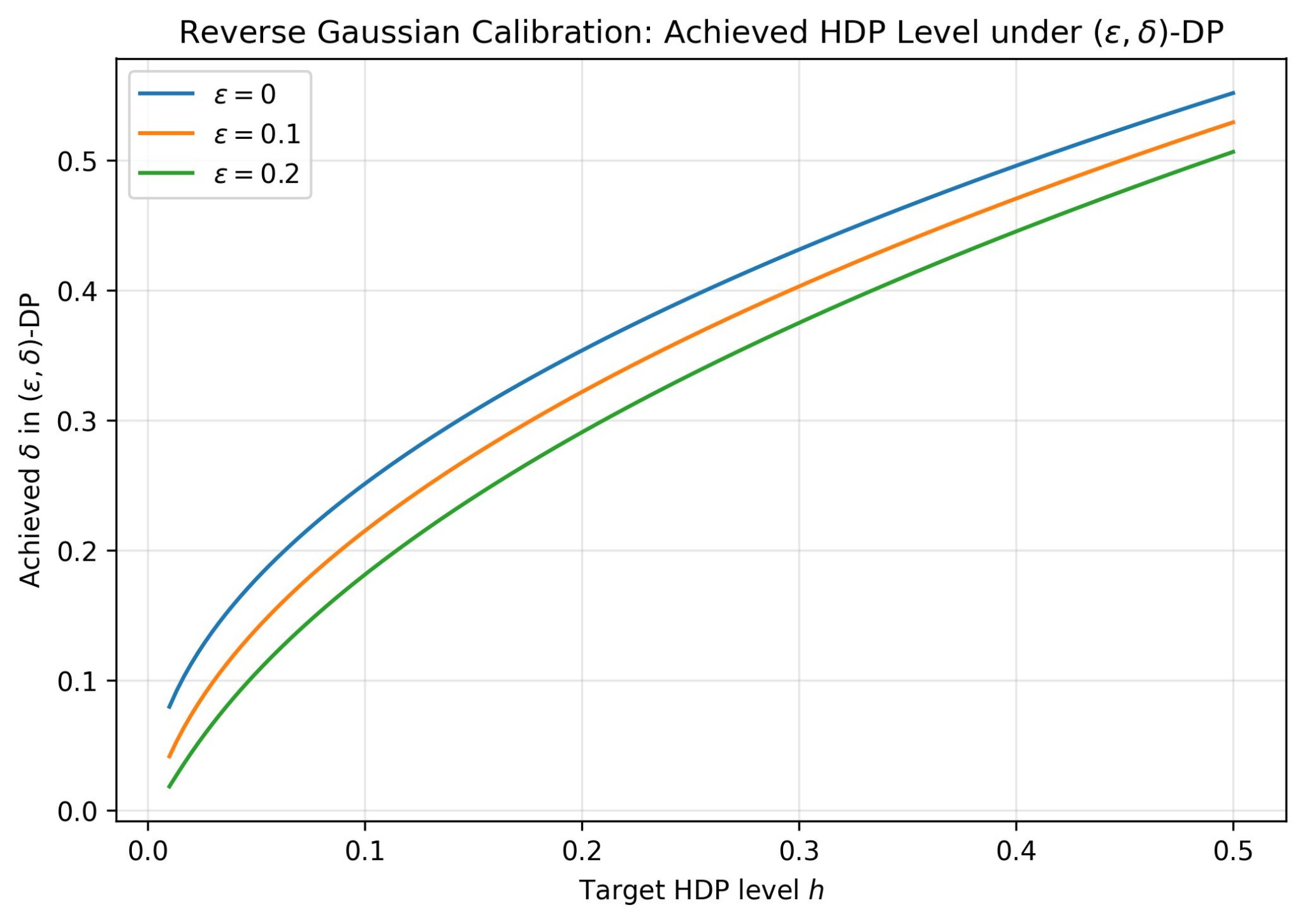}
\par\smallskip{\small (b) Gaussian: HDP to $(\epsilon,\delta)$-DP}
\end{minipage}
\caption{Exact Gaussian forward and reverse calibration. The left
panel gives the attained Hellinger level under approximate-DP
calibration; the right gives the attained approximate-DP boundary at
a fixed Hellinger level. The complete Gaussian and scalar Laplace
comparisons are in Supplementary Figure 11.4.}
\label{fig:calibration-main}
\end{figure}

Figure~\ref{fig:calibration-main} illustrates both Gaussian calibration
directions, evaluating the attained privacy quantities at the variance
chosen for the initial requirement. Additional numerical comparisons
and the exact scalar Laplace calibration (Proposition 23) are in the
supplementary material.

\addtocounter{theorem}{1}
\begin{remark}[Divergence scale and noise calibration]
\label{rem:cal-divergence-scale}
For Gaussian additive mechanisms, the perturbation variance in
Theorem~\ref{thm:PDP_mech} is minimized at $\lambda=-1/2$ within the
power-divergence family at a common numerical PDP budget $0<\epsilon<4$. The exact
Gaussian calibrations above identify the attained operational privacy
guarantees; further comparisons are in the supplementary material.
\end{remark}

\paragraph{Calibration under composition.}
For $K$ Gaussian releases with common $L_2$-sensitivity bound $\Delta$,
let $\sigma_0^2=\sigma_H^2(h_0)$ be the reference variance at
prescribed HDP level $h_0$, and put
$\mu_0=\Delta/\sigma_0$ and
$\delta_0=\delta_G(\epsilon_0;\mu_0)$.
Under equal HDP allocation, the per-step variance is $K\sigma_0^2$.
Basic approximate-DP composition instead allocates
$(\epsilon_0/K,\delta_0/K)$ at each step, requiring
\begin{equation}
\sigma_{H,K}^2=K\sigma_0^2,\qquad
\sigma_{\mathrm{basic},K}^2
=\sigma_{\mathrm{BW}}^2(\epsilon_0/K,\delta_0/K).
\label{eq:cal-equal-variance}
\end{equation}

For example, with $\Delta=1$, $h_0=0.50$, $K=10$, and
$\epsilon_0=1$, the corresponding final $\delta_0\approx0.33$.
HDP allocation requires per-step Gaussian variance $4.35$, compared
with $32.36$ under basic approximate-DP composition, approximately
$7.45$ times as large. The corresponding final approximate-DP
guarantee in this Gaussian setting follows from the Gaussian testing
and composition bound detailed in the supplementary material.
The complete composition calculations and other comparisons are given
there as well.

\section{Private minimum Hellinger distance estimation} 

In this section, we will briefly discuss minimum Hellinger distance estimation for continuous
i.i.d. data and modify the estimation method to satisfy HDP using PGD and PNR algorithms. We also study the consistency and efficiency of the PMHDEs.

\subsection{Minimum Hellinger distance estimation} \label{sec:MHDE}

The minimum Hellinger distance estimation method for i.i.d. observations, proposed in Beran (\cite{Beran1977}), has been extended to several statistical models, including dependent data (see \cite{Basu2011, Cheng2006, Li2019}). A useful feature of these estimators is that they are, like maximum likelihood estimators (MLEs), efficient when the posited parametric model is true. However, unlike MLEs, they are also robust with a ``high breakdown point''. In other words, the MHDE achieves the dual goal of robustness and efficiency in the true model.  For a comprehensive discussion of minimum divergence theory, see \cite{Basu2011}. 

Let $\{X_1, X_2,\cdots, X_n\}$ be i.i.d. real-valued random variables with density $g(\cdot)$, and postulated to belong to a parametric family $\{f_{\bta} : \bta \in \Theta\subset\mathbb{R}^m\}$. The minimum Hellinger distance estimator in the population, $\bta_g$,  if it exists, is the minimizer of the
$||f_{\bta}^{\frac{1}{2}}-g^{\frac{1}{2}}||_2$; that is,
\begin{eqnarray*}
\bta_g=\argmin ||f_{\bta}^{\frac{1}{2}}-g^{\frac{1}{2}}||_2 = \argmin HD(f_{\bta}, g) .
\end{eqnarray*}
\cite{Beran1977} and \cite{Cheng2006} establish the existence of $\bta_g$ under family regularity, described in Section 7.
Replacing $g(\cdot)$ by $g_n(\cdot)$, where $g_n(\cdot)$ is a nonparametric estimator of $g(\cdot)$,
one obtains the MHDE. In this paper, we use the kernel density estimator (KDE) of $g(\cdot)$; namely,
\begin{align*}
    g_n(x) = \frac{1}{n\cdot c_n}\sum_{i=1}^{n} K\left(\frac{x-X_i}{c_n}\right),
\end{align*}
where $K(\cdot)$ is a kernel density with support $(-\beta,\beta)$ for $\beta\in(0,\infty)$, and $c_n$ (referred to as bandwidth) is a sequence of constants converging to 0 such that $n c_n \ra \ff$. Thus, the loss function of the MHDE is given by
\begin{eqnarray}\label{eq:loss}
    L_n(\bta) = 2 HD^2(f_{\bta}, g_n) = 2 \int_{\mathbb{R}}(\sqrt{f_{\bta}(x)}-\sqrt{g_n(x)})^2 \mathrm{d}x,
\end{eqnarray}
where we include factor 2 to draw connections to the general power divergence family described above. We notice here that other non-parametric density estimators, such as wavelet-based density estimators, can be used since they possess similar $L_1$ properties like the KDE (see \cite{chacon2005l1}).
Statistical properties such as consistency and asymptotic normality of the MHDE have been established under the assumptions {\bf(A1)-(A8)} 
in Section 7.
In the rest of the paper we assume that these conditions hold.

Computationally, the estimators are typically derived using optimization algorithms such as GD and the NR method. Using an ``additive mechanism'' of the HDP described in Section \ref{sec:HDP}, we derive private versions of these estimators. The mechanism involves adding \emph{appropriate noise} at every iteration of the optimization algorithm, referred to as the PGD and PNR algorithms. The resulting optimization is called private optimization (also referred to in the literature as noisy optimization). The variability induced by noise addition depends on the $L_2-$sensitivity of the gradient and Hessian of $L_n(\bta)$. Analysis of this
is much more subtle, unlike the M-estimator, and requires some new technical ideas (see Theorem \ref{thm:sensitivity}  below), and when incorporated into the algorithms, allows an improved numerical performance. We now describe PGD and PNR algorithms and study their statistical properties.

\subsection{Almost sure local convexity}
In this section, we leverage the properties of Hellinger distance,  assumptions {\bf(A1)-(A8)},
and additional moment conditions to establish almost surely locally strongly convex (ASLSC) properties of the loss function $L_n(\bta)$.  To this end, we need a few additional notations.

\begin{assumptionU}\label{asp:U1}
Let $u_{\bta,i}(x)= \frac{\partial}{\partial \ta_i} \log f_{\bta}(x)$.
Assume that for all ~ $1 \le j \le m, 0 \le k_j \le 6 $ and  $k_1+ k_2+ \cdots k_m \le 6$,
$$\mathbf{E}_{\bta}\left[\prod_{i=1}^m |u_{\bta,i}(X)|^{k_i}\right] < \ff.$$
Additionally, assume that the expectation above is continuous in $\bta$.
\end{assumptionU}

\begin{assumptionU}\label{asp:U2}
Assume that all the partial and cross-partial derivatives of $f_{\bta}$ up to order three exist and are continuous. Set $u_{\bta,i,j}(x)= \frac{\partial}{\partial \ta_i} u_{\bta,j}(x), ~\text{and} \quad 
u_{\bta,i,j,l}(x)= \frac{\partial}{\partial \ta_i} u_{\bta,j,l}(x)$.
Assume that for all $ 1 \le i, j,l \le m$,
$\mathbf{E}_{\bta}\left[|u_{\bta, i,j}(X)|^2\right] < \ff$ and $\mathbf{E}_{\bta}\left[|u_{\bta, i,j,l}(X)|^2\right] < \ff$. Also, the Fisher information matrix is positive definite for all $\bta \in \Theta$ and, in particular,
$I(\bta)=((I_{ij}(\bta) \coloneqq -\mathbf{E}_{\bta}\left[u_{\bta, i,j}(X)\right])) < \ff$ for all $\bta \in \Theta$.
Additionally, for all $1\leq i,j,l\leq m$,
$\mathbf{E}_{\bta}
[u_{\bta,l}^{2}(X)u_{\bta,i,j}^{2}(X)]<\infty$.
The expectations in this assumption are continuous in $\bta$
on $\Theta$.
\end{assumptionU}

Let $\nabla L_n(\bta)$ and $H_n(\bta)$ denote the gradient and Hessian of $L_n$; that is,
\begin{align*}
    \nabla L_n(\bta) =& -2 \int g_n^{\frac{1}{2}}(x)f_{\bta}^{\frac{1}{2}}(x) \mathbf{u}_{\bta}(x) \mathrm{d}x \quad \text{and}\\
    H_n(\bta) =& -\int g_n^{\frac{1}{2}}(x)f_{\bta}^{\frac{1}{2}}(x) [\mathbf{u}_{\bta}(x)\mathbf{u}^T_{\bta}(x)+2\mathbf{\dot u}_{\bta}(x)]\mathrm{d}x.
\end{align*}  
In the above, $\mathbf{u}_{\bta}(x)= \nabla_{\bta}\log f_{\bta}(x)=[u_{\bta,1}(x),\cdots,u_{\bta,m}(x)]^T$ is the score vector and $\mathbf{\dot u}_{\bta}(x)$ is the matrix of \AV{second} partials of \AV{$\log f_{\bta}(x)$} with respect to components of $\bta$. That is,
\begin{align*}
\mathbf{\dot u}_{\bta}(x)=\begin{bmatrix}
        u_{\bta,1,1}(x) & \cdots &  u_{\bta,1,m}(x)\\
        \vdots & \ddots & \vdots\\
        u_{\bta,m,1}(x) & \cdots &  u_{\bta,m,m}(x)
    \end{bmatrix}.
\end{align*}
Our next proposition establishes the uniform boundedness of the gradient and the Hessian of the loss function.
\begin{proposition} \label{prop:GD_bound}
    With probability 1, $\sup\limits_{\bta\in\Theta}||\nabla L_n(\bta)||_2\leq B_1$, $\sup\limits_{\bta\in\Theta}||H_n(\bta)||_2\leq B_2$, for some constants $B_1,B_2\in(0,\infty)$.
\end{proposition}

\begin{proof}
First, note, using the Cauchy-Schwarz inequality, that
\[
||\nabla L_n(\bta)||_2 \le  2\sqrt{\operatorname{tr}I(\bta)}.
\]
Hence, by taking the supremum on both sides of the above inequality and using assumption {\bf(A3)}
and the compactness of $\Theta$, it follows  for some $0 < B_1 < \ff$
\[
\sup_{\bta \in \Theta} ||\nabla L_n(\bta)||_2 \le  B_1.
\]
Turning to $H_n(\bta)$, we apply Cauchy-Schwarz inequality to every component of the Hessian matrix and use assumption~{\bf\ref{asp:U1}} and the compactness of $\Ta$ to verify that there exists a constant $B_2$  such that $\sup_{\bta \in \Ta}||H_n(\bta)|| \le B_2$.
\end{proof}

Next, we turn to almost sure convergence of the Hessian matrices. We note that the sequence $A_n$ of $m \times m$ matrices converge to $A$ if the $(i,j)^{th}$ element of $A_n$ converges to $(i,j)^{th}$ element of $A$. 
\begin{proposition}{\label{prop:loss_converge}}
Under the assumptions {\bf(A1)-(A8)}
and {\bf\ref{asp:U1}}-{\bf\ref{asp:U2}}, the Hessian matrix $H_n(\bta)$ converges almost surely to $H_{\ff}(\bta)$ for all $\bta \in \Ta$ and $H_{\ff}(\bta)= I(\bta)- D(\bta)$, where the $(i,j)^{\text{th}}$ element of $D(\bta)$ is given by
\begin{eqnarray*}
D_{i,j}(\bta) = \int_{\mathbb{R}}\left(g^{\frac{1}{2}}(x)-f_{\bta}^{\frac{1}{2}}(x)\right)f_{\bta}^{\frac{1}{2}}(x) [u_{\bta,i}(x)u_{\bta,j}(x)+2 u_{\bta,i,j}(x)]\mathrm{d}x.
\end{eqnarray*}
Furthermore, $H_{\ff}(\bta)$ is continuous in $\bta$.
\end{proposition}

\begin{proof}
We will first establish that $H_{n, i,j}(\bta)$ converges to $H_{\ff,i,j}(\bta)$. To this end, using equation (42) 
in Section 10,
notice that $H_{n, i,j}(\bta)= I_{i,j}(\bta) - D_{n,i,j}(\bta) ~\text{where}$
\begin{eqnarray*}
D_{n,i,j}(\bta) = \int_{\mathbb{R}}\left(g_n^{\frac{1}{2}}(x)-f_{\bta}^{\frac{1}{2}}(x)\right)f_{\bta}^{\frac{1}{2}}(x) [u_{\bta,i}(x)u_{\bta,j}(x)+2 u_{\bta,i,j}(x)]\mathrm{d}x
\end{eqnarray*}
Now, adding and subtracting $g^{\frac{1}{2}}(\cdot)$ to the RHS of the above equation, we obtain
\begin{align*}
    D_{n,i,j}(\bta)
    = D_{i,j}(\bta) \AV{+} \int_{\mathbb{R}}\left(g_n^{\frac{1}{2}}(x)-g^{\frac{1}{2}}(x)\right)f_{\bta}^{\frac{1}{2}}(x) [u_{\bta,i}(x)u_{\bta,j}(x)+2 u_{\bta,i,j}(x)]\mathrm{d}x.
\end{align*}
Next, applying the Cauchy-Schwarz inequality to the second term on the RHS of the above equation and using assumptions {\bf\ref{asp:U1}} and {\bf\ref{asp:U2}}
it follows that $D_{n, i,j}(\bta)$ converges almost surely to $D_{i,j}(\bta)$. This implies convergence of $H_n(\bta)$ to $H_{\ff}(\bta)$. Now, combining the above equations, the expression for $H_{\ff}(\bta)$ follows. Turning to the continuity of $H_{\ff}(\bta)$, it follows from Cauchy-Schwarz inequality,  Scheffe's Theorem, and Assumption {\bf(A4)} 
that $D(\bta)$ is continuous in $\bta$. Finally, continuity of $H_{\ff}(\bta)$ follows from {\bf(A3)} 
and the continuity of $D(\bta)$.
\end{proof}

\noindent Below, we use $\la_{min}(A)$ and $\la_{max}(A)$ to denote the minimum and maximum eigenvalue of a square matrix $A$.
\begin{proposition} \label{prop:Hessian_positive}
Under assumptions {\bf(A1)-(A8)}
and {\bf\ref{asp:U1}}-{\bf\ref{asp:U2}}, there exists an $\ep >0$ such that if $HD(g,f_{\bta_g})<\ep$, then there exists an open ball of radius $r_\ep$, centered at $\bta_g$, $B_{r_\ep}(\bta_g)$, such that for all $\bta \in B_{r_\ep}(\bta_g)$, $H_{\infty}(\bta)$ is strictly positive definite. Furthermore, $\la_{max}(H_{\ff}(\bta)) \le C$, where $0 <C < \ff$ is independent of $\bta$.
\end{proposition}

\begin{proof}
First, notice that using equation (43)
in Section 10
that $D_{i,j}(\bta) \le c\cdot HD(g, f_{\bta})$ where $c>0$ is independent of $\bta$, which implies that $D(\bta) \le C'HD(g, f_{\bta}) \mathbf{J}_m $ where $\mathbf{J}_m $ is a $m \times m$ matrix of ones and $0 <C' < \ff$. If
$HD(g,f_{\bta_g}) < \ep$ is small, then by Proposition \ref{prop:loss_converge}  and Weyl's inequality, it follows that the minimal eigenvalue of $H_{\ff}(\bta_g)$, $\la_{min}(H_{\ff}(\bta_g))$, is close to that of $I(\bta_g)$; that is, there exists $\ep'$ such that $|\la_{min}(H_{\ff}(\bta_g))- \la_{min}(I(\bta_g))| < \ep^{\prime}$. Since $H_{\ff}(\bta)$ is continuous in $\bta$ by Proposition \ref{prop:loss_converge} and  $\la_{min}(I(\bta_g)) >0$, it follows that there exists a neighborhood $B_{r_{\ep}}(\bta_g)$ such that $\la_{min}(H_{\ff}(\bta)) >0$ for all $\bta \in B_{r_{\ep}}(\bta_g)$. The proof regarding $\la_{max}(H_{\ff}(\bta))$ is similar.
\end{proof}

\noindent Our next proposition is concerned with the ASLSC of $H_n(\bta).$

\begin{proposition} \label{prop:convex}
Let assumptions {\bf(A1)-(A8)}
and {\bf\ref{asp:U1}}-{\bf\ref{asp:U2}} hold.
Also, assume that $HD(g,f_{\bta_g})<\ep$, where $\ep$ is
as in Proposition~\ref{prop:Hessian_positive}.
Then there exists an open ball of radius $r$, centered at $\bta_g$, $B_r(\bta_g)$, and $0<\tau_1\leq \tau_2<\infty$ (independent of $\bta$) and an almost surely finite random integer $N$ such that for all $\bta \in B_r(\bta_g)$ and $n \ge N$,
\begin{align*}
    \tau_1\leq \la_{min}(H_n(\bta))\leq \la_{max}(H_n(\bta))\leq \tau_2
\end{align*}
with probability one. That is, $L_n(\bta)$ is almost surely locally strongly convex and $\tau_2-$smooth.
\end{proposition}

\begin{proof}
By Proposition~\ref{prop:Hessian_positive}, choose $r>0$ so that
$K=\bar B_r(\bta_g)$ is contained in its positivity neighborhood.
Continuity of $H_\infty(\bta)$ and compactness of $K$ give
\[
a_K=\inf_{\bta\in K}\la_{min}(H_\infty(\bta))>0,
\qquad
b_K=\sup_{\bta\in K}\la_{max}(H_\infty(\bta))<\infty.
\]
From the Hessian formulas and Cauchy--Schwarz, for $\bta\in K$,
\begin{align*}
|H_{n,i,j}(\bta)-H_{\infty,i,j}(\bta)|^2
&\leq HD^2(g_n,g)\,
\mathbf E_{\bta}\left[
\{u_{\bta,i}(X)u_{\bta,j}(X)+2u_{\bta,i,j}(X)\}^2\right]\\
&\leq HD^2(g_n,g)\left\{
2\mathbf E_{\bta}[u_{\bta,i}^2(X)u_{\bta,j}^2(X)]
+8\mathbf E_{\bta}[u_{\bta,i,j}^2(X)]\right\}.
\end{align*}
The moment terms are finite and continuous in $\bta$ by
{\bf\ref{asp:U1}}--{\bf\ref{asp:U2}}, and hence uniformly bounded on $K$.
Summing the squared entrywise bounds yields a deterministic
$C_K<\infty$ such that
\[
\sup_{\bta\in K}\|H_n(\bta)-H_\infty(\bta)\|_2
\leq C_K HD(g_n,g)\xrightarrow{a.s.}0,
\]
The convergence follows from $HD^2(g_n,g)\leq\|g_n-g\|_1$ and the
almost sure $L_1$ convergence of $g_n$ to $g$.
There is therefore an almost surely finite random integer $N$ such that
the uniform error is at most $a_K/2$ for every $n\geq N$.
Weyl's inequality gives the asserted bounds with
$\tau_1=a_K/2$ and $\tau_2=b_K+a_K/2$.
\end{proof}

\noindent Our next proposition summarizes some useful properties of $L_n(\bta)$ and is based on the definition of almost sure $\tau_1$ strong convexity and $\tau_2$ smoothness. The proof is similar to the discussion in \cite{Boyd2004} Section 9.1.2.
\begin{proposition} \label{prop:convex_ineq}
The following inequalities hold for $\bta_1,\bta_2\in B_r(\bta_g)$:
\begin{enumerate}
\item $L_n(\bta_1)\geq L_n(\bta_2)+\langle\nabla L_n(\bta_2),(\bta_1-\bta_2)\rangle+\frac{\tau_{1}}{2}||\bta_1-\bta_2||_2^2$.
        \item $\langle\nabla L_n(\bta_1)-\nabla L_n(\bta_2),\bta_1-\bta_2\rangle\geq \tau_{1}||\bta_1-\bta_2||_2^2$.
        \item $L_n(\bta_1)\leq L_n(\bta_2)+\langle\nabla L_n(\bta_2),(\bta_1-\bta_2)\rangle+\frac{\tau_{2}}{2}||\bta_1-\bta_2||_2^2$.
        \item $\langle\nabla L_n(\bta_1)-\nabla L_n(\bta_2),\bta_1-\bta_2\rangle\leq \tau_{2}||\bta_1-\bta_2||_2^2$.
    \end{enumerate}
\end{proposition}
\noindent Our next Proposition is concerned with the almost sure Lipschitz property of the Hessian of the Hellinger loss function, (\ref{eq:loss}) which is required to establish certain utility properties of our proposed algorithms in Section \ref{sec:private_opti} below.  The proof is in Section 6.
\begin{proposition} \label{prop:Lipschitz}
Under the assumptions {\bf(A1)-(A8)}
and {\bf\ref{asp:U1}}-{\bf\ref{asp:U2}},
the Hessian matrix $H_n(\bta)$ is almost surely Lipschitz; that is, if $\bta_1,\bta_2\in\Ta$, then there exists $\al \in (0, \ff)$ such that
\begin{eqnarray*}
    ||H_n(\bta_1)-H_n(\bta_2)||_2 \le \alpha ||\bta_1-\bta_2||_2
\end{eqnarray*}
holds for any $n$ with probability one.
\end{proposition}

\subsection{Private optimization} \label{sec:private_opti}

As explained above, in this section, we systematically develop private versions of the GD and NR algorithms. We begin by observing that the estimator is a solution to the Hellinger-score equation
\begin{eqnarray} \label{eq:GD}
    \nabla{L_n(\bta)}=0, \quad \text{where } \nabla{L_n(\bta)}=\left[ \frac{\partial L_n(\bta)}{\partial \theta_1},\cdots, \frac{\partial L_n(\bta)}{\partial \theta_m} \right]^T.
\end{eqnarray}
We first consider the GD algorithm. To obtain the solution to (\ref{eq:GD}), given a potential root $\hat{\bta}_n^{(k)}$ of the equation, we obtain an updated solution by minimizing the objective function,
\[
Q(\bta)=L_n(\hat\bta_n^{(k)}) + \langle\nabla L_n(\hat\bta_n^{(k)}),\bta-\hat\bta_n^{(k)}\rangle + \frac{1}{2\eta}||\bta-\hat\bta_n^{(k)}||_2^2.
\]
Taking the derivative and setting it equal to zero, one obtains $\hat\bta_{n}^{(k+1)} = \hat\bta_{n}^{(k)} - \eta \nabla L_n(\hat\bta_{n}^{(k)})$, where $\eta$ is a pre-determined step-size and frequently referred to as the learning rate.  The idea is to update the estimator $\hat\bta^{(k)}_n$ until it reaches the zero of $\nabla L_n(\bta)$. Letting $k$ increase without bound ensures that $\hat\bta_{n}^{(k)}$ is close to $\hat\bta_{n}$, where $\hat\bta_{n}$ is the stationary point of $L_n(\bta)$. It is known that $\hat\bta_n$ is not guaranteed to be the global minimizer of the loss function (see \cite{agarwal2009}). However, under family regularity and a large sample size, the algorithm will converge to the global minimizer by choosing the starting point appropriately. We focus on the private version of the above algorithm, and as explained previously, we introduce an appropriate amount of noise in each iteration of the optimization algorithm. Specifically, using the additive mechanism described in the previous section, we introduce the noise $N_k$ to obtain the private version of $Q(\cdot)$. That is,
\[Q_k(\bta)=L_n(\hat\bta_n^{(k)}) + \langle\nabla L_n(\hat\bta_n^{(k)})+N_k,\bta-\hat\bta_n^{(k)}\rangle + \frac{1}{2\eta}||\bta-\hat\bta_n^{(k)}||_2^2.
\]
While the convergence properties of the non-private sequence $\bta_n^{(k)}$ are typically obtained using the convexity and smoothness properties of the loss function, we, on the other hand, leverage the properties of Hellinger distance and convergence of kernel densities to establish \emph{ASLSC} of $L_n(\bta)$ as in Proposition \ref{prop:convex}. Additionally,  we establish convergence rates (Theorem~\ref{thm:sensitivity}), which are required to establish the efficiency of the estimators.
Thus, perhaps more importantly, the private estimator obtained via our private optimization algorithms is not only efficient but also satisfies the privacy levels under some practical conditions described in section \ref{sec:efficiency_PMHDE} below.

To obtain $\hat\bta_{n}^{(K)}$ to be $\ep-$HDP, we start from $\hat\bta_{n}^{(0)}$. For all $k\geq 1$, we design a mechanism $M_k(\cdot,\cdot)$ to obtain $\hat\bta_{n}^{(k)}$ from $\hat\bta_{n}^{(k-1)}$ (treated as a plug-in constant vector), which is $\ep'-$HDP. Then using Corollary~\ref{cor:composition}, it will follow that the $K-$composition mechanism $M^{(K)}$ applied to the starting point $\hat\bta_{n}^{(0)}$ and the dataset $D$ to obtain $\hat\bta_{n}^{(K)}$ satisfies $h_K(\ep')-$HDP. If $h_K(\ep')\leq \ep$, then $\hat\bta_{n}^{(K)}$ will satisfy $\ep-$HDP. 
Turning to the mechanism $M_k(\cdot, \cdot)$, let
$M_k(w,D) = w(D) - \eta\cdot M_k'(w',D)$, where $M_k'(w',D) = w'(D) +\sigma\cdot Z_{k}$ and $w'(D) = \nabla L_n(w(D))$,
and $Z_k$ is the perturbing random vector which are i.i.d. for $k=1,2, \cdots K$. We design $M_k'(w,D)$ to be $\ep'-$HDP and hence using the post-processing property of the mechanism, $M_k(w,D)$ satisfies $\ep'-$HDP. Finally, using the 
Corollary~\ref{cor:composition} we conclude that $M^{(K)}(w,D)$ satisfies $\ep-$HDP. 

We next describe the mechanism $M_k'(w,D)$. By the previous description, it is an additive mechanism and we take $Z_k=[Z_{k,1},\cdots,Z_{k,m}]^T \sim N(\bm{0},\mathbf{I})$. Next to determine $\sigma$ we use Proposition~\ref{prop:HDP_mech} to obtain
\begin{align*}
    \sigma_{n,\ep'}= \Delta_n\sqrt{\frac{1}{-8\log(1-0.5\epsilon')}} \coloneqq \Delta_n c_{\ep'},
\end{align*}
where $\Delta_n$ is the upper-bound of the $L_2$ sensitivity of the query function $w'$ on dataset $D$ which is $\nabla L_n(\bta)$. Note that $\nabla L_n(\bta)$ is a function of the dataset $D$ for fixed $\bta=\hat\bta_{n}^{(k-1)}$.  Our next proposition describes a \emph{weak upper bound} on the $L_r$ sensitivity of the $\nabla L_n(\bta)$ and $H_n(\bta)$.

\begin{proposition} \label{prop:weak_sensitivity}
Suppose that assumptions {\bf(A1)-(A8)}
in Section 7
and assumptions {\bf\ref{asp:U1}}-{\bf\ref{asp:U2}} hold. 
Then for $r=1,2$,
\begin{align}\label{eq:sensitivity_1}
\Delta_{L_r} [\nabla L_n(\bta)]= O(n^{-\frac{1}{2}}), ~~
\Delta_{L_r}[H_n(\bta)]= O(n^{-\frac{1}{2}}), \quad \text {as} ~ n \ra \ff.
\end{align}
\end{proposition}
\noindent Behavior of the sensitivity of the gradient of the loss function, $\nabla L_n(\bta)$, is essential to study the convergence rate and the asymptotic efficiency of the PMHDE. As explained before, sensitivity is defined on a pair of adjacent datasets with an unbounded range. In the HD setting, the sensitivity appears through the integrals of kernel densities of adjacent datasets, yielding a weak upper bound. The disadvantage of this weak-upper bound is that it does not yield asymptotic normality of the private estimator. Under additional privacy constraints, we provide in 
Theorem~\ref{thm:sensitivity} below,  \emph{a sharper upper bound} for the sensitivity. We first turn to the algorithm for private gradient descent.
\begin{definition} [Private gradient descent (PGD)]
~
\begin{enumerate}
    \item {\bf via Gaussian noise:}
    \begin{align} \label{eq:PGD}
    \hat\bta_{n}^{(k+1)} = \hat\bta_{n}^{(k)} - \eta \left(\nabla L_n(\hat\bta_{n}^{(k)})+N_{n,k}\right),
\end{align}
$N_{n,k}=\Delta_n c_{\ep'}\mathbf{Z}_k$ where $\Delta_n$ is a valid $L_2$ sensitivity bound for $\nabla L_n(\bta)$:
global when calibrated uniformly over all datasets, or conditional on
the neighborhood of a dataset satisfying {\bf(U3)} when the sharp rate
of Theorem~\ref{thm:sensitivity} is used. $\ep'$ is the privacy level in each iteration. $\mathbf{Z}_k=[Z_{k,1},\cdots,Z_{k,m}]^T \sim N(\bm{0},\mathbf{I})$.

~

\item {\bf via Laplace noise:}
\begin{align}
    \hat\bta_{n}^{(k+1)} = \hat\bta_{n}^{(k)} \AV{-} \eta \left(\nabla L_n(\hat\bta_{n}^{k})+Y_k\right),
\end{align}
where $Y_k\sim Lap(0,b)$. If the parameter space is one dimension, $b$ is obtained by solving
\begin{align*}
    -2\left[e^{-\frac{\Delta_n^{(1)}}{2b}}+\frac{\Delta_n^{(1)}}{2b}e^{-\frac{\Delta_n^{(1)}}{2b}}-1\right]=\epsilon'.
\end{align*}
If the parameter space is $m-$dimension, $Y_k=(Y_{1,k},\cdots,Y_{m,k})$, $Y_{i,k}\sim Lap(0,b)$, $b=\frac{\Delta_n^{(1)}}{2\log(\frac{1}{1-0.5\epsilon'})}$, where $\Delta_n^{(1)}$ is the $L_1$ sensitivity of $\nabla L_n(\bta)$. $\ep'$ is the privacy level in each iteration. 
In Proposition~\ref{prop:ep_choice} below, we describe a method to choose $\ep'$ for both the mechanisms.
\end{enumerate}
\end{definition}

We summarize the iterations as an algorithm in the Gaussian case.

\begin{algorithm}[H]
\caption{Private gradient descent (PGD)}
\begin{algorithmic}
\STATE \textbf{Input:} MHDE loss function $L_n(\bta)$, number of iteration $K$, learning rate $\eta$, MHDE privacy level $\epsilon$, each iteration privacy level $\epsilon'$ from Proposition~\ref{prop:ep_choice}, initial point $\hat\bta_n^{(0)}$.
\STATE \textbf{Output:} Private MHDE $\hat\bta_n^{(K)}$.
\STATE $k=1$.
\WHILE{$k\leq K$}
    \STATE Generate $Z_k$ from $N(0,I)$, with same dimension of $\hat\bta_n^{(0)}$.
    \STATE Calculate $\Delta_n$, the $L_2$ sensitivity of $\nabla L_n(\bta)$.
    \STATE Update $\hat\bta_n^{(k)}$ by $\hat\bta_n^{(k)} = \hat\bta_n^{(k-1)} - \eta \left(\nabla L_n(\hat\bta_n^{(k-1)})+\Delta_n\cdot c_{\ep'}\cdot Z_{k}\right)$.
    \STATE $k \gets k+1$.
\ENDWHILE
\RETURN $\hat\bta_n^{(K)}$.
\end{algorithmic}
\end{algorithm}
We will show that for a fixed $\ep$ and $K$ (depending on $n$ and pre-determined), the algorithm returns PMHDE, which satisfies $\ep-$HDP.
While the GD is useful, its convergence rate can be arbitrarily slow. Hence, frequently, in applications,  the NR method is used to obtain MHDE, which guarantees a quadratic convergence rate. For this reason, we now describe the private NR algorithm. We recall that the standard NR algorithm follows the iteration
\begin{align*}
    \hat\bta_{n}^{(k+1)} = \hat\bta_{n}^{(k)} -H_n^{-1}(\hat\bta_{n}^{(k)}) \nabla L_n(\hat\bta_{n}^{(k)}).
\end{align*}
Next, to obtain the private versions of the Hessian and Hellinger score, we use Corollary~\ref{cor:composition}, Proposition~\ref{prop:HDP_mech}, and Proposition~\ref{prop:symmetric_HDP} to determine the appropriate noise in the additive mechanism.
\begin{definition}[Private Newton-Raphson (PNR) via Gaussian noise]
    The private Newton-Raphson iterates are
\begin{align}\label{eq:PNR_iteration}
    \hat\bta_{n}^{(k+1)} = \hat\bta_{n}^{(k)} -\eta\left(H_n(\hat\bta_{n}^{(k)})+W_{n,k}\right)^{-1} \left(\nabla L_n(\hat\bta_{n}^{(k)})+N_{n,k}\right),
\end{align}
where $W_{n,k}\in\mathbb{R}^{m\times m}$ and $N_{n,k}\in\mathbb{R}^{m\times1}$ are the noise added to satisfy $\ep-$HDP. That is,
\begin{align*}
    N_{n,k} = \Delta_n\cdot c_{\ep'/2}\cdot Z_{k} ~~\text{and}~~
    W_{n,k} = \Delta_n^{(H)}\cdot c_{\ep'/2}\cdot\Tilde Z_{k},
\end{align*}
where $\Delta_n$ and  $\Delta_n^{(H)}$ are appropriate estimates of $L_2$ sensitivities of $\nabla L_n(\bta)$ and $H_n(\bta)$. Also, $Z_{k}$ is $m-$dimensional Gaussian vector with independent standard normal components; $\Tilde Z_{k}$ is an $m\times m$   symmetric matrix with i.i.d. standard Gaussian upper-triangular entries, including the diagonal. $\ep'$ is the privacy level in each iteration. 
\end{definition}

We summarize the iterations as an algorithm for the Gaussian case.

\begin{algorithm}[H]
\caption{Private Newton-Raphson (PNR)}
\begin{algorithmic}
\STATE \textbf{Input:} MHDE loss function $L_n(\bta)$, number of iteration $K$, learning rate $\eta$, MHDE privacy level $\epsilon$, each iteration privacy level $\epsilon'$ from Proposition~\ref{prop:ep_choice}, initial point $\hat\bta_n^{(0)}$.
\STATE \textbf{Output:} Private MHDE $\hat\bta_n^{(K)}$.
\STATE $k=1$.
\WHILE{$k\leq K$}
    \STATE Generate $Z_k$ from $N(0,I)$, with same dimension of $\hat\bta_n^{(0)}$.
    \STATE Calculate Hessian matrix of $L_n(\bta)$ at $\hat\bta_n^{(k-1)}$: $H_n(\hat\bta_n^{(k-1)})$.
    \STATE Generate $\tilde Z_k$ from $N(0,I)$, with same dimension of $H_n(\hat\bta_n^{(k-1)})$.
    \STATE Calculate $\Delta_n$, the $L_2$ sensitivity of $\nabla L_n(\bta)$.
    \STATE Calculate $\Delta_n^{(H)}$, the $L_2$ sensitivity of $H_n(\bta)$.
    \STATE Calculate $N_{n,k} = \Delta_n\cdot c_{\ep'/2}\cdot Z_{k}$.
    \STATE Calculate $W_{n,k} = \Delta_n^{(H)}\cdot c_{\ep'/2}\cdot \Tilde Z_{k}$.
    \STATE Update $\hat\bta_n^{(k)}$ by $\hat\bta_n^{(k)} = \hat\bta_n^{(k-1)} -\eta\cdot\left(H_n(\hat\bta_n^{(k-1)})+W_{n,k}\right)^{-1} \left(\nabla L_n(\hat\bta_n^{(k-1)})+N_{n,k}\right)$.
  \STATE $k \gets k+1$.
\ENDWHILE
\RETURN $\hat\bta_n^{(K)}$.
\end{algorithmic}
\end{algorithm}

\noindent It is possible to use the additive Laplace mechanism to obtain the PMHDE, where one replaces $N_{n,k}$ and $W_{n,k}$ by a Laplace distribution with variance as described in PGD, and replaces $\ep'$ with $\ep'/2$.

\begin{proposition} \label{prop:ep_choice}
Let the iteration number $K$ and the privacy budget $\ep$ be given. Let $\ep'$ denote the privacy level at every iteration of the PGD and PNR algorithms and let $\hat\bta_n^{(K)}$ denote the PMHDE. If $\ep'$ satisfies $\ep=h_K(\ep')$, then $\hat\bta_n^{(K)}$ satisfies $\ep=h_K(\ep')-$HDP. In particular, if $\ep'=\frac{\ep}{K}$ then $\ep[1-\ep(K-1)(4K)^{-1}] \leq h_K(\ep') \le \ep$.
\end{proposition}

We end this subsection with a brief discussion about similar algorithms studied in the literature, namely for the M-estimators. First, unlike the M-estimators, the difficulty in our problem is that the loss function  $L_n(\bta)$ is usually not convex in $\bta$. We address this issue by leveraging the properties at the optimal point,  which is required for statistical analysis of the estimator in non-private settings. Next, the gradient $\nabla_{\bta} \log f_{\bta}(x)$ (score function)
is not always a bounded function of $x$. This is an important issue since the sensitivity of the estimator depends on the score function, and in the M-estimator case, they are assumed to be bounded. However, this assumption leads to a loss of statistical efficiency.

\subsection{Utility of PMHDE} \label{sec:utility_PMHDE}

In this section, we describe the utility properties of the PMHDE. In addition to the assumptions in Section 7
we need additional weak family regularity conditions to study the convergence properties of the algorithms.
We emphasize here that our proof method also yields the convergence of the non-private GD and NR algorithms, which have not been studied in the literature before. In summary, establishing the properties of private algorithms only requires weak family regularity conditions. 

As explained above,  the weak upper bound does not yield asymptotic normality of the private estimator. However, in practice, extreme points in a dataset are not typically revealed to safeguard privacy. 
Under this consideration, we assume a range of the data that increases as the sample size increases.  Specifically, we consider the kernel density estimator 
\begin{eqnarray} \label{eq:modified_gn}
\bar{g}_n(x)= \frac{1}{n\cdot c_n}\sum_{i=1}^{n}K\left(\frac{x-X_i}{c_n}\right)\bm{1}_{(X_i\in B_n)},
\end{eqnarray}
where $b_n \nearrow \ff$ and $B_n=(-b_n, b_n)$. In the rest of the paper, the loss function $L_n(\bta)$ is based on $\bar{g}_n(\cdot)$. With this choice of the query function, we derive a sharper upper bound for the sensitivity, which plays an essential role in the proof of the asymptotic normality of the PMHDE. We need an additional regularity condition on the postulated family of densities.  We recall that $c_n$ is the bandwidth associated with $\bar g_n(\cdot)$. 
\begin{assumptionU}\label{asp:U3}
Let $p\in(1,2)$ and satisfy $\frac{1}{p}+\frac{1}{q}=1$. Suppose $K\in L_p(\mathbb R)$ and write
$S_n=(-b_n-\beta c_n,b_n+\beta c_n)$, which contains the support of $\bar g_n$. Choose the compact interior region $A_n=[-a_n,a_n]\subset S_n$ by the
weighted-tail construction in the proof of Theorem~\ref{thm:sensitivity},
with $b_n$ large enough to contain this region, and let
$\delta_n=\inf_{x\in A_n}\bar g_n(x)$. For the observed dataset, assume
that $\delta_n>0$, $\delta_n\to0$, and
$(nc_n)^{-(1-1/p)}\leq C_\delta\sqrt{\delta_n}$ for some $C_\delta<\infty$. Additionally, assume that 
\begin{eqnarray*}
\mathbf{E}_{\bta}\left[ \Vert\mathbf{u}_{\bta}(X)f_{\bta}^{\frac{1}{2}-\frac{1}{q}}(X) \Vert_1^q \right],~\mathbf{E}_{\bta}\left[\Vert\mathbf{u}_{\bta}(X) \mathbf{u}^T_{\bta}(X) f_{\bta}^{\frac{1}{2}-\frac{1}{q}}(X)\Vert_1^q \right], \text{and } \mathbf{E}_{\bta}\left[\Vert\mathbf{\dot u}_{\bta}(X) f_{\bta}^{\frac{1}{2}-\frac{1}{q}}(X)\Vert_1^q \right] 
\end{eqnarray*}
are all finite and continuous in $\bta$ for all $\bta \in \Theta$.
\end{assumptionU}

For fixed $D$, write $\Delta_{L_r}(h;D)=
\sup_{D':\,\|D-D'\|_H=1}\|h(D)-h(D')\|_r$, $r=1,2$.
Under {\bf(U3)}, the following bound is conditional on $D$ and uniform
in its one-record neighbors; it also holds in Frobenius norm for the Hessian.
\begin{theorem}[Sensitivity for MHDE] \label{thm:sensitivity}
Suppose that assumptions of Section 7
and assumptions {\bf\ref{asp:U1}}-{\bf\ref{asp:U3}} hold. 
Then for $r=1,2$ and $p\in(1,2)$,
\begin{align} \label{eq:sensitivity_2}
\Delta_{L_r} [\nabla L_n(\bta);D]= O(n^{-\frac{1}{p}}), ~~
\Delta_{L_r}[H_n(\bta);D]= O(n^{-\frac{1}{p}}), \quad \text{as} ~ n \ra \ff.
\end{align}
The constants in the above expressions depend on $f_{\bta}(\cdot)$ and the dimension $m$.
\end{theorem}
This is not an unrestricted global-sensitivity assertion for unbounded
data. The supplementary local-mass and covering conditions give an
event $\mathcal E_n$ with
$\mathbf P_g(\mathcal E_n^c)\le M_n\exp\{-c(nc_n)^{2/p-1}\}$;
when $\log M_n=o((nc_n)^{2/p-1})$, its probability tends to one,
as in random differential privacy \citep{HallRinaldoWasserman2013}.
The sharp PGD/PNR privacy claims below have this conditional meaning;
unrestricted HDP requires global calibration.

Our next result is concerned with the utility of the method, measured using the $L_2$ distance between the private and non-private estimators.
\begin{theorem}[Utility of PGD via Gaussian noise] \label{thm:utility_GD}
Let assumptions {\bf(A1)-(A8)}
in Section 7
and assumptions {\bf\ref{asp:U1}}-{\bf\ref{asp:U3}} hold. Then, the PMHDE, $\hat{\bta}_n^{(K_n)}$, obtained via the PGD algorithm satisfies $\ep-$HDP. Furthermore,
there exists a strictly positive learning rate $\eta$, initial value $\hat{\bta}_n^{(0)}$, $N$ such that for all $n>N$, there exist $p\in(1,2]$ and $K_n$ satisfying $K_n\geq c_1\log n$ for some $c_1\in(0,\infty)$
\begin{align} \label{eq:GD_utility}
||\hat\bta_n^{(K_n)}-\hat\bta_n||_2\leq c_2 n^{-\frac{1}{p}}({K_n\log(K_n/\xi)})^{\frac{1}{2}}
\end{align}
with probability at least $1-\xi$, where $c_2\in(0,\infty)$ is a constant depending on $f_{\bta}$ and $m$. That is, $||\hat\bta_n^{(K_n)}-\hat\bta_n||_2=O_p\left(n^{-\frac{1}{p}}({K_n\log(K_n)})^{\frac{1}{2}}\right)$.
\end{theorem}

\begin{remark}
~
\begin{enumerate}
\item The calculations show that the upper bound in the above theorem is $||\hat\bta^{(K_n)}_n-\hat\bta_n||_2\leq C r_{noi}$, where $$r_{noi}=\Delta_n\cdot \frac{4m^{\frac{1}{2}}+2(2\log\frac{K_n}{\xi})^{\frac{1}{2}}}{\left(-\log(1-\frac{\epsilon}{2K_n})\right)^{\frac{1}{2}}}.$$ 
The dominant term in the numerator is $\Delta_n(\log K_n/\xi)^{\frac{1}{2}}$ while the denominator is asymptotic to $K_n^{-\frac{1}{2}}$; that is, $\left(-\log(1-\frac{\epsilon}{2K_n})\right)^{\frac{1}{2}}\sim K_n^{-\frac{1}{2}}$. This yields the approximate upper bound in the theorem. 
\item As explained previously, the iteration number $K$ and privacy level $\epsilon$ are predetermined. In practice, $K$ is chosen based on the sample size $n$, $K_n \sim \log n$. 

\item Using the weak upper bound for sensitivity in Proposition~\ref{prop:weak_sensitivity}, namely, $\Delta_n\sim  n^{-\frac{1}{2}}$, it follows that $r_{noi}\to 0$ as $n\to\infty$. However, with this choice, one obtains consistency and not the limit distribution.
However, invoking the additional assumption {\bf\ref{asp:U3}}, one can use the sharper upper bound, namely $\Delta_n \sim n^{-\frac{1}{p}}$ for $ p \in (1, 2)$. While this choice continues to yield consistency, it also yields the asymptotic distribution that coincides with the asymptotic distribution of the non-private MHDE, as we shall see in Section \ref{sec:efficiency_PMHDE} below.  
\item For both private and non-private versions of the algorithm, the values of $\eta$, $c_2$, and global optimization property in the Theorem depend on the asymptotic properties of  $L_n(\bta)$ in a neighborhood around $\bta_g$. Also, the proof of the utility of the PMHDE relies on the above-mentioned properties of the loss function.
\item For the privacy guarantees, the initial value is chosen independently of the confidential sample.

\end{enumerate}
\end{remark}
In the Theorem \ref{thm:utility_GD}, $\ep$ can be replaced by $\ep_n\to0$ as $n\to\ff$. The next corollary describes this case.

\begin{corollary}
    Under the same conditions as in Theorem \ref{thm:utility_GD}, for $p\in(1,2)$, $\ep_n\to0$, and $\ep_n\geq c_0n^{-\beta}$ for some $c_0>0$, where $\beta\in(0,\frac{2}{p}-1)$, as $n\to\ff$
    \begin{align*}
        ||\hat\bta_n^{(K_n)}-\hat\bta_n||_2\leq c_3 n^{-\frac{1}{p}+\frac{\beta}{2}}({K_n\log(K_n/\xi)})^{\frac{1}{2}} \to 0,
    \end{align*}
    where $c_3\in(0,\ff)$ depends on $f_{\bm \theta}$, $m$, and $c_0$.
\end{corollary}

\begin{theorem}[Utility of PNR via Gaussian noise] \label{thm:utility_NT}
    Let assumptions {\bf(A1)-(A8)}
    in Section 7
    and assumptions {\bf\ref{asp:U1}}-{\bf\ref{asp:U3}} hold. Then the PMHDE, $\hat\bta_n^{(K_n)}$, obtained via the PNR algorithm satisfies $\ep-$HDP. Furthermore, there exists a learning rate $\eta>0$, initial value $\hat\bta_n^{(0)}$, and $N$ such that for all $n >N$, there exist $p \in (1, 2]$ and $K_n$ satisfying $K_n\geq c_1 \log(\log n)$ for some $c_1\in(0,\infty)$ 
    \begin{align} \label{eq:NR_utility}
        ||\hat\bta_n^{(K_n)}-\hat\bta_n||_2\leq c_2 n^{-\frac{1}{p}}(K_n\log (K_n/\xi))^{\frac{1}{2}}
    \end{align}
    with probability at least $1-\xi$, where $c_2\in(0,\infty)$ is a constant depending on $f_{\bta}$ and $m$. That is, $||\hat\bta_n^{(K_n)}-\hat\bta_n||_2=O_p\left(n^{-\frac{1}{p}}(K_n\log (K_n))^{\frac{1}{2}}\right)$.
\end{theorem}

\begin{remark}
~
\begin{enumerate}
\item The calculations show that the upper bound in the above theorem is $||\hat\bta_n^{(K_n)}-\hat\bta_n||_2\leq C\cdot r_{noi}$, where 
\begin{align*}
    r_{noi}\sim& \Delta_n^{(H)}\cdot\frac{\left(2m\log\frac{4K_nm}{\xi}\right)^{\frac{1}{2}}}{\left(-8\log(1-\frac{\epsilon}{4K_n})\right)^{\frac{1}{2}}}.
\end{align*}
The specific expression is shown in the proof.
The dominant term is $\Delta_n^{(H)}\cdot\frac{\left(2m\log\frac{4K_nm}{\xi}\right)^{\frac{1}{2}}}{\left(-8\log(1-\frac{\epsilon}{4K_n})\right)^{\frac{1}{2}}}$. The denominator $\left(-8\log(1-\frac{\epsilon}{4K_n})\right)^{\frac{1}{2}}\sim K_n^{-\frac{1}{2}}$ and the numerator $\left(2m\log\frac{4K_nm}{\xi}\right)^{\frac{1}{2}}\sim (\log K_n/\xi)^{\frac{1}{2}}$. This yields the approximate upper bound in the theorem. 
\item Arguing as in the PGD algorithm, $K$ is chosen as $K_n\sim \log(\log n)$ for reducing computational complexity and obtaining consistency of the estimator. We notice here that with fewer iterations, compared to the PGD algorithm, one obtains the PMHDE with $\ep-$HDP guarantees. 
\item Using the same arguments as in item (3) of the remark following Theorem~\ref{thm:utility_GD}, the asymptotic properties are now determined by $\Delta_n^{(H)}$. Specifically, using Theorem~\ref{thm:sensitivity}, the choice of $\Delta_n^{(H)}\sim n^{-\frac{1}{2}}$ leads to consistency alone, and under additional Assumption~{\bf\ref{asp:U3}}, one also obtains the asymptotic distribution by using $\Delta_n^{(H)}\sim n^{-\frac{1}{p}}$ for $p\in(1,2)$.
\end{enumerate}
\end{remark}

As before, in Theorem \ref{thm:utility_NT}, fixed $\ep >0$ can be replaced by $\ep_n\to0$ as $n\to\ff$. Next corollary states the case.

\begin{corollary}
    Under same conditions as in Theorem \ref{thm:utility_NT}, for $p\in(1,2)$, $\ep_n\to0$, and $\ep_n\geq c_0n^{-\beta}$ for some $c_0>0$, where $\beta\in(0,\frac{2}{p}-1)$, as $n\to\ff$
    \begin{align*}
        ||\hat\bta_n^{(K_n)}-\hat\bta_n||_2\leq c_3 n^{-\frac{1}{p}+\frac{\beta}{2}}({K_n\log(K_n/\xi)})^{\frac{1}{2}} \to 0,
    \end{align*}
    where $c_3\in(0,\ff)$ depends on $f_{\bm \theta}$, $m$, and $c_0$.
\end{corollary}

\subsection{Efficiency of PMHDEs} \label{sec:efficiency_PMHDE}

We now discuss the statistical properties of the PMHDE. Noting that our loss function is obtained using $\bar{g}_n(\cdot)$ and its almost sure $L_1$ convergence (using generalized dominated convergence Theorem) to $g(\cdot)$, we apply PGD and PNR algorithms for obtaining PMHDE. We recall that $K$ is the number of pre-determined iterations of the gradient descent or Newton-Raphson algorithm. In the Theorem below, we use $K_n$ for $K$ to emphasize its dependence on $n$. We note here that the efficiency proof will rely on the sharper bound in Theorem \ref{thm:sensitivity}. Before we state the Theorem, we recall that $H_{\ff}(\bta)=\lim\limits_{n\to\infty}H_n(\bta)$, where $H_n(\cdot)$ is the Hessian matrix. That is,
\begin{align*}
    H_{\ff}(\bta)= -\int g^{1/2}(x)f_{\bta}^{1/2}(x) [\mathbf{u}_{\bta}(x)\mathbf{u}^T_{\bta}(x)+2\mathbf{\dot u}_{\bta}(x)]\mathrm{d}x.
\end{align*}
Also, set $\Sigma_g= 4^{-1} \int_{\Real} \rho_{\bta_g}(x)\rho_{\bta_g}^T(x)\mathrm{d}x$, where $\rho_{\bta}(x)=4 H^{-1}_{\ff}(\bta)\nabla f^{\frac{1}{2}}_{\bta}(x)$.

\begin{theorem} \label{thm:asym}
Let assumptions {\bf(A1)-(A8)}
in Section 7
and assumptions {\bf\ref{asp:U1}}-{\bf\ref{asp:U3}} hold. Let $ K_n > C \log(n)$ for the gradient-descent algorithm and $K_n \ge C \log(\log n)$ for the Newton-Raphson algorithm, where $0 < C < \ff$. Let $\hat{\bta}^{(K_n)}_n$ denote the private Hellinger distance estimator of $\bta_g$ evaluated using one of gradient-descent or Newton-Raphson algorithms.  Then the following hold: 
\begin{enumerate}
\item $\lim_{n \ra \ff} \sqrt{n}||\hat\bta^{(K_n)}_{n}-\hat\bta_n||_2 =0$, in probability.
\item $\lim_{n \ra \ff} ||\hat\bta^{(K_n)}_{n}-\bta_{g}||_2 =0$, in probability.
\item $\sqrt{n}(\hat\bta^{(K_n)}_{n}-\bta_g)\overset{d}{\to} N(0, \Sigma_g)$, as $n \ra \ff$. Furthermore, if $g=f_{\bta_0}$,  then $\Sigma_g=I^{-1}(\bta_0)$.
\end{enumerate}
\end{theorem}
It is worth emphasizing here that non-private estimator $\hat\bta_n$ obtained by minimizing $L_n(\cdot)$ (derived using $\bar{g}_n$) is asymptotically normal with mean vector $\bm{0}$ and covariance matrix $\Sigma_g$. That is, PMHDE and MHDE have the same asymptotic distribution implying that PMHDE is efficient. The above theorem also holds if $\ep$ is replaced by $\ep_n\to0$ satisfying $\ep_n\geq c_0n^{-\beta}$ for some $c_0>0$, where $\beta\in(0,\frac{2}{p}-1)$ and $p\in(1,2)$ is the constant in Assumption {\bf\ref{asp:U3}}.

\subsection{Private confidence interval} \label{sec:CI}
From Theorem \ref{thm:asym} above, one obtains that $n^{\frac{1}{2}}(\hat{\bta}_n -\bta_g)$ converges in distribution to a Gaussian distribution with mean vector $\mathbf{0}$ and covariance matrix $\Sigma_g$ and when the model is correctly specified, $\Sigma_g=I^{-1}(\bta_0)$. To construct the confidence interval, we use the private estimator $\hat\bta_n^{(K_n)}$, a private version of $H_n(\bta)$, and the model Fisher information $I(\hat\bta_n^{(K_n)})$. The last quantity is computed from the private estimator by post-processing and requires no additional privacy budget.
The Hessian is privatized using the symmetric matrix mechanism in Proposition \ref{prop:symmetric_HDP}. By composition and post-processing (see \cite{Wang2018}), the previously stated $3\ep$-HDP bound for the confidence interval remains valid, though conservative. 

To derive the private version of $\Sigma_g$, it is convenient to use an alternative expression frequently referred to as the sandwich formula in the literature. Towards this derivation, recalling the loss function and using the first-order Taylor approximation of the gradient, and $\nabla L_n(\hat\bta_n) = \bm{0}$, we obtain
\begin{align*}
\sqrt{n}(\hat\bta_n-\bta_g) =& -\sqrt{n} H_n^{-1}(\bta^*_n)\cdot \nabla L_n(\bta_g),  
\end{align*}
where $H_n(\cdot)$, as before, is the Hessian of $L_n(\bta)$, and $\bta^*_n=\alpha\hat\bta_n+(1-\alpha)\bta_g\in \Theta$ for some $\alpha\in[0,1]$.
Now, under assumptions {\bf(A1)-(A8)},
it follows that $\sqrt{n}\nabla L_n(\bta_g)\stackrel{d}{\to} N(0,I(\bta_g))$ as $n\to\infty$ (similar to \cite{Cheng2006} Lemma 4.4 and Lemma 4.5). Notice that $I(\bta_g)$ can be expressed as
$\mathbf{E}_{\bta_g}[\mathbf{u}_{\bta_g}(X)\mathbf{u}^T_{\bta_g}(X)])$.
Hence, using the almost sure  convergence  of $H_n(\bta^*_n)$ to  $H_{\ff}(\bta_g)$, the limiting covariance matrix of 
$\sqrt{n}(\hat\bta_n-\bta_g)$ is
$H_\infty^{-1}(\bta_g)\cdot I(\bta_g) \cdot H_\infty^{-1}(\bta_g)$. Evaluating the model Fisher information at the private estimator gives
\begin{align*}
\mathbf{\hat V}\coloneqq \hat{\mathbf{V}}[\sqrt{n}(\hat\bta_n^{(K_n)}-\bta_g)]= H_n^{-1}(\hat\bta_n^{(K_n)})\cdot I(\hat\bta_n^{(K_n)})\cdot H_n^{-1}(\hat\bta_n^{(K_n)}),
\end{align*}
where $I(\bta)=\int_{\Real} f_{\bta}(x)\mathbf u_{\bta}(x)\mathbf u_{\bta}^T(x)\mathrm{d}x$.
This is commonly referred to as the sandwich formula. Now, using the symmetric matrix mechanism, as explained above, we obtain the private version of $H_n(\cdot)$ as $H_n(\cdot)+ \Delta_n^{(H)}\cdot c_{\ep}\cdot\tilde Z$, where $\tilde Z$ is $m\times m$ random upper triangular symmetric matrix including diagonals, whose components are i.i.d. standard Gaussian. The private version of $\nabla L_n(\cdot)$ is given by $\nabla L_n(\cdot)+ \Delta_n\cdot c_{\ep}\cdot \mathbf{Z}$, where $\mathbf{Z}=[Z_1,\cdots,Z_m]\sim N(\bm{0},\mathbf{I})$. The $1-\alpha$ confidence interval for the $j^{th}$ element of $\bta_g$ is given by $\hat\theta_{n,j}^{(K_n)}\pm z_{1-\alpha/2}\cdot\left(\frac{\mathbf{\hat V}_{jj}}{n}\right)^{\frac{1}{2}}$, where $\mathbf{\hat V}_{jj}$ is the $(j,j)^{th}$ component of $\mathbf{\hat V}$ and $\hat\theta_{n,j}^{(K_n)}$ is the $j^{th}$ element of $\hat\bta_{n}^{(K_n)}$.

Since the plug-in method for the construction of the confidence interval does not take into account the perturbation in the last step, a correction is required. Using the perturbation random variables introduced for PGD and PNR algorithms, the corrected confidence interval for $j^{th}$ element of $\bta_g$, for the PGD algorithm, is $\hat\theta_{n,j}^{(K_n)}\pm z_{1-\alpha/2}\cdot\left(\frac{\mathbf{\hat V}_{jj}}{n}+\eta^2(\Delta_n\cdot c_{\ep'})^2\right)^{\frac{1}{2}}$. Similarly, using the perturbation random variable, $\tilde N_{n,K_n}$, of the last iteration defined in Lemma 6.4 in Section 6, 
an approximation for the variance of 
$\tilde N_{n,K_n}$ is $C^{NR}=\eta^2 \left(H_n^{-1}(\hat\bta_n^{(K_n)})+ 
W_{n,K_n} \right)^{-1}\cdot (\Delta_n\cdot c_{\ep'/2})\cdot\left(H_n^{-1}(\hat\bta_n^{(K_n)})+ 
W_{n,K_n} \right)^{-1}$. The corrected confidence interval for $j^{th}$ element of $\bta_g$, for the PNR algorithm, is $\hat\theta_{n,j}^{(K_n)}\pm z_{1-\alpha/2}\cdot\left(\frac{\mathbf{\hat V}_{jj}}{n}+ C^{NR}_{j,j} \right)^{\frac{1}{2}}$, where $C^{NR}_{j,j}$ is the $(j,j)^{th}$ component of $C^{NR}$. Similar ideas were also considered in \cite{avella2023}, where the limiting variance results from the M-estimation explicitly involves the bound of the M-estimating score function.

\subsection{Robustness under contamination}
\label{sec:robustness_PMHDE}

We conclude this section by summarizing the robustness properties
preserved by the private estimators. The efficiency and approximation
results provide the connection with the contamination response of MHDE. Write $T(g)=\bta_g$ for the MHDE functional, and let
$g_{\alpha,z}=(1-\alpha)f_{\bta_0}+\alpha q_z$, where $q_z$ is a
contaminating density and $0<\alpha<1$. Under the conditions of
Theorem~\ref{thm:asym} at the relevant sampling laws, PMHDE inherits
the asymptotic linear representation of MHDE and hence its influence
function. For each fixed $\alpha,z$,
it also has the same limiting $\alpha$-influence function,
$\{T(g_{\alpha,z})-T(f_{\bta_0})\}/\alpha$, in the finite-contamination
formulation of \citet{Beran1977}.
Propositions 14.1 and 14.2
in Supplementary Section 14 give the statements and proofs
directly from the efficiency theorem.

The common limiting functional also gives a common asymptotic breakdown
point on the same admissible contamination family, provided the
consistency conditions of Supplementary Proposition 14.3 hold at each
contaminated law on the full open parameter space. In particular, the
value is $1/2$ under these conditions and the classical Hellinger
separation conditions in Supplementary Corollary 14.5. This statement
takes sample size to infinity first at each fixed contaminated law. When sample size and
contamination severity vary together, the PGD and PNR approximation
bounds yield a two-sided probability comparison with MHDE and the
breakdown-transfer result in
Supplementary Proposition 14.7. The appendix includes
these bounds, the local algorithmic error radii, and their consequences
for varying contamination fractions. Breakdown is considered on the
full model parameter space, rather than enforced by a fixed compact
constraint on the output.

Robustness concerns departures of the sampling distribution from the
model. HDP concerns changes in the law of the released output under a
neighboring-dataset replacement. The approximation and efficiency
results connect the private algorithms to the statistical robustness
of MHDE.

The following section examines the finite-sample behavior of these
procedures, including their response to gross-error contamination.

\section{Numerical experiments} \label{sec:numerical}
In this section, we present results from several numerical experiments. The experiments compare the outputs from both private and non-private optimization algorithms. We also study the coverage of the private confidence intervals.  We begin by describing the simulation design.

\textbf{Simulation design}:
The data are simulated from a normal distribution with a mean of five and a variance of four. The kernel density estimator $g_n(\cdot)$ is given by
\begin{align*}
        g_n(x) = \frac{1}{n\cdot c_n}\sum_{i=1}^{n}K\left(\frac{x-X_i}{c_n}\right),
\end{align*}
where $K(x) = \frac{3}{4} \cdot (1 - x^2)\cdot 1_{|x|\leq1}(x)$ is the Epanichikov kernel. The bandwidth $c_n$ is chosen using Silverman's bandwidth selection (\cite{Silverman1986}) and then fixed for all replications. The loss function is approximated using the Monte-Carlo approach; that is,
\begin{align*}
    L_n(\bta) =& 2\cdot \int_{\mathbb{R}}\left( f_{\bta}^{1/2}(x)-g_n^{1/2}(x) \right)^2 \mathrm{d}x
    \approx 2\left[ 2-2 \frac{1}{\AV{r_n}}\sum_{i=1}^{r_n} \sqrt{\frac{f_{\bta}(X_{n,i})}{g_n(X_{n,i})}} \right],
\end{align*}
where $r_n$ is the number of Monte Carlo samples and $\{X_{n,i}\cdots X_{n, r_n}\}|(X_1, \cdots X_n)\stackrel{i.i.d.}\sim g_n(\cdot).$ The algorithm in \cite{Cheng2006} is used to generate data from $g_n(\cdot)$. We apply Algorithm 1 (PGD) and Algorithm 2 (PNR) with start point $(1,1)$ and choose $K=50$ for the PGD algorithm and $K=5$ for the PNR algorithm to obtain the PMHDE.  One can also choose smaller values of $K$ that ensure convergence of the algorithm. The learning rate for both the PGD and PNR algorithms is taken to be $\eta=0.5$. The sharp sensitivities, $\Delta_n$ and $\Delta_n^{(H)}$, are derived with $p=1.7$. We also tried other choices of $p$, yielding similar results (data not presented). Using standard calculations, we approximate $\Delta_n$ and $\Delta_n^{(H)}$ by 
\begin{align*}
    \Delta_n = \frac{2\sqrt{6}}{\sigma}\cdot\left(\frac{1}{n}\right)^{1/p},\quad 
    \Delta_n^{(H)} = \frac{\sqrt{118}}{\sigma^2}\cdot\left(\frac{1}{n}\right)^{1/p},
\end{align*}
$\sigma$ is chosen as the private estimator from the previous iteration of the algorithm. The 95\% confidence interval is calculated for both the parameters as described in Section \ref{sec:CI}. All the simulation results are based on datasets with sample sizes varying from 50 to 1000 and 5000 replications (not all data are presented).

Table \ref{T1} contains results for PMHDE and MHDE for the mean $\mu$ (std. error) and standard deviation $\sigma$ (std. error) with sample size 1000, privacy levels $\ep= 0.2, 0.6, 2.0$, and bandwidth $c_n=0.448$.  $\ep=2$ corresponds to a non-private estimator, namely the MHDE. The confidence intervals are $2\ep-$HDP. As $\ep$ decreases (corresponding to increased privacy), we notice that the estimator smoothly deviates from the non-private estimator. Similarly, the standard error increases, implying that the perturbation is at work. The uncorrected coverage of the confidence interval deteriorates with increased privacy. Specifically, for $\ep=0.20$, even though the average estimate of the mean is $4.996$, the confidence interval fails to capture the true value $67.3\%$ of times. However, the finite sample correction, as outlined in Section \ref{sec:CI}, improves the coverage to $82.4\%$. A similar interpretation also holds for the PNR algorithm even though, in this case, the computational complexity is reduced due to a ten-fold reduction in the values of $K$. Results for the PNR algorithm are summarized in Table 11.1
in Section 11.

\begin{figure}[ht]
    \includegraphics[width=0.7\linewidth]{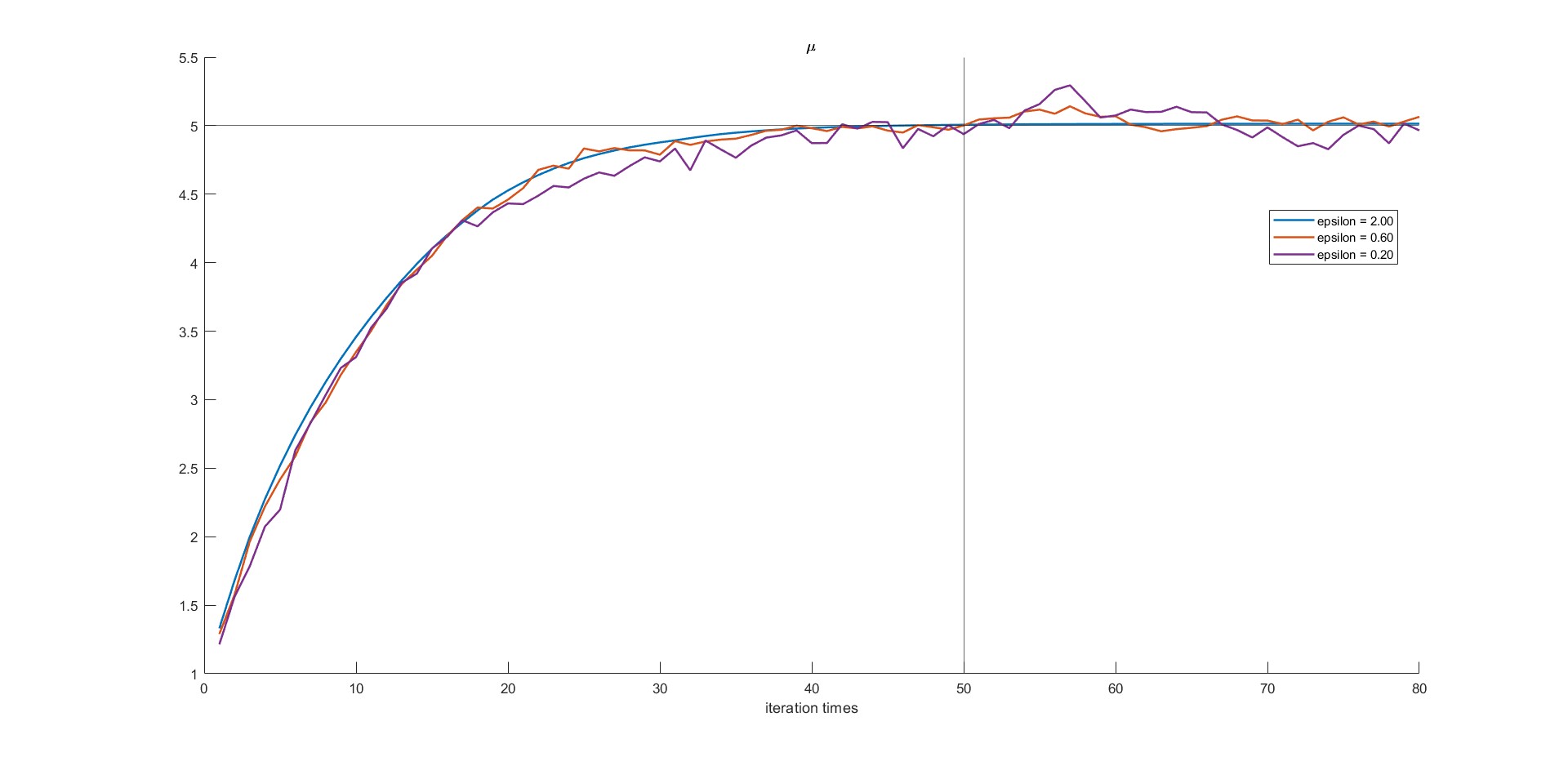}
    \includegraphics[width=0.7\linewidth]{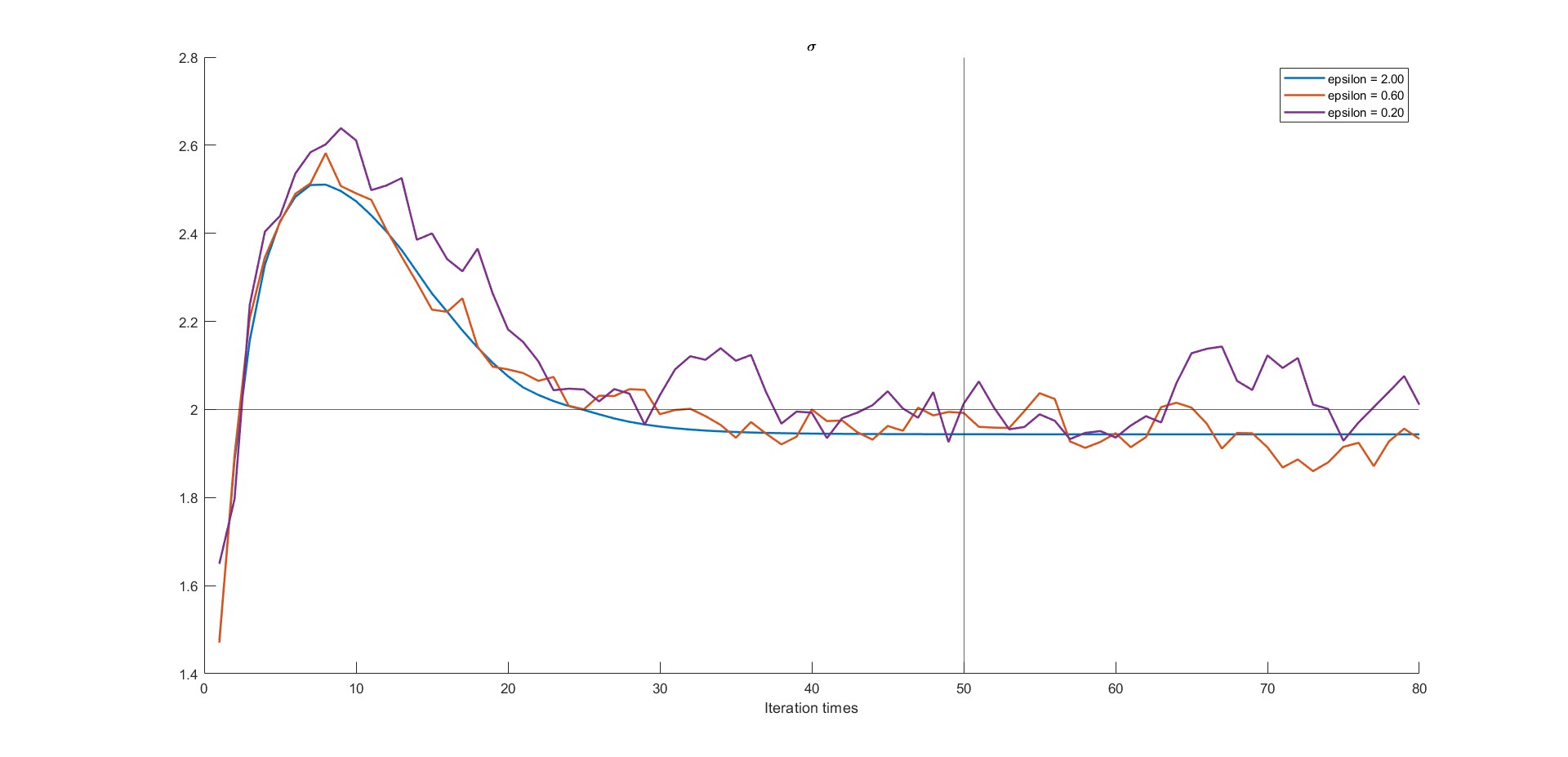}
    \caption{Private gradient descent path}
    \label{F1}
\end{figure}

\begin{figure}[ht]
    \includegraphics[width=0.7\linewidth]{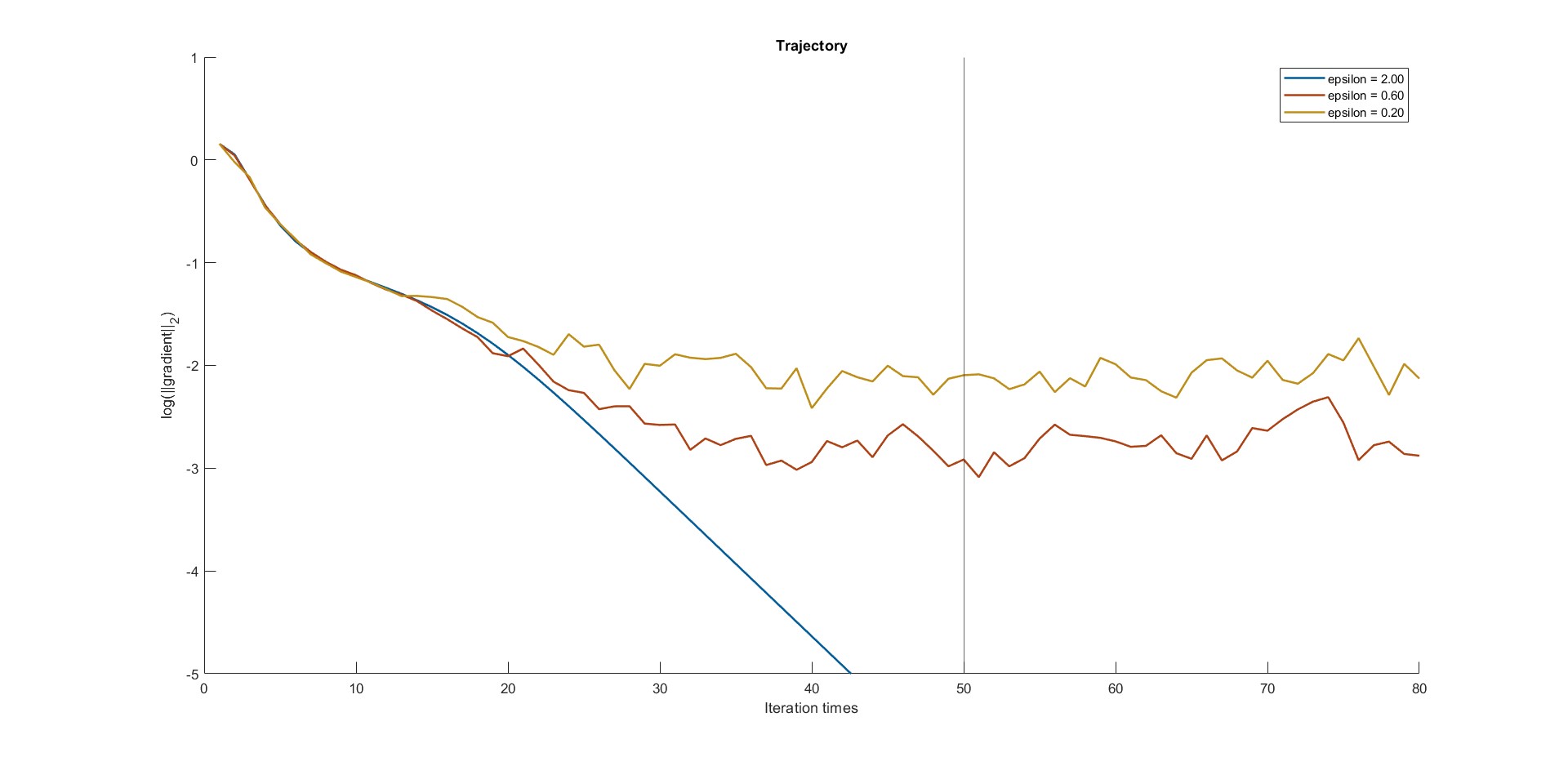}
    \caption{Private gradient descent trajectory}
    \label{F2}
\end{figure}

\begin{figure}[ht]
    \includegraphics[width=0.7\linewidth]{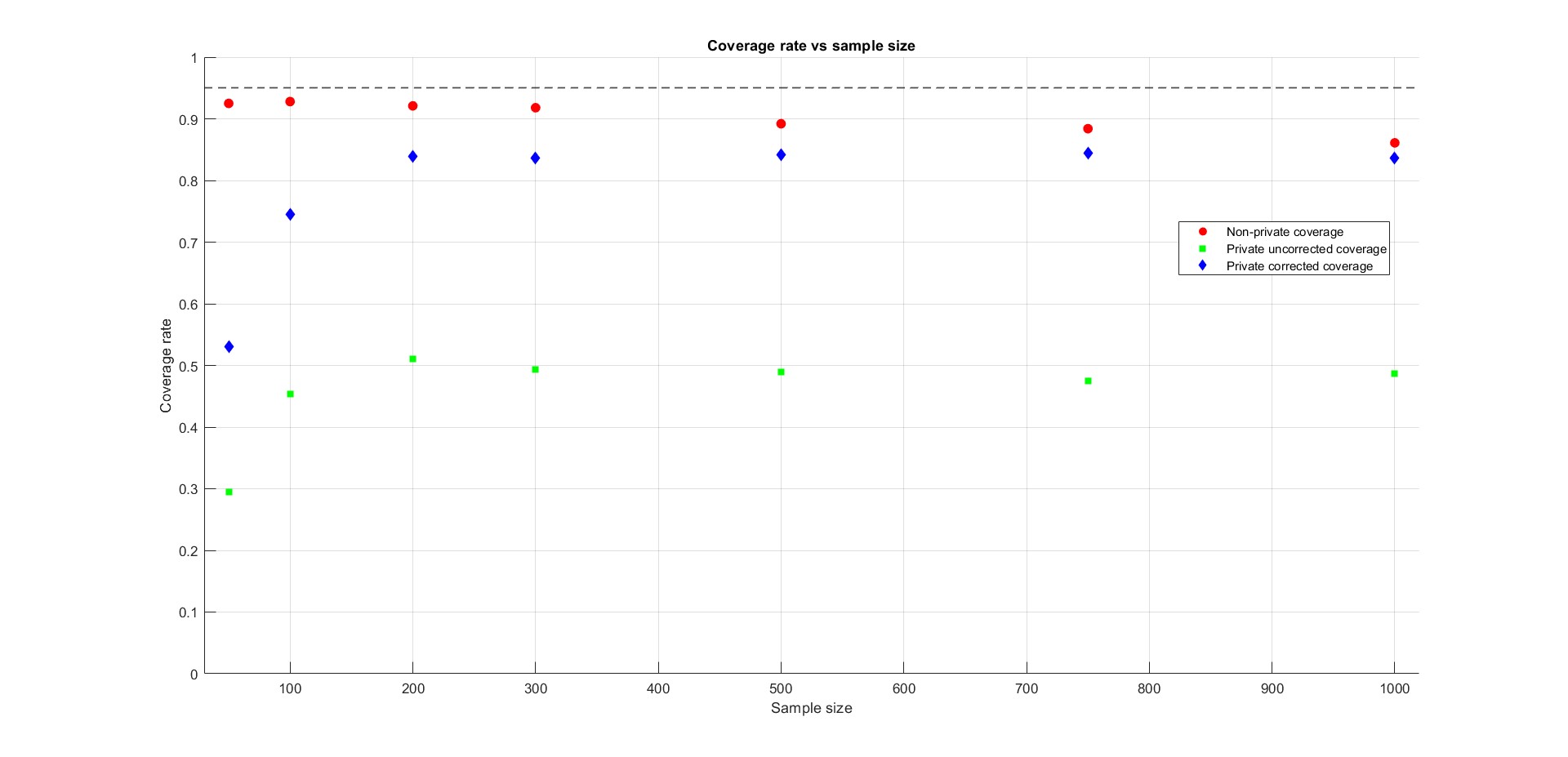}
    \caption{Private and non-private gradient descent 95\% confidence interval coverage}
    \label{F5}
\end{figure}

The behavior of the solution $\hat\bta^{(j)}_n$ for PGD and PNR algorithms across iterations, representing the solution trajectory, are given in Figure \ref{F1}, Figure \ref{F2}, and Figure 11.1,
Figure 11.2
respectively.  The figures are based on 20 repetitions. 
The coverage rate of the 95\% confidence interval for $\mu$ against the sample size, for $\ep=0.6$,  is given in Figure \ref{F5} for PGD and Figure 11.3
for PNR respectively.

In Table 11.2,
we present results comparing our HDP to PDP for values of $\lambda$ away from $\frac{-1}{2}$. The noise variance is now derived using (\ref{eq:additive_noise}) Theorem \ref{thm:PDP_mech}. We observe that the standard error of the estimates using the PGD and PNR algorithm increases as one deviates from the optimal value of $\frac{-1}{2}$. Inspection of the coverage shows that despite an increase in the standard error, the coverage rate of the CI is poor in contrast to the HDP setting.

\begin{table}[ht]
\centering
\begin{tabular}{l|l|c|c|c}
\toprule
& & \multicolumn{3}{c}{$\epsilon$} \\
\cmidrule(lr){3-5}
& &  2.00 & 0.60 & 0.20 \\
\midrule
\multirow{2}{*}{Estimator} 
& $\mu$: Mean (Std. Error) & 4.991 (0.083) & 4.989 (0.2) & 4.996 (0.349) \\ \cmidrule(lr){2-5}
& $\sigma$: Mean (Std. Error) & 1.984 (0.058) & 2.002 (0.144) & 2.043 (0.256) \\
\midrule
\multirow{2}{*}{CI coverage for $\mu$} 
& Corrected & 0.861 & 0.836 & 0.824 \\ \cmidrule(lr){2-5}
& Uncorrected & 0.861 & 0.487 & 0.327 \\
\midrule
\multirow{2}{*}{CI coverage for $\sigma$} 
& Corrected & 0.819 & 0.933 & 0.927 \\ \cmidrule(lr){2-5}
& Uncorrected & 0.819 & 0.459 & 0.284 \\
\bottomrule
\end{tabular}
\caption{Results for different values of $\epsilon$. Sample size is 1000, $K = 50$.}
\label{T1}
\end{table}

\textbf{Robustness}:
We now describe the robustness properties of PMHDE by investigating the behavior under a gross-error contamination model.  Denote by $f_{\bta,\alpha}$ the contamination model,
\begin{align*}
    f_{\bta,\alpha}(x) = (1-\alpha)f_{\bta}(x) + \alpha U_{z},
\end{align*}
where $\AV{U_{z}}$ is the uniform density on the interval $[z-\kappa,z+\kappa]$ for a small $\kappa>0$. Note that $f_{\bta,\alpha}(x)$ represents $\AV{100}\alpha\%$ contamination with distant outliers. In our experiments, $U_{z}$ is the uniform density on $[q(0.985,f_{\bta}),q(0.995,f_{\bta})]$, where $q(0.985,f_{\bta})$ is the 98.5\% quantile of $f_{\bta}$ and $q(0.995,f_{\bta})$ is the 99.5\% quantile of $f_{\bta}$; that is, $[9.34,10.15]$. We apply varying contamination levels with $\alpha=0,0.05,0.1,0.2,0.3$ for both PGD and PNR Algorithms.  The result for the PGD algorithm is shown in Table \ref{T5}, while that for the PNR algorithm is shown in Table 11.4.

We note that at high privacy levels, the iterative algorithm will yield estimates with high variability or tend to deviate from the true value since noise with larger variance is introduced in each iteration. This phenomenon is more common when the sample sizes are smaller. In these cases, it is common to use a thresholding strategy, which results in data in the extreme tails being suppressed. One approach to establishing the thresholds is by looking at the extreme tails of non-private estimators, while other methods are also feasible. In Section 11,
we provide numerical experiments illustrating the behavior of private estimators for different sample sizes (ranging from 200 to 500) and privacy levels.

Turning to Table 11.4
last row ($\ep=0.2$), we note that the standard error for PMHDE is larger due to aberrant values of the private estimate of $\mu$ in certain data sets. It turns out, in the no contamination case, 35 out of 5000 experiments yield estimates of private $\mu$ that are much larger than 10 or smaller than 0.15, while the true value is 5. In these cases, the usefulness of the estimate is in question. In applied settings, it is common not to release such values, and ad-hoc measures are adopted to circumvent this problem. We used the lower 0.7\% and upper 99.5\% percentiles of a Gaussian distribution with non-private $\hat{\mu}_n$ and $ \hat{\sigma}_n$ to threshold the private estimates. This resulted in the following estimates for the case $\ep=0.2$:
4.863(0.807), 5.012(0.838), 5.135(0.811), 5.3340(0.850), 5.4430(0.921).
The results illustrate that PMHDE retains robustness (compared to MLE) even under contamination.

\begin{table}[ht]
\centering
\scriptsize
\begin{tabular}{c|c|c|c|c|c}
\toprule
& \multicolumn{5}{c}{Contamination percentage $\alpha$} \\
\cmidrule(lr){2-6}
& 0\% & 5\% & 10\% & 20\% & 30\% \\
\midrule
MLE (Std. Error) & 5.001 (0.063) & 5.241 (0.061) & 5.476 (0.059) & 5.952 (0.056) & 6.422 (0.052) \\
\midrule
PMHDE $\epsilon=2$ (Std. Error) & 4.991 (0.083) & 5.159 (0.089) & 5.291 (0.09) & 5.52 (0.095) & 5.715 (0.102) \\
\midrule
PMHDE $\epsilon=0.6$ (Std. Error) & 4.986 (0.199) & 5.158 (0.203) & 5.289 (0.204) & 5.516 (0.207) & 5.712 (0.214) \\
\midrule
PMHDE $\epsilon=0.2$ (Std. Error) & 4.992 (0.353) & 5.15 (0.349) & 5.288 (0.355) & 5.494 (0.367) & 5.675 (0.523) \\
\bottomrule
\end{tabular}
\caption{Contamination results, gradient descent, sample size is 1000, $\mu=5$.}
\label{T5}
\end{table}

\section{Extensions and concluding remarks} \label{sec:ext and cr}
In this paper, we developed new differential privacy concepts called $\ep-$HDP and $(\la,\ep)-$PDP and illustrated the optimality of the Hellinger choice within the power-divergence family for minimizing Gaussian noise variance at a common privacy budget. We used these concepts to develop PMHDE estimators that are robust, efficient, and private. These estimators are derived by privatizing the classical gradient descent and Newton-Raphson algorithms via a Gaussian mechanism. We analyzed the convergence properties of these algorithms and established that the resulting estimators are private, efficient, and robust. Since the models do not satisfy the strong convexity properties, we utilize ASLSC and smoothness derived from standard assumptions to analyze the resulting estimators. 

Our methods also work when the Gaussian mechanism is replaced by the Laplace mechanism. Almost all properties in Theorem \ref{thm:utility_GD}, Theorem \ref{thm:utility_NT}, and Theorem \ref{thm:asym} go through if one uses concentration bounds for Laplace random variables. Our initial analysis suggests that, in Theorem \ref{thm:utility_GD} and Theorem \ref{thm:utility_NT},  the utility takes the following form: $||\hat\bta_n^{(K_n)}-\hat\bta_n||_2=O_p\left(n^{-\frac{1}{p}}K_n^{\frac{1}{2}}\log (K_n)\right)$ which also yields efficiency without any change. A detailed analysis of this case with a concentration inequality for the Laplace random variables and other probabilistic properties of compositions, especially when the number of queries diverges, will be discussed elsewhere.

It is also possible to extend the results to other minimum divergence estimators, such as minimum negative exponential disparity estimators, blended weight Hellinger distance estimators, and recently developed $S-$estimators (see \cite{ghosh2017minimum}). We have not carried out all the technical details carefully; however, an initial heuristic analysis shows analogous versions of our results will continue to hold for each case. A unified approach for all these estimators under minimal conditions would be useful and is being studied by the authors.

Finally, replacing the power divergence with a general convex function to define an extended notion of privacy is also useful. However, obtaining closed-form expressions for the variance in the additive mechanism presents certain technical challenges.





\vspace{0.2in}

\noindent\textbf{Supplementary material and source code}

The online supplementary material contains the proofs of the main results in Section 6 and the regularity conditions, additional proofs and calculations, and numerical results in Sections 7--11. It is available at \url{https://github.com/Frederick00D/HDP/blob/main/Hellinger_Differential_Privacy_Supplementary.pdf}.

The repository at \url{https://github.com/Frederick00D/HDP} contains implementations of the main algorithms and the source code to generate tables and figures in the manuscript and the supplementary material.

\clearpage
\numberwithin{theorem}{section}
\numberwithin{figure}{section}
\numberwithin{table}{section}
\setcounter{theorem}{0}
\setcounter{equation}{0}
\setcounter{figure}{0}
\setcounter{table}{0}
\setcounter{assumption}{0}
\setcounter{assumptionU}{0}
\renewcommand{\theassumptionU}{{\bf (U\arabic{assumptionU})}}
\renewcommand{\theassumption}{{\bf (A\arabic{assumption})}}
\renewcommand{\theHequation}{supp.\arabic{equation}}
\def\theHtheorem{supp.\arabic{section}.\arabic{theorem}}
\def\theHfigure{supp.\arabic{section}.\arabic{figure}}
\def\theHtable{supp.\arabic{section}.\arabic{table}}
\def\theHassumptionU{supp.\arabic{assumptionU}}
\def\theHassumption{supp.\arabic{assumption}}
\setcounter{section}{5}

\section{Proofs} \label{sec:proof}
In this section we provide proofs and detailed results in main-paper order. Unprefixed result numbers identify results in the main paper or statements relocated here; section-prefixed result numbers identify new supplementary results. Section numbers and equation references within the supplement refer to this document unless the main paper is named. Sections 7--11 contain regularity conditions, mechanism calculations, privacy comparisons, supporting proofs, and additional numerical results, respectively.
\subsection{Proof of Proposition 6}
\label{app:pd-partition-proof}
\begin{proof}
Put $\nu=P_1+P_2$ and $q=dP_1/d\nu$, so $dP_2/d\nu=1-q$.
For the proof, write $D_\lambda(P_1,P_2)=\int J_\lambda(q)\,d\nu$, where
\[
J_\lambda(u)=
\frac{u^{\lambda+1}(1-u)^{-\lambda}
-(\lambda+1)u+\lambda(1-u)}
{\lambda(\lambda+1)},
\qquad 0<u<1.
\]
At $\lambda=0,-1$, respectively, use
$J_0(u)=u\log\{u/(1-u)\}-2u+1$ and
$J_{-1}(u)=(1-u)\log\{(1-u)/u\}+2u-1$.
Give these functions their lower-semicontinuous extensions to $[0,1]$.
They are nonnegative and convex: on $(0,1)$,
$J_\lambda(1/2)=J_\lambda'(1/2)=0$ and
$J_\lambda''(u)=u^{\lambda-1}(1-u)^{-\lambda-2}>0$.
The linear terms in the numerator integrate to $-1$, so the integral
is the original power divergence, with the stated boundary conventions.

For a cell $A$ of positive $\nu$-measure, put
$\bar q_A=P_1(A)/\nu(A)$. Jensen's inequality gives
\[
D_\lambda(P_1^\pi,P_2^\pi)
=\sum_{A\in\pi:\,\nu(A)>0}\nu(A)J_\lambda(\bar q_A)
\leq\int J_\lambda(q)\,d\nu.
\]
Applying this inequality within the cells of a coarser partition proves
monotonicity under refinement.

Let $(I_{r,j})_j$ be the nested dyadic partitions of $[0,1]$, with
interval length $2^{-r}$ and $1$ included in the last interval, and set
$A_{r,j}=\{q\in I_{r,j}\}$. Omit empty cells to obtain $\pi_r$, and define
\[
s_r=\sum_j
\min\left\{
r,\inf_{u\in\overline{I_{r,j}}}J_\lambda(u)
\right\}\mathbf1_{A_{r,j}}.
\]
These are nonnegative simple functions. Nestedness and lower
semicontinuity give $s_r\uparrow J_\lambda(q)$, including points where
the limit is infinite. Also, $\bar q_{A_{r,j}}$ belongs to
$\overline{I_{r,j}}$, and hence
\[
\int s_r\,d\nu
\leq D_\lambda(P_1^{\pi_r},P_2^{\pi_r})
\leq D_\lambda(P_1,P_2).
\]
The monotone convergence theorem proves the finite-partition identity, also when the divergence is infinite.
The privacy equivalence follows by applying the identity to each
neighboring pair.
\end{proof}

\subsection{Proof of Theorem 8}
We start with the proof of Theorem 8.
We begin with the case $\la(\la+1) \neq 0$. Since $M_1(w, D)$ is $(\la, \ep)$-PDP, the power divergence between $M_1(w, D)$ and $M_1(w, D')$ is atmost $\ep_1$. For brevity, we denote the random variables $M_1(w, D)$ and $M_1(w, D')$ by $X_1$ and $X_2$ respectively. Let $p_1(\cdot)$ and $p_2(\cdot)$ denote their densities. Thus, by the $(\la, \ep)-$PDP property it follows that,
\begin{align} \label{eq:1}
    \frac{1}{\la(\la+1)}\mathbf{E}_{p_2}\left[\left(\frac{p_1(X_2)}{p_2(X_2)}\right)^{\la+1}\right] \le\ep_1 +\frac{1}{\la(\la+1)}
\end{align}
Next, let the random variables $Y_1|X_1$ and $Y_2|X_2$ represent the compositions of mechanisms $M_2(M_1(w,D),D)$ given $M_1(w, D)$ and $M_2(M_1(w,D'),D')$ given $M_1(w, D')$ with conditional densities $q_{Y_1|X_1}(\cdot)$ and $q_{Y_2|X_2}(\cdot)$ respectively.
Again using the PDP property, it follows that for a generic random variable $V\sim q_{Y_2|X_2}(\cdot)$,
\begin{align} \label{eq:2}
    \frac{1}{\la(\la+1)}\mathbf{E}_{q_{Y_2|X_2}}\left[\left(\frac{q_{Y_1|X_1}(V)}{q_{Y_2|X_2}(V)}\right)^{\la+1}\right] \le\ep_2 +\frac{1}{\la(\la+1)}
\end{align}
Now, to calculate power divergence of $M^{(2)}(w, D)$ and $M^{(2)}(w, D')$, we need to calculate the joint power divergence of $(X_1,Y_1)$ and $(X_2, Y_2)$. Now, using that the joint density is the product of conditional density and the marginal density, it follows that
{\footnotesize
\begin{eqnarray*}
    D_{\la}(M^{(2)}(w, D),M^{(2)}(w, D'))= \frac{1}{\la(\la+1)}\mathbf{E}_{p_2}\left\{\left(\frac{p_1(X_2)}{p_2(X_2)}\right)^{\la+1} \mathbf{E}_{q_{Y_2|X_2}}\left[\left(\frac{q_{Y_1|X_1}(V)}{q_{Y_2|X_2}(V)}\right)^{\la+1}\right]-1\right\}.
\end{eqnarray*}}
Now, suppose $\la(\la+1) >0$. Then, using (\ref{eq:1}) and (\ref{eq:2}) it follows that
\begin{eqnarray*}
  D_{\la}(M^{(2)}(w, D),M^{(2)}(w, D')) &\le& \frac{1}{\la(\la+1)} [\ep_1(\la(\la+1))+1]  [\ep_2(\la(\la+1))+1]- \frac{1}{\la(\la+1)}\\
  &=& \epsilon_1+\epsilon_2 + \lambda(\lambda+1)\epsilon_1\epsilon_2.
\end{eqnarray*}
Next, if $\la(\la+1) <0$, then using the condition $0<\ep <-[\la(\la+1)]^{-1}$ the above inequality continues to hold.
Finally, consider the case $\la=0$. In this case, the power divergence reduces to the Kullback Leibler (KL) divergence between the densities. Hence, using the PDP property, it follows that
\begin{align*}
    D_0(p_1, p_2)=KL(p_1, p_2) \le \ep_1, ~~\text{and}~~ D_0(q_{Y_1|X_1}, q_{Y_2|X_2})=KL(q_{Y_1|X_1}, q_{Y_2|X_2}) \le \ep_2.
\end{align*}
Hence, with $V \sim q_{Y_1|X_1}(\cdot)$, we have
\begin{eqnarray*}
    D_{\la}(M^{(2)}(w, D),M^{(2)}(w, D'))= \mathbf{E}_{p_1}\left[\log \frac{p_1(X_2)}{p_2(X_2)} \right] + \mathbf{E}_{q_{Y_1|X_1}}\left[\log\frac{q_{Y_1|X_1}(V)}{q_{Y_2|X_2}(V)}\right] \leq \ep_1+\ep_2.
\end{eqnarray*}
The proof of the case $\la=-1$ is similar. This completes the proof of (1). Next, the proof of (2) follows exactly as in (1) except that $M_2(X_1,D)|X_1$ and $M_2(X_2,D')|X_2$ are now replaced by $M_2(w,D)$ and $M_2(w,D')$. 
That is, we replace $\frac{q_{Y_1|X_1}(V)}{q_{Y_2|X_2}(V)}$ in (\ref{eq:2})  by $\frac{q_{Y_1}(V)}{q_{Y_2}(V)}$, where $q_{Y_1}$ and $q_{Y_2}$ are the unconditional distributions of $M_2(w,D)$ and $M_2(w,D')$ respectively and $V$ has the density $q_{Y_2}$.          

To start the proof of part (3), we first notice that adjacent $\mathbf{D}^{(2)}$ and $\mathbf{D}^{(2)'}$ can be decomposed into two distinct cases. By definition, 
\begin{align*}
    ||\mathbf{D}^{(2)}-\mathbf{D}^{(2)'}||_H = \sum_{i=1}^2 ||D_i-D_i'||_H=1.
\end{align*}
Since the Hamming distance is a non-negative integer, the above equation holds if either: Case (1): $||D_1-D_1'||_H=1$ and $||D_2-D_2'||_H=0$; or Case (2): $||D_1-D_1'||_H=0$ and $||D_2-D_2'||_H=1$. Let, as before, $p_1(\cdot)$, $p_2(\cdot)$ denote the distributions of $M_1(w_1,D_1)$, $M_1(w_1,D_1')$. Also, let $q_1(\cdot)$, and $q_2(\cdot)$ are density functions of $M_2(w_2,D_2)$ and $M_2(w_2,D_2')$ respectively. The joint density of $M^{(2)}(\mathbf{W}^{(2)},\mathbf{D}^{(2)})$ and $M^{(2)}(\mathbf{W}^{(2)},\mathbf{D}^{(2)'})$ are therefore given by $h_1(x,y)$ and $h_2(x,y)$ respectively, where $h_1(x,y)=p_1(x)q_1(y)$ and $h_2(x,y)=p_2(x)q_2(y)$.

If case (1) happens, $q_1(\cdot)=q_2(\cdot)$ holds since $D_2=D_2'$. This leads to $\frac{h_1(x,y)}{h_2(x,y)}=\frac{p_1(x)}{p_2(x)}$. The PD between $M^{(2)}(\mathbf{W}^{(2)},\mathbf{D}^{(2)})$ and $M^{(2)}(\mathbf{W}^{(2)},\mathbf{D}^{(2)'})$ is reduced to the PD between $M_1(w_1,D_1)$ and $M_1(w_1,D_1')$, since
\begin{align*}
    D_\la(M^{(2)}(\mathbf{W}^{(2)},\mathbf{D}^{(2)}),M^{(2)}(\mathbf{W}^{(2)},\mathbf{D}^{(2)'})) =& \frac{1}{\la(\la+1)}\mathbf{E}_{h_2}\left[\left(\frac{h_1(X,Y)}{h_2(X,Y)}\right)^{\la+1} -1\right]\\
    =& \frac{1}{\la(\la+1)}\mathbf{E}_{p_2}\left[\left(\frac{p_1(X)}{p_2(X)}\right)^{\la+1} -1\right]\\
    =& D_\la(M_1(w_1,D_1),M_1(w_1,D_1')) \leq \ep_1.
\end{align*}
The last inequality follows from $M_1(w_1,D_1)$ is $(\la,\ep_1)-$PDP. Similarly in Case (2),
\begin{align*}
    D_\la(M^{(2)}(\mathbf{W}^{(2)},\mathbf{D}^{(2)}),M^{(2)}(\mathbf{W}^{(2)},\mathbf{D}^{(2)'}))=D_\la(M_2(w_2,D_2),M_2(w_2,D_2')) \leq \ep_2.
\end{align*}
Combining Case (1) and Case (2), we get\[D_\la(M^{(2)}(\mathbf{W}^{(2)},\mathbf{D}^{(2)}),M^{(2)}(\mathbf{W}^{(2)},\mathbf{D}^{(2)'}))\leq \max\{\ep_1,\ep_2\}\] which implies that the parallel composition 
is $(\la,\max\{\ep_1,\ep_2\})-$PDP. This completes the proof of Theorem 8. \QED

\subsection{Proof of Theorem 9} 
We begin with the case when $Y_i$'s are $N(0, \si^2)$. In this case, using Lemma \ref{appb:lem1} we have that
\begin{eqnarray*}
    D_{\la}(M(w,D), M(w, D')) &=& \frac{1}{\lambda(\lambda+1)} \left[  e^{\frac{\lambda(\lambda+1)||\mathbf{v}||_2^2}{2\sigma^2}}   -1 \right],
\end{eqnarray*}
where $\mathbf{v}$ is the difference between the mean of $M(w,D)$ and the mean of $M(w,D')$, and for $r=1,2$, $||\mathbf{v}||_r\AV{\leq}\Delta_{L_r}W$. Thus, \AV{a sufficient condition for} $D_{\la}(M(w,D), M(w, D')) \le \ep$ \AV{is} 
\begin{eqnarray*}
\sigma \geq ||\mathbf{v}||_2\sqrt{\frac{\lambda(\lambda+1)}{2\log (1+\lambda(\lambda+1)\epsilon)}},
\end{eqnarray*}
which is well-defined for all $\ep>0$ when $\la(\la+1)>0$, and for $0<\ep<-[\la(\la+1)]^{-1}$ when $\la(\la+1)<0$. Thus, we choose
\[
\si^2_{\la, \ep}=(\Delta_{L_2}W)^2\frac{\lambda(\lambda+1)}{2\log (1+\lambda(\lambda+1)\epsilon)}.
\]
Next, when $\la=0$,
\begin{eqnarray*}
    D_0(M(w, D), M(w, D'))=\AV{\frac{||w_1-w_2||_2^2}{2\sigma^2}}
\end{eqnarray*}
Again, $D_0{M(w, D), M(w, D')} \le \ep$ \AV{is ensured by choosing}
\begin{eqnarray*}
    \si^2_{0, \ep}=  \frac{(\Delta_{L_2} W)^2 }{2\epsilon}.
\end{eqnarray*}
The case for $\la=-1$ is similar. Next, turning to the Laplace case, using Lemma \ref{appb:lem2}, notice that
\begin{eqnarray*}
D_{\la}(M(w,D), M(w, D')) &=&\frac{1}{\lambda(\lambda+1)} \left[ \left(\prod_{i=1}^{m} \int_{\mathbb{R}} \frac{1}{2b} e^{-\frac{(\lambda+1)|y_i-v_i|-\lambda |y_i|}{b}}\mathrm{d}y_i \right) -1 \right]\\
&\leq&\frac{1}{\lambda(\lambda+1)} \left[  e^{\frac{sign(\lambda)(\lambda+1)||v||_1}{b}} -1 \right].
\end{eqnarray*}
Similarly,
\begin{align*}
    D_{\la}(M(w,D'), M(w, D)) \leq \frac{1}{\lambda(\lambda+1)} \left[  e^{\frac{sign(\lambda+1)(\lambda)||v||_1}{b}} -1 \right].
\end{align*}
Thus, $D_{\la}(M(w,D), M(w, D')) \le \ep$ \AV{is ensured by choosing}
\begin{eqnarray*}
    b\geq \max\left\{ \frac{sign(\lambda)(\lambda+1)||v||_1}{\log(\lambda(\lambda+1)\epsilon+1)},  \frac{sign(\lambda+1)(\lambda)||v||_1}{\log(\lambda(\lambda+1)\epsilon+1)} \right\}.
\end{eqnarray*}
Thus, we choose $b_{\la, \ep}$ to be
\begin{align*}
    b_{\la, \ep}= \max\left\{ \frac{sign(\lambda)(\lambda+1)\Delta_{L_1} W}{\log(\lambda(\lambda+1)\epsilon+1)},  \frac{sign(\lambda+1)(\lambda)\Delta_{L_1} W}{\log(\lambda(\lambda+1)\epsilon+1)} \right\}.
\end{align*}
Turning to the case $\la(\la+1)=0$, we note that
\begin{eqnarray*}
D_0(M(w,D), M(w, D'))&=&\sum_{i=1}^m \int_{\mathbb{R}}\frac{1}{2b}e^{-\frac{|y_i|}{b}}\cdot \frac{|y_i-v_i|-|y_i|}{b}\mathrm{d}y_i\\
&\leq& \frac{||w_1-w_2||_1}{b}
\end{eqnarray*}
Hence, $D_0(M(w,D), M(w, D')) \le \ep$  \AV{holds if} $b\geq  \frac{\AV{\Delta_{L_1} W}}{\epsilon}$. Thus,
\begin{eqnarray*}
    b_{0, \ep}=\frac{\Delta_{L_1} W}{\epsilon}.
\end{eqnarray*}
The proof for the case $\la=-1$ is similar. \QED

\subsection{Proof of Theorem 10}
The result follows from Theorem 8 by taking $\la=-\frac12$ and
using $D_{-1/2}=2D_{HD}$. \QED

\subsection{Proof of Corollary 11}

The proof is based on the following iterative argument for adaptive and sequential compositions. Setting $\ep_1= \ep$ and $\ep_2=h_1(\ep)$ it follows from Theorem 10
part 1. and part 2., that $h_2(\ep)=\ep+h_1(\ep)-\frac{1}{2}\ep h_1(\ep)$. Now iterating, we obtain $h_{j+1}(\ep)= \ep + h_j(\ep) - \frac{1}{2}\ep h_j(\ep)$. The proof for the parallel compositions follows from part 3. of Theorem 10. \QED

\subsection{Proof of Proposition 12
}
The proof of the Proposition follows using a comparison argument. Recall that the total variation distance between two densities can be expressed as one-half the $L_1-$norm, which is bounded above by the Hellinger distance between the densities. That is, 
\begin{align*}
     TV(p_1, p_2)= \frac{1}{2}||p_1-p_2||_1 \leq HD(p_1,p_2).
\end{align*}
Now, if $HD^2(p_1,p_2) \le \ep$, then $TV(p_1, p_2) \le \sqrt{\ep}$ which implies that $M(\cdot,\cdot)$ satisfies $\sqrt{\ep}-$TV privacy. Hence, using \cite{Ghazi2024} page 209, it follows that  $M$ also satisfies $(0,\sqrt{\ep})$ differential privacy.
Turning to $\mu$-GDP, suppose that $M$ satisfies either of the
additional mechanism assumptions in Proposition 12. Let
$\eta=\sup_{D\sim D'}\operatorname{TV}(P_D,P_{D'})$, where $P_D$
denotes the output law at $D$. For the Gaussian mechanism, the
exact Gaussian testing formula in \cite{Dong2022} gives the
smallest GDP parameter as
$2\Phi^{-1}\{(1+\eta)/2\}$. The same identity holds for the
scalar Laplace mechanism by
Proposition~\ref{prop:laplace-gdp-exact}.
Since the first part gives $\eta\leq\sqrt{\epsilon}$,
monotonicity of $\Phi^{-1}$ implies that $M$ satisfies
$\mu$-GDP with
\(
\mu=2\Phi^{-1}\left(\frac{1+\sqrt{\epsilon}}{2}\right)
\).
\QED

\subsection{Proof of Theorem 13}

For $k$-neighbor datasets $D$ and $D'$, we need to bound\[D_{HD}(M(w,D),M(w,D')).\]By the definition of $k$-neighbor datasets, there exist datasets $B_0,\ldots,B_k$ with $B_0=D$ and $B_k=D'$, such that 
\[||B_i-B_{i+1}||_H=1,\qquad i=0,1,\ldots,k-1.\] Also, since 
$M(w, B_i)$ is $\ep-$HDP for all $0 \le i\le k$, we get that
\[D_{HD}(M(w,B_i), M(w, B_{i+1}))\le \ep.\]
Now, using that $D_{HD}^{\frac{1}{2}}(\cdot, \cdot)=HD(\cdot, \cdot)$ is a metric, using the triangle inequality 
\begin{eqnarray}{\label{eq:3}}
D_{HD}^{\frac{1}{2}}(M(w,D), M(w, D')) \le \sum_{i=1}^k D_{HD}^{\frac{1}{2}}(M(w,B_{i-1}), M(w, B_i)) \le k \sqrt{\ep}.
\end{eqnarray}
The result follows by squaring both sides of (\ref{eq:3}). \QED

\subsection{Proof of Proposition 14}
The HDP mechanism bounds follow from Theorem 9 by taking
$\la=-\frac12$ and replacing the PDP budget by $2\ep$,
using $D_{-1/2}=2D_{HD}$. \QED

\subsection{Proof of Proposition 16}
\begin{proof}
Writing $a=P_1(A)$ and $b=P_2(A)$, the partition identity implies
\begin{equation}
(\sqrt a-\sqrt b)^2+
(\sqrt{1-a}-\sqrt{1-b})^2\leq h.
\label{eq:hdp-op-binary}
\end{equation}
Set $a=\sin^2\theta$ and $b=\sin^2\phi$, with
$\theta,\phi\in[0,\pi/2]$. The inequality becomes
$\cos(\theta-\phi)\geq c_h$. With $\omega_h=\arccos c_h$, it follows that
\[
\max\{0,\phi-\omega_h\}\leq\theta
\leq\min\{\pi/2,\phi+\omega_h\}.
\]
Taking squared sines gives the event-probability bounds of Proposition 16.
Bernoulli laws with probabilities at the two endpoints attain the
bounds. Setting $b=0$ proves the zero-probability bound of Proposition 16.
Dropping the second term in \eqref{eq:hdp-op-binary} gives
$|\sqrt a-\sqrt b|\leq\sqrt h$. Since
$\sqrt a\leq\sqrt b+\sqrt h$, we obtain
$|a-b|\leq2\sqrt{hb}+h$, proving the relative-probability bound of Proposition 16.
\end{proof}

\subsection{Proof of Proposition 17}
\begin{proof}
Let $R,R_i$ be the actual and counterfactual joint laws of
$(D,Y)$. Since both have prior $\Pi$ on $D$,
\[
D_{HD}(R,R_i)
=\int_{\mathcal D}D_{HD}(P_u,P_{u^{(-i)}})\,\Pi(du)\leq h.
\]
Thus their affinity $\rho$ is at least $c_h$. Let $m,m_i$ be the two
output densities relative to a common dominating measure $\zeta$, and
let $a(y)$ be the affinity between the two posterior laws. Disintegration
and Cauchy--Schwarz give
\[
\rho=\int\sqrt{m(y)m_i(y)}\,a(y)\,d\zeta(y),
\qquad
\rho^2\leq\int m(y)a(y)^2\,d\zeta(y).
\]
For any probability laws $Q_1,Q_2$, another application of
Cauchy--Schwarz gives
$\operatorname{TV}(Q_1,Q_2)^2\leq1-\rho(Q_1,Q_2)^2$.
Hence the left side of the claimed expectation bound is at most
$1-\rho^2\leq1-c_h^2=s_h^2$.
Markov's inequality proves the asserted probability bound.
Where the counterfactual marginal density is zero, any posterior version
may be assigned; the inequalities remain valid.
\end{proof}

\subsection{HDP and approximate-DP conversions (Propositions 19--20)}
\paragraph{From HDP to differential privacy.}
For $t\geq0$, define the directional pairwise privacy profiles
\begin{equation}
\delta_{12}(t)=\sup_A\{P_1(A)-e^tP_2(A)\}
=\int(p_1-e^tp_2)_+\,d\nu,
\qquad
\delta_{21}(t)=\int(p_2-e^tp_1)_+\,d\nu.
\label{eq:hdp-op-profiles}
\end{equation}
An $(\epsilon,\delta)$-DP guarantee bounds both quantities at
$t=\epsilon$ for every neighboring pair. The profile formulation records
all such guarantees; see \citet{BalleBartheGaboardi2020}.

\noindent\textbf{Proposition 19 (HDP to approximate DP).}
An $h$-HDP mechanism satisfies $(\epsilon,\Delta_h(\epsilon))$-DP for
every $\epsilon\geq0$, where
\begin{equation}
\Delta_h(\epsilon)=
\frac{\sqrt{(e^\epsilon-1)^2+4e^\epsilon s_h^2}
-(e^\epsilon-1)}{2}.
\label{eq:hdp-op-to-dp}
\end{equation}
This conversion is sharp over the Hellinger class for each
$\epsilon$. In particular,
$\Delta_h(0)=s_h$ and
$\lim_{\epsilon\to\infty}\Delta_h(\epsilon)=s_h^2$.
\par\medskip

\begin{proof}
Put $a=e^\epsilon$ and $\rho=\rho(P_1,P_2)\geq c_h$.
Cauchy--Schwarz gives
\[
\int|p_1-ap_2|\,d\nu
\leq\left\{(1+a)^2-4a\rho^2\right\}^{1/2}.
\]
Since
$\delta_{12}(\epsilon)
=\{\int|p_1-ap_2|\,d\nu+1-a\}/2$,
the asserted bound follows.
Interchanging $P_1,P_2$ proves the reverse direction.
For sharpness, suppose $0<h<2$, set
$S=\{(a-1)^2+4as_h^2\}^{1/2}$ and $d=\Delta_h(\epsilon)$,
and let
$b=(d-s_h^2)/S$ and $u=d+ab$.
These numbers belong to $[0,1]$, and direct substitution gives
\[
D_{HD}\{\operatorname{Ber}(u),\operatorname{Ber}(b)\}=h,
\qquad u-ab=d.
\]
Thus the bound is attained on the success event. The cases $h=0,2$
follow from identical and mutually singular laws, respectively.
The endpoint formulas follow from \eqref{eq:hdp-op-to-dp}.
\end{proof}

The function $\Delta_h$ is a universal envelope over the HDP class,
not necessarily the attained profile of a specified mechanism.
Exact Gaussian and Laplace calibrations retain that mechanism-specific
information. For example, different divergence parameters can describe
the same neighboring output laws; comparing their numerical values
alone does not compare disclosure probabilities.

There is also a direct relation between power divergences and the
pairwise profiles:
\begin{equation}
D_\lambda(P_1,P_2)
=\int_0^\infty
\left\{
e^{\lambda t}\delta_{12}(t)
+e^{-(\lambda+1)t}\delta_{21}(t)
\right\}\,dt,
\qquad \lambda\in\mathbb R.
\label{eq:hdp-op-power-profile}
\end{equation}
Here the integral is allowed to be infinite. To verify the identity,
when $a\geq b>0$ and $\lambda\ne-1,0$, integrate pointwise:
\[
\int_0^{\log(a/b)}e^{\lambda t}(a-e^tb)\,dt
=
\frac{
a^{\lambda+1}b^{-\lambda}-(\lambda+1)a+\lambda b
}{\lambda(\lambda+1)}.
\]
For $a<b$, use the second term in
\eqref{eq:hdp-op-power-profile}. At $\lambda=0,-1$, integrate directly
to obtain the corresponding logarithmic expression; zero-density
cases use the same extended integrals. Tonelli's theorem and
$\int p_1\,d\nu=\int p_2\,d\nu=1$ prove the identity.
In particular,
\begin{equation}
D_{HD}(P_1,P_2)
=\frac12\int_0^\infty e^{-t/2}
\{\delta_{12}(t)+\delta_{21}(t)\}\,dt.
\label{eq:hdp-op-hd-profile}
\end{equation}
Thus HDP constrains a weighted integral of the pairwise disclosure
profiles, whereas a specified approximate-DP guarantee constrains
those profiles at a specified threshold. These identities are
pairwise; no interchange of a supremum over neighboring datasets and
an integral is being asserted.

\paragraph{From differential privacy to HDP.}
The converse conversion also follows directly from the neighboring
output laws.

\noindent\textbf{Proposition 20 (Approximate DP to HDP).}
Under symmetric adjacency, an $(\epsilon,\delta)$-DP mechanism, with
$\epsilon\geq0$ and $0\leq\delta\leq1$, satisfies $h_*$-HDP, where
\begin{equation}
h_*=2\left\{1-\frac{1-\delta}{\cosh(\epsilon/2)}\right\}.
\label{eq:hdp-op-from-dp}
\end{equation}
The bound is sharp. In particular, pure $\epsilon$-DP implies
$2\{1-1/\cosh(\epsilon/2)\}$-HDP.
\par\medskip

\begin{proof}
Put $a=e^\epsilon$, $u=\min\{p_1,ap_2\}$, and
$v=\min\{p_2,ap_1\}$. The two directional DP inequalities give
$\int u\,d\nu\geq1-\delta$ and $\int v\,d\nu\geq1-\delta$.
Where they are positive, $a^{-1}\leq u/v\leq a$, and hence
\[
\sqrt{p_1p_2}\geq\sqrt{uv}
\geq\frac{\sqrt a}{1+a}(u+v).
\]
Integration yields
$\rho(P_1,P_2)\geq2(1-\delta)\sqrt a/(1+a)$ and proves
\eqref{eq:hdp-op-from-dp}.
Equality is attained by
\begin{align*}
P_1&=\left(
\delta,\frac{(1-\delta)a}{1+a},
\frac{1-\delta}{1+a},0
\right),\\
P_2&=\left(
0,\frac{1-\delta}{1+a},
\frac{(1-\delta)a}{1+a},\delta
\right).
\end{align*}
Indeed,
$\sum_j(P_{1j}-aP_{2j})_+
=\sum_j(P_{2j}-aP_{1j})_+=\delta$,
and the affinity equals $2(1-\delta)\sqrt a/(1+a)$.
This is the four-output extremal experiment used by
\citet{Kairouz2015}. Replacing the four points by uniform laws on
disjoint unit intervals gives an attaining example with Lebesgue
densities.
\end{proof}

\paragraph{Distinguishing events and comparison of protections.}
Let $U_0,U_1$ be uniform laws on disjoint unit intervals. For $0<h<2$,
set $r=h-h^2/4$ and take
\[
P_2=U_0,\qquad P_1=(1-r)U_0+rU_1.
\]
Then $D_{HD}(P_1,P_2)=2(1-\sqrt{1-r})=h$.
The event $A=\operatorname{supp}(U_1)$ has
$P_2(A)=0$ and $P_1(A)=r>0$, so it violates every finite pure-DP
guarantee. Observing this event identifies the first of the two
neighboring datasets with certainty. Its probability is nevertheless
bounded sharply by the zero-probability bound of Proposition 16.
The same pair satisfies $(0,r)$-DP, so this example is not excluded
by approximate DP at all parameter values.
It also rules out a general implication from positive-budget HDP to
finite-parameter GDP: testing $P_2$ against $P_1$ has positive power
at zero type-I error, whereas a finite Gaussian trade-off has value
one at that endpoint \citep{Dong2022}.

The example can be realized as a mechanism on a binary database domain,
or by selecting $P_1$ or $P_2$ according to a binary database attribute.
It identifies a distinction between the constraints, not a universal
ordering of all disclosure risks. In particular, mechanisms with the
same attained Hellinger distance can have crossing privacy profiles;
Appendix~\ref{app:hdp-crossing} gives an example.
Relative protection must therefore be assessed at specified parameters
and for specified disclosure quantities. The exact calibrations
developed in Section 2.4 of the main paper provide the corresponding comparisons
within particular mechanism families.

\subsection{Proof and additional calculations for Proposition 21}
For the Gaussian mechanism in Section 2.4 of the main paper,
consider neighboring output laws $P_i=N(w_i,\sigma^2I_m)$ with $d>0$.
Their log-likelihood ratio is
\begin{equation}
\ell_{12}(y)=\log\frac{p_1(y)}{p_2(y)}
=\frac{(w_1-w_2)^T\{y-(w_1+w_2)/2\}}{\sigma^2}.
\label{eq:cal-gaussian-lr}
\end{equation}

Its distribution under the two neighboring output laws is

\begin{equation}
\ell_{12}(Y)\sim
\begin{cases}
N(\mu_d^2/2,\mu_d^2),&Y\sim P_1,\\
N(-\mu_d^2/2,\mu_d^2),&Y\sim P_2.
\end{cases}
\label{eq:cal-gaussian-loss-laws}
\end{equation}

\begin{proof}
Substitution of $Y=w_i+\sigma Z$ into the preceding likelihood-ratio identity
proves \eqref{eq:cal-gaussian-loss-laws}. Under $P_2$, for every
$s\in\mathbb R$,
\begin{equation}
\mathbf E_{P_2}e^{s\ell_{12}(Y)}
=\exp\{s(s-1)\mu_d^2/2\}.
\label{eq:cal-gaussian-moments}
\end{equation}
Taking $s=1/2$ gives the affinity $e^{-\mu_d^2/8}$ and proves
the Hellinger formula of Proposition 21.
For $\epsilon\geq0$, the event maximizing
$P_1(A)-e^\epsilon P_2(A)$ is
$A_\epsilon=\{\ell_{12}>\epsilon\}$. Therefore,
\begin{align*}
P_1(A_\epsilon)
&=\Phi(\mu_d/2-\epsilon/\mu_d),\\
P_2(A_\epsilon)
&=\Phi(-\mu_d/2-\epsilon/\mu_d).
\end{align*}
Their difference, with the second term multiplied by $e^\epsilon$,
is the Gaussian privacy-profile formula in Proposition 21. Interchanging the laws gives the
same profile. Setting $\epsilon=0$ gives the TV distance.
Finally, for $\mu>0$,
\begin{equation}
\frac{\partial}{\partial\mu}\delta_G(\epsilon;\mu)
=\phi(\mu/2-\epsilon/\mu)>0.
\label{eq:cal-gaussian-monotonicity}
\end{equation}
Here one uses
$e^\epsilon\phi(-\mu/2-\epsilon/\mu)
=\phi(\mu/2-\epsilon/\mu)$.
The Hellinger distance and TV are also increasing and continuous in
$\mu$. Taking suprema over neighboring pairs proves the final assertion,
even when the sensitivity supremum is approached but not attained.
\end{proof}

The likelihood-ratio calculation and exact approximate-DP condition are
those underlying the analytic Gaussian mechanism of
\citet{BalleWang2018}. In particular, it is the difference of the two
event probabilities above, rather than a bound on the first one alone,
that gives the exact calibration.

The same likelihood-ratio moments yield the power-divergence
identity stated after Proposition 21 in the main paper, including
$D_{-1/2}=2D_{HD}$ and the Gaussian R\'enyi divergence formula.

The Gaussian testing trade-off is
$G_{\mu_d}(a)=\Phi\{\Phi^{-1}(1-a)-\mu_d\}$, where $a$ is the
type-I error probability. This follows by projecting the observation
onto the direction $w_1-w_2$ and applying the likelihood-ratio test.
Thus the sharp GDP parameter of $M_\sigma$ is $\mu=\Delta/\sigma$
\citep{Dong2022}. Within this Gaussian family, the GDP parameter,
the attained Hellinger level, and the entire approximate-DP profile
are determined by the same sensitivity-to-noise ratio.

\subsection{Exact scalar Laplace calibration (Proposition 23)}
Consider $M_b(D)=W(D)+Z_b$ for a scalar query, where $Z_b$ has density
$(2b)^{-1}e^{-|z|/b}$, $b>0$, and variance $2b^2$.
In this paragraph, $\Delta=\sup_{D\sim D'}|W(D)-W(D')|\in(0,\infty)$.
For a neighboring pair, translation and reflection allow us to take
$P_1=\operatorname{Laplace}(0,b)$ and
$P_2=\operatorname{Laplace}(d,b)$, with $d\geq0$.
Their log-likelihood ratio is
\begin{equation}
\ell_{12}(y)=\frac{|y-d|-|y|}{b}
=\begin{cases}
d/b,&y\leq0,\\
(d-2y)/b,&0<y<d,\\
-d/b,&y\geq d.
\end{cases}
\label{eq:cal-laplace-lr}
\end{equation}

\noindent\textbf{Proposition 23 (Exact scalar Laplace privacy quantities).}
For the scalar Laplace mechanism, put $r=\Delta/(2b)$. Its exact
worst-case Hellinger level and approximate-DP profile are
\begin{equation}
h_L(b)=2\{1-(1+r)e^{-r}\},
\qquad
\delta_b^L(\epsilon)=[1-e^{\epsilon/2-r}]_+,
\qquad \epsilon\geq0.
\label{eq:cal-laplace-quantities}
\end{equation}
Its exact TV distance is $1-e^{-r}$, and its smallest pure-DP
parameter is $2r=\Delta/b$.
\par\medskip

\begin{proof}
For a pair separated by $d$, splitting the affinity integral at $0,d$
gives
\begin{align*}
\rho(P_1,P_2)
&=\frac1{2b}\int_{\mathbb R}
\exp\left\{-\frac{|y|+|y-d|}{2b}\right\}\,dy\\
&=\frac12e^{-d/(2b)}+\frac d{2b}e^{-d/(2b)}
+\frac12e^{-d/(2b)}
=\left(1+\frac d{2b}\right)e^{-d/(2b)}.
\end{align*}
For $0\leq\epsilon<d/b$, the positive region of
$p_1-e^\epsilon p_2$ is $(-\infty,a_\epsilon)$, where
$a_\epsilon=(d-\epsilon b)/2$. Direct integration gives
\begin{align*}
P_1(Y<a_\epsilon)&=1-\tfrac12e^{\epsilon/2-d/(2b)},\\
e^\epsilon P_2(Y<a_\epsilon)&=\tfrac12e^{\epsilon/2-d/(2b)}.
\end{align*}
Subtracting proves the profile formula in this range.
For $\epsilon\geq d/b$, the likelihood ratio never exceeds
$e^\epsilon$, and the profile is zero. Reflection gives the same
reverse-direction profile. Both discrepancies increase continuously
with $d$, so taking the sensitivity supremum replaces $d$ by $\Delta$.
The assertions about TV and pure DP follow from the profile.
\end{proof}

The Laplace privacy profile is also given in
\citet[Theorem~3]{BalleBartheGaboardi2020}; the affinity calculation
supplies its attained Hellinger level.
For $\epsilon\geq0$, $0\leq\delta<1$, and
$(\epsilon,\delta)\ne(0,0)$, the minimum admissible scale is
\begin{equation}
b_*(\epsilon,\delta)
=\frac{\Delta}{\epsilon-2\log(1-\delta)},
\qquad
\operatorname{Var}(Z_{b_*})=2b_*^2.
\label{eq:cal-laplace-forward}
\end{equation}
Indeed, the exact condition $[1-e^{\epsilon/2-r}]_+\leq\delta$
is equivalent to $r\leq\epsilon/2-\log(1-\delta)$.
At the minimum scale, put $r_*=\epsilon/2-\log(1-\delta)$.
Keeping $b_*$ fixed, the achieved Hellinger level is
\begin{equation}
h_{\mathrm{achieved}}^L(\epsilon,\delta)
=2\{1-(1+r_*)e^{-r_*}\}.
\label{eq:cal-laplace-achieved-h}
\end{equation}

Conversely, for $0<h<2$, let $r_h>0$ be the unique solution of
$(1+r_h)e^{-r_h}=1-h/2$. Uniqueness follows because the left side
has derivative $-r e^{-r}<0$ for $r>0$ and decreases from $1$ to $0$.
The minimum scale and variance for $h$-HDP are
\begin{equation}
b_H(h)=\frac{\Delta}{2r_h},
\qquad
\operatorname{Var}(Z_{b_H})=\frac{\Delta^2}{2r_h^2},
\label{eq:cal-laplace-reverse}
\end{equation}
and the exact profile at this fixed scale is
\begin{equation}
\delta_h^L(\epsilon)=[1-e^{\epsilon/2-r_h}]_+.
\label{eq:cal-laplace-h-profile}
\end{equation}
For $0\leq\epsilon\leq2r_h$, substituting
$(\epsilon,\delta_h^L(\epsilon))$ into
\eqref{eq:cal-laplace-forward} returns $b_H(h)$.
At $\epsilon=2r_h$, this is the tight pure-DP guarantee.
For $\epsilon>2r_h$, the profile is still zero, but the requirement
$(\epsilon,0)$ permits the smaller minimum scale $\Delta/\epsilon$;
these looser pure-DP guarantees do not recover the fixed scale.
This distinguishes the Laplace reverse correspondence from the
strictly positive Gaussian profile at every finite $\epsilon$.

\subsection{Gaussian forward and reverse calibration notes}
In particular, at $\epsilon=0$ the exact parameter is
\begin{equation}
\delta_h^G(0)=2\Phi(\mu_h/2)-1,
\qquad
\sigma_{\mathrm{BW}}^2(0,\delta)
=\frac{\Delta^2}
{4[\Phi^{-1}\{(1+\delta)/2\}]^2}.
\label{eq:cal-gaussian-zero-epsilon}
\end{equation}
The corresponding worst-case $L_1$ distance is $2\delta_h^G(0)$.
At any finite $\epsilon$, a nonconstant Gaussian query with finite
positive variance has $\delta_G(\epsilon;\Delta/\sigma)>0$, so it
does not satisfy pure $\epsilon$-DP.
The reverse correspondence is not a single pair: it is the curve defining
the Gaussian reverse-calibration boundary in Corollary 22. Pairs above this curve are valid but
looser descriptions of the same release. The connection with the
Gaussian testing profile is given in \citet{Dong2022}.
If only an upper bound on $\Delta$ is used, these formulas provide
certified levels; pairwise attained values use the actual separation.


\subsection{Additional comparisons for Remark 24}
\noindent\textbf{Additional details for Remark 24 (Divergence scale and noise calibration).}
The dependence of a noise requirement on the divergence scale is not
specific to Gaussian mechanisms. For a fixed mechanism family, let
$d_j(v)$ be its worst-case discrepancy under criterion $j$ at variance
$v$. If $d_1(v)\leq d_2(v)$ throughout the relevant range, then
$\{v:d_2(v)\leq t\}\subseteq\{v:d_1(v)\leq t\}$.
Thus the more compressed scale requires no larger minimum variance
at the same numerical threshold. This elementary constraint ordering
is different from comparison at corresponding operational guarantees.

For example, when $\Delta=1$ and $h=0.5$,
\[
\sigma_H^2(h)=\frac1{8\log(4/3)}=0.434507\ldots,
\qquad \delta_h^G(0)=0.551865\ldots.
\]
These are attained parameters of the same Gaussian release.
Imposing TV at most $0.5$ instead gives minimum variance
$1/\{4[\Phi^{-1}(0.75)]^2\}=0.549527\ldots$.
At the corresponding TV value $0.551865\ldots$, the exact TV
calibration returns $0.434507\ldots$.
The difference in the numerical coordinates is therefore distinct from
a difference in protection at matched parameters.
The general scale effect is also visible in the scalar Laplace family,
where eliminating $r$ gives
$h=2\{\eta+(1-\eta)\log(1-\eta)\}$ for
$\eta=\operatorname{TV}\in[0,1)$.
\par\medskip
\subsection{Gaussian calibration under composition}
For a fixed neighboring pair, consider independent Gaussian releases
with mean separations $d_j$ and noise variances $\sigma_j^2$.
Their joint log-likelihood ratio is the sum of the individual ones.
Writing $\mu_j=d_j/\sigma_j$ and
$\mu_{1:K}^2=\sum_{j=1}^K\mu_j^2$, this sum has distribution
$N(\mu_{1:K}^2/2,\mu_{1:K}^2)$ under the first joint law and
$N(-\mu_{1:K}^2/2,\mu_{1:K}^2)$ under the second. Hence
\begin{equation}
D_{HD}(P_{1,1:K},P_{2,1:K})
=2\left\{1-\exp\left(-\frac18\sum_{j=1}^K\mu_j^2\right)\right\}
=2\left\{1-\prod_{j=1}^K(1-h_j/2)\right\},
\label{eq:cal-gaussian-product}
\end{equation}
where $h_j=2(1-e^{-\mu_j^2/8})$.
The joint profile is exactly $\delta_G(\epsilon;\mu_{1:K})$.
These are pairwise identities. For a mechanism-wide product identity,
the supremum is over $\sum_j d_j(D,D')^2/\sigma_j^2$; separate
sensitivity suprema need not be attained by the same neighboring pair.

For adaptive releases, suppose that at every fixed preceding transcript
$z$, the conditional Gaussian mechanism at step $j$ has standardized
separation at most a deterministic $\bar\mu_j$ for every neighboring
pair. Its conditional affinity is at least
$e^{-\bar\mu_j^2/8}$.
If $R_{1,j},R_{2,j}$ denote the laws of the first $j$ releases, successive
integration gives
\[
\rho(R_{1,j},R_{2,j})
\geq e^{-\bar\mu_j^2/8}\rho(R_{1,j-1},R_{2,j-1}).
\]
Thus the transcript is $\bar h$-HDP, where
\begin{equation}
\bar h=2\left\{1-\exp\left(-\frac18\sum_{j=1}^K\bar\mu_j^2\right)\right\}.
\label{eq:cal-adaptive-hdp}
\end{equation}
The GDP composition theorem gives parameter
$\bar\mu_{1:K}=(\sum_j\bar\mu_j^2)^{1/2}$ and hence the
approximate-DP guarantees
$(\epsilon,\delta_G(\epsilon;\bar\mu_{1:K}))$
\citep{Dong2022}. These are transcript guarantees in the adaptive case;
the final iterate need not be Gaussian. Releasing only that iterate
is post-processing. This is the form used for the iterative estimation
procedures below.

To compare equal allocations, fix a reference Gaussian mechanism with
Hellinger level $h_0\in(0,2)$, variance
$\sigma_0^2=\sigma_H^2(h_0)$, and parameter $\mu_0=\Delta/\sigma_0$.
For $K$ steps with the same sensitivity $\Delta$, HDP and GDP
allocations are respectively
\begin{equation}
h_{\mathrm{step}}=2\{1-(1-h_0/2)^{1/K}\},
\qquad
\mu_{\mathrm{step}}=\mu_0/\sqrt K.
\label{eq:cal-equal-allocation}
\end{equation}
Both give the same per-step Gaussian variance:
\begin{equation}
\sigma_{H,K}^2=\sigma_{\mathrm{GDP},K}^2=K\sigma_0^2.
\label{eq:supp-cal-equal-variance}
\end{equation}
Indeed, substitute \eqref{eq:cal-equal-allocation} into
the reverse Gaussian calibration in Corollary 22, or use
$\Delta^2/\mu_{\mathrm{step}}^2$.

For comparison with basic approximate-DP composition, fix
$\epsilon_0\geq0$ and let
$\delta_0=\delta_G(\epsilon_0;\mu_0)$.
Equal allocation uses $(\epsilon_0/K,\delta_0/K)$ at each step, with
per-step variance
$\sigma_{\mathrm{BW}}^2(\epsilon_0/K,\delta_0/K)$.
For the TV rule $\eta_{1:K}\leq1-\prod_j(1-\eta_j)$, let
$\eta_0=2\Phi(\mu_0/2)-1$ and use
\begin{equation}
\eta_{\mathrm{step}}=1-(1-\eta_0)^{1/K},
\qquad
\sigma_{\mathrm{TV},K}^2
=\frac{\Delta^2}
{4[\Phi^{-1}\{(1+\eta_{\mathrm{step}})/2\}]^2}.
\label{eq:cal-tv-allocation}
\end{equation}
The TV composition rule follows by coupling each pair of conditional
releases while their histories agree: the probability that all steps
agree is at least $\prod_j(1-\eta_j)$.

Every per-step Gaussian variance in this comparison is calibrated
exactly to its stated per-step requirement. Basic approximate-DP and
the displayed TV rule provide the specified composition guarantees;
they need not equal the attained profile of the joint Gaussian laws.
For independent releases with unit separation and a common per-step
variance $v_K$, that attained Hellinger level is
$2\{1-e^{-K/(8v_K)}\}$ by \eqref{eq:cal-gaussian-product}.
Thus the composition comparison can report both the parameters delivered
by the stated rule and the attained quantities when the joint law is
available. For $h_0=0.5$, the reference constants are, approximately,
$\sigma_0^2=0.434507$, $\mu_0=1.517055$, and
$\eta_0=0.551865$; choosing $\epsilon_0=1$ gives
$\delta_0=0.327168$. The smallest GDP parameter of the scalar Laplace mechanism is
$2\Phi^{-1}\{1-\tfrac12e^{-\Delta/(2b)}\}$; the full-profile
proof and the inverse calibration in the noise scale and Hellinger
level are given in Appendix~\ref{app:laplace-gdp}.


\subsection{Proof of Proposition 30
}

To show $H_n(\bta)$ is $\alpha-$Lipschitz continuous, it is enough to show that the $(i,j)^{th}$ component of $H_n(\bta)$ is Lipschitz continuous for all $(i,j)$, where the Lipschitz constant depends only on $m$ and the upper bounds in assumptions {\bf(U1)} and {\bf(U2)}.
That is, we will show that for any $\bta^{(1)},\bta^{(2)}\in \Ta$,
\begin{align*}
    |H_{n,i,j}(\bta^{(1)})-H_{n,i,j}(\bta^{(2)})|\leq \alpha_{i,j} ||\bta^{(1)}-\bta^{(2)}||_2,
\end{align*}
where $H_{n,i,j}(\bta)$ is the $(i,j)^{th}$ component of $H_n(\bta)$.  Recall that
\begin{align*}
    H_{n,i,j}(\bta) 
    =& -\int_{\mathbb{R}}g_n^{\frac{1}{2}}(x)f_{\bta}^{\frac{1}{2}}(x) u_{\bta,j}(x)u_{\theta,i}(x)\mathrm{d}x -2\int_{\mathbb{R}}g_n^{\frac{1}{2}}(x)f_{\bta}^{\frac{1}{2}}(x) u_{\bta,i,j}(x)\mathrm{d}x.
\end{align*}
Then for any $\bta^{(1)},\bta^{(2)}\in \Ta$,
\begin{eqnarray}\label{eq:prop_Lip_1}
    |H_{n,i,j}(\bta^{(1)})-H_{n,i,j}(\bta^{(2)})| &\leq& 
    \int_{\mathbb{R}}g_n^{\frac{1}{2}}(x) |T_{1,i,j}(\bta^{(1)},x)-T_{1,i,j}(\bta^{(2)},x)|\mathrm{d}x \nonumber\\
    & & + 2\int_{\mathbb{R}}g_n^{\frac{1}{2}}(x) |T_{2,i,j}(\bta^{(1)},x)-T_{2,i,j}(\bta^{(2)},x)|\mathrm{d}x,
\end{eqnarray}
where
\begin{align*}
    T_{1,i,j}(\bta,x)=f_{\theta}^{\frac{1}{2}}(x) u_{\bta,j}(x)u_{\bta,i}(x),\quad
    T_{2,i,j}(\bta,x)=f_{\bta}^{\frac{1}{2}}(x) u_{\bta,i,j}(x).
\end{align*}
Notice that for $\bta\in\Ta$, $T_{1,i,j}(\bta,x)$ and $T_{2,i,j}(\bta,x)$ are differentiable in $\bta$ by Assumption {\bf(U2)}. 
Using Cauchy-Schwarz inequality and the upper bounds in Assumption {\bf(U2)}, for $\bta\in\Ta$ we obtain that $g_n^{1/2}(x)||\nabla T_{1,i,j}(\bta,x)||_2$ and $g_n^{1/2}(x)||\nabla T_{2,i,j}(\bta,x)||_2$ are integrable with respect to $x$, where the gradient is taken with respect to $\bta$. Set $\bta_t=\bta^{(2)}+t(\bta^{(1)}-\bta^{(2)})$.
By the convexity of $\Ta$ (see Assumption {\bf\ref{asp:A1}}),
$\bta_t\in\Ta$ for $0\leq t\leq1$. The fundamental theorem of
calculus and Cauchy--Schwarz give
\begin{align}
|T_{1,i,j}(\bta^{(1)},x)-T_{1,i,j}(\bta^{(2)},x)|
&\leq \|\bta^{(1)}-\bta^{(2)}\|_2
\int_0^1\|\nabla T_{1,i,j}(\bta_t,x)\|_2\,\mathrm{d}t,
\label{eq:prop_Lip_2}\\
|T_{2,i,j}(\bta^{(1)},x)-T_{2,i,j}(\bta^{(2)},x)|
&\leq \|\bta^{(1)}-\bta^{(2)}\|_2
\int_0^1\|\nabla T_{2,i,j}(\bta_t,x)\|_2\,\mathrm{d}t.
\label{eq:prop_Lip_3}
\end{align}
For $\ell=1,2$, Cauchy--Schwarz and
$\int_{\mathbb R}g_n(x)\,\mathrm{d}x=1$ yield
\[
\int_{\mathbb R}g_n^{1/2}(x)
\|\nabla T_{\ell,i,j}(\bta,x)\|_2\,\mathrm{d}x
\leq
\left(\int_{\mathbb R}
\|\nabla T_{\ell,i,j}(\bta,x)\|_2^2\,\mathrm{d}x\right)^{1/2}.
\]
The right-hand side is bounded uniformly in $\bta\in\Ta$
by the moment conditions in {\bf(U1)}--{\bf(U2)} and
compactness of $\Ta$. Thus, setting
\begin{align*}
\alpha_{i,j}=\sup_{\bta\in\Ta}\biggl\{
&\int_{\mathbb R}g_n^{1/2}(x)
\|\nabla T_{1,i,j}(\bta,x)\|_2\,\mathrm{d}x\\
&+2\int_{\mathbb R}g_n^{1/2}(x)
\|\nabla T_{2,i,j}(\bta,x)\|_2\,\mathrm{d}x\biggr\},
\end{align*}
we obtain constants bounded uniformly in $n$.
Multiplying (\ref{eq:prop_Lip_2}) and (\ref{eq:prop_Lip_3})
by $g_n^{1/2}(x)$, integrating, and applying Tonelli's theorem
in (\ref{eq:prop_Lip_1}), we obtain
\begin{align*}
|H_{n,i,j}(\bta^{(1)})-H_{n,i,j}(\bta^{(2)})|
\leq \alpha_{i,j}\|\bta^{(1)}-\bta^{(2)}\|_2.
\end{align*}
This completes the proof. \QED

Before we prove the theorem, we recall that
\begin{eqnarray*}
    g_n(x)= \frac{1}{nc_n}\sum_{i=1}^nK\left(\frac{x-X_i}{c_n}\right).
\end{eqnarray*}
and for the neighboring i.i.d. observations $\{X'_1,X_2,\cdots, X_n\}$, the corresponding density estimator is 
\begin{eqnarray*}
    \tilde{g}_n(x)= \frac{1}{nc_n}\sum_{i=2}^nK\left(\frac{x-X_i}{c_n}\right)+\frac{1}{nc_n}K\left(\frac{x-X'_1}{c_n}\right).
\end{eqnarray*}
The corresponding loss functions are given by 
\begin{eqnarray*}
    L_n(\bta)=2HD^2(g_n,f_{\bta})  ~~\text{and}~~   \tilde L_n(\bta) = 2HD^2(\tilde g_n,f_{\bta}),
\end{eqnarray*}
and the Hessian of the loss functions are given by $H_n(\bta)$ and $\tilde H_n(\bta)$.

\subsection{Proof of Proposition 31
and Theorem 35
}
We begin with the proof of Proposition 31.  First, notice that for all $1 \le i \le m$
\begin{eqnarray*}
    \frac{\partial}{\partial\theta_i}L_n(\bta) = -2\int_{\mathbb{R}}g_n^{\frac{1}{2}}(x)f_{\bta}^{\frac{1}{2}}(x) u_{\bta,i}(x)\mathrm{d}x
\end{eqnarray*}
Hence, 
\begin{align*}
    \bar\Delta^{(i)} \coloneqq \frac{\partial}{\partial\theta_i}(L_n(\bta)-\tilde L_n(\bta)) =& \AV{-}2\int_{\mathbb{R}}(g_n^{\frac{1}{2}}(x)-\tilde g_n^{\frac{1}{2}}(x))f_{\bta}^{\frac{1}{2}}(x) u_{\bta,i}(x)\mathrm{d}x
\end{align*}
where we have suppressed $n$ in the notation $\bar{\Delta}^{(i)}$. Using Cauchy-Schwarz inequality, the $HD^2(g_n, \tilde g_n)$ is bounded above by $||g_n-\tilde g_n||_1$ and, assumption {\bf\ref{asp:A2}}  it follows that $$\AV{|\bar\Delta^{(i)}|}\leq 2 HD(g_n,\tilde g_n) \cdot \left[ \int_{\mathbb{R}} f_{\bta}(x) u_{\bta,i}^2(x)\mathrm{d}x\right]^{1/2} \le C_i||g_n-\tilde g_n||_1^{\frac{1}{2}}.$$
Now, 
\begin{align}\label{eq:sensitivity_3}
    ||g_n-\tilde g_n||_1^{1/2} = \left[ \frac{1}{n\cdot c_n}\int_{\mathbb{R}}|K\left(\frac{x-X_1}{c_n}\right)-K\left(\frac{x-X_1'}{c_n}\right)| \mathrm{d}x \right]^{1/2} \leq \left( \frac{2}{n} \right)^{1/2}.
\end{align}
Hence $\AV{|\bar\Delta^{(i)}|} \leq 2C_i \left( \frac{2}{n} \right)^{1/2}$. Let $\mathbf{\bar\Delta}=[\bar\Delta^{(1)},\cdots,\bar\Delta^{(m)}]$. Then $\Delta_{L_1}(\nabla L_n(\bta))\leq \AV{\sup_{D\sim D^{\prime}}}||\mathbf{\bar\Delta}||_1 \leq  C\cdot m\cdot n^{-\frac{1}{2}}$ and $\Delta_{L_2}(\nabla L_n(\bta))\leq \AV{\sup_{D\sim D^{\prime}}}||\mathbf{\bar\Delta}||_2 \leq  C\cdot \sqrt{m}\cdot n^{-\frac{1}{2}}$, where $C=\max_{i}\{2\sqrt{2}C_i\}$.

We now turn to the Hessian. Recall the definition $u_{\bta,i}(x)=\frac{1}{f_{\bta}(x)}\cdot \frac{\partial}{\partial\theta_i} f_{\bta}(x)$, $u_{\bta,i,j}=\frac{\partial}{\partial\theta_j}u_{\bta,i}$, and
\begin{align*}
    H_{n,i,j}(\bta)
    = -\int_{\mathbb{R}}g_n^{\frac{1}{2}}(x)f_{\bta}^{\frac{1}{2}}(x) [u_{\bta,j}(x)u_{\bta,i}(x)+2 u_{\bta,i,j}(x)]\mathrm{d}x.
\end{align*}
Hence,
\begin{align*}
    \bar\Delta^{(i,j)} \coloneqq H_{n,i,j}(\theta)- \tilde H_{n,i,j}(\theta) = \AV{-}\int_{\mathbb{R}} (g_n^{\frac{1}{2}}(x)-{\tilde g_n}^{\frac{1}{2}}(x))f_{\bta}^{\frac{1}{2}}(x) [u_{\bta,j}(x)u_{\bta,i}(x)+2 u_{\bta,i,j}(x)]\mathrm{d}x
\end{align*}
Using Cauchy-Schwarz inequality, the $HD^2(g_n, \tilde g_n)$ is bounded above by $||g_n-\tilde g_n||_1$ and, assumptions {\bf(U1)}-{\bf(U2)} it follows that 
\begin{align*}
    |\bar\Delta^{(i,j)}|\leq& 2 HD(g_n,\tilde g_n) \cdot \left[ \int_{\mathbb{R}} \AV{f_{\bta}(x)} [u_{\bta,j}(x)u_{\bta,i}(x)+2 u_{\bta,i,j}(x)]^{\AV{2}}\AV{\mathrm{d}x}\right]^{1/2}\\
    \le& (C_{i,1}+C_{i,2})||g_n-\tilde g_n||_1^{\frac{1}{2}},
\end{align*}
where $C_{i,1}$ and $C_{i,2}$ are upper bounds given in assumptions {\bf(U1)} and {\bf(U2)}. Using (\ref{eq:sensitivity_3}), it follows that $\AV{|\bar\Delta^{(i,j)}|}\leq C_{i,j}\cdot n^{-\frac{1}{2}}$. Let $\mathbf{\bar\Delta}$ be a $m\times m$ matrix with $(i,j)^{\text{th}}$ element $\bar\Delta^{(i,j)}$. Then $\Delta_{L_1}(H_n(\bta))\leq \AV{\sup_{D\sim D^{\prime}}}||\mathbf{\bar\Delta}||_1 \leq  C\cdot m\cdot n^{-\frac{1}{2}}$ and $\Delta_{L_2}(H_n(\bta))\leq \AV{\sup_{D\sim D^{\prime}}}||\mathbf{\bar\Delta}||_2 \leq  C\cdot m\cdot n^{-\frac{1}{2}}$ for some $0 <C < \ff$.
This completes the proof of Proposition 31.

We next turn to the Proof of Theorem 35.
We first construct the interior region used in assumption {\bf(U3)}. Put
\[
\begin{split}
V_{\bta}(x)=f_{\bta}(x)\biggl\{&\sum_{i=1}^m u_{\bta,i}(x)^2
+2\sum_{i,j=1}^m u_{\bta,i}(x)^2u_{\bta,j}(x)^2\\
&+8\sum_{i,j=1}^m u_{\bta,i,j}(x)^2\biggr\},
\qquad
T(a)=\sup_{\bta\in\Theta}\int_{|x|>a}V_{\bta}(x)\,\mathrm{d}x.
\end{split}
\]
By {\bf(U1)}--{\bf(U2)}, $V_{\bta}$ is integrable and its
integral is continuous in $\bta$. Pointwise continuity and
Scheff\'e's lemma imply $L_1$-continuity in $\bta$; compactness of
$\Theta$ gives
\[
T(a)\longrightarrow0\qquad(a\longrightarrow\infty).
\]
For $p\in(1,2)$, choose $a_n\uparrow\infty$ with
$T(a_n)\le n^{1-2/p}$, then choose $b_n>a_n+\beta c_n$ and set
\[
A_n=[-a_n,a_n]\subset S_n,
\qquad S_n=(-b_n-\beta c_n,b_n+\beta c_n).
\]
The construction is independent of the dataset. The bound
$\|u_{\bta}u_{\bta}^T+2\dot u_{\bta}\|_F^2
\le2\|u_{\bta}u_{\bta}^T\|_F^2+8\|\dot u_{\bta}\|_F^2$
controls both exterior weights. The KDE lower bound in {\bf(U3)}
is imposed only for the given dataset on this $A_n$; also,
$0<\delta_n\le(2a_n)^{-1}\to0$ since $\int\bar g_n\le1$.

To reduce notational complexity, \emph{redefine}
$\tilde{g}_n(x)$ as follows:
\begin{eqnarray*}
    \tilde{g}_n(x)= \frac{1}{nc_n}\sum_{i=2}^nK\left(\frac{x-X_i}{c_n}\right)\bm{1}_{(X_i\in B_n)}+\frac{1}{nc_n}K\left(\frac{x-X'_1}{c_n}\right)\bm{1}_{(X_1'\in B_n)},
\end{eqnarray*}
where \(B_n=(-b_n,b_n)\), as defined in the main paper. Let $\bar\Delta^{(i)}=\bar\Delta^{(i,1)}+\bar\Delta^{(i,2)}$, where
\begin{align*}
\bar\Delta^{(i,1)}
&=-2\int_{A_n}(\sqrt{\bar g_n(x)}-\sqrt{\tilde g_n(x)})
f_{\bta}^{1/2}(x)u_{\bta,i}(x)\,\mathrm{d}x,\\
\bar\Delta^{(i,2)}
&=-2\int_{\mathbb R\setminus A_n}
(\sqrt{\bar g_n(x)}-\sqrt{\tilde g_n(x)})
f_{\bta}^{1/2}(x)u_{\bta,i}(x)\,\mathrm{d}x.
\end{align*}
For the given dataset satisfying {\bf(U3)},
$\bar g_n(x)\geq\delta_n>0$ on $A_n$, and hence
\[
|\sqrt{\bar g_n(x)}-\sqrt{\tilde g_n(x)}|
=\frac{|\bar g_n(x)-\tilde g_n(x)|}
{\sqrt{\bar g_n(x)}+\sqrt{\tilde g_n(x)}}
\leq\frac{|\bar g_n(x)-\tilde g_n(x)|}{\sqrt{\delta_n}}.
\]
Next, the triangle inequality and a
change of variables give, for every one-record replacement,
\[
\|\bar g_n-\tilde g_n\|_{L_p(\mathbb R)}
\leq\frac{2\|K\|_{L_p(\mathbb R)}}{n c_n^{1-1/p}}.
\]
Set $R_{\bta}(x)=|f_{\bta}^{1/2}(x)u_{\bta,i}(x)|$ and
$C_R=\sup_{\bta\in\Theta}\|R_{\bta}\|_{L_q(\mathbb R)}<\infty$,
where finiteness follows from {\bf(U3)} and compactness of $\Theta$.
H\"older's inequality and the condition on $\delta_n$ yield
\begin{align*}
|\bar\Delta^{(i,1)}|
&\leq\frac{2C_R}{\sqrt{\delta_n}}
\|\bar g_n-\tilde g_n\|_{L_p(\mathbb R)}\\
&\leq\frac{4C_R\|K\|_{L_p(\mathbb R)}}
{n c_n^{1-1/p}\sqrt{\delta_n}}
\leq4C_R C_\delta\|K\|_{L_p(\mathbb R)}\,n^{-1/p}.
\end{align*}
This bound is uniform in $\bta$ and over the given dataset's
one-record neighbors. The same argument applies to the Hessian
entries using $|f_{\bta}^{1/2}(u_{\bta,i}u_{\bta,j}+2u_{\bta,i,j})|$
in place of $R_{\bta}$; its $L_q$ norm is also uniformly bounded by {\bf(U3)}.
Next, we turn to $\bar\Delta^{(i,2)}$. The two truncated KDEs vanish outside $S_n$, so only $S_n\setminus A_n$ contributes. By Cauchy--Schwarz and the inequality $(\sqrt a-\sqrt b)^2\leq|a-b|$,
\begin{align*}
|\bar\Delta^{(i,2)}|^2
&\leq4\left[\int_{S_n\setminus A_n}(\sqrt{\bar g_n(x)}-\sqrt{\tilde g_n(x)})^2\,\mathrm{d}x\right]
\left[\int_{S_n\setminus A_n}f_{\bta}(x)u_{\bta,i}(x)^2\,\mathrm{d}x\right]\\
&\leq\frac8n\int_{S_n\setminus A_n}f_{\bta}(x)u_{\bta,i}(x)^2\,\mathrm{d}x,
\end{align*}
where $\|\bar g_n-\tilde g_n\|_1\leq2/n$ for every neighboring pair. For the
region constructed above, $S_n\setminus A_n\subset\{|x|>a_n\}$.
Consequently, by the definition of $V_{\bta}$ and the choice of $a_n$,
\[
\sup_{\bta\in\Theta}\max_i
\int_{S_n\setminus A_n}f_{\bta}(x)u_{\bta,i}(x)^2\,\mathrm{d}x
\le T(a_n)\le n^{1-2/p}.
\]
The same bound applies to every Hessian exterior weight
$f_{\bta}(u_{\bta,i}u_{\bta,j}+2u_{\bta,i,j})^2$, since it is
bounded by $V_{\bta}$.
It follows that $|\bar\Delta^{(i,2)}|\leq\sqrt{8}\,n^{-1/p}$. Combining the two parts, for $i=1,\ldots,m$,
\[
|\bar\Delta^{(i)}|\leq C_6 n^{-1/p}.
\]
Setting $\mathbf{\bar\Delta}=[\bar\Delta^{(1)},\ldots,\bar\Delta^{(m)}]^T$, we therefore have, for the neighboring pair under consideration,
\[
\|\mathbf{\bar\Delta}\|_1\leq mC_6n^{-1/p},\qquad
\|\mathbf{\bar\Delta}\|_2\leq\sqrt m\,C_6n^{-1/p}.
\]
The same argument gives the Hessian bound, completing the proof of
Theorem 35. \QED

\paragraph{A random-data extension.}
The construction is deterministic, but {\bf(U3)} is conditional on
the sample. Suppose $K(u)\ge\kappa>0$ on $|u|\le\rho$ and cover
$A_n$ by $M_n$ intervals $I_{n,j}\subset B_n$ of length $\rho c_n$ with
\[
\min_{j\le M_n}\mathbf P_g\{X\in I_{n,j}\}
\ge c_0c_n(nc_n)^{-2(1-1/p)}.
\]
Set $v_n=(nc_n)^{2/p-1}\to\infty$. Each interval has expected count
at least $c_0v_n$. Chernoff's inequality and a union bound imply that,
with probability at least $1-M_n\exp(-c_0v_n/8)$, all counts are at
least $c_0v_n/2$. A one-record replacement leaves at least $c_0v_n/4$
in each interval for large $n$. Since each $x\in A_n$ belongs to a
covering interval, $K\ge\kappa$ yields, simultaneously for the sample
and all its one-record replacements,
\[
\inf_{x\in A_n}\bar g_n(x)\ \wedge\
\inf_{x\in A_n}\tilde g_n(x)
\ \ge\ \frac{\kappa c_0}{4}\frac{v_n}{nc_n}
=\frac{\kappa c_0}{4}(nc_n)^{-2(1-1/p)}.
\]
Thus the interior rate requirement in {\bf(U3)} holds on an event
$\mathcal E_n$ with
\[
\mathbf P_g(\mathcal E_n^c)
\le M_n\exp\{-c_0(nc_n)^{2/p-1}/8\}.
\]
If $\log M_n=o((nc_n)^{2/p-1})$, then
$\mathbf P_g(\mathcal E_n)\to1$. The $n^{-1/p}$ bound thus holds
conditionally and uniformly over the observed dataset's one-record
neighbors, as in random differential privacy
\citep{HallRinaldoWasserman2013}, not for every unbounded dataset.
The local-mass and covering conditions require verification for the
chosen family and sequences.

\subsection{Proof of Proposition 34
}

For PGD, recalling that the mechanism $M_k(w,D)=w - \eta\left(\nabla L_n(w)+\Delta_n\cdot c_{\ep'} Z_k\right)$. $M_k$ satisfies $\ep'-$HDP by Proposition 14
and the post processing property. Now starting with the initial estimate $w=\hat\bta_n^{(0)}$,  we obtain $\hat\bta_n^{(k)}$, for $k \ge 1$ using the iteration
\begin{align*}
    \hat\bta_n^{(k)} = M_k(\hat\bta_n^{(k-1)},D) = \hat\bta_n^{(k-1)} - \eta\left(\nabla L_n(\hat\bta_n^{(k-1)})+\Delta_n\cdot c_{\ep'} Z_k\right).
\end{align*}
Hence, by Corollary 11,
$\hat\bta_n^{(K)}$ satisfies $h_K(\epsilon')-$HDP. Finally, by the choice of $\epsilon'$ satisfying $h_K(\epsilon')=\epsilon$, it follows that $\hat\bta_n^{(K)}$ satisfies $\epsilon-$HDP. Next, turning to PNR, the mechanism is $M_k(w,D)=w -\left(H_n(w)+W_{n,k}\right)^{-1} \left(\nabla L_n(w)+N_{n,k}\right)$, where $W_{n,k}\in\mathbb{R}^{m\times m}$ and $N_{n,k}\in\mathbb{R}^{m\times1}$ are the independent random variables added to satisfy the HDP property and 
\begin{align*}
    N_{n,k} = \Delta_n\cdot c_{\ep'/2}\cdot Z_{k}, \quad
    W_{n,k} = \Delta_n^{(H)}\cdot c_{\ep'/2}\cdot \Tilde Z_{k}.
\end{align*}
Let $M_{k,1}(w,D)=\left(H_n(w)+W_{n,k}\right)^{-1}$, $M_{k,2}(w,D)=\left(\nabla L_n(w)+N_{n,k}\right)$. Then $M_{k,1}$ and $M_{k,2}$ satisfies $\frac{\ep'}{2}-$HDP by Proposition 14,
Proposition 15,
and the post processing property. Hence, by Corollary 11,
it follows that $M_k(w,D)=w-M_{k,1}(w,D)\cdot M_{k,2}(w,D)$ satisfies $\epsilon'-$HDP. Finally, starting with the initial estimate $w=\hat\bta_n^{(0)}$, we obtain $\hat\bta_n^{(k)}$  for $ k \ge 1$ by iterating
\begin{align*}
    \hat\bta_n^{(k)} = M_k(\hat\bta_n^{(k-1)},D) = \hat\bta_n^{(k-1)} - M_{k,1}(\hat\bta_n^{(k-1)},D)\cdot M_{k,\AV{2}}(\hat\bta_n^{(k-1)},D)
\end{align*}
Hence, by Corollary 11,
$\hat\bta_n^{(K)}$ satisfies $h_K(\epsilon')-$HDP. Also, by the choice of $\epsilon'$ satisfying $h_K(\epsilon')=\epsilon$, it follows that $\hat\bta_n^{(K)}$ satisfies $\epsilon-$HDP. Finally, to obtain the bounds for $h_K(\ep K^{-1})$, first notice that $h_2(\ep K^{-1})=h_1(\ep K^{-1})+\epsilon'-\frac{1}{2}h_1(\ep K^{-1})\epsilon' \leq h_1(\ep K^{-1})+\epsilon'=2\ep K^{-1}$. Iterating, it follows that $h_K(\ep K^{-1})\leq \ep$. Now, turning to the lower bound, by iterating Corollary 11
we obtain
\begin{align*}
    h_K(\ep K^{-1}) =& \epsilon' + \epsilon' \left[(K-1)-\frac{1}{2}\sum_{j=1}^{K-1}h_j(\ep K^{-1})\right] = K\epsilon' - \epsilon'\frac{1}{2}\sum_{j=1}^{K-1}h_j(\ep K^{-1}).
\end{align*}
Next, using the upper bound, $h_j(\ep K^{-1})\leq \epsilon\cdot j K^{-1}$, we obtain
\begin{align*}
    h_K(\ep K^{-1})\geq  \ep[1-\ep(K-1)(4K)^{-1}].    \quad \QED
\end{align*}

\subsection{Proof of Theorem 36
}
The proof of the theorem relies on the behavior of the Hellinger loss function at \emph{private} estimates. Intuitively, we show that under ASLSC and $\tau_2-$smoothness, the closeness of the loss functions implies the closeness of the parameter estimates and vice-versa. This is achieved via Lemma \ref{lem:completness}-Lemma \ref{lem:continuity}. We recall that $N$ is defined in Proposition 28. 
In this proof, for the ease of exposition, we set $\tau_1$ and $\tau_2$ to be $2\tau_1$ and $2\tau_2$.

\begin{lemma} \label{lem:completness}
Assume that assumptions {\bf\ref{asp:A1}}-{\bf\ref{asp:A8}} in Section~\ref{app:A} and {\bf(U1)}-{\bf(U2)} hold and that $||\hat\bta_n-\bta_g||_2\leq \frac{1}{2}r$. Also, assume that for all $k=1,\cdots,K$, and $n \ge N$, $||\hat\bta_n-\hat\bta_n^{(k)}||_2\leq \frac{1}{2}r$ with probability $1-\frac{k\xi}{K}$ and $||\hat\bta_n^{(k+1)}-\hat\bta_n^{(k)}||_2\leq r$ with probability $1-\frac{(k+1)\xi}{K}$, where $r$ is as defined in Proposition 28.
Let $N_{n,k} = \Delta_n c_{\ep'}Z_k$, where $Z_k \sim N(\bm{0},\mathbf{I})$. Then, with probability $1-\frac{k\xi}{K}$,
\begin{align} \label{eq:completness}
    L_n(\hat\bta_n^{(k)})-L_n(\hat\bta_n)\leq (1-\gamma)^k (L_n(\hat\bta_n^{(0)})-L_n(\hat\bta_n)) + \frac{3r\cdot ||N_{n,k}||_2}{2\gamma},
\end{align}
where $\eta$ and $\gamma$ are chosen such that $0 < \gamma\leq 2\eta\tau_{1}\le 2\eta\tau_{2}<1$.
\end{lemma}

\noindent\textbf{Proof}: Recall that
\begin{align*}
    Q_k(\bta)=L_n(\hat\bta_n^{(k)}) + \langle\nabla L_n(\hat\bta_n^{(k)})+N_k,\bta-\hat\bta_n^{(k)}\rangle + \frac{1}{2\eta}||\bta-\hat\bta_n^{(k)}||_2^2.
\end{align*}
Since $\hat\bta_n^{(k+1)}$ minimizes $Q_k(\bta)$, it follows that
by setting $\bta_\gamma=\gamma\hat\bta_n + (1-\gamma)\hat\bta_n^{(k)}$, that
\begin{align} \label{eq:Qub1}
Q_k(\hat\bta_n^{(k+1)}) \leq Q_k(\bta_\gamma)=& 
    L_n(\hat\bta_n^{(k)}) + \gamma\langle \nabla L_n(\hat\bta_n^{(k)}),\hat\bta_n-\hat\bta_n^{(k)}\rangle\\
    & + \frac{\gamma^2}{2\eta}||\hat\bta_n-\hat\bta_n^{(k)}||_2^2 + \gamma \langle N_{n,k},\hat\bta_n-\hat\bta_n^{(k)}\rangle. \nonumber
\end{align}
Now using Proposition 29 
part (1) we obtain 
\begin{align*}
    \langle\nabla L_n(\hat\bta_n^{(k)}),\hat\bta_n-\hat\bta_n^{(k)}\rangle \leq& L_n(\hat\bta_n)-L_n(\hat\bta_n^{(k)}) - \tau_{1}||\hat\bta_n-\hat\bta_n^{(k)}||_2^2.
\end{align*}
Now, using this bound in the inequality (\ref{eq:Qub1}), we obtain
\begin{align}{\label{eq:Qub2}}
Q_k(\hat{\bta}_n^{(k+1)}  )  \leq& L_n(\hat\bta_n^{(k)}) - \gamma[L_n(\hat\bta_n^{(k)})-L_n(\hat\bta_n)]  + \left(\frac{\gamma^2}{2\eta}-\gamma\tau_{1}\right)||\hat\bta_n-\hat\bta_n^{(k)}||_2^2 \\
& + \gamma\langle N_{n,k},\hat\bta_n-\hat\bta_n^{(k)}\rangle, \nonumber
\end{align}
yielding the upper bound of $Q_k(\hat{\bta}_n^{(k+1)})$.
We next obtain a lower bound for $Q_k(\hat{\bta}_n^{(k+1)})$. To this end, we use part (3) of Proposition 29.
Specifically, using $ L_n(\hat\bta_n^{(k)}) \geq L_n(\hat\bta_n^{(k+1)})- \langle\nabla L_n(\hat\bta_n^{(k)}),\hat\bta_n^{(k+1)}-\hat\bta_n^{(k)}\rangle - \tau_{2} ||\hat\bta_n^{(k)}-\hat\bta_n^{(k+1)}||_2^2$, we obtain that
\begin{align*}
    Q_k(\hat\bta_n^{(k+1)})
    \geq & L_n(\hat\bta_n^{(k+1)}) + (\frac{1}{2\eta}- \tau_{2})||\hat\bta_n^{(k+1)}-\hat\bta_n^{(k)}||_2^2 + \langle N_{n,k},\hat\bta_n^{(k+1)}-\hat\bta_n^{(k)}\rangle.
\end{align*}
Now since $2\eta\tau_{2}\leq1$, it follows that 
\begin{align} \label{eq:Qlb}
    Q_k(\hat\bta_n^{(k+1)}) \geq& L_n(\hat\bta_n^{(k+1)}) + \langle N_{n,k},\hat\bta_n^{(k+1)}-\hat\bta_n^{(k)}\rangle.
\end{align}
Now using (\ref{eq:Qub2}) and (\ref{eq:Qlb}) it follows that
\begin{align*}
L_n(\hat\bta_n^{(k+1)})-L_n(\hat\bta_n)
\leq& (1-\gamma)(L_n(\hat\bta_n^{(k)})-L_n(\hat\bta_n)) + \left(\frac{\gamma^2}{2\eta}-\gamma\tau_{1}\right)||\hat\bta_n-\hat\bta_n^{(k)}||_2^2 \\
&+ \gamma\langle N_{n,k},\hat\bta_n-\hat\bta_n^{(k)}\rangle- \langle N_{n,k},\hat\bta_n^{(k+1)}-\hat\bta_n^{(k)}\rangle.
\end{align*}
Now choosing $\gamma$ so that $0<\gamma\leq 2\eta\tau_{1}<1$ and applying Cauchy-Schwarz inequality, it follows that
\begin{align}
L_n(\hat\bta_n^{(k+1)})-L_n(\hat\bta_n)
\leq& (1-\gamma)(L_n(\hat\bta_n^{(k)})-L_n(\hat\bta_n)) + ||N_{n,k}||_2(||\hat\bta_n-\hat\bta_n^{(k)}||_2 \nonumber\\
& + ||\hat\bta_n^{(k+1)}-\hat\bta_n^{(k)}||_2) \label{eq:lem_completeness}\\
\leq& (1-\gamma)(L_n(\hat\bta_n^{(k)})-L_n(\hat\bta_n))+ \frac{3 r ||N_{n,k}||_2}{2}, \nonumber
\end{align}
where the last inequality follows from the assumptions $||\hat\bta_n-\hat\bta_n^{(k)}||_2\leq \frac{1}{2}r$ and $||\hat\bta_n^{(k+1)}-\hat\bta_n^{(k)}||_2\leq r$. Now iterating the above inequality, it follows that
\begin{align*}
    L_n(\hat\bta_n^{(k)})-L_n(\hat\bta_n) \leq& (1-\gamma)^k (L_n(\hat\bta_n^{(0)})-L_n(\hat\bta_n)) + \frac{3r\cdot ||N_{n,k}||_2}{2}\cdot\frac{(1-(1-\gamma)^k)}{\gamma}\\
    \leq& (1-\gamma)^k (L_n(\hat\bta_n^{(0)})-L_n(\hat\bta_n)) + \frac{3r\cdot ||N_{n,k}||_2}{2\gamma}.~~ \QED
\end{align*}
Our next key result is Lemma \ref{lem:continuity} below, which verifies that under the assumptions of Lemma \ref{lem:continuity} the private and non-private estimators are close for large $n$ and for every iteration $k=0,1, \cdots K$. The proof of this lemma relies on the notion that, under the assumptions in Section~\ref{app:A} and {\bf(U1)}-{\bf(U3)},  if the loss functions are ``close'', then arguments of the loss functions are also ``close''. This is the content of our next lemma and the proof is based on almost sure local strong convexity and is provided in Section~\ref{app:D}.
\begin{lemma} \label{lem:continuity_relation}
Let assumptions {\bf\ref{asp:A1}}-{\bf\ref{asp:A8}} in Section~\ref{app:A} and {\bf(U1)}-{\bf(U2)} hold. Then for $\bta \in B_r(\bta_g)$ and $n \ge N$, if $L_n(\bta)-L_n(\hat\bta)\leq \frac{r^2}{4}\tau_{1} $ then $ ||\bta-\hat\bta||_2\leq\frac{r}{2}$. Furthermore, if $||\bta-\hat\bta||_2\leq\frac{r}{2} $ for $\bta \in B_r(\bta_g)$, then  for $n \ge N$, $  L_n(\bta)-L_n(\hat\bta)\leq \frac{r^2}{4}\tau_{2}.$
\end{lemma}

We next turn to the key result verifying the validity of the conditions in Lemma \ref{lem:completness} above.

\begin{lemma} \label{lem:continuity}
Under assumptions {\bf\ref{asp:A1}}-{\bf\ref{asp:A8}} and {\bf(U1)}-{\bf(U2)}, for $\eta\leq\frac{1}{\tau_{2}}$, assume that for $n \ge N$, $\hat\bta_n\in B_{r/c}(\bta_g)\subset B_{r/2}(\bta_g)$, where $c>2\left(\frac{\tau_2}{\tau_1}\right)^{\frac{1}{2}}$, 
then there exists $\hat\bta_n^{(0)}$, such that $L_n(\hat\bta_n^{(k)})-L_n(\hat\bta_n)\leq \tau_{1}\frac{r^2}{4}$ and $||\hat\bta_n^{(k)}-\hat\bta_n||_2\leq \frac{r}{2}$ hold with probability $1-\frac{k\xi}{K}$ for all $k=0,\cdots,K$.

\end{lemma}
The proof of this lemma is similar to the proof of Lemma 18 in \cite{avella2023}. A mildly different proof is given in Section~\ref{app:D}.

We now turn to the proof of Theorem 36. 

\textbf{Proof of Theorem 36
}:
Using Proposition 34 
with $K$ replaced by $K_n$, it follows that $\hat{\bta}_n^{(K_n)}$ satisfies $\ep-$HDP. We next turn to verification of the utility.
The key idea is to use Proposition 29
(1)  and Lemma  \ref{lem:completness} and iterate until the required bound is reached. Towards this, using
$\nabla L_n(\hat\bta_n)=(0,\cdots,0)$ and taking $\gamma \in (0, 2 \eta \tau_1)$, it follows that
\begin{align*}
||\hat\bta_n^{(k)}-\hat\bta_n||_2^2
    \leq& \frac{(1-\gamma)^k (L_n(\hat\bta_n^{(0)})-L_n(\hat\bta_n))}{\tau_{1}} + \frac{3r\cdot ||N_{n,k}||_2} {2\gamma \tau_{1}}.
\end{align*}
Using concentration inequality for $L_2-$norm of the Gaussian vector (see \cite{Rigollet2023}), namely,
\begin{align*}
    P\left( ||Z_k||_2 \geq 4\sqrt{m}+2\sqrt{2[\log K-\log \xi]}\right)\leq \frac{\xi}{K},
\end{align*}
it follows by setting $\ep'= \frac{\ep}{K}$ that
\begin{align} \label{eq:noise_bound}
P(||N_{n,k}||_2  \leq \Delta_n c_{\ep'}\left(4\sqrt{m}+2\sqrt{2[\log K - \log \xi]}\right) \coloneqq  r_{noi}) > 1-\frac{\xi}{K}.
\end{align}
We emphasize here that $r_{noi}$ depends on $n, K_n$ and $\xi$. Now, first consider the case $L_n(\hat\bta_n^{(0)})\neq L_n(\hat\bta_n)$. By choosing $k >k_0 \coloneqq \max \{0,\frac{-\log(L_n(\hat\bta_n^{(0)})-L_n(\hat\bta_n))+\log(3r\cdot r_{noi})-\log (2\gamma)}{\log(1-\gamma)}\}$, it follows that with probability $(1-\frac{ k \xi}{K})$, and for all $k \ge k_0$
\begin{align}\label{eq:GD_step1}
||\hat\bta_n^{(k)}-\hat\bta_n||_2^2
    \leq& \frac{3r\cdot r_{noi}} {\gamma \tau_{1}}\coloneqq C_0^2 r_{noi}.
\end{align}
We notice here that this bound  is of order $\Delta_n^{\frac{1}{2}}(K_n \log K_n)^{\frac{1}{4}}$. However, this will not yield efficiency. Our goal is to remove the square root from $\Delta_n$. This suggests one needs larger $k$ in the above bound. This is accomplished by additional iterations (see Theorem 2 in \cite{avella2023}).
To this end, we need the following claim, whose proof is given below.

\vspace{0.05in}
\noindent{\bf Claim:} For $k>k_0$, choose $n,K$ such that $r_{noi}^{\frac{1}{2}}<\frac{1}{2\eta}C_0$. Then
$$L_n(\hat\bta_n^{(k+1)})-L_n(\hat\bta_n) \le   (1-\gamma)(L_n(\hat\bta_n^{(k)})-L_n(\hat\bta_n)) + (2\eta\tau_{2}+\frac{3}{2})C_0 r_{noi}^{\frac{3}{2}}.$$

Using the claim with  $k=k_0+j-1$ and iterating we get
\begin{align*}
    L_n(\hat\bta_n^{(k_0+j)})-L_n(\hat\bta_n)
    \leq& (1-\gamma)^j(L_n(\hat\bta_n^{(k_0)})-L_n(\hat\bta_n)) + \left[\sum_{i=0}^{j}(1-\gamma)^i\right](2\eta\tau_{2}+\frac{3}{2})C_0r_{noi}^{\frac{3}{2}}\\
    \leq& (1-\gamma)^j(L_n(\hat\bta_n^{(k_0)})-L_n(\hat\bta_n)) + \frac{1}{\gamma}(2\eta\tau_{2}+\frac{3}{2})C_0 r_{noi}^{\frac{3}{2}}.
\end{align*}
Now, using Proposition 29
(1) and utilizing $\nabla L_n(\hat\bta_n)=(0,\cdots,0)$ , it follows that $||\hat\bta_n^{(k_0+j)}-\hat\bta_n||_2^2\leq \frac{L_n(\hat\bta_n^{(k_0+j)})-L_n(\hat\bta_n)}{\tau_{1}}$ and hence
\begin{align*}
||\hat\bta_n^{(k_0+j)}-\hat\bta_n||_2^2\leq \frac{(1-\gamma)^j (L_n(\hat\bta_n^{(k_0)})-L_n(\hat\bta_n))}{\tau_{1}} + (2\eta\tau_{2}+\frac{3}{2})\frac{C_0}{\gamma \tau_{1}} r_{noi}^{\frac{3}{2}}.
\end{align*}
Next, we choose
\begin{align*}
    j\geq k_1\coloneqq \max\left\{0, \frac{-\log(L_n(\hat\bta_n^{(k_0)})-L_n(\hat\bta_n))+\log\left(\frac{1}{2}\frac{C_0}{\gamma \tau_{1}} r_{noi}^{\frac{3}{2}}\right)}{\log(1-\gamma)}\right\},
\end{align*}
and setting $C_1=({\frac{2\eta\tau_{2}+2}{\gamma \tau_{1}}C_0})^{\frac{1}{2}}$, we obtain for $k>k_0+k_1$,
\begin{align*}
    ||\hat\bta_n^{(k)}-\hat\bta_n||_2^2\leq C_1^2 r_{noi}^{\frac{3}{2}}.
\end{align*}
We notice that the power of $r_{noi}$ is now $3/2$ and is still below the required power of 2. Hence, continuing the iterations and using the Claim with starting value $k_0+k_1+\cdots k_i$, we obtain for $k \ge k_0+k_1+\cdots k_i$,
\begin{eqnarray*}
    ||\hat\bta^{(k)}_n-\hat \bta||_2^2 \le C^2_i r_{noi}^{2-\frac{1}{2^i}} \quad \text{if} \quad r_{noi}^{\frac{1}{2^i}}\leq\frac{C_{i-1}}{2\eta},
\end{eqnarray*}
where $C_i=\left(\frac{2\eta\tau_{2}+2}{\gamma \tau_{1}}C_{i-1}\right)^{\frac{1}{2}}=\left(\frac{2\eta\tau_{2}+2}{\gamma \tau_{1}}\right)^{1-\frac{1}{2^i}}\cdot C_0^{\frac{1}{2^i}}$. Finally, taking $i=\log_2(n)$ we get $k>k_0+\cdots+k_{\log_2(n)}$
\begin{eqnarray*}
    ||\hat\bta_n^{(k)}-\hat\bta_n||_2^2\leq C_{\log_2(n)}^2 r_{noi}^{2-\frac{1}{n}},~~\text{if} \quad r_{noi}^{\frac{1}{n}}\leq\frac{C_{\log_2(n)-1}}{2\eta},
\end{eqnarray*}
where $C_{\log_2(n)}=\left(\frac{2\eta\tau_{2}+2}{\gamma \tau_{1}}\right)^{1-\frac{1}{n}}\cdot C_0^{\frac{1}{n}}$. 
Now, letting $ n \ra \ff$, notice that $C_{\log_2(n)}$ converges to $ C_{\ff}(\ga) \coloneqq 2(\eta \tau_2+1)(\gamma \tau_1)^{-1}$. Also, notice that $r_{noi}^{\frac{1}{n}}$ converges to 1. Now choosing $ \gamma \in (0, 2\eta \tau_1)$ and $C_{\ff}>1$ (such a $\gamma$ exists) it follows that
\begin{eqnarray}\label{eq:GD_utility2}
    \limsup_{n \ra \ff} r_{noi}^{-2}||\hat\bta_n^{(k_n)}-\hat\bta_n||_2 = C_{\ff}.
\end{eqnarray}
This requires $K\geq k_0+k_1+\cdots+k_{\log_2(n)}\sim (\log n)\cdot(\log r_{noi})$ which implies $K\geq c\log n$, since $r_{noi}$ is bounded by a constant  by choice of $n$ and $K$. Next, we notice that $r_{noi}=\Delta_n c_{\ep'}\left(4\sqrt{m}+2\sqrt{2[\log K - \log \xi]}\right)$ by Theorem 35
and  $\Delta_n \sim c\cdot n^{-\frac{1}{p}}$ for $p \in (1, 2]$.
Hence, (\ref{eq:GD_utility2}) becomes 
\begin{eqnarray*}
    ||\hat\bta_n^{(K_n)}-\hat\bta_n||_2 \le c\cdot n^{-\frac{1}{p}}({K_n\log(K_n/\xi)})^{\frac{1}{2}},
\end{eqnarray*}
for large $n$ with high probability. Thus, to complete the proof of the Theorem, we now establish the claim.\\
{\bf Proof of the Claim:}
Notice that by Proposition 29
inequality 4 that
\begin{eqnarray}\label{eq:GD_claim}
||\nabla L_n(\hat\bta_n^{(k)})||_2 = ||\nabla L_n(\hat\bta_n^{(k)})-\nabla L_n(\hat\bta_n)||_2 \leq 2\tau_{2}||\hat\bta_n-\hat\bta_n^{(k)}||_2.
\end{eqnarray}
Now,  first using equation (24) in Definition 32 of the main paper
and the expression above and applying (\ref{eq:GD_step1}) it follows that 
\begin{align*}
||\hat\bta_n^{(k+1)}-\hat\bta_n^{(k)}||_2 \leq 2\eta\tau_{2} C_0 r_{noi}^{\frac{1}{2}}+ \eta r_{noi}.
\end{align*}
From the inequality (\ref{eq:lem_completeness}) in the proof of Lemma \ref{lem:completness},  using (\ref{eq:GD_step1}) and from (\ref{eq:GD_claim}) it follows that
\begin{align*}
L_n(\hat\bta_n^{(k+1)})-L_n(\hat\bta_n) \leq&
(1-\gamma)(L_n(\hat\bta_n^{(k)})-L_n(\hat\bta_n)) + (2\eta\tau_{2}+1)C_0 r_{noi}^{\frac{3}{2}} +\eta r_{noi}^2. 
\end{align*}
Next, choosing $n,K$ such that $r_{noi}^{\frac{1}{2}}<\frac{1}{2\eta}C_0$ it follows that
\begin{align*}
L_n(\hat\bta_n^{(k+1)})-L_n(\hat\bta_n)\leq& (1-\gamma)(L_n(\hat\bta_n^{(k)})-L_n(\hat\bta_n)) + (2\eta\tau_{2}+\frac{3}{2})C_0 r_{noi}^{\frac{3}{2}}.
\end{align*}
This completes the proof of the claim and the Theorem. \QED

\subsection{Proof of Theorem 39
}
The proof of the Theorem relies on the Lemma \ref{lem:NR_decomp}-Lemma \ref{lem:NR_continuity}
whose proofs use matrix concentration inequality and is similar to the idea of proof of Theorem 36.
We recall that the concentration inequality for $L_2-$norm of the Gaussian vector and matrix,  (see \cite{Rigollet2023,Tropp2015}), is given by
\begin{align*}
    & P\left(||N_{n,k}||_2 \leq \Delta_n \cdot c_{\ep'/2}\cdot [4\sqrt{m}+2\sqrt{2(\log 2K - \log \xi)}]\right)\geq 1-\frac{\xi}{2K},~~ \text{and}~~ \\
    & P\left(||W_{n,k}||_2 \leq \Delta_n^{(H)}\cdot c_{\ep'/2}\cdot \sqrt{2m\log(4Km/\xi)}\right) \geq 1-\frac{\xi}{2K}.
\end{align*}
We use these upper bounds on the norms with probability $1-\frac{\xi}{K}$ in the following lemmas and proofs. In this proof, for the ease of exposition, we set $\tau_1$ and $\tau_2$ to be $2\tau_1$ and $2\tau_2$, and choose $\eta=1$. Our first lemma provides a useful alternative expression for $\hat\bta_n^{(k+1)}-\hat\bta_n^{(k)}$. 

\begin{lemma}\label{lem:NR_decomp}
\begin{align*}
    \left(H_n(\hat\bta_n^{(k)})+W_{n,k}\right)^{-1} \left(\nabla L_n(\hat\bta_n^{(k)})+N_{n,k}\right)= H_n^{-1}(\hat\bta_n^{(k)})\nabla L_n(\hat\bta_n^{(k)}) + \tilde N_{n,k}.
\end{align*}
Under assumptions {\bf\ref{asp:A1}}-{\bf\ref{asp:A8}} of Section~\ref{app:A} and {\bf(U1)}-{\bf(U2)}, if $\hat\bta_n^{(k)}\in B_r(\bta_g)$, then for large $n$, $||\tilde N_{n,k}||_2 \leq \frac{\left\Vert N_{n,k}\right\Vert_2}{2\tau_{1}} + \frac{B_1\cdot||W_{n,k}||_2}{2\tau_{1}^2} + \frac{||N_{n,k}||_2\cdot||W_{n,k}||_2}{2\tau_{1}^2}$ holds with probability $1-\frac{\xi}{2K}$. Additionally, $\kappa \sim n^{-\frac{1}{p}}(K\log (K/\xi))^{\frac{1}{2}}$, there exists $N_\kappa$ such that for all $n>N_\kappa$ and $k=1,2, \cdots K$  $P(||\tilde N_{n,k}||_2\leq \kappa)>1-\frac{\xi}{K}$.
\end{lemma}

\noindent\textbf{Proof}: Using Neumann series formula, note that
\begin{align*}
    &\left(H_n(\hat\bta_n^{(k)})+W_{n,k}\right)^{-1} \left(\nabla L_n(\hat\bta_n^{(k)})+N_{n,k}\right) \\
    =&  H_n^{-1}(\hat\bta_n^{(k)})\left[\mathbf{I} + \sum_{j=1}^\infty (-W_{n,k} H_n^{-1}(\hat\bta_n^{(k)}))^j\right] \left(\nabla L_n(\hat\bta_n^{(k)})+N_{n,k}\right)\\
    \coloneqq & H_n^{-1}(\hat\bta_n^{(k)})\nabla L_n(\hat\bta_n^{(k)}) + \tilde N_{n,k},
\end{align*}
where
\begin{align*}
    \tilde N_{n,k}=H_n^{-1}(\hat\bta_n^{(k)}) \left\{ N_{n,k} + \left[\sum_{j=1}^\infty (-W_{n,k} H_n^{-1}(\hat\bta_n^{(k)}))^j\right]\left(\nabla L_n(\hat\bta_n^{(k)})+N_{n,k}\right) \right\}.
\end{align*}
Now, applying the properties of matrix norms, Proposition 28,
and Proposition 25,
we obtain 
\begin{align*}
    ||\Tilde N_{n,k}||_2 
    \leq \frac{1}{2\tau_{1}} \cdot \left[ \left\Vert N_{n,k}\right\Vert_2 + \left[\sum_{j=1}^\infty \left(\frac{||W_{n,k}||_2}{2\tau_{1}}\right)^j\right]\left(B_1+||N_{n,k}||_2\right) \right].
\end{align*}
Let $n$ be large enough such that $||W_{n,k}||_2\leq\tau_{1}$ with probability $1-\frac{\xi}{2K}$. Then it follows that
\begin{align*}
    ||\Tilde N_{n,k}||_2 \leq \frac{\left\Vert N_{n,k}\right\Vert_2}{2\tau_{1}} + \frac{B_1\cdot||W_{n,k}||_2}{2\tau_{1}^2} + \frac{||N_{n,k}||_2\cdot||W_{n,k}||_2}{2\tau_{1}^2}.
\end{align*}
Notice that as $n\to\infty$, both $||N_{n,k}||_2$ and $||W_{n,k}||_2$ converge to $0$ in probability at  rate $n^{-\frac{1}{p}}(K\log (K/\xi))^{\frac{1}{2}}$. \QED

\begin{lemma} \label{lem:NR_completeness}
    Under assumptions {\bf\ref{asp:A1}}-{\bf\ref{asp:A8}} and {\bf(U1)}-{\bf(U2)}, if $\hat\bta_n^{(k)}\in B_{r}(\bta_g)$, then\[||\nabla L_n(\hat\bta_n^{(k+1)})||_2\leq \frac{\alpha}{2\tau_{1}^2}||\nabla L_n(\hat\bta_n^{(k)})||_2^2 + C ||\tilde N_{n,k}||_2.\]This bound holds with probability $1-\frac{\xi}{K}$.
\end{lemma}

\noindent\textbf{Proof}: Recall that from the PNR iteration, namely, $\hat\bta_n^{(k+1)}=\hat\bta_n^{(k)}-H_n^{-1}(\hat\bta_n^{(k)})\nabla L_n(\hat\bta_n^{(k)}) + \tilde N_{n,k}$, that $\nabla L_n(\hat\bta_n^{(k)})+H_n(\hat\bta_n^{(k)})\cdot[\hat\bta_n^{(k)}-\hat\bta_n^{(k+1)}-\tilde N_{n,k}]=0$. We now rewrite $||\nabla L_n(\hat\bta_n^{(k+1)})||_2$ as 
\begin{gather*}
    ||\nabla L_n(\hat\bta_n^{(k+1)})||_2 = ||T_1-T_2+H_n(\hat\bta_n^{(k)})\tilde N_{n,K}||_2,~~ \text{where}\\
    T_1= \nabla L_n(\hat\bta_n^{(k+1)})-\nabla L_n(\hat\bta_n^{(k)})~~\text{and}~~ T_2= H_n(\hat\bta_n^{(k)})(\hat\bta_n^{(k+1)}-\hat\bta_n^{(k)}).
\end{gather*}
Notice that $T_1-T_2$ can be written as
\begin{align*}
    T_1-T_2 =& \int_{0}^{1} H_n(\hat\bta_n^{(k)}+t(\hat\bta_n^{(k+1)}-\hat\bta_n^{(k)}))\cdot (\hat\bta_n^{(k+1)}-\hat\bta_n^{(k)})\mathrm{d}t - \int_{0}^{1}H_n(\hat\bta_n^{(k)})(\hat\bta_n^{(k+1)}-\hat\bta_n^{(k)})\mathrm{d}t\\
    =& (\hat\bta_n^{(k+1)}-\hat\bta_n^{(k)}) \int_{0}^{1}\left[H_n(\hat\bta_n^{(k)}+t(\hat\bta_n^{(k+1)}-\hat\bta_n^{(k)}))- H_n(\hat\bta_n^{(k)})\right]\mathrm{d}t.
\end{align*}
Using Proposition 30
(namely the Lipschitz property of the Hessian), it follows that with probability 1,
\begin{align*}
    ||T_1-T_2||_2 \leq ||\hat\bta_n^{(k+1)}-\hat\bta_n^{(k)}||_2 \cdot \int_{0}^{1} \alpha\cdot t\cdot ||\hat\bta_n^{(k+1)}-\hat\bta_n^{(k)}||_2 \mathrm{d}t = \frac{\alpha}{2}||\hat\bta_n^{(k+1)}-\hat\bta_n^{(k)}||_2^2.
\end{align*}
Using the upper bound of Proposition 25,
$\hat\bta_n^{(k+1)}-\hat\bta_n^{(k)}=-H_n^{-1}(\hat\bta_n^{(k)})\nabla L_n(\hat\bta_n^{(k)}) + \tilde N_{n,k}$,  Proposition 28,
and for large $n$ that $||\tilde N_{n,k}||_2\leq 1$ with probability $1-\frac{\xi}{K}$, we obtain
\begin{align*}
    ||\nabla L_n(\hat\bta_n^{(k+1)})||_2\leq& \frac{\alpha}{2}||\hat\bta_n^{(k+1)}-\hat\bta_n^{(k)}||_2^2 + B_2 ||\tilde N_{n,k}||_2 \leq \frac{\alpha}{2\tau_{1}^2}||\nabla L_n(\hat\bta_n^{(k)})||_2^2 + C ||\tilde N_{n,k}||_2,
\end{align*}
where the constant $C\in(0,\infty)$ only depends on $\alpha,\tau_1,B_1,B_2$. \QED

The next lemma concerns the ``distance'' between the private and non-private estimators at every iteration, and the proof is based on induction. The choice of
$\hat\bta_n^{(0)}$, verifies the assumption that the assumptions in Lemma \ref{lem:NR_completeness} hold; that's is, for all $k=1,2 \cdots K$, $\hat\bta_n^{(k)}\in B_{r}(\bta_g)$.

\begin{lemma} \label{lem:NR_continuity}
Under assumption {\bf\ref{asp:A1}}-{\bf\ref{asp:A8}} and {\bf(U1)}-{\bf(U2)}, if $\hat\bta_n\in B_{r/2}(\bta_g)$, and\[||\nabla L_n(\hat\bta_n^{(0)})||_2\leq\min\{\frac{\tau_{1}r}{2},\frac{\tau_{1}^2}{\alpha}\},\]then for $k=0,1,\cdots,K$, $||\hat\bta_n^{(k)}-\hat\bta_n||_2\leq \frac{r}{2}$ holds with probability $1-\frac{k\xi}{K}$.
\end{lemma}

\noindent\textbf{Proof}:
We prove the lemma using the following claim:

\noindent\textbf{Claim}: If $\hat\bta_n\in B_{r/2}(\bta_g)$ and  $||\nabla L_n(\bta)||_2\leq 2\tau_{1}r$, then  $||\bta-\hat\bta_n||_2\leq r$.

First, we finish the proof of the lemma using the Claim and then prove the Claim.
We prove the lemma by induction. First notice that by assumption $||\nabla L_n(\hat\bta_n^{(0)})||_2\leq\min\{\frac{\tau_{1}r}{2},\frac{\tau_{1}^2}{\alpha}\}\leq \tau_{1}r$ and hence from the claim it follows, with $\bta=\hat\bta_n^{(0)}$ and $r$ replaced by $\frac{r}{2}$, that $||\hat\bta_n^{(0)}-\hat\bta_n||_2\leq \frac{r}{2}$. We start the inductive hypothesis with $k=k_0$. That is,  assume for $k=k_0$, $||\nabla L_n(\hat\bta_n^{(k_0)})||_2\leq\min\{\tau_{1}r,\frac{\tau_{1}^2}{\alpha}\}$ and $||\hat\bta_n^{(k_0)}-\hat\bta_n||_2\leq \frac{r}{2}$, and $\hat\bta_n^{(k_0)}\in B_r(\bta_g)$. Also from Lemma \ref{lem:NR_completeness}, and  for large $n$ such that $||\tilde N_{n,k}||_2\leq \min\{\tau_{1}r,\frac{\tau_{1}^2}{\alpha}\}$ with probability $1-\frac{\xi}{K}$,  we obtain
\begin{align*}
    ||\nabla L_n(\hat\bta_n^{(k_0+1)})||_2\leq& \frac{\alpha}{2\tau_{1}^2}||\nabla L_n(\hat\bta_n^{(k_0)})||_2^2 + C ||\tilde N_{n,k}||_2\\
    \leq& \frac{\alpha}{2\tau_{1}^2}\cdot \left(\frac{\tau_{1}^2}{\alpha}\right)^2 + C ||\tilde N_{n,k}||_2
    \leq \min\{\frac{\tau_{1}^2}{\alpha},\tau_{1}r\}.
\end{align*}
Now, applying the claim with $\bta=\hat\bta_n^{(k_0+1)}$ and replacing $r$ by $\frac{r}{2}$, it follows that $||\hat\bta_n^{(k_0+1)}-\hat\bta_n||_2\leq\frac{r}{2}$. This completes the induction. Now, we turn to the proof of the claim.

\noindent\textbf{Proof of the claim}: The proof uses the ASLSC property and is similar to the one used in \cite{avella2023}. Specifically, we establish the proof using contradiction. To this end, suppose $||\bta-\hat\bta_n||_2> r$; let $\tilde \bta$ denote the point on the boundary of $\mathcal{B}_r(\hat\bta)$. By Proposition 29,
\begin{align*}
    \nabla L_n(\tilde\bta) ^T\cdot (\tilde\bta-\hat\bta) \geq 2\tau_{1}||\tilde\bta-\hat\bta||_2^2.
\end{align*}
Define $\bm{v}=\frac{\tilde\bta-\hat\bta}{||\tilde\bta-\hat\bta||_2}$; then we have
\begin{align*}
    \nabla L_n(\tilde\bta) ^T\cdot \bm{v} \geq 2\tau_{1}||\tilde\bta-\hat\bta||_2 = 2\tau_{1}r.
\end{align*}
Set $f(t)=\nabla L_n(\hat\bta+t\cdot v) ^T\cdot v$ for $t\geq0$, then $f'(t)=\bm{v}^T H_n(\hat\bta+t\cdot \bm{v})\cdot \bm{v}\geq0$, since Hessian matrix is positive definite by Proposition 28.
Hence $f(t)$ is increasing in $t$ and this implies that
\begin{align*}
    ||\nabla L_n(\bta)||_2 \geq \nabla L_n(\bta)^T\cdot \bm{v} \geq \nabla L_n(\tilde\bta)^T\cdot \bm{v} \geq 2\tau_{1}r
\end{align*}
which is a contradiction since  $||\nabla L_n(\bta)||_2\leq 2\tau_{1}r$. Therefore, it follows that
\begin{align*}
    ||\nabla L_n(\bta)||_2\leq 2\tau_{1}r \Longrightarrow ||\bta-\hat\bta||_2\leq r.
\end{align*}
This completes the proof of the claim and the lemma. \QED

We now turn to the proof of the Theorem. 

\textbf{Proof of Theorem 39
}: Using Proposition 34
with $K$ replaced by $K_n$, it follows that $\hat{\bta}_n^{(K_n)}$ satisfies $\ep-$HDP. We next turn to verification of the utility.
We assume $N$ is large enough to satisfy the conditions in Lemma \ref{lem:NR_decomp}. That is for $n>N$ such that $P(||\tilde N_{n,k}||_2\leq r_{noi})\geq 1-\frac{\xi}{K}$ for $r_{noi}\sim n^{-\frac{1}{p}}(K\log (K/\xi))^{\frac{1}{2}}$. We will use Lemma \ref{lem:NR_completeness} and Lemma \ref{lem:NR_continuity} to obtain the following claim:

\noindent\textbf{Claim}: For $\hat\bta_n\in B_{r/2}(\bta_g)$ and $||\nabla L_n(\hat\bta_n^{(0)})||_2\leq\min\{\frac{\tau_{1}r}{2},\frac{\tau_{1}^2}{\alpha}\}$, the inequality
\begin{align*}
    \frac{\alpha}{2\tau_{1}^2}||\nabla L_n(\hat\theta_n^{(K)})||_2\leq \left( \frac{\alpha}{2\tau_{1}^2}||\nabla L_n(\hat\theta_n^{(0)})||_2 \right)^{2^K}+3C\cdot r_{noi}
\end{align*}
holds for some constant $C\in(0,\infty)$ with probability $1-\xi$.

Using Proposition 29
(2), and multiplying both side by $\frac{\alpha}{\tau_1}$, we obtain
\begin{align*}
    \frac{\alpha}{\tau_1}||\hat\bta_n^{(K)}-\hat\bta_n||_2
    \leq& \frac{\alpha}{2\tau_{1}^2} ||\nabla L_n(\hat\bta_n^{(K)})-\nabla L_n(\hat\bta_n)||_2.
\end{align*}
Now using the fact that $\nabla L_n(\hat\bta_n)=0$, the claim, we obtain (since $||\nabla L_n(\hat\bta_n^{(0)})||_2\leq\min\{\frac{\tau_{1}r}{2},\frac{\tau_{1}^2}{\alpha}\}$)
\begin{align*}
    \frac{\alpha}{\tau_1}||\hat\bta_n^{(K)}-\hat\bta_n||_2 \leq& \left(\frac{1}{2} \right)^{2^{K}} + 3Cr_{noi}.
\end{align*}

Choose $K$ large such that $\left(\frac{1}{2} \right)^{2^{K}}\leq Cr_{noi}$, that is $K\geq\frac{1}{\log 2}\log\frac{\log Cr_{noi}}{\log(1/2)}$, then
\begin{align*}
    ||\hat\bta_n^{(K)}-\hat\bta_n||_2\leq 4Cr_{noi}.
\end{align*}

By Lemma \ref{lem:NR_decomp}, $r_{noi}\sim n^{-\frac{1}{p}}(K\log (K/\xi))^{\frac{1}{2}}$. Using the sharp bound of $\Delta_n^{(H)}$ in Theorem 35,
we obtain the utility bound.
This also implies that $K\geq C'\log\log n$ for some $C'\in(0,\infty)$. We complete the proof by establishing the claim.

\noindent\textbf{Proof the the claim}: We prove the claim by induction. Notice that for $k=1$, the claim is true by Lemma \ref{lem:NR_completeness} and Lemma \ref{lem:NR_continuity}. Assume that the claim holds for $k=k_0$. Then for $k=k_0+1$, using Lemma \ref{lem:NR_completeness} and the choice of $\hat\bta_n^{(0)}$ such that $||\nabla L_n(\hat\bta_n^{(0)})||_2\leq\min\{\frac{\tau_{1}r}{2},\frac{\tau_{1}^2}{\alpha}\}$, it follows that
\begin{align*}
    \frac{\alpha}{2\tau_{1}^2}||\nabla L_n(\hat\bta_n^{(k_0+1)})||_2 \leq& \left(\frac{\alpha}{2\tau_{1}^2}||\nabla L_n(\hat\bta_n^{(k_0)})||_2\right)^2 + C r_{noi}.   
\end{align*}
Now by inductive hypothesis, it follows that
\begin{align*}
    \frac{\alpha}{2\tau_{1}^2}||\nabla L_n(\hat\bta_n^{(k_0+1)})||_2\leq& \left[\left( \frac{\alpha}{2\tau_{1}^2}||\nabla L_n(\hat\bta_n^{(0)})||_2 \right)^{2^{k_0}}+3C r_{noi}\right]^2 +  C r_{noi}\\
    \leq& \left( \frac{\alpha}{2\tau_{1}^2}||\nabla L_n(\hat\bta_n^{(0)})||_2 \right)^{2^{k_0+1}} + \frac{3}{2}Cr_{noi} + 9C^2r_{noi}^2 + Cr_{noi}.
\end{align*}
Let $n$ be large  such that $9C^2r_{noi}^2\leq \frac{1}{2}Cr_{noi}$. It then follows that
\begin{align*}
    \frac{\alpha}{2\tau_{1}^2}||\nabla L_n(\hat\bta_n^{(k_0+1)})||_2 
    \leq& \left( \frac{\alpha}{2\tau_{1}^2}||\nabla L_n(\hat\bta_n^{(0)})||_2 \right)^{2^{k_0+1}} + 3Cr_{noi}.
\end{align*}
This completes the induction and the proof of the Claim and the Theorem. \QED

We next turn to the proof of Theorem 42.
First, we recall that $Q$ is the distribution associated with the mechanism, representing the noise distribution.
\subsection{Proof of Theorem 42
}
We begin with part (1). Suppose $\bta_n^{(K_n)}$ is obtained using the PGD or PNR algorithm. Then, using Theorem 36
or Theorem 39
with $p\in (1,2)$ it follows that
$n^{\frac{1}{2}}||\hat\bta_n^{(K_n)}-\hat\bta_n||_2$ converges to zero in probability (with respect to the joint distribution of $P_g \times Q$) since $K_n  \sim \log n $ for PGD algorithm and $K_n \sim \log \log n$ for the PNR algorithm. Turning to part (2), observe that
\begin{eqnarray} {\label{eq:efficiency}}
    (\hat\bta_n^{(K_n)}-\bta_g) = (\hat\bta_n^{(K_n)}-\hat\bta_n) + (\hat\bta_n-\bta_g).
\end{eqnarray}
Now, taking the norm, the first term on the RHS of the above equation converges to 0 in probability by part (1), and the second term converges to zero almost surely under the assumptions {\bf\ref{asp:A1}}-{\bf\ref{asp:A8}} in Section~\ref{app:A}. Finally, turning to part (3), by multiplying both sides of (\ref{eq:efficiency}) by $\sqrt{n}$, the first term converges to zero in $P_g \times Q$ probability by part (1). The second term converges to a normal distribution under the assumptions {\bf\ref{asp:A1}}-{\bf\ref{asp:A8}} in Section~\ref{app:A} under $P_g$, by Theorem \ref{thm:app_A1} in Section~\ref{app:A}. Hence, $n^{\frac{1}{2}}(\hat\bta_n^{(K_n)}-\bta_g)$ converges in distribution (under $P_g \times Q$) to a multivariate normal distribution; that is,
\begin{eqnarray*}
  \lim_{n \ra \ff}  P_g\times Q\left(n^{\frac{1}{2}}(\hat\bta_{n}^{(K_n)}-\bta_g)\leq \bm{x} \right) = P(\mathbf{Z}\leq \bm{x}),
\end{eqnarray*}  
where $\mathbf{Z} \sim N(0,\Sigma_g)$.

\clearpage

\section{Regularity conditions and asymptotic properties of MHDE}\label{app:A}
Let $f(x)$ and $g(x)$ be any two probability density functions. The Hellinger distance between $f(x)$ and $g(x)$ is defined as the $L_2$-norm of the difference between the square root of density functions, that is,
\begin{align*}
    HD^2(f,g) = ||f^{\frac{1}{2}}(x)-g^{\frac{1}{2}}(x)||_2^2 = \int\left[f^{\frac{1}{2}}(x)-g^{\frac{1}{2}}(x)\right]^2 \mathrm{d}x.
\end{align*}
Let $\{X_1, X_2,\cdots, X_n\}$ be i.i.d. real-valued random variables with density $g(\cdot)$, and postulated to belong to a parametric family $\{f_{\bta} : \bta \in \Theta\subset\mathbb{R}^m\}$. The minimum Hellinger distance estimator in the population, $\bta_g$,  if it exists, is the minimizer of the
$||f_{\bta}^{\frac{1}{2}}-g^{\frac{1}{2}}||_2$; that is,
\begin{eqnarray*}
\bta_g=\argmin ||f_{\bta}^{\frac{1}{2}}-g^{\frac{1}{2}}||_2 = \argmin HD(f_{\bta}, g) .
\end{eqnarray*}
When $g(\cdot)=f_{\bta_0}(\cdot)$, $\bta_g=\bta_0$. We also assume that $\bta_g$ and $\bta_0$ belong to the interior of $\Theta$. \cite{Beran1977} and  \cite{Cheng2006} 
establish that  under the assumption,
\begin{assumption} \label{asp:A1}
    $\Theta \subset \Real^m$ is compact and convex and the family $\{f_{\bta}(\cdot): \bta \in \Theta\}$ is identifiable; that is, if 
 $\bta_1\neq\bta_2$ then $f_{\bta_1}(\cdot)\neq f_{\bta_2}(\cdot)$ on a set of positive Lebesgue measure.
\end{assumption}
\noindent that $\bta_g$ exists and is unique. We will assume this condition holds. 
In practice, one replaces $g(\cdot)$ by $g_n(\cdot)$, where $g_n(\cdot)$ is a nonparametric estimate of of $g(\cdot)$; specifically, a kernel density estimator, defined below.
\begin{align*}
    g_n(x) = \frac{1}{n\cdot c_n}\sum_{i=1}^{n} K\left(\frac{x-X_i}{c_n}\right).
\end{align*}
The MHDE is obtained by minimizing the loss function
\begin{align*}
    \hat\bta_{n} = \arg\min_{\bta\in\Theta} L_n(\bta), \quad\text{where } L_n(\bta) = \int_{\mathbb{R}}(\sqrt{f_{\bta}(x)}-\sqrt{g_n(x)})^2 \mathrm{d}x.
\end{align*}
Asymptotic properties of $\hat\bta_{n}$ rely on the bandwidth $c_n$ and additional regularity assumptions on the parametric family. We  provide the assumptions below: 
\begin{assumption}\label{asp:A2}
    The kernel function $K(\cdot)$ is symmetric (about 0) density with compact support. The bandwidth $c_n$ satisfies $c_n\to0$, $n^{\frac{1}{2}}c_n^2\to0$, $n^{\frac{1}{2}}c_n\to\infty$.
\end{assumption}
\begin{assumption}\label{asp:A3}
$f_{\bta}(x)$ is twice continuously differentiable in $\bta$. Also, the Fisher information matrix $I(\bta)$ is 
positive definite and continuous in $\bta$ with finite maximum eigenvalue.
\end{assumption}
\begin{assumption}\label{asp:A4}
    $||\mathbf{u}_{\bta}(\cdot)f_{\bta}^{\frac{1}{2}}(\cdot)||_2$, $||\mathbf{\dot u}_{\bta}(\cdot)f_{\bta}^{\frac{1}{2}}(\cdot)||_2$, $||\mathbf{u}_{\bta}(\cdot)\mathbf{u}^T_{\bta}(\cdot)f_{\bta}^{\frac{1}{2}}(\cdot)||_2$  exist and are continuous in $\bta$. 
\end{assumption}
\begin{assumption}\label{asp:A5}
    Let $\{a_n,n\geq1\}$ be a sequence diverging to infinity. Assume $\lim\limits_{n\to\infty} n \sup\limits_{t\in supp(K)} \mathbf{P}(|X-c_n t|>a_n)=0$, where $supp(K)$ is the support of the kernel density $K(\cdot)$ and $X$ is a generic random variable with density $f_{\bta_g}(\cdot)$.
\end{assumption}
\begin{assumption}\label{asp:A6}
    Let $M(n) = \sup\limits_{|x|\leq a_n} \sup\limits_{t\in supp(K)} |f_{\bta_g}^{-1}(x)f_{\bta_g}(x+t c_n)|$. Assume $\sup\limits_{n\geq1} M(n) < \infty$. 
\end{assumption}
\begin{assumption}\label{asp:A7}
    The score function has a regular central behavior,
    \begin{align*}
        \lim_{n\to\infty}(n^{\frac{1}{2}}c_n)^{-1} \int_{-a_n}^{a_n} \mathbf{u}_{\bta_g}(x)\mathrm{d}x = \bm{0}; ~\text{also, assume that } \lim_{n\to\infty}(n^{\frac{1}{2}}c_n^4) \int_{-a_n}^{a_n} \mathbf{u}_{\bta_g}(x)\mathrm{d}x = \bm{0}.        
    \end{align*}
\end{assumption}
\begin{assumption}\label{asp:A8}
    The score function is smooth in an $L_2$ sense; i.e.
    \begin{align*}
        \lim_{n\to\infty}\sup\limits_{t\in supp(K)} \int_{\R} [\mathbf{u}_{\bta_g}(x+t c_n)-\mathbf{u}_{\bta_g}(x)]^2 f_{\bta_g}(x) \mathrm{d}x = \bm{0}. 
    \end{align*} 
\end{assumption}

It is known that, under the above conditions, ${\hat\bta_{n}}$ is known to be unique,  consistent, and asymptotically efficient 
(see \cite{Beran1977}, \cite{Cheng2006}). Write $g(x)=\lim\limits_{n\to\infty}g_n(x)$ (which exists by {\bf\ref{asp:A2}}) and set
\begin{align*}
    &\rho_{\bta}(x) = -4 \left[\int g^{\frac{1}{2}}(x) f^{\frac{1}{2}}_{\bta}(x) [\mathbf{u}_{\bta}(x)\mathbf{u}^T_{\bta}(x)+2\mathbf{\dot u}_{\bta}(x)] \mathrm{d}x\right]^{-1} \cdot \nabla f^{\frac{1}{2}}_{\bta}(x) ~~\text{and}\\
    & \Sigma_g= 4^{-1} \int \rho_{\bta_g}(x)\rho_{\bta_g}^T(x)\mathrm{d}x.
\end{align*}
The next theorem is concerned with the limit distribution of MHDE and is similar to the proof in \cite{Cheng2006}
when the true model is $g(\cdot)$.
\begin{theorem} \label{thm:app_A1}
Under the assumptions~{\bf\ref{asp:A1}}-{\bf\ref{asp:A8}}, 
$\rho_{\bta}(\cdot)$ is 
continuous at $\bta_g$. Furthermore,
\begin{enumerate}
    \item $||\hat\bta_{n} - \bta_g||_2 \overset{P}{\to} 0$,
    \item $\sqrt{n}(\hat\bta_{n}-\bta_g)\overset{d}{\to} N(0,\Sigma_g)$.
    \item In particular, if $g(\cdot)=f_{\bta_0}(\cdot)$, then 
    $\sqrt{n}(\hat\bta_{n}-\bta_0)\overset{d}{\to} N(0,I^{-1}(\bta_0))$.
\end{enumerate}
\end{theorem}

\clearpage

\section{Calculations for additive mechanisms}\label{app:B}
\subsection{Gaussian mechanism}

\begin{lemma}\label{appb:lem1}
    For two $m$ dimensional random variable $\mathbf{X}\sim N(\mathbf{w}_1,\sigma^2\cdot\mathbf{I})$ and $\mathbf{Y}\sim N(\mathbf{w}_2,\sigma^2\cdot\mathbf{I})$, the power divergence with parameter $\la$ is given by
    \begin{align*}
        D_\la(\mathbf{X},\mathbf{Y})= \begin{cases} 
        \frac{1}{\lambda(\lambda+1)} \left[  e^{\frac{\lambda(\lambda+1)||\mathbf{v}||_2^2}{2\sigma^2}}  -1 \right], & \la(\la+1)\neq0\\
        \frac{||\mathbf{v}||_2^2}{2\sigma^2},& \la(\la+1)=0,
        \end{cases}
    \end{align*}
    where $\mathbf{v}=\mathbf{w}_1-\mathbf{w}_2$. In particular, if $\la=-\frac{1}{2}$ then $D_\la(\mathbf{X},\mathbf{Y})=-4\left[e^{-\frac{||\mathbf{v}||_2^2}{8\AV{\sigma^2}}}-1\right]$.
\end{lemma}

\textbf{Proof}:
Denote the density function for $\mathbf{X}$ and $\mathbf{Y}$ by $p(\cdot)$ and $q(\cdot)$ correspondingly, that is
\begin{align*}
    p(\mathbf{x}) = \frac{1}{(\sqrt{2\pi}\sigma)^m} e^{-\frac{||\mathbf{x}-\mathbf{w}_1||_2^2}{2\sigma^2}} ~~\text{and}~~ q(\mathbf{x}) = \frac{1}{(\sqrt{2\pi}\sigma)^m} e^{-\frac{||\mathbf{x}-\mathbf{w}_2||_2^2}{2\sigma^2}}.
\end{align*}
Let $\mathbf{y}=\mathbf{x}-\mathbf{w}_2$, $\mathbf{v}=\mathbf{w}_1-\mathbf{w}_2$, and denote by $y_i,v_i$ the $i^{th}$ element of $\mathbf{y}$ and $\mathbf{v}$. For the case $\la(\la+1)\neq0$, the power divergence with parameter $\lambda$ between $\mathbf{X}$ and $\mathbf{Y}$ is given by
\begin{align*}
    D_\lambda(\mathbf{X},\mathbf{Y}) =& \frac{1}{\lambda(\lambda+1)}\int_{\mathbb{R}^m}\left[ \frac{p^{\lambda+1}(\mathbf{x})}{q^{\lambda+1}(\mathbf{x})}\cdot q(\mathbf{x}) -q(\mathbf{x})\right] \mathrm{d}\mathbf{x}\\
    =& \frac{1}{\lambda(\lambda+1)} \left[ \int_{\mathbb{R}^m} \frac{1}{(\sqrt{2\pi}\sigma)^m} e^{-\frac{(\lambda+1)||\mathbf{x}-\mathbf{w}_1||_2^2-\lambda ||\mathbf{x}-\mathbf{w}_2||_2^2}{2\sigma^2}}\mathrm{d}\mathbf{x} -1 \right]\\
    =& \frac{1}{\lambda(\lambda+1)} \left[ \int_{\mathbb{R}^m} \frac{1}{(\sqrt{2\pi}\sigma)^m} e^{-\frac{(\lambda+1)||\mathbf{y}-\mathbf{v}||_2^2-\lambda ||\mathbf{y}||_2^2}{2\sigma^2}}\mathrm{d}\mathbf{y} -1 \right]\\
    =& \frac{1}{\lambda(\lambda+1)} \left[ \left(\prod_{i=1}^{m} \int_{\mathbb{R}} \frac{1}{\sqrt{2\pi}\sigma} e^{-\frac{(\lambda+1)(y_i-v_i)^2-\lambda y_i^2}{2\sigma^2}}\mathrm{d}y_i \right) -1 \right]\\
    =& \frac{1}{\lambda(\lambda+1)} \left[ \left(\prod_{i=1}^{m} e^{\frac{\lambda(\lambda+1)v_i^2}{2\sigma^2}} \int_{\mathbb{R}} \frac{1}{\sqrt{2\pi}\sigma} e^{-\frac{(y_i-(1+\lambda)v_i)^2 }{2\sigma^2}}\mathrm{d}y_i \right) -1 \right]\\
    =& \frac{1}{\lambda(\lambda+1)} \left[  e^{\frac{\lambda(\lambda+1)||\mathbf{v}||_2^2}{2\sigma^2}}   -1 \right].
\end{align*}
Next consider the case $\la(\la+1)=0$. Denote the $i^{th}$ element of $\mathbf{x}$, $\mathbf{w}_1$, and $\mathbf{w}_2$ by $x_i,w_{1,i},w_{2,i}$ correspondingly.  First we study the case $\la=0$.  
\begin{align*}
    D_{0}(\mathbf{X},\mathbf{Y}) =& \int_{\mathbb{R}^m}p(\mathbf{x})\log(\frac{p(\mathbf{x})}{q(\mathbf{x})})\mathrm{d}x\\
    =& \int_{\mathbb{R}^m}p(\mathbf{x})\frac{-(||\mathbf{x}-\mathbf{w}_1||_2^2-||\mathbf{x}-\mathbf{w}_2||_2^2)}{2\sigma^2}\mathrm{d}\mathbf{x}\\
    =& \sum_{i=1}^m \int_{\mathbb{R}}\frac{1}{\sqrt{2\pi}\sigma}e^{-\frac{(x_i-w_{1,i})^2}{2\sigma^2}}\cdot \frac{(w_{1,i}-w_{2,i})(2x-w_{1,i}-w_{2,i})}{2\sigma^2}\mathrm{d}x_i\\
    =& \sum_{i=1}^m \left( \frac{w_{1,i}-w_{2,i}}{2\sigma^2} \mathbf{E}_{X\sim N(w_{1,i},\sigma^2)} [2X-w_{1,i}-w_{2,i}] \right)\\
    =& \frac{||\mathbf{w}_1-\mathbf{w}_2||_2^2}{2\sigma^2} = \frac{||\mathbf{v}||_2^2}{2\sigma^2}
\end{align*}
The case $ \la=-1 $ is similar and this completes the proof. \QED

\subsection{Laplace mechanism}

\begin{lemma}\label{appb:lem2}
    For two $m$ dimensional random variable $\mathbf{X}$ and $\mathbf{Y}$, where $X_i\sim Lap(w_{1,i},b)$ and $Y_i\sim Lap(w_{2,i},b)$, the power divergence between them,  with parameter $\la$ is given by
    \begin{align*}
        D_\la(\mathbf{X},\mathbf{Y})= \begin{cases}
        \frac{1}{\lambda(\lambda+1)} \left[ \left(\prod_{i=1}^{m} \int_{\mathbb{R}} \frac{1}{2b} e^{-\frac{(\lambda+1)|y_i-v_i|-\lambda |y_i|}{b}}\mathrm{d}y_i \right) -1 \right], & \la(\la+1)\neq0\\
        \sum_{i=1}^m \int_{\mathbb{R}} \frac{1}{2b}e^{-\frac{|y_j|}{b}}\cdot \frac{|y_i-v_i|-|y_i|}{b}\mathrm{d}y_i, & \la(\la+1)=0,
        \end{cases}
    \end{align*}
    where $v_i=w_{1,i}-w_{2,i}$. Furthermore
    \begin{align*}
        D_\la(\mathbf{X},\mathbf{Y})\leq \begin{cases}
        \frac{1}{\lambda(\lambda+1)} \left[  e^{\frac{sign(\lambda)(\lambda+1)||\mathbf{v}||_1}{b}} -1 \right], & \la(\la+1)\neq0\\
        \frac{||\mathbf{v}||_1}{b}, & \la(\la+1)=0,
        \end{cases}
    \end{align*}
    In particular, if $\la=-\frac{1}{2}$, then\begin{align*}D_\la(\mathbf{X},\mathbf{Y})=&-4 \left[ \left(\prod_{i=1}^{m} \int_{\mathbb{R}} \frac{1}{2b} e^{-\frac{|y_i-v_i| + |y_i|}{2b}}\mathrm{d}y_i \right) -1 \right]\\ \leq& -4 \left[  e^{\frac{-||\mathbf{v}||_1}{2b}} -1 \right].\end{align*}
\end{lemma}

\textbf{Proof}:
Denote the density function for $\mathbf{X}$ and $\mathbf{Y}$ by $p(\cdot)$ and $q(\cdot)$ correspondingly, that is
\begin{align*}
    p(\mathbf{x}) = \frac{1}{(2b)^m} e^{-\frac{||\mathbf{x}-\mathbf{w}_1||_1}{b}} ~~\text{and}~~ q(\mathbf{x}) = \frac{1}{(2b)^m} e^{-\frac{||\mathbf{x}-\mathbf{w}_2||_1}{b}}.
\end{align*}
Let $\mathbf{y}=\mathbf{x}-\mathbf{w}_2$, $\mathbf{v}=\mathbf{w}_1-\mathbf{w}_2$, and denote $y_i,v_i$ the $i^{th}$ element of $\mathbf{y}$ and $\mathbf{v}$. For the case $\la(\la+1)\neq0$, the power divergence with parameter $\lambda$ between $\mathbf{X}$ and $\mathbf{Y}$ is given by
\begin{align*}
    D_\lambda(\mathbf{X},\mathbf{Y}) =& \frac{1}{\lambda(\lambda+1)}\int_{\mathbb{R}^m} \frac{p^{\lambda+1}(\mathbf{x})}{q^{\lambda+1}(\mathbf{x})}\cdot q(\mathbf{x}) -q(\mathbf{x}) \mathrm{d}\mathbf{x}\\
    =& \frac{1}{\lambda(\lambda+1)} \left[ \int_{\mathbb{R}^m} \frac{1}{(2b)^m} e^{-\frac{(\lambda+1)||\mathbf{y}-\mathbf{v}||_1-\lambda ||\mathbf{y}||_1}{b}}\mathrm{d}\mathbf{y} -1 \right]\\
    =& \frac{1}{\lambda(\lambda+1)} \left[ \left(\prod_{i=1}^{m} \int_{\mathbb{R}} \frac{1}{2b} e^{-\frac{(\lambda+1)|y_i-v_i|-\lambda |y_i|}{b}}\mathrm{d}y_i \right) -1 \right].
\end{align*}
Furthermore,
\begin{align*}
    D_\lambda(\mathbf{X},\mathbf{Y}) \leq& \frac{1}{\lambda(\lambda+1)} \left[ \left(\prod_{i=1}^{m} \int_{\mathbb{R}} \frac{1}{2b} e^{-\frac{(\lambda+1)|y_i|- sign(\lambda)(\lambda+1)|v_i|-\lambda |y_i|}{b}}\mathrm{d}y_i \right) -1 \right]\\
    =& \frac{1}{\lambda(\lambda+1)} \left[ \left(\prod_{i=1}^{m} e^{\frac{sign(\lambda)(\lambda+1)|v_i|}{b}} \int_{\mathbb{R}} \frac{1}{2b} e^{\AV{-}\frac{|y_i|}{b}}\mathrm{d}y_i \right) -1 \right]\\
    =& \frac{1}{\lambda(\lambda+1)} \left[  e^{\frac{sign(\lambda)(\lambda+1)||\mathbf{v}||_1}{b}} -1 \right].
\end{align*}
Next we consider the case $\la(\la+1)=0$. Denote the $i^{th}$ element of $\mathbf{x}$, $\mathbf{w}_1$, and $\mathbf{w}_2$ by $x_i,w_{1,i},w_{2,i}$ correspondingly.  We first study the case $\la=0$.To this end,
\begin{align*}
    D_{0}(\mathbf{X},\mathbf{Y}) =& \int_{\mathbb{R}^m}p(\mathbf{x})\log(\frac{p(\mathbf{x})}{q(\mathbf{x})})\mathrm{d}\mathbf{x}\\
    =& \int_{\mathbb{R}^m}p(\mathbf{x})\frac{-(||\mathbf{x}-\mathbf{w}_1||_1-||\mathbf{x}-\mathbf{w}_2||_1)}{b}\mathrm{d}\mathbf{x}\\
    =& \sum_{i=1}^m \int_{\mathbb{R}}\frac{1}{2b}e^{-\frac{|x_i-w_{1,i}|}{b}}\cdot \frac{|x_i-w_{2,i}|-|x_i-w_{1,i}|}{b}\mathrm{d}x_i\\
    =& \sum_{i=1}^m \int_{\mathbb{R}}\frac{1}{2b}e^{-\frac{|y_i|}{b}}\cdot \frac{|y_i-v_i|-|y_i|}{b}\mathrm{d}y_i.
\end{align*}
Furthermore,
\begin{align*}
    D_{0}(\mathbf{X},\mathbf{Y}) \leq& \sum_{i=1}^m \int_{\mathbb{R}}\frac{1}{2b}e^{-\frac{|x_i-w_{1,i}|}{b}}\cdot \left|\frac{|x_i-w_{2,i}|-|x_i-w_{1,i}|}{b}\right|\mathrm{d}x_i\\
    \leq& \sum_{i=1}^m \int_{\mathbb{R}}\frac{1}{2b}e^{-\frac{|x_i-w_{1,i}|}{b}}\cdot \frac{|w_{1,i}-w_{2,i}|}{b}\mathrm{d}x_i = \frac{||\mathbf{w}_1-\mathbf{w}_2||_1}{b}.
\end{align*}
The case $\la=-1$ is similar, and this completes the proof. \QED

\subsection{Exact Laplace mechanism}

\begin{lemma} \label{appb:lem3}
    For two $m$ dimensional random variable $\mathbf{X}$ and $\mathbf{Y}$, where $X_i\sim Lap(w_{1,i},b)$ and $Y_i\sim Lap(w_{2,i},b)$, the $\la-$power divergence between, $D_\lambda(\mathbf{X},\mathbf{Y})$, for $\la(\la+1) \neq 0$ \AV{and $\la\neq-\frac12$} given by
\begin{align*}
    \frac{1}{\lambda(\lambda+1)} \left[ \left(\prod_{i=1}^{m} \frac{1}{2b} \left[e^{\frac{\lambda |v_i|}{b}} \left(b + \frac{b}{2\lambda+1}\right) + e^{-\frac{(\lambda+1) |v_i|}{b}} \left(b-\frac{b}{2\lambda+1}\right) \right] \right) -1 \right].
\end{align*}
If $\la(\la+1)=0$ or  \AV{$\la=-\frac12$}, then     
\begin{align*}
D_\lambda(\mathbf{X},\mathbf{Y})= 
    \begin{cases}
    -4 \left[ \left(\prod_{i=1}^{m} \AV{\biggl(}e^{\frac{-|v_i|}{2b}}+\frac{|v_i|}{2b}e^{-\frac{|v_i|}{2b}}\AV{\biggr)} \right) -1 \right], & \la= -\frac{1}{2}\\
    \frac{1}{2b}\sum_{i=1}^m  \left[2|v_i|-2b+2be^{-\frac{|v_i|}{b}}\right], & \la(\la+1)=0.
    \end{cases}
\end{align*}
\end{lemma}

\textbf{Proof}:
In case $\la(\la+1)\neq0, \la\neq -\frac{1}{2}$, using Lemma~\ref{appb:lem2}, it follows that 
\begin{align*}
    D_\lambda(\mathbf{X},\mathbf{Y})
    =& \frac{1}{\lambda(\lambda+1)} \left[ \left(\prod_{i=1}^{m} \int_{\mathbb{R}} \frac{1}{2b} e^{-\frac{(\lambda+1)|y_i-v_i|-\lambda |y_i|}{b}}\mathrm{d}y_i \right) -1 \right].
\end{align*}
For each $i$, we remove the absolute sign by studying case $v_i<0$ and $v_i\geq0$. If $v_i<0$,
\begin{align*}
    &\int_{\mathbb{R}} \frac{1}{2b} e^{-\frac{(\lambda+1)|y_i-v_i|-\lambda |y_i|}{b}}\mathrm{d}y_i\\
    =& \frac{1}{2b} \left[\int_{-\infty}^{v_i} e^{-\frac{-(\lambda+1)(y_i-v_i)+\lambda y_i}{b}}\mathrm{d}y_i + \int_{v_i}^{0} e^{-\frac{(\lambda+1)(y_i-v_i)+\lambda y_i}{b}}\mathrm{d}y_i + \int_{0}^{\infty} e^{-\frac{(\lambda+1)(y_i-v_i)-\lambda y_i}{b}}\mathrm{d}y_i \right]\\
    =& \frac{1}{2b} \left[e^{\frac{-\lambda v_i}{b}} \left(b + \frac{b}{2\lambda+1}\right) + e^{\frac{(\lambda+1) v_i}{b}} \left(b-\frac{b}{2\lambda+1}\right) \right].
\end{align*}
If $v_i\geq 0$,
\begin{align*}
    &\int_{\mathbb{R}} \frac{1}{2b} e^{-\frac{(\lambda+1)|y_i-v_i|-\lambda |y_i|}{b}}\mathrm{d}y_i\\
    =& \frac{1}{2b} \left[\int_{-\infty}^{0} e^{-\frac{-(\lambda+1)(y_i-v_i)+\lambda y_i}{b}}\mathrm{d}y_i + \int_{0}^{v_i} e^{-\frac{-(\lambda+1)(y_i-v_i)-\lambda y_i}{b}}\mathrm{d}y_i + \int_{v_i}^{\infty} e^{-\frac{(\lambda+1)(y_i-v_i)-\lambda y_i}{b}}\mathrm{d}y_i \right]\\
    =& \frac{1}{2b} \left[e^{\frac{\lambda v_i}{b}} \left(b + \frac{b}{2\lambda+1}\right) + e^{\frac{-(\lambda+1) v_i}{b}} \left(b-\frac{b}{2\lambda+1}\right) \right].
\end{align*}
Combining the cases $v_i<0$ and $v_i\geq0$, we get
\begin{align*}
    \int_{\mathbb{R}} \frac{1}{2b} e^{-\frac{(\lambda+1)|y_i-v_i|-\lambda |y_i|}{b}}\mathrm{d}y_i
    = \frac{1}{2b} \left[e^{\frac{\lambda |v_i|}{b}} \left(b + \frac{b}{2\lambda+1}\right) + e^{-\frac{(\lambda+1) |v_i|}{b}} \left(b-\frac{b}{2\lambda+1}\right) \right].
\end{align*}
Therefore,
\begin{align*}
    D_\lambda(\mathbf{X},\mathbf{Y})  
    = \frac{1}{\lambda(\lambda+1)} \left[ \left(\prod_{i=1}^{m} \frac{1}{2b} \left[e^{\frac{\lambda |v_i|}{b}} \left(b + \frac{b}{2\lambda+1}\right) + e^{-\frac{(\lambda+1) |v_i|}{b}} \left(b-\frac{b}{2\lambda+1}\right) \right] \right) -1 \right].
\end{align*}
We now turn to the case $\la= -\frac{1}{2}$. If $v_i<0$, by the same calculation of the integral, it follows that
\begin{align*}
    \int_{\mathbb{R}} \frac{1}{2b} e^{-\frac{(\lambda+1)|y_i-v_i|-\lambda |y_i|}{b}}\mathrm{d}y_i
    = e^{\frac{v_i}{2b}}-\frac{v_i}{2b}e^{\frac{v_i}{2b}}.
\end{align*}
If $v_i\geq0$,
\begin{align*}
    \int_{\mathbb{R}} \frac{1}{2b} e^{-\frac{(\lambda+1)|y_i-v_i|-\lambda |y_i|}{b}}\mathrm{d}y_i
    = e^{\frac{-v_i}{2b}}+\frac{v_i}{2b}e^{-\frac{v_i}{2b}}.
\end{align*}
Combining $v_i<0$ and $v_i\geq0$, we get
\begin{align*}
    \int_{\mathbb{R}} \frac{1}{2b} e^{-\frac{(\lambda+1)|y_i-v_i|-\lambda |y_i|}{b}}\mathrm{d}y_i
    = e^{\frac{-|v_i|}{2b}}+\frac{|v_i|}{2b}e^{-\frac{|v_i|}{2b}}.
\end{align*}
Therefore
\begin{align*}
    D_\lambda(\mathbf{X},\mathbf{Y}) = -4 \left[ \left(\prod_{i=1}^{m} e^{\frac{-|v_i|}{2b}}+\frac{|v_i|}{2b}e^{-\frac{|v_i|}{2b}} \right) -1 \right].
\end{align*}
Finally, we turn to the case $\la(\la+1)=0$. We study the case $\la=0$. Using Lemma~\ref{appb:lem2}, it follows that 
\begin{align*}
    D_0(\mathbf{X},\mathbf{Y})=
    \sum_{i=1}^m \int_{\mathbb{R}} \frac{1}{2b}e^{-\frac{|y_j|}{b}} \cdot \frac{|y_i-v_i|-|y_i|}{b}\mathrm{d}y_i.   
\end{align*}
If $v_i\geq0$,
\begin{align*}
    & \int_{\mathbb{R}} e^{-\frac{|y_i|}{b}}\cdot \frac{|y_i-v_i|-|y_i|}{b}\mathrm{d}y_i\\
    =& \int_{-\infty}^{0} e^{-\frac{-y_i}{b}} \frac{-(y_i-v_i)+y_i}{b}\mathrm{d}y_i + \int_{0}^{v_i} e^{-\frac{y_i}{b}} \frac{-(y_i-v_i)-y_i}{b}\mathrm{d}y_i + \int_{v_i}^{\infty} e^{-\frac{y_i}{b}} \frac{(y_i-v_i)-y_i}{b}\mathrm{d}y_i\\
    =& 2v_i-2b+2be^{-\frac{v_i}{b}}.
\end{align*}
If $v_i<0$,
\begin{align*}
    & \int_{\mathbb{R}} e^{-\frac{|y_i|}{b}}\cdot \frac{|y_i-v_i|-|y_i|}{b}\mathrm{d}y_i\\
    =& \int_{-\infty}^{v_i} e^{-\frac{-y_i}{b}} \frac{-(y_i-v_i)+y_i}{b}\mathrm{d}y_i + \int_{v_i}^{0} e^{-\frac{-y_i}{b}} \frac{(y_i-v_i)+y_i}{b}\mathrm{d}y_i + \int_{0}^{\infty} e^{-\frac{y_i}{b}} \frac{(y_i-v_i)-y_i}{b}\mathrm{d}y_i\\
    =& -2v_i-2b+2be^{\frac{v_i}{b}}.
\end{align*}
Combining the cases $v_i<0$ and $v_i\geq0$, we get
\begin{align*}
    \int_{\mathbb{R}} e^{-\frac{|y_i|}{b}}\cdot \frac{|y_i-v_i|-|y_i|}{b}\mathrm{d}y_i = 2|v_i|-2b+2be^{-\frac{|v_i|}{b}}.
\end{align*}
Therefore,
\begin{align*}
    D_0(\mathbf{X},\mathbf{Y}) 
    = \frac{1}{2b}\sum_{i=1}^m  \left[2|v_i|-2b+2be^{-\frac{|v_i|}{b}}\right].
\end{align*}
The case $\la=-1$ is similar, and this completes the proof. \QED

\clearpage

\section{Relationships between privacy frameworks}\label{app:C}

\subsection{Proof of Remark 5
}

\textbf{Link to $\rho-$zCDP}: Suppose $M(D)=W(D)+\mathbf Y$, where $\mathbf Y\sim N_m(0,\sigma^2I_m)$ with a dataset-independent covariance, satisfies $(\la,\ep)-$PDP for $\la>0$. Put $\Delta=\Delta_{L_2}W$ and $\rho=\Delta^2/(2\sigma^2)$. For every order $a>1$ and every adjacent $D,D'$, the Gaussian R\'enyi divergence satisfies
\[
D_a^{(R)}(P_D,P_{D'})
=\frac{a\|W(D)-W(D')\|_2^2}{2\sigma^2}\leq a\rho.
\]
At the given PDP order, the Gaussian power-divergence identity yields
\[
\frac{e^{\lambda(\lambda+1)\rho}-1}{\lambda(\lambda+1)}\leq\epsilon,
\qquad
\rho\leq\frac{\log\{1+\lambda(\lambda+1)\epsilon\}}{\lambda(\lambda+1)}\leq\epsilon.
\]
The Gaussian bound holds at \emph{every} R\'enyi order $a>1$, and hence $M$ satisfies $\rho-$zCDP, in particular $\epsilon-$zCDP.

\noindent\textbf{Link to $(\epsilon,\delta)-$differential privacy}: Suppose a mechanism $M$ satisfies $(\la,\ep)-$PDP, then by definition, $D_{\lambda}(f_1,f_2)\leq \ep$, where $f_1$ is the density of $M(w,D)$ and $f_2$ is the density of $M(w,D')$. We now determine the relationship to $(\ep, \delta)-$DP.

If $\lambda>0$, then
\begin{align*}
    & D_{\lambda}(f_1,f_2)=\frac{1}{\lambda(\lambda+1)} \left[ \int_{\mathbb{R}^m} \frac{f_1^{\lambda+1}(x)}{f_2^{\lambda}(x)}\mathrm{d}x -1 \right]\leq \epsilon\\
    \Longleftrightarrow & \int_{\mathbb{R}^m} \frac{f_1^{\lambda+1}(x)}{f_2^{\lambda}(x)}\mathrm{d}x \leq \lambda(\lambda+1)\epsilon + 1 = e^{\log(\lambda(\lambda+1)\epsilon + 1)}.
\end{align*}
For any set $A\subset \mathbb{R}^m$, applying Holder inequality for $p=\lambda+1$ and $q=\frac{\lambda+1}{\lambda}$, it follows that
\begin{align*}
    \mathbf{P}_{X\sim f_1}(X\in A) =& \int_{A} f_1(x) \mathrm{d}x = \int_{A} \frac{f_1(x)}{[f_2(x)]^{\frac{\lambda}{\lambda+1}}} \cdot [f_2(x)]^{\frac{\lambda}{\lambda+1}} \mathrm{d}x\\
    \leq&  \left(\int_{A}\left(\frac{f_1(x)}{[f_2(x)]^{\frac{\lambda}{\lambda+1}}} \right)^{p} \mathrm{d}x\right)^{\frac{1}{p}} \cdot \left(\int_{A}\left([f_2(x)]^{\frac{\lambda}{\lambda+1}}\right)^{q} \mathrm{d}x\right)^{\frac{1}{q}}\\
    =& \left(\int_{A}\frac{[f_1(x)]^{\lambda+1}}{[f_2(x)]^{\lambda}} \mathrm{d}x\right)^{\frac{1}{\lambda+1}} \cdot \left(\int_{A} f_2(x)\mathrm{d}x\right)^{\frac{\lambda}{\lambda+1}}\\
    \leq& \left(\int_{\mathbb{R}^m}\frac{[f_1(x)]^{\lambda+1}}{[f_2(x)]^{\lambda}} \mathrm{d}x\right)^{\frac{1}{\lambda+1}} \cdot \left(\int_{A} f_2(x)\mathrm{d}x\right)^{\frac{\lambda}{\lambda+1}}\\
    \leq& e^{\frac{1}{\lambda+1}\log(\lambda(\lambda+1)\epsilon + 1)} \cdot \left[\mathbf{P}_{X\sim f_2}(X\in A)\right]^{\frac{\lambda}{\lambda+1}}\\
    =& \left[ e^{\frac{1}{\lambda}\log(\lambda(\lambda+1)\epsilon + 1)} \cdot \mathbf{P}_{X\sim f_2}(X\in A)\right]^{\frac{\lambda}{\lambda+1}}.
\end{align*}
If $e^{\frac{1}{\lambda}\log(\lambda(\lambda+1)\epsilon + 1)} \cdot \mathbf{P}_{X\sim f_2}(X\in A)>\delta^{\frac{\lambda+1}{\lambda}}$, then
\begin{align}
    \mathbf{P}_{X\sim f_1}(X\in A) \leq& \left[ e^{\frac{1}{\lambda}\log(\lambda(\lambda+1)\epsilon + 1)} \cdot \mathbf{P}_{X\sim f_2}(X\in A)\right]^{\frac{\lambda}{\lambda+1}} \nonumber\\
    =& e^{\frac{1}{\lambda}\log(\lambda(\lambda+1)\epsilon + 1)} \cdot \mathbf{P}_{X\sim f_2}(X\in A)\cdot \left[ e^{\frac{1}{\lambda}\log(\lambda(\lambda+1)\epsilon + 1)} \cdot \mathbf{P}_{X\sim f_2}(X\in A)\right]^{\frac{-1}{\lambda+1}} \nonumber\\
    \leq& e^{\frac{1}{\lambda}\log(\lambda(\lambda+1)\epsilon + 1)} \cdot \mathbf{P}_{X\sim f_2}(X\in A)\cdot \delta^{\frac{-1}{\lambda}} \nonumber\\
    =& e^{\frac{1}{\lambda}\log(\frac{\lambda(\lambda+1)\epsilon + 1}{\delta})} \cdot \mathbf{P}_{X\sim f_2}(X\in A)  \label{appc:eq1}
\end{align}
If $e^{\frac{1}{\lambda}\log(\lambda(\lambda+1)\epsilon + 1)} \cdot \mathbf{P}_{X\sim f_2}(X\in A)\leq \delta^{\frac{\lambda+1}{\lambda}}$, then
\begin{align*}
    \mathbf{P}_{X\sim f_1}(X\in A) \leq& \left[ e^{\frac{1}{\lambda}\log(\lambda(\lambda+1)\epsilon + 1)} \cdot \mathbf{P}_{X\sim f_2}(X\in A)\right]^{\frac{\lambda}{\lambda+1}} =\delta.
\end{align*}
Therefore
\begin{align*}
    \mathbf{P}_{X\sim f_1}(X\in A) \leq& e^{\frac{1}{\lambda}\log(\frac{\lambda(\lambda+1)\epsilon + 1}{\delta})} \cdot \mathbf{P}_{X\sim f_2}(X\in A) + \delta.
\end{align*}
This implies that $M$ satisfies $(\frac{1}{\lambda}\log(\frac{\lambda(\lambda+1)\epsilon + 1}{\delta}),\delta)-$DP.

If $\lambda<-1$, write $\lambda'=-\lambda-1>0$, 
\begin{align*}
    & D_{\lambda}(f_2,f_1)=\frac{1}{\lambda(\lambda+1)} \left[ \int_{\mathbb{R}^m} \frac{f_2^{\lambda+1}(x)}{f_1^{\lambda}(x)}\mathrm{d}x -1 \right]\leq \epsilon\\
    \Longleftrightarrow & \int_{\mathbb{R}^m} \frac{f_2^{\lambda+1}(x)}{f_1^{\lambda}(x)}\mathrm{d}x \leq \lambda(\lambda+1)\epsilon + 1 = e^{\log(\lambda(\lambda+1)\epsilon + 1)}\\
    \Longleftrightarrow & \int_{\mathbb{R}^m} \frac{f_1^{\lambda'+1}(x)}{f_2^{\lambda'}(x)}\mathrm{d}x \leq \lambda'(\lambda'+1)\epsilon + 1 = e^{\log(\lambda'(\lambda'+1)\epsilon + 1)}.
\end{align*}
Applying Holder inequality for $p=\lambda'+1>1$, $q = \frac{\lambda'+1}{\lambda'}>1$, by the same method, we obtain
\begin{align}
    \mathbf{P}_{X\sim f_1}(X\in A) \leq \left[ e^{\frac{1}{\lambda'}\log(\lambda'(\lambda'+1)\epsilon + 1)} \cdot \mathbf{P}_{X\sim f_2}(X\in A)\right]^{\frac{\lambda'}{\lambda'+1}}. \label{appc:eq2}
\end{align}
Using the same $\delta$, it follows that
\begin{align*}
    \mathbf{P}_{X\sim f_1}(X\in A) \leq&  e^{\frac{-1}{\lambda+1}\log(\frac{\lambda(\lambda+1)\epsilon + 1}{\delta})} \cdot \mathbf{P}_{X\sim f_2}(X\in A) + \delta.
\end{align*}
This implies that $M$ satisfies $(\frac{-1}{\lambda+1}\log(\frac{\lambda(\lambda+1)\epsilon + 1}{\delta}),\delta)-$DP.

~

\noindent\textbf{Link to $\mu-$GDP for Gaussian and scalar Laplace mechanisms}: The GDP conclusions of Remark 5 of the main paper concern the additive mechanisms with common dataset-independent noise scales, rather than arbitrary PDP mechanisms.
For the Gaussian mechanism $M(D)=W(D)+\mathbf Y$, with $\mathbf Y\sim N_m(0,\sigma^2I_m)$, let $\Delta_2=\Delta_{L_2}W$. The exact Gaussian testing trade-off is $G_{\mu_G}$, where $\mu_G=\Delta_2/\sigma$. For $\lambda(\lambda+1)\neq0$, the Gaussian power divergence of an adjacent pair with shift $v$ satisfies
\[
D_\lambda\{N_m(0,\sigma^2I_m),N_m(v,\sigma^2I_m)\}
=\frac{\exp\{\lambda(\lambda+1)\|v\|_2^2/(2\sigma^2)\}-1}{\lambda(\lambda+1)}.
\]
Thus $(\lambda,\epsilon)-$PDP implies
\[
\mu_G^2\leq
\frac{2}{\lambda(\lambda+1)}\log\{1+\lambda(\lambda+1)\epsilon\},
\]
with the continuous bound $\mu_G^2\leq2\epsilon$ at $\lambda=0,-1$.

For the scalar Laplace mechanism $M(D)=W(D)+Z_b$, where $Z_b\sim\operatorname{Lap}(0,b)$ and $b$ is common to all datasets, let $\Delta_1=\Delta_{L_1}W$ and $r=\Delta_1/(2b)$. Proposition~\ref{prop:laplace-gdp-exact} establishes from the complete testing profile that its smallest GDP parameter is
\[
\mu_L(r)=2\Phi^{-1}\left(1-\frac12e^{-r}\right).
\]
In particular, the scalar Laplace PDP bound restricts $r$ through the exact power-divergence formula in Lemma~\ref{appb:lem3}. To make the dependence on $(\lambda,\epsilon)$ precise, put
\[
d_\lambda(t)=D_\lambda\{\operatorname{Lap}(0,1),\operatorname{Lap}(t,1)\},\qquad t\geq0.
\]
The formulas in Lemma~\ref{appb:lem3} show that $d_\lambda$ is strictly increasing for $t>0$. In the nontrivial admissible PDP range, $\Delta_1/b\leq d_\lambda^{-1}(\epsilon)$, so the mechanism is also $\mu$-GDP with
\[
\mu=2\Phi^{-1}\left(1-\frac12
\exp\left\{-\frac{d_\lambda^{-1}(\epsilon)}2\right\}\right).
\]
This uses the full family-specific GDP comparison, not a general conversion from a single power-divergence order. \QED

\clearpage
\section{Supporting results for Hessian convergence and utility}\label{app:D}

\subsection[Hessian convergence in Proposition 26]{Details on the convergence of $H_n(\cdot)$ to $H_{\ff}(\cdot)$ in Proposition 26
}
We write $H_{n,i,j}(\bta)$ as the $i-th$ row and $j-th$ column element of $H_{n}(\bta)$, and $I_{i,j}(\bta)$ as the $i-th$ row and $j-th$ column element of $I(\bta)$. Then we only need to show $H_{n,i,j}(\bta_0)\to I_{i,j}(\bta_0)$ for any $i,j=1,\cdots,m$.
Recall that
\begin{align*}
    H_{n,i,j}(\bta) =& \frac{\partial^2}{\partial\bta_i\partial\bta_j}L_n(\bta)= -T_{1,n,i,j}(\bta) -2 T_{2,n,i,j}(\bta),
\end{align*}
where
\begin{align*}
    T_{1,n,i,j}(\bta) = \int_{\mathbb{R}}g_n^{\frac{1}{2}}(x)f_{\bta}^{\frac{1}{2}}(x) u_{\bta,j}(x)u_{\bta,i}(x)\mathrm{d}x,\quad
    T_{2,n,i,j}(\bta) = \int_{\mathbb{R}}g_n^{\frac{1}{2}}(x)f_{\bta}^{\frac{1}{2}}(x) u_{\bta,i,j}(x)\mathrm{d}x.
\end{align*}
We decompose $T_{1,n,i,j}(\bta)$ and $T_{2,n,i,j}(\bta)$ as follows:
\begin{align*}
    T_{1,n,i,j}(\bta) 
    = T_{1,n,i,j}^{(1)}(\bta) + I_{i,j}(\bta), \quad \text{and} \quad
    T_{2,n,i,j}(\bta)
    = T_{2,n,i,j}^{(1)}(\bta) - I_{i,j}(\bta),
\end{align*}
where
\begin{align*}
    &T_{1,n,i,j}^{(1)}(\bta)=\int_{\mathbb{R}}\left(g_n^{\frac{1}{2}}(x)-f_{\bta}^{\frac{1}{2}}(x)\right)f_{\bta}^{\frac{1}{2}}(x) u_{\bta,j}(x)u_{\bta,i}(x)\mathrm{d}x, \quad\text{and}\\
    &T_{2,n,i,j}^{(1)}(\bta)=\int_{\mathbb{R}}\left(g_n^{\frac{1}{2}}(x)-f_{\bta}^{\frac{1}{2}}(x)\right)f_{\bta}^{\frac{1}{2}}(x) u_{\bta,i,j}(x)\mathrm{d}x.
\end{align*}
Then $H_{n,i,j}(\bta)$ can be written as follows:
\begin{align} \label{eq:appD_1}
    H_{n,i,j}(\bta)= I_{i,j}(\bta) - D_{n,i,j}(\bta),
\end{align}
where $D_{n,i,j}(\bta)= T_{1,n,i,j}^{(1)}(\bta) + 2 T_{2,n,i,j}^{(1)}(\bta)$. We are going to show $D_{n,i,j}(\bta)\to D_{i,j}(\bta)$ almost surely as $n\to\infty$, where $D_{i,j}(\bta) = \int_{\mathbb{R}}\left(g^{\frac{1}{2}}(x)-f_{\bta}^{\frac{1}{2}}(x)\right)f_{\bta}^{\frac{1}{2}}(x) [u_{\bta,i}(x)u_{\bta,j}(x)+2 u_{\bta,i,j}(x)]\mathrm{d}x$. Notice that
\begin{align*}
    D_{n,i,j}(\bta)
    = D_{i,j}(\bta) \AV{+}\int_{\mathbb{R}}\left(g_n^{\frac{1}{2}}(x)-g^{\frac{1}{2}}(x)\right)f_{\bta}^{\frac{1}{2}}(x) [u_{\bta,i}(x)u_{\bta,j}(x)+2 u_{\bta,i,j}(x)]\mathrm{d}x.
\end{align*}
We are going to show $\int_{\mathbb{R}}\left(g_n^{\frac{1}{2}}(x)-g^{\frac{1}{2}}(x)\right)f_{\bta}^{\frac{1}{2}}(x) [u_{\bta,i}(x)u_{\bta,j}(x)+2 u_{\bta,i,j}(x)]\mathrm{d}x\to 0$ almost surely as $n\to\infty$. Using the Cauchy–Schwarz inequality and the upper bounds in assumptions {\bf(U1)}-{\bf(U2)}, it follows that as $n\to\infty$,
\begin{align*}
    & \left|\int_{\mathbb{R}}\left(g_n^{\frac{1}{2}}(x)-g^{\frac{1}{2}}(x)\right)f_{\bta}^{\frac{1}{2}}(x) [u_{\bta,i}(x)u_{\bta,j}(x)+2 u_{\bta,i,j}(x)]\mathrm{d}x\right|\\
    \leq& HD(g_n,g)\cdot \left(\mathbf{E}_{\bta}\left[u^2_{\bta,i}(X)u^2_{\bta,j}(X)\right]\right)^{\AV{\frac12}} + 2 HD(g_n,g)\cdot \left(\mathbf{E}_{\bta}\left[u^2_{\bta,i,j}(X)\right]\right)^{\AV{\frac12}}\\
    \leq& c\cdot HD(g_n,g)
    \xlongrightarrow{a.s.}  0.
\end{align*}
The convergence follows from $HD(g_n,g)\xlongrightarrow{a.s.} 0 ~\text{when}~ n\to\infty$, since by assumption {\bf\ref{asp:A2}}, $||g_n-g||_1$ converges to  zero almost surely. Thus,
$H_{\ff}(\bta)= I(\bta)- D(\bta)$. This completes the proof.

\subsection[Upper bound in Proposition 27]{Establish the upper bound for $D_{i,j}(\cdot)$ in Proposition 27
}
Using Cauchy- Schwarz inequality and {\bf(U1)}-{\bf(U2)} the result follows. To see this, notice that
\begin{align*}
    |D_{i,j}(\bta)| =& \left| \int_{\mathbb{R}}\left(g^{\frac{1}{2}}(x)-f_{\bta}^{\frac{1}{2}}(x)\right)f_{\bta}^{\frac{1}{2}}(x) [u_{\bta,i}(x)u_{\bta,j}(x)+2 u_{\bta,i,j}(x)]\mathrm{d}x \right|.
\end{align*}
Now, splitting the RHS of the above equation, we see that it is bounded above by
\begin{align*}
 \left| \int_{\mathbb{R}}\left(g^{\frac{1}{2}}(x)-f_{\bta}^{\frac{1}{2}}(x)\right)f_{\bta}^{\frac{1}{2}}(x) u_{\bta,i}(x)u_{\bta,j}(x)\mathrm{d}x \right| + 2 \left| \int_{\mathbb{R}}\left(g^{\frac{1}{2}}(x)-f_{\bta}^{\frac{1}{2}}(x)\right)f_{\bta}^{\frac{1}{2}}(x) u_{\bta,i,j}(x)\mathrm{d}x \right|.
\end{align*} 
Now, applying Cauchy-Schwarz inequality, we get
\begin{align} \label{eq:appD_2}
   |D_{i,j}^2(\bta)| \leq& \AV{2}HD^2(g,f_{\bta})\cdot \mathbf{E}_{\bta}\left[u^2_{\bta,i}(X)u^2_{\bta,j}(X)\right] + \AV{8}HD^2(g,f_{\bta})\cdot \mathbf{E}_{\bta}\left[u^2_{\bta,i,j}(X)\right] \nonumber\\
    \leq& c\cdot HD^2(g,f_{\bta}),
\end{align}
where $0<c=\AV{1+\sup_{\bta\in\Theta}\left\{2\mathbf{E}_{\bta}\left[u^2_{\bta,i}(X)u^2_{\bta,j}(X)\right]+8\mathbf{E}_{\bta}\left[u^2_{\bta,i,j}(X)\right]\right\}}<\ff$.

\subsection{Proof of Lemma \ref{lem:continuity_relation}}

\textbf{Statement:}
Let assumptions {\bf\ref{asp:A1}}-{\bf\ref{asp:A8}} and {\bf(U1)}-{\bf(U2)} hold. Then for $\bta \in B_r(\bta_g)$ and $n \ge N$, if $L_n(\bta)-L_n(\hat\bta)\leq \frac{r^2}{4}\tau_{1} $ then $ ||\bta-\hat\bta||_2\leq\frac{r}{2}$. Furthermore, if $||\bta-\hat\bta||_2\leq\frac{r}{2} $ for $\bta \in B_r(\bta_g)$, then  for $n \ge N$, $  L_n(\bta)-L_n(\hat\bta)\leq \frac{r^2}{4}\tau_{2}.$

\noindent\textbf{Proof}: Let $n \ge N$ and $\bta,\hat\bta \in B_r(\bta_g)$. Suppose $L_n(\bta)-L_n(\hat\bta)\leq \frac{r^2}{4}\tau_{1}$.
Then using Proposition 29
(1), it follows that  $L_n(\bta)\geq L_n(\hat\bta) + \langle\nabla L_n(\hat\bta),\bta-\hat\bta\rangle +\tau_{1}||\bta-\hat\bta||_2^2$. Since $\nabla L_n(\hat\bta)=(0,\cdots,0)$ it follows that 
$\frac{r^2 \tau_1}{4} \ge L_n(\bta)-L_n(\hat\bta) \geq \tau_{1}||\bta-\hat\bta||_2^2$, the result follows. The rest of the proof follows similarly, using Proposition 29
(3); that is, if $||\bta-\hat\bta||_2\leq\frac{r}{2} $, then
$L_n(\bta)-L_n(\hat\bta) \leq \langle\nabla L_n(\hat\bta),\bta-\hat\bta\rangle+\tau_{2,n}||\bta-\hat\bta||_2^2 =\tau_{2}||\bta-\hat\bta||_2^2\leq \frac{r^2}{4}\tau_{2}$. \QED

\subsection{Proof of Lemma \ref{lem:continuity}}

\textbf{Statement:}
Under assumptions {\bf\ref{asp:A1}}-{\bf\ref{asp:A8}} and {\bf(U1)}-{\bf(U2)}, for $\eta\leq\frac{1}{\tau_{2}}$, assume that for $n \ge N$, $\hat\bta_n\in B_{r/c}(\bta_g)\subset B_{r/2}(\bta_g)$, where $c>2\left(\frac{\tau_2}{\tau_1}\right)^{\frac{1}{2}}$, 
then there exists $\hat\bta_n^{(0)}$, such that $L_n(\hat\bta_n^{(k)})-L_n(\hat\bta_n)\leq \tau_{1}\frac{r^2}{4}$ and $||\hat\bta_n^{(k)}-\hat\bta_n||_2\leq \frac{r}{2}$ hold with probability $1-\frac{k\xi}{K}$ for all $k=0,\cdots,K$.

\noindent\textbf{Proof}: The Lemma states that under the stated conditions, the estimators from each iteration, $\hat\bta_n^{(k)}$, belong to the ball $B_{r/2}(\hat\bta_n) \subset B_{r}(\bta_g)$.  We prove the result by induction. First, for $k=0$, we choose the initial estimator to be a consistent estimator of $\bta_g$. Hence for large $n$, $||\hat\bta^{(0)}_n-\bta_g||\le \left(\frac{\tau_1}{\tau_2}\right)^{\frac{1}{2}}\frac{r}{2}-\frac{r}{c}$. Hence, for large $n$, $||\hat\bta_n^{(0)}-\hat\bta_n||_2\leq \left(\frac{\tau_1}{\tau_2}\right)^{\frac{1}{2}}\frac{r}{2}$. By Lemma \ref{lem:continuity_relation}, $L_n(\hat\bta_n^{(0)})-L_n(\hat\bta_n)\leq \tau_{1}\frac{r^2}{4}$ and $||\hat\bta_n^{(0)}-\hat\bta_n||_2\leq \frac{r}{2}$ hold.
Hence, by induction hypothesis, let 
$L_n(\hat\bta_n^{(k)})-L_n(\hat\bta_n)\leq \tau_{1}\frac{r^2}{4}$ and $||\hat\bta_n^{(k)}-\hat\bta_n||_2\leq \frac{r}{2}$ hold. We will establish that  $L_n(\hat\bta_n^{(k+1)})-L_n(\hat\bta_n)\leq \tau_{1}\frac{r^2}{4}$ and $||\hat\bta_n^{(k+1)}-\hat\bta_n||_2\leq \frac{r}{2}$. The proof of this relies on the behavior of
$||\nabla L_n(\hat\bta_n^{(k)})||_2$, $||N_{n,k}||_2$, and their relationships which is described in the following claim whose proof is relegated to the end.

\noindent\textbf{Claim}: If $||\nabla L_n(\hat\bta_n^{(k)})||_2 \geq \sqrt{\frac{(1+2\eta\tau_{2})B||N_{n,k}||_2 + \eta\tau_{2}||N_{n,k}||_2^2}{1-\eta\tau_{2}}}$, where $B$ is the upper bound of $||\nabla L_n(\hat\bta_n^{(k)})||_2$ from Proposition 25,
then the following inequality holds.
\begin{align*}
    L_n(\hat\bta_n^{(k+1)})-L_n(\hat\bta_n) \leq L_n(\hat\bta_n^{(k)})-L_n(\hat\bta_n).
\end{align*}
Using the claim for $k$, and the assumption $L_n(\hat\bta_n^{(k)})-L_n(\hat\bta_n)\leq \tau_{1}\frac{r^2}{4}$, we obtain that $L_n(\hat\bta_n^{(k+1)})-L_n(\hat\bta_n)\leq \tau_{1}\frac{r^2}{4}$. Next, applying Lemma \ref{lem:continuity_relation} we get $||\hat\bta_n^{(k+1)}-\hat\bta_n||_2\leq \frac{r}{2}$. This completes the proof under the condition of the claim. 

Next, we turn to the case $||\nabla L_n(\hat\bta_n^{(k)})||_2 < \sqrt{\frac{(1+2\eta\tau_{2})B||N_{n,k}||_2 + \eta\tau_{2}||N_{n,k}||_2^2}{1-\eta\tau_{2}}}$. By (\ref{eq:noise_bound}), it follows that with probability $1-\frac{k\xi}{K}$ (since we have $k$ iterations here)
\begin{align*}
    ||\nabla L_n(\hat\bta_n^{(k)})||_2 < \sqrt{\frac{(1+2\eta\tau_{2})B r_{noi} + \eta\tau_{2}r_{noi}^2}{1-\eta\tau_{2}}} \coloneqq \bar r_{noi}.
\end{align*}
Using Proposition 29
(1), $\tau_{1}||\hat\bta_n-\hat\bta_n^{(k)}||_2^2\leq L_n(\hat\bta_n)-L_n(\hat\bta_n^{(k)})-\langle\nabla L_n(\hat\bta_n^{(k)}),\hat\bta_n-\hat\bta_n^{(k)}\rangle$. Since $\hat{\bta}_n$ is the minimizer of $L_n(\bta)$, it follows that $\tau_{1}||\hat\bta_n-\hat\bta_n^{(k)}||_2^2\leq |\langle\nabla L_n(\hat\bta_n^{(k)}),\hat\bta_n-\hat\bta_n^{(k)}\rangle|$. Now, applying the Cauchy-Schwarz inequality, it follows that  $\tau_{1}||\hat\bta_n-\hat\bta_n^{(k)}||_2^2\leq||\nabla L_n(\hat\bta_n^{(k)})||_2\cdot||\hat\bta_n-\hat\bta_n^{(k)}||_2$. Hence, we obtain
\begin{align*}
    ||\hat\bta_n-\hat\bta_n^{(k)}||_2 \leq \frac{||\nabla L_n(\hat\bta_n^{(k)})||_2}{\tau_1} \leq \frac{\bar r_{noi}}{\tau_1}.
\end{align*}
Now using $\hat\bta_n^{(k)}-\hat\bta_n^{(k+1)}=\eta(\nabla L_n(\hat\bta_n^{k})+N_k)$, it follows that
\begin{align*}
    ||\hat\bta_n-\hat\bta_n^{(k+1)}||_2\leq& ||\hat\bta_n-\hat\bta_n^{k}||_2 + \eta||\nabla L_n(\hat\bta_n^{(k)})||_2 + \eta||N_k||_2 \leq \frac{\Bar{r}_{noi}}{\tau_{1}} + \eta \Bar{r}_{noi} + \eta r_{noi}\\
    \le& \left(\frac{\tau_1}{\tau_2}\right)^{\frac{1}{2}}\frac{r}{2} \leq \frac{r}{2},
\end{align*}
where the last inequality follows by taking $n$ large. This is equivalent to
choosing $n$ such that 
\begin{align*}
\Delta_n\frac{\left(4\sqrt{m}+2\sqrt{2\log(\frac{K}{\xi})}\right)}{\sqrt{\AV{-}8\log(1-0.5\frac{\epsilon}{K})}} \leq r_u.
\end{align*}
Finally, the inequality $L_n(\hat\bta^{(k+1)})-L_n(\hat\bta_n)\leq \tau_{1}\frac{r^2}{4}$ follows using Lemma \ref{lem:continuity_relation} and $||\hat\bta_n-\hat\bta_n^{(k+1)}||_2\leq \left(\frac{\tau_1}{\tau_2}\right)^{\frac{1}{2}}\frac{r}{2}$. This completes the induction. To complete the proof of the Lemma, we now establish the claim.

\noindent\textbf{Proof of the claim}: Using equation (24) in Definition 32 of the main paper, let $\bta^*=\gamma\hat\bta_n^{(k+1)}+(1-\gamma)\hat\bta_n^{(k)}$ for some $\gamma\in[0,1]$ in the Taylor expansion of $L_n(\hat\bta_n^{(k+1)})$ up-to second order, and apply Cauchy–Schwarz inequality to get
\begin{align*}
    L_n(\hat\bta_n^{(k+1)})-L_n(\hat\bta_n) =& L_n(\hat\bta_n^{(k)}-\eta(\nabla L_n(\hat\bta_n^{(k)})+N_k))-L_n(\hat\bta_n)\\
    \leq& L_n(\hat\bta_n^{(k)})-L_n(\hat\bta_n)-\eta||\nabla L_n(\hat\bta_n^{(k)})||_2^2 + \eta||N_{n,k}||_2\cdot||\nabla L_n(\hat\bta_n^{(k)})||_2\\
    & + \frac{\eta^2}{2}\nabla L_n(\hat\bta_n^{(k)})^T H_n(\bta^*)\nabla L_n(\hat\bta_n^{(k)}) +\frac{\eta^2}{2}N_{n,k}^T H_n(\bta^*)N_{n,k}\\
    & + \eta^2 N_{n,k}^T H_n(\bta^*)\nabla L_n(\hat\bta_n^{(k)}).
\end{align*}
Furthermore, use Proposition 29
(3) and Cauchy–Schwarz inequality to get
\begin{align*}
    &\nabla L_n(\hat\bta_n^{(k)})^T H_n(\theta^*)\nabla L_n(\hat\bta_n^{(k)})\leq 2\tau_{2}||\nabla L_n(\hat\bta_n^{(k)})||_2^2,\\
    &N_{n,k}^T H_n(\bta^*)N_k \leq 2\tau_{2}||N_{n,k}||_2^2, ~~\text{and}\\
    &|N_{n,k}^T H_n(\bta^*)\nabla L_n(\hat\bta_n^{(k)})|\leq ||N_{n,k}^T||_2\cdot||H_n(\bta^*)\nabla L_n(\hat\bta_n^{(k)})||_2 \leq 2\tau_{2}||N_{n,k}^T||_2\cdot||\nabla L_n(\hat\bta_n^{(k)})||_2.
\end{align*}
These give the upper bound of $L_n(\hat\bta_n^{(k+1)})-L_n(\hat\bta_n)$ as follows:
\begin{align*}
    &L_n(\hat\bta_n^{(k+1)})-L_n(\hat\bta_n)\leq L_n(\hat\bta_n^{(k)})-L_n(\hat\bta_n) + Term_{n,k},
\end{align*}
where $Term_{n,k}$ is given by
\begin{align*}
    &-\eta||\nabla L_n(\hat\bta_n^{(k)})||_2^2 + \eta||N_{n,k}||_2\cdot||\nabla L_n(\hat\bta_n^{(k)})||_2 + \eta^2\tau_{2}||\nabla L_n(\hat\bta_n^{(k)})||_2^2 + \eta^2\tau_{2}||N_{n,k}||_2^2 \\
    & +2\eta^2 \tau_{2}||N_{n,k}^T||_2\cdot||\nabla L_n(\hat\bta_n^{(k)})||_2.
\end{align*}
Hence the condition $||\nabla L_n(\hat\bta_n^{(k)})||_2 \geq \sqrt{\frac{(1+2\eta\tau_{2})B||N_{n,k}||_2 + \eta\tau_{2}||N_{n,k}||_2^2}{1-\eta\tau_{2}}}$ implies that $Term_{n,k}\leq0$, that is
\begin{align*}
    &-\eta||\nabla L_n(\hat\bta_n^{(k)})||_2^2 + \eta||N_{n,k}||_2\cdot||\nabla L_n(\hat\bta_n^{(k)})||_2 + \eta^2\tau_{2}||\nabla L_n(\hat\bta_n^{(k)})||_2^2\\
    & + \eta^2\tau_{2}||N_{n,k}||_2^2 +2\eta^2 \tau_{2}||N_{n,k}^T||_2\cdot||\nabla L_n(\hat\bta_n^{(k)})||_2 \leq 0,
\end{align*}
and furthermore
\begin{align*}
    L_n(\hat\bta_n^{(k+1)})-L_n(\hat\bta_n) \leq L_n(\hat\bta_n^{(k)})-L_n(\hat\bta_n).
\end{align*}
This completes the proof. \QED

\clearpage

\section{Numerical calculations and additional experiments} \label{app:E}

In this section, we provide a Monte-Carlo approximation to the loss function and give calculation details for the Normal distribution used in numerical experiments.
\begin{align*}
\tilde L_n(\bta) =& 2\int_{\mathbb{R}} (\sqrt{f_{\bta}(x)}-\sqrt{g_n(x)})^2\mathrm{d}x = 2 \left[ 2-2\int_{\mathbb{R}} g_n(x) \left(\frac{f_{\bta}(x)}{g_n(x)}\right)^{\frac{1}{2}}\mathrm{d}x \right]\\
    \approx&  2 \left[ 2- \frac{2}{r_n}\sum_{i=1}^{r_n} \left(\frac{f_{\bta}(X_{n,i})}{g_n(X_{n,i})}\right)^{\frac{1}{2}} \right],
\end{align*}
where $\{X_{n,1}, \cdots X_{n, r_n}\}|(X_1, \cdots X_n)$ are i.i.d. $g_n(\cdot)$. Next, the gradient is given by
\begin{align*}
    \nabla \tilde L_n(\bta)
    = - \frac{2}{r_n} \sum_{i=1}^{r_n} \left(\frac{f_{\bta}(X_{n,i})}{g_n(X_{n,i})}\right)^{\frac{1}{2}} u_{\bta}(X_{n,i}),
\end{align*}
where $u_{\bta}(x)=\frac{\nabla f_{\bta}(x)}{f_{\bta}(x)}$, while the Hessian is given by
\begin{align*}
    \tilde H_n(\bta)
    =& \frac{1}{r_n} \sum_{i=1}^{r_n} \left(\frac{f_{\bta}(X_{n,i})}{g_n(X_{n,i})}\right)^{\frac{1}{2}} u_{\bta}(X_{n,i})\cdot u^T_{\bta}(X_{n,i}) \\
    & - \frac{2}{r_n} \sum_{i=1}^{r_n} \left(\frac{f_{\bta}(X_{n,i})}{g_n(X_{n,i})}\right)^{\frac{1}{2}} \frac{1}{f_{\bta}(X_{n,i})} H_f(X_{n,i}),
\end{align*}
where $H_f$ is Hessian of $f_{\bta}(\cdot)$. For the normal distribution, the gradient is given by
\begin{align*}
\nabla f_{\bta}(x) = f_{\bta}(x)\cdot 
\begin{bmatrix}
\frac{x-\mu}{\sigma^2} \\
        \frac{(x-\mu)^2-\sigma^2}{\sigma^3}
    \end{bmatrix}
,\quad u_{\bta}(x) = 
    \begin{bmatrix}
        \frac{x-\mu}{\sigma^2} \\
        \frac{(x-\mu)^2-\sigma^2}{\sigma^3}
    \end{bmatrix},
\end{align*}
and the Hessian of $f_{\bta}(\cdot)$ is given by
\begin{align*}
    H_f(x) = f_{\bta}(x)\cdot 
    \begin{bmatrix}
        \frac{(x-\mu)^2-\sigma^2}{\sigma^4} & \frac{(x-\mu)((x-\mu)^2-3\sigma^2)}{\sigma^5} \\
        \frac{(x-\mu)((x-\mu)^2-3\sigma^2)}{\sigma^5} & \frac{(x-\mu)^4-5\sigma^2(x-\mu)^2+2\sigma^4}{\sigma^6}
    \end{bmatrix}.
\end{align*}

\clearpage

\subsection{Additional figures and tables for PGD and PNR}

\begin{figure}[ht]
    \includegraphics[width=0.9\linewidth]{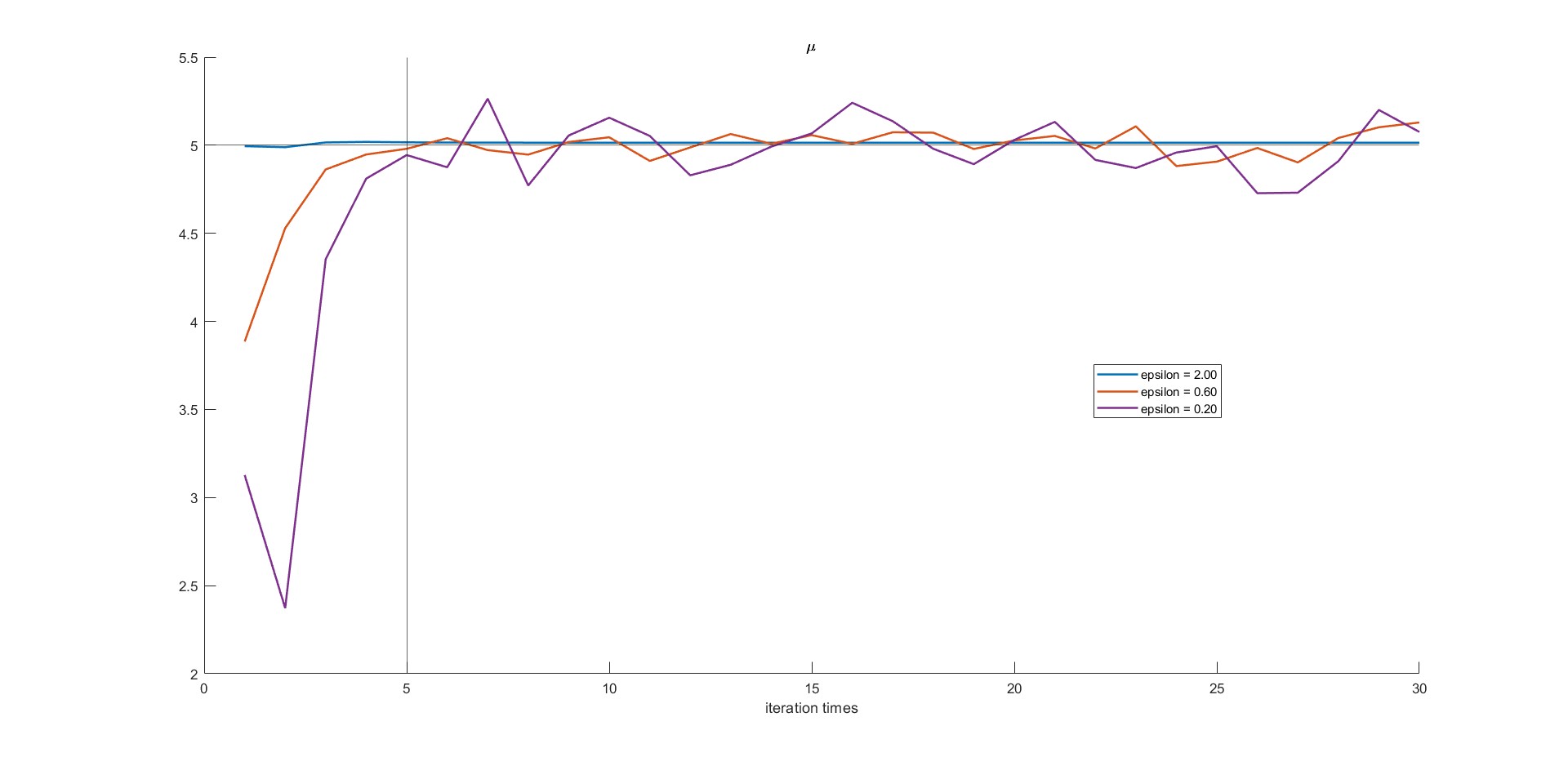}
    \includegraphics[width=0.9\linewidth]{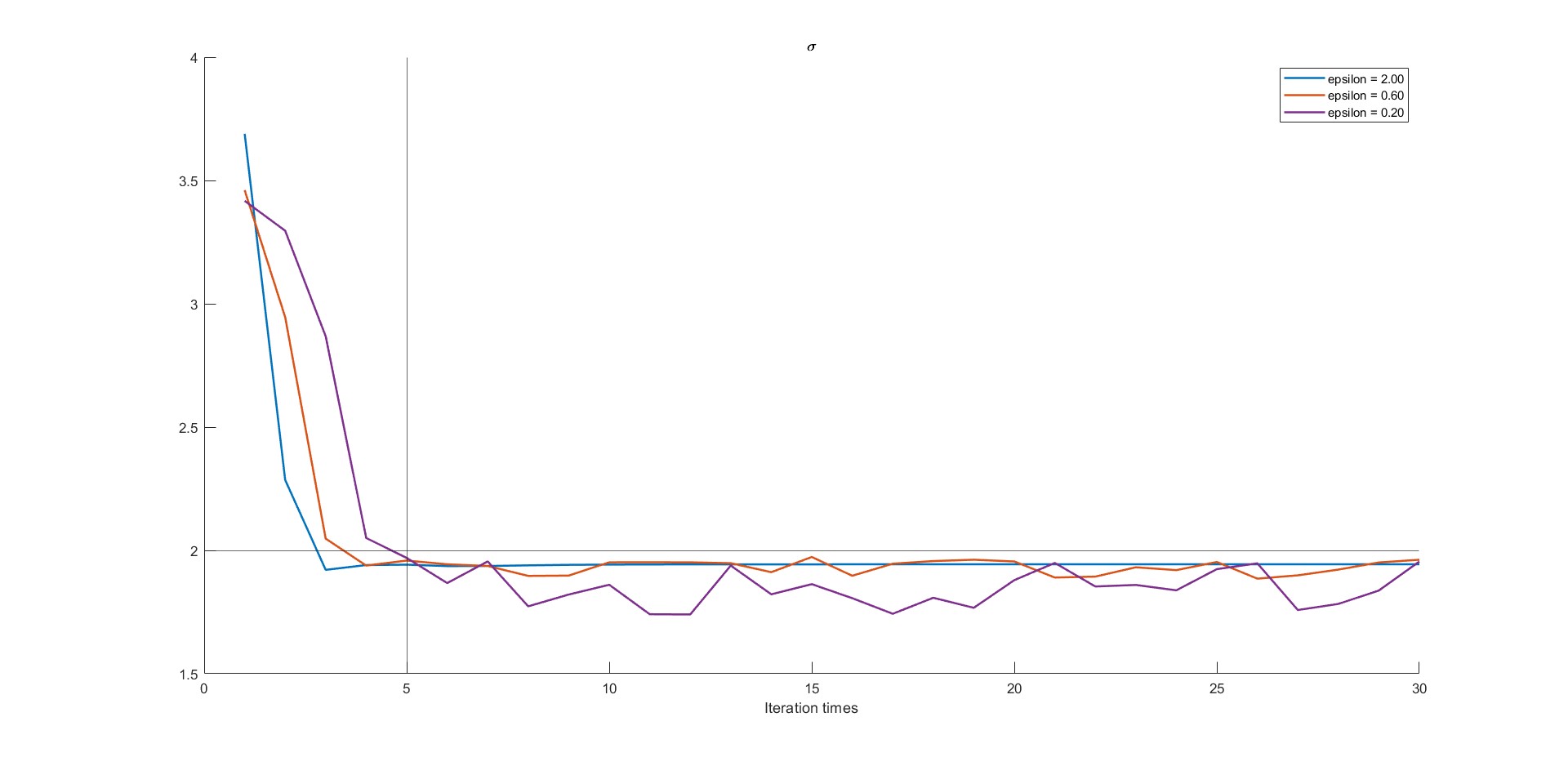}
    \caption{Private Newton's method path}
    \label{F3}
\end{figure}

\begin{figure}[ht]
    \includegraphics[width=0.9\linewidth]{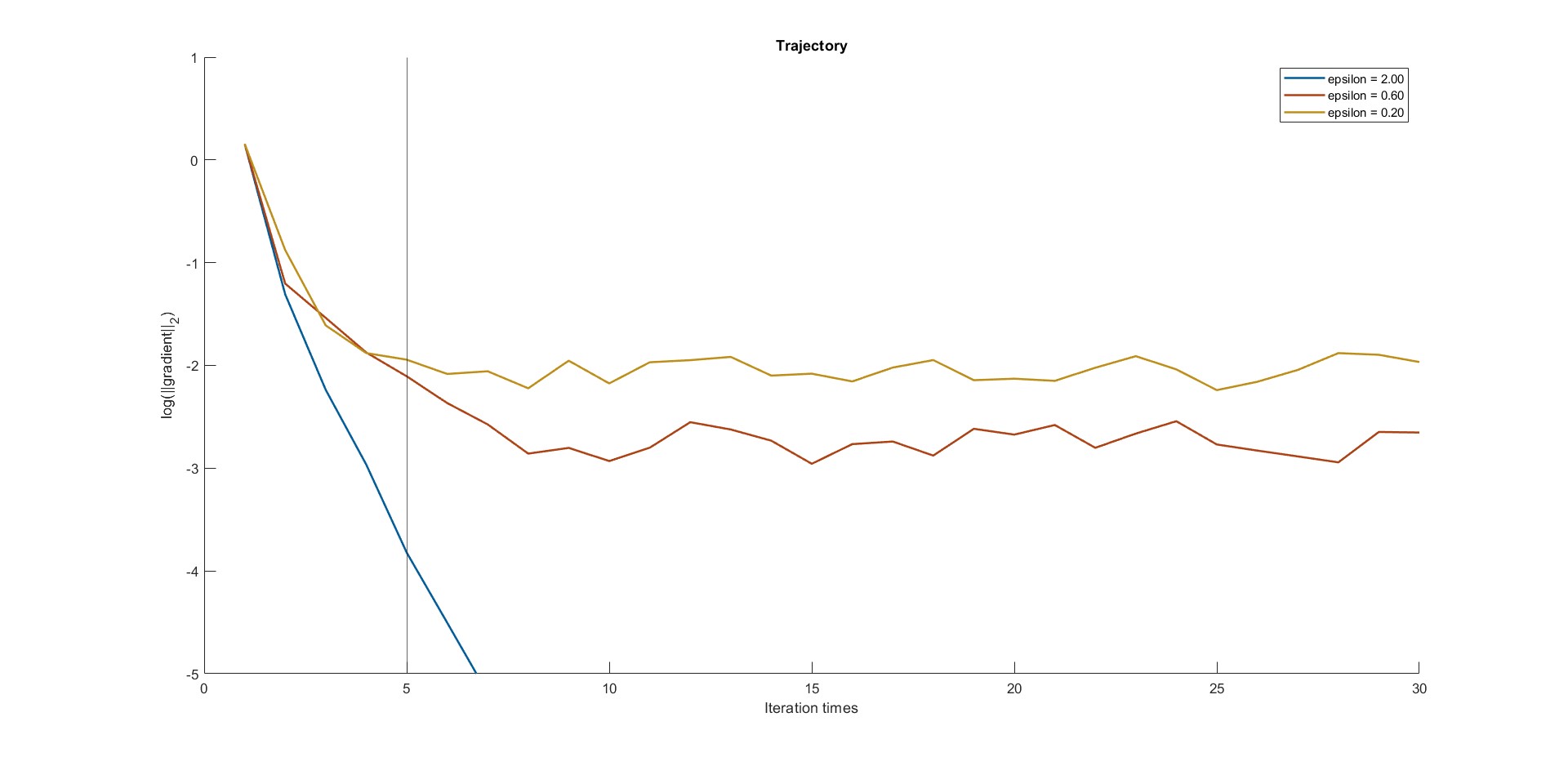}
    \caption{Private Newton's method trajectory}
    \label{F4}
\end{figure}

\begin{figure}[ht]
    \includegraphics[width=0.9\linewidth]{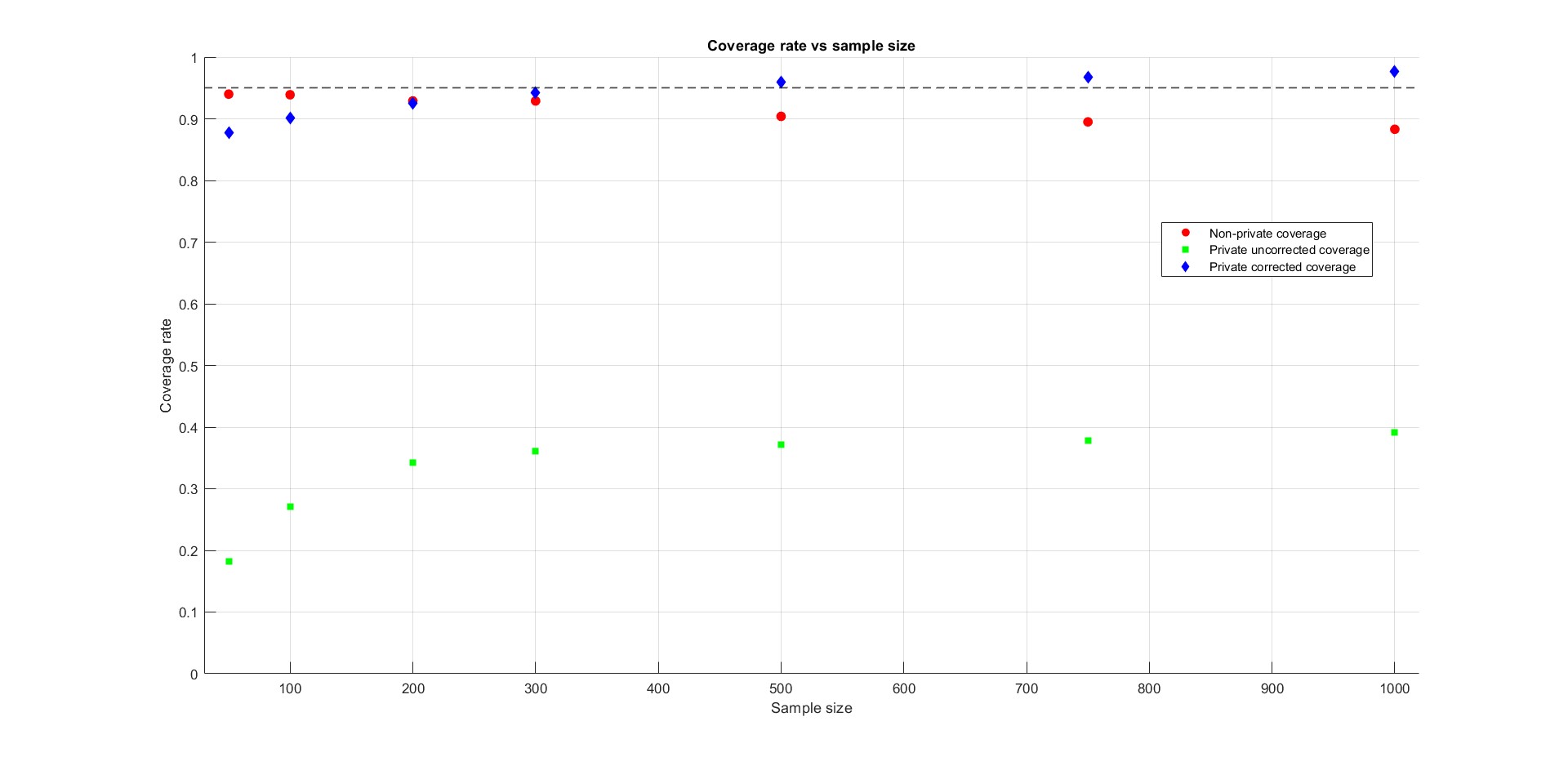}
    \caption{Private and non-private Newton's method 95\% confidence interval coverage}
    \label{F6}
\end{figure}

\begin{table}[ht]
\centering
\begin{tabular}{l|l|c|c|c}
\toprule
& & \multicolumn{3}{c}{$\epsilon$} \\
\cmidrule(lr){3-5}
& &  2.00 & 0.60 & 0.20 \\
\midrule
\multirow{2}{*}{Estimator} 
& $\mu$: Mean (Std. Error) & 5 (0.08) & 4.948 (0.332) & 4.868 (1.756) \\ \cmidrule(lr){2-5}
& $\sigma$: Mean (Std. Error) & 1.975 (0.076) & 1.987 (0.349) & 2.196 (1.99) \\
\midrule
\multirow{2}{*}{CI coverage for $\mu$} 
& Corrected & 0.883 & 0.977 & 0.95 \\ \cmidrule(lr){2-5}
& Uncorrected & 0.883 & 0.391 & 0.247 \\
\midrule
\multirow{2}{*}{CI coverage for $\sigma$} 
& Corrected & 0.739 & 0.913 & 0.904 \\ \cmidrule(lr){2-5}
& Uncorrected & 0.739 & 0.442 & 0.264 \\
\bottomrule
\end{tabular}
\caption{Results for different values of $\epsilon$ (Newton-Raphson). Sample size is 1000, $K=5$.}
\label{T2}
\end{table}

\begin{table}[ht]
\centering
\begin{tabular}{l|l|c|c|c}
\toprule
& & \multicolumn{3}{c}{$\lambda=1$, $\epsilon$} \\
\cmidrule(lr){3-5}
& &  Non-private & 1.20 & 0.40 \\
\midrule
\multirow{2}{*}{Estimator} 
& $\mu$: Mean (Std. Error) & 4.991 (0.083) & 4.986 (0.295) & 4.96 (2.009) \\ \cmidrule(lr){2-5}
& $\sigma$: Mean (Std. Error) & 1.984 (0.058) & 2.023 (0.21) & 2.038 (1.809) \\
\midrule
\multirow{2}{*}{CI coverage for $\mu$} 
& Corrected & 0.861 & 0.823 & 0.817 \\ \cmidrule(lr){2-5}
& Uncorrected & 0.861 & 0.371 & 0.292 \\
\midrule
\multirow{2}{*}{CI coverage for $\sigma$} 
& Corrected & 0.819 & 0.93 & 0.913 \\ \cmidrule(lr){2-5}
& Uncorrected & 0.819 & 0.332 & 0.266 \\
\bottomrule
\end{tabular}
\caption{Results for different values of $\epsilon$ (Gradient descent). Sample size is 1000, $K=50$.}
\label{T3}
\end{table}

\begin{table}[ht]
\centering
\begin{tabular}{l|l|c|c|c}
\toprule
& & \multicolumn{3}{c}{$\lambda=1$, $\epsilon$} \\
\cmidrule(lr){3-5}
& &  Non-private & 1.20 & 0.40 \\
\midrule
\multirow{2}{*}{Estimator} 
& $\mu$: Mean (Std. Error) & 5 (0.08) & 4.927 (0.727) & 4.76 (5.962) \\ \cmidrule(lr){2-5}
& $\sigma$: Mean (Std. Error) & 1.975 (0.076) & 2.107 (1.265) & 2.513 (7.652) \\
\midrule
\multirow{2}{*}{CI coverage for $\mu$} 
& Corrected & 0.883 & 0.965 & 0.934 \\ \cmidrule(lr){2-5}
& Uncorrected & 0.883 & 0.288 & 0.208 \\
\midrule
\multirow{2}{*}{CI coverage for $\sigma$} 
& Corrected & 0.739 & 0.918 & 0.894 \\ \cmidrule(lr){2-5}
& Uncorrected & 0.739 & 0.309 & 0.235 \\
\bottomrule
\end{tabular}
\caption{Results for different values of epsilon (Newton-Raphson). Sample size is 1000, $K=5$.}
\label{T4}
\end{table}

\begin{table}[ht]
\centering
\scriptsize
\begin{tabular}{c|c|c|c|c|c}
\toprule
& \multicolumn{5}{c}{Contamination percentage $\alpha$} \\
\cmidrule(lr){2-6}
& 0\% & 5\% & 10\% & 20\% & 30\% \\
\midrule
MLE (Std. Error) & 5.001 (0.063) & 5.241 (0.061) & 5.476 (0.059) & 5.952 (0.056) & 6.422 (0.052) \\
\midrule
PMHDE $\epsilon=2$ (Std. Error) & 5 (0.08) & 5.174 (0.085) & 5.309 (0.087) & 5.555 (0.091) & 5.778 (0.096) \\
\midrule
PMHDE $\epsilon=0.6$ (Std. Error) & 4.952 (0.326) & 5.119 (0.333) & 5.252 (0.341) & 5.472 (0.38) & 5.646 (0.42) \\
\midrule
PMHDE $\epsilon=0.2$ (Std. Error) & 4.942 (13.391) & 4.905 (5.895) & 5.054 (2.78) & 5.349 (4.023) & 5.169 (15.035) \\
\bottomrule
\end{tabular}
\caption{Contamination results, Newton-Raphson, sample size is 1000, $\mu=5$}
\label{T6}
\end{table} 

\clearpage
\subsection[Estimation and coverage rate for PDP]{Estimation and coverage rate for $(\la,\ep)-$PDP}
In this section, we provide PMHDE and coverage rates for PGD and PNR algorithms for different values of $\lambda$.

\begin{table}[ht]
\centering
\begin{tabular}{l|l|c|c|c}
\toprule
& & \multicolumn{3}{c}{$\lambda=-0.1$, $\epsilon$} \\
\cmidrule(lr){3-5}
& &  Non-private & 1.20 & 0.40 \\
\midrule
\multirow{2}{*}{Estimator} 
& $\mu$: Mean (Std. Error) & 4.991 (0.083) & 4.992 (0.212) & 4.979 (0.454) \\ \cmidrule(lr){2-5}
& $\sigma$: Mean (Std. Error) & 1.984 (0.058) & 2.001 (0.152) & 2.045 (0.291) \\
\midrule
\multirow{2}{*}{CI coverage for $\mu$} 
& Corrected & 0.861 & 0.836 & 0.82 \\ \cmidrule(lr){2-5}
& Uncorrected & 0.861 & 0.468 & 0.33 \\
\midrule
\multirow{2}{*}{CI coverage for $\sigma$} 
& Corrected & 0.819 & 0.931 & 0.927 \\ \cmidrule(lr){2-5}
& Uncorrected & 0.819 & 0.428 & 0.296 \\
\bottomrule
\end{tabular}
\caption{Results for different values of $\epsilon$ (Gradient descent). Sample size is 1000, $K = 50$.}
\label{T7}
\end{table}

\begin{table}[ht]
\centering
\begin{tabular}{l|l|c|c|c}
\toprule
& & \multicolumn{3}{c}{$\lambda = -0.1$, $\epsilon$} \\
\cmidrule(lr){3-5}
& &  Non-private & 1.20 & 0.40 \\
\midrule
\multirow{2}{*}{Estimator} 
& $\mu$: Mean (Std. Error) & 5 (0.08) & 4.955 (0.348) & 4.827 (3.731) \\ \cmidrule(lr){2-5}
& $\sigma$: Mean (Std. Error) & 1.975 (0.076) & 1.992 (0.353) & 2.223 (1.554) \\
\midrule
\multirow{2}{*}{CI coverage for $\mu$} 
& Corrected & 0.883 & 0.977 & 0.948 \\ \cmidrule(lr){2-5}
& Uncorrected & 0.883 & 0.382 & 0.234 \\
\midrule
\multirow{2}{*}{CI coverage for $\sigma$} 
& Corrected & 0.739 & 0.917 & 0.902 \\ \cmidrule(lr){2-5}
& Uncorrected & 0.739 & 0.41 & 0.264 \\
\bottomrule
\end{tabular}
\caption{Results for different values of $\epsilon$ (Newton). Sample size is 1000, $K = 5$.}
\label{T8}
\end{table}

\begin{table}[ht]
\centering
\begin{tabular}{l|l|c|c|c}
\toprule
& & \multicolumn{3}{c}{$\lambda = 0.5$, $\epsilon$} \\
\cmidrule(lr){3-5}
& &  Non-private & 1.20 & 0.40 \\
\midrule
\multirow{2}{*}{Estimator} 
& $\mu$: Mean (Std. Error) & 4.991 (0.083) & 4.991 (0.256) & 4.976 (0.549) \\ \cmidrule(lr){2-5}
& $\sigma$: Mean (Std. Error) & 1.984 (0.058) & 2.013 (0.18) & 2.056 (0.322) \\
\midrule
\multirow{2}{*}{CI coverage for $\mu$} 
& Corrected & 0.861 & 0.826 & 0.817 \\ \cmidrule(lr){2-5}
& Uncorrected & 0.861 & 0.396 & 0.306 \\
\midrule
\multirow{2}{*}{CI coverage for $\sigma$} 
& Corrected & 0.819 & 0.931 & 0.922 \\ \cmidrule(lr){2-5}
& Uncorrected & 0.819 & 0.379 & 0.277 \\
\bottomrule
\end{tabular}
\caption{Results for different values of $\epsilon$ (Gradient descent). Sample size is 1000, $K = 50$.}
\label{T9}
\end{table}

\begin{table}[ht]
\centering
\begin{tabular}{l|l|c|c|c}
\toprule
& & \multicolumn{3}{c}{$\lambda = 0.5$, $\epsilon$} \\
\cmidrule(lr){3-5}
& &  Non-private & 1.20 & 0.40 \\
\midrule
\multirow{2}{*}{Estimator} 
& $\mu$: Mean (Std. Error) & 5 (0.08) & 4.942 (0.483) & 4.808 (2.064) \\ \cmidrule(lr){2-5}
& $\sigma$: Mean (Std. Error) & 1.975 (0.076) & 2.03 (0.523) & 2.297 (2.313) \\
\midrule
\multirow{2}{*}{CI coverage for $\mu$} 
& Corrected & 0.883 & 0.972 & 0.94 \\ \cmidrule(lr){2-5}
& Uncorrected & 0.883 & 0.322 & 0.222 \\
\midrule
\multirow{2}{*}{CI coverage for $\sigma$} 
& Corrected & 0.739 & 0.915 & 0.899 \\ \cmidrule(lr){2-5}
& Uncorrected & 0.739 & 0.355 & 0.254 \\
\bottomrule
\end{tabular}
\caption{Results for different values of $\epsilon$ (Newton). Sample size is 1000, $K = 5$.}
\label{T10}
\end{table}

\clearpage
\subsection{Additional results for HDP and robustness evaluation}
In this section, we provide additional results for PMHDE for sample sizes 200, 300, and 500 for both PGD and PNR algorithms. As explained in Section 4 of the main paper,
when $n$ and $\ep$ are both small, the algorithms can produce aberrant values, reducing their usefulness. For this reason, we use only estimates within the lower 0.7\% and upper 99.5\% percentiles of a Gaussian distribution with non-private $\hat{\mu}_n$ and $ \hat{\sigma}_n$. All the Tables in this section are based on such a thresholding strategy. Since the confidence intervals are unaffected by thresholding, we retain all the simulation experiments for constructing the confidence intervals.

Tables \ref{T11} and \ref{T12}  provide the estimators and the coverage rates for sample size 200, while Tables \ref{T13} and \ref{T14} provide the behavior of PMHDE under contamination for the sample size 200. The corresponding Tables for sample size 300 and  500 are given in Tables \ref{T15}, \ref{T16}, \ref{T17}, \ref{T18}, \ref{T19}, \ref{T20}, \ref{T21}, \ref{T22} respectively.

{\bf{Sample size 200:}}

\begin{table}[ht]
\centering
\begin{tabular}{l|l|c|c|c}
\toprule
& & \multicolumn{3}{c}{$\epsilon$} \\
\cmidrule(lr){3-5}
& &  2.00 & 0.60 & 0.20 \\
\midrule
\multirow{2}{*}{Estimator} 
& $\mu$: Mean (Std. Error) & 4.992 (0.153) & 4.978 (0.538) & 4.844 (1.22) \\ \cmidrule(lr){2-5}
& $\sigma$: Mean (Std. Error) & 1.952 (0.104) & 2.036 (0.588) & 1.744 (6.164) \\
\midrule
\multirow{2}{*}{CI coverage for $\mu$} 
& Corrected & 0.921 & 0.839 & 0.626 \\ \cmidrule(lr){2-5}
& Uncorrected & 0.921 & 0.51 & 0.27 \\
\midrule
\multirow{2}{*}{CI coverage for $\sigma$} 
& Corrected & 0.846 & 0.921 & 0.6 \\ \cmidrule(lr){2-5}
& Uncorrected & 0.846 & 0.484 & 0.216 \\
\bottomrule
\end{tabular}
\caption{Results for different values of $\epsilon$ (Gradient descent). Sample size is 200, $K = 50$.}
\label{T11}
\end{table}

\begin{table}[ht]
\centering
\begin{tabular}{l|l|c|c|c}
\toprule
& & \multicolumn{3}{c}{$\epsilon$} \\
\cmidrule(lr){3-5}
& &  2.00 & 0.60 & 0.20 \\
\midrule
\multirow{2}{*}{Estimator} 
& $\mu$: Mean (Std. Error) & 4.993 (0.148) & 4.819 (1.242) & 4.703 (1.972) \\ \cmidrule(lr){2-5}
& $\sigma$: Mean (Std. Error) & 1.93 (0.139) & 2.465 (2.452) & 3.118 (3.792) \\
\midrule
\multirow{2}{*}{CI coverage for $\mu$} 
& Corrected & 0.929 & 0.925 & 0.888 \\ \cmidrule(lr){2-5}
& Uncorrected & 0.929 & 0.343 & 0.157 \\
\midrule
\multirow{2}{*}{CI coverage for $\sigma$} 
& Corrected & 0.771 & 0.879 & 0.829 \\ \cmidrule(lr){2-5}
& Uncorrected & 0.771 & 0.359 & 0.169 \\
\bottomrule
\end{tabular}
\caption{Results for different values of $\epsilon$ (Newton). Sample size is 200, $K = 5$.}
\label{T12}
\end{table}

\begin{table}[ht]
\centering
\scriptsize
\begin{tabular}{c|c|c|c|c|c}
\toprule
& \multicolumn{5}{c}{Contamination percentage $\alpha$} \\
\cmidrule(lr){2-6}
& 0\% & 5\% & 10\% & 20\% & 30\% \\
\midrule
MLE (Std. Error) & 5 (0.142) & 5.238 (0.138) & 5.483 (0.134) & 5.948 (0.127) & 6.418 (0.119) \\
\midrule
PMHDE $\epsilon=2$ (Std. Error) & 4.99 (0.153) & 5.169 (0.155) & 5.301 (0.158) & 5.562 (0.165) & 5.783 (0.172) \\
\midrule
PMHDE $\epsilon=0.6$ (Std. Error) & 4.977 (0.521) & 5.141 (0.534) & 5.283 (0.54) & 5.508 (0.548) & 5.721 (0.576) \\
\midrule
PMHDE $\epsilon=0.2$ (Std. Error) & 4.856 (1.228) & 5.024 (1.16) & 5.13 (1.22) & 5.313 (1.215) & 5.463 (1.306) \\
\bottomrule
\end{tabular}
\caption{Contamination results, gradient descent, sample size is 200, $\mu=5$.}
\label{T13}
\end{table}

\begin{table}[ht]
\centering
\scriptsize
\begin{tabular}{c|c|c|c|c|c}
\toprule
& \multicolumn{5}{c}{Contamination percentage $\alpha$} \\
\cmidrule(lr){2-6}
& 0\% & 5\% & 10\% & 20\% & 30\% \\
\midrule
MLE (Std. Error) & 5 (0.142) & 5.238 (0.138) & 5.483 (0.134) & 5.948 (0.127) & 6.418 (0.119) \\
\midrule
PMHDE $\epsilon=2$ (Std. Error) & 4.991 (0.148) & 5.174 (0.15) & 5.308 (0.151) & 5.578 (0.155) & 5.823 (0.156) \\
\midrule
PMHDE $\epsilon=0.6$ (Std. Error) & 4.809 (1.275) & 4.974 (1.221) & 5.021 (1.216) & 5.17 (1.308) & 5.304 (1.375) \\
\midrule
PMHDE $\epsilon=0.2$ (Std. Error) & 4.744 (1.964) & 4.815 (2.071) & 4.866 (2.116) & 4.954 (2.21) & 5.002 (2.206) \\
\bottomrule
\end{tabular}
\caption{Contamination results, Newton, sample size is 200, $\mu=5$.}
\label{T14}
\end{table}

\clearpage
{\bf{Sample size 300:}}

\begin{table}[ht]
\centering
\begin{tabular}{l|l|c|c|c}
\toprule
& & \multicolumn{3}{c}{$\epsilon$} \\
\cmidrule(lr){3-5}
& &  2.00 & 0.60 & 0.20 \\
\midrule
\multirow{2}{*}{Estimator} 
& $\mu$: Mean (Std. Error) & 4.989 (0.128) & 4.979 (0.405) & 4.926 (0.866) \\ \cmidrule(lr){2-5}
& $\sigma$: Mean (Std. Error) & 1.962 (0.086) & 2.023 (0.338) & 1.952 (1.971) \\
\midrule
\multirow{2}{*}{CI coverage for $\mu$} 
& Corrected & 0.918 & 0.837 & 0.733 \\ \cmidrule(lr){2-5}
& Uncorrected & 0.918 & 0.493 & 0.317 \\
\midrule
\multirow{2}{*}{CI coverage for $\sigma$} 
& Corrected & 0.85 & 0.944 & 0.763 \\ \cmidrule(lr){2-5}
& Uncorrected & 0.85 & 0.477 & 0.29 \\
\bottomrule
\end{tabular}
\caption{Results for different values of $\epsilon$ (Gradient descent). Sample size is 300, $K = 50$.}
\label{T15}
\end{table}

\begin{table}[ht]
\centering
\begin{tabular}{l|l|c|c|c}
\toprule
& & \multicolumn{3}{c}{$\epsilon$} \\
\cmidrule(lr){3-5}
& &  2.00 & 0.60 & 0.20 \\
\midrule
\multirow{2}{*}{Estimator} 
& $\mu$: Mean (Std. Error) & 4.993 (0.123) & 4.863 (0.934) & 4.697 (1.722) \\ \cmidrule(lr){2-5}
& $\sigma$: Mean (Std. Error) & 1.946 (0.113) & 2.233 (1.442) & 2.785 (6.688) \\
\midrule
\multirow{2}{*}{CI coverage for $\mu$} 
& Corrected & 0.929 & 0.943 & 0.897 \\ \cmidrule(lr){2-5}
& Uncorrected & 0.929 & 0.361 & 0.182 \\
\midrule
\multirow{2}{*}{CI coverage for $\sigma$} 
& Corrected & 0.777 & 0.9 & 0.836 \\ \cmidrule(lr){2-5}
& Uncorrected & 0.777 & 0.377 & 0.204 \\
\bottomrule
\end{tabular}
\caption{Results for different values of $\epsilon$ (Newton). The Sample size is 300, $K = 5$.}
\label{T16}
\end{table}

\begin{table}[ht]
\centering
\scriptsize
\begin{tabular}{c|c|c|c|c|c}
\toprule
& \multicolumn{5}{c}{Contamination percentage $\alpha$} \\
\cmidrule(lr){2-6}
& 0\% & 5\% & 10\% & 20\% & 30\% \\
\midrule
MLE (std.error) & 4.999 (0.115) & 5.238 (0.113) & 5.474 (0.11) & 5.938 (0.103) & 6.433 (0.096) \\
\midrule
PMHDE $\epsilon=2$ (Std. Error) & 4.989 (0.128) & 5.164 (0.132) & 5.309 (0.135) & 5.542 (0.14) & 5.756 (0.149) \\
\midrule
PMHDE $\epsilon=0.6$ (Std. Error) & 4.992 (0.408) & 5.161 (0.398) & 5.291 (0.408) & 5.514 (0.416) & 5.717 (0.434) \\
\midrule
PMHDE $\epsilon=0.2$ (Std. Error) & 4.925 (0.858) & 5.095 (0.868) & 5.241 (0.844) & 5.396 (0.905) & 5.562 (0.932) \\
\bottomrule
\end{tabular}
\caption{Contamination results, gradient descent, sample size is 300, $\mu=5$.}
\label{T17}
\end{table}

\begin{table}[ht]
\centering
\scriptsize
\begin{tabular}{c|c|c|c|c|c}
\toprule
& \multicolumn{5}{c}{Contamination percentage $\alpha$} \\
\cmidrule(lr){2-6}
& 0\% & 5\% & 10\% & 20\% & 30\% \\
\midrule
MLE (Std. Error) & 4.999 (0.115) & 5.238 (0.113) & 5.474 (0.11) & 5.938 (0.103) & 6.433 (0.096) \\
\midrule
PMHDE $\epsilon=2$ (Std. Error) & 4.993 (0.124) & 5.173 (0.127) & 5.32 (0.129) & 5.562 (0.132) & 5.806 (0.135) \\
\midrule
PMHDE $\epsilon=0.6$ (Std. Error) & 4.856 (0.947) & 5.009 (0.922) & 5.114 (0.882) & 5.296 (0.991) & 5.42 (1.055) \\
\midrule
PMHDE $\epsilon=0.2$ (Std. Error) & 4.724 (1.724) & 4.863 (1.729) & 4.904 (1.782) & 5.023 (1.849) & 5.063 (1.894) \\
\bottomrule
\end{tabular}
\caption{Contamination results, Newton, sample size is 300, $\mu=5$.}
\label{T18}
\end{table}

\clearpage
{\bf{Sample size 500:}}

\begin{table}[ht]
\centering
\begin{tabular}{l|l|c|c|c}
\toprule
& & \multicolumn{3}{c}{$\epsilon$} \\
\cmidrule(lr){3-5}
& &  2.00 & 0.60 & 0.20 \\
\midrule
\multirow{2}{*}{Estimator} 
& $\mu$: Mean (Std. Error) & 4.989 (0.106) & 4.986 (0.296) & 4.951 (0.559) \\ \cmidrule(lr){2-5}
& $\sigma$: Mean (Std. Error) & 1.973 (0.072) & 2.016 (0.213) & 2.058 (0.799) \\
\midrule
\multirow{2}{*}{CI coverage for $\mu$} 
& Corrected & 0.892 & 0.842 & 0.798 \\ \cmidrule(lr){2-5}
& Uncorrected & 0.892 & 0.489 & 0.326 \\
\midrule
\multirow{2}{*}{CI coverage for $\sigma$} 
& Corrected & 0.848 & 0.941 & 0.888 \\ \cmidrule(lr){2-5}
& Uncorrected & 0.848 & 0.44 & 0.301 \\
\bottomrule
\end{tabular}
\caption{Results for different values of $\epsilon$ (gradient descent). Sample size is 500, $K = 50$.}
\label{T19}
\end{table}

\begin{table}[ht]
\centering
\begin{tabular}{l|l|c|c|c}
\toprule
& & \multicolumn{3}{c}{$\epsilon$} \\
\cmidrule(lr){3-5}
& &  2.00 & 0.60 & 0.20 \\
\midrule
\multirow{2}{*}{Estimator} 
& $\mu$: Mean (Std. Error) & 4.995 (0.102) & 4.929 (0.631) & 4.785 (1.322) \\ \cmidrule(lr){2-5}
& $\sigma$: Mean (Std. Error) & 1.959 (0.094) & 2.089 (1.514) & 2.554 (2.009) \\
\midrule
\multirow{2}{*}{CI coverage for $\mu$} 
& Corrected & 0.904 & 0.96 & 0.916 \\ \cmidrule(lr){2-5}
& Uncorrected & 0.904 & 0.371 & 0.213 \\
\midrule
\multirow{2}{*}{CI coverage for $\sigma$} 
& Corrected & 0.767 & 0.908 & 0.874 \\ \cmidrule(lr){2-5}
& Uncorrected & 0.767 & 0.391 & 0.243 \\
\bottomrule
\end{tabular}
\caption{Results for different values of $\epsilon$. Sample size is 500 (Newton), $K = 5$.}
\label{T20}
\end{table}

\begin{table}[ht]
\centering
\scriptsize
\begin{tabular}{c|c|c|c|c|c}
\toprule
& \multicolumn{5}{c}{Contamination percentage $\alpha$} \\
\cmidrule(lr){2-6}
& 0\% & 5\% & 10\% & 20\% & 30\% \\
\midrule
MLE (Std. Error) & 4.998 (0.09) & 5.24 (0.087) & 5.473 (0.086) & 5.946 (0.08) & 6.424 (0.076) \\
\midrule
PMHDE $\epsilon=2$ (Std. Error) & 4.989 (0.107) & 5.157 (0.109) & 5.293 (0.115) & 5.526 (0.119) & 5.738 (0.126) \\
\midrule
PMHDE $\epsilon=0.6$ (Std. Error) & 4.986 (0.293) & 5.157 (0.294) & 5.291 (0.299) & 5.514 (0.302) & 5.719 (0.312) \\
\midrule
PMHDE $\epsilon=0.2$ (Std. Error) & 4.987 (0.556) & 5.146 (0.553) & 5.256 (0.571) & 5.464 (0.583) & 5.65 (0.636) \\
\bottomrule
\end{tabular}
\caption{Contamination results, gradient descent, sample size is 500, $\mu=5$.}
\label{T21}
\end{table}

\begin{table}[ht]
\centering
\scriptsize
\begin{tabular}{c|c|c|c|c|c}
\toprule
& \multicolumn{5}{c}{Contamination percentage $\alpha$} \\
\cmidrule(lr){2-6}
& 0\% & 5\% & 10\% & 20\% & 30\% \\
\midrule
MLE (Std. Error) & 4.998 (0.09) & 5.24 (0.087) & 5.473 (0.086) & 5.946 (0.08) & 6.424 (0.076) \\
\midrule
PMHDE $\epsilon=2$ (Std. Error) & 4.995 (0.103) & 5.169 (0.105) & 5.309 (0.11) & 5.553 (0.113) & 5.793 (0.117) \\
\midrule
PMHDE $\epsilon=0.6$ (Std. Error) & 4.915 (0.631) & 5.069 (0.606) & 5.2 (0.629) & 5.383 (0.659) & 5.546 (0.716) \\
\midrule
PMHDE $\epsilon=0.2$ (Std. Error) & 4.777 (1.316) & 4.923 (1.306) & 4.995 (1.384) & 5.151 (1.386) & 5.259 (1.419) \\
\bottomrule
\end{tabular}
\caption{Contamination results, Newton, sample size is 500, $\mu=5$.}
\label{T22}
\end{table}

\clearpage
\subsection{Forward and reverse privacy calibrations}
\begin{figure}[H]
\centering
\begin{minipage}[t]{0.485\textwidth}
\centering
\includegraphics[width=\linewidth,trim=0 0 0 76bp,clip]{gaussian_forward_calibration.jpg}
\par\smallskip{\small (a) Gaussian forward: attained $h$ vs. target $\delta$}
\end{minipage}\hfill
\begin{minipage}[t]{0.485\textwidth}
\centering
\includegraphics[width=\linewidth,trim=0 0 0 76bp,clip]{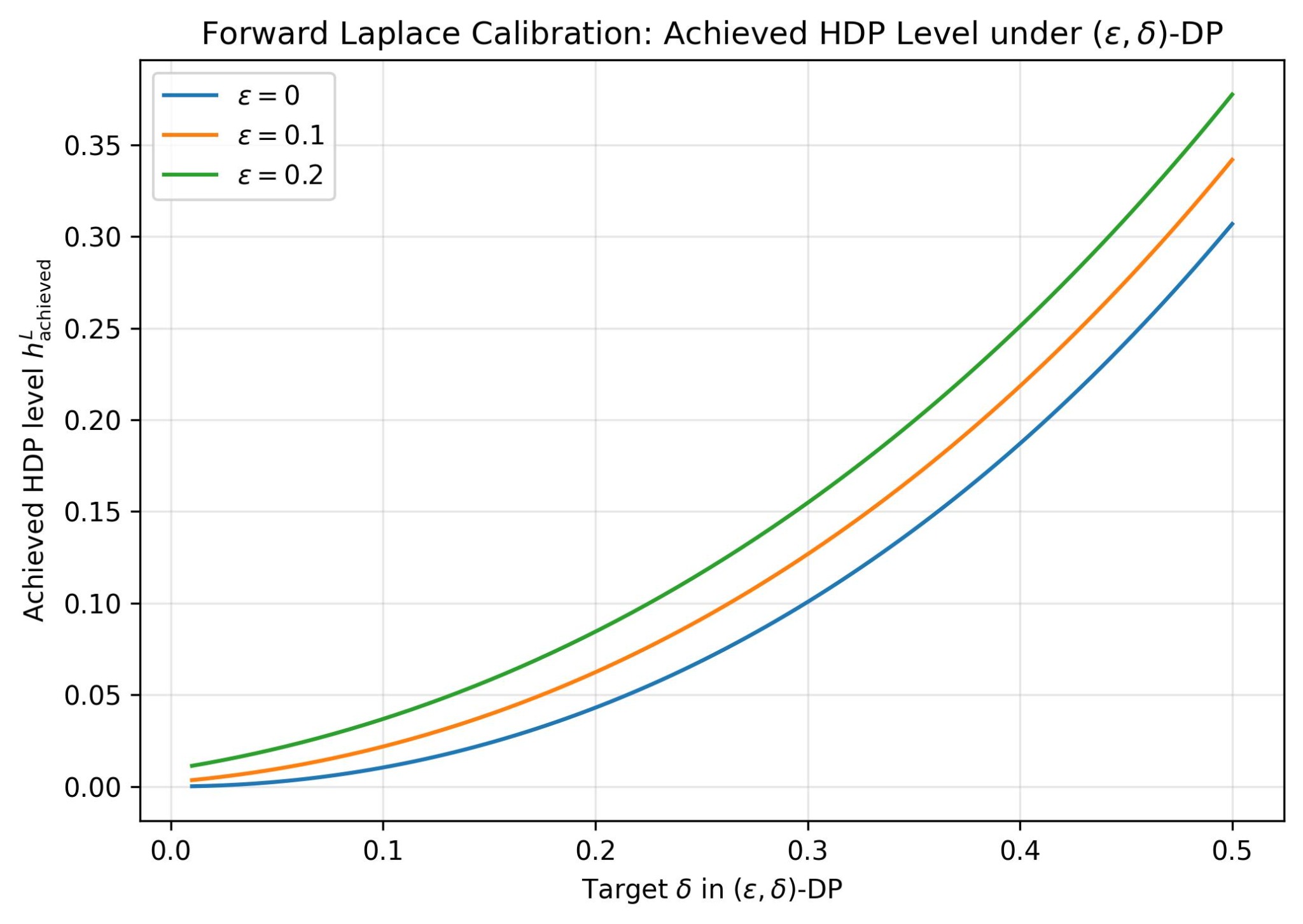}
\par\smallskip{\small (b) Laplace forward: attained $h$ vs. target $\delta$}
\end{minipage}
\par\medskip
\begin{minipage}[t]{0.485\textwidth}
\centering
\includegraphics[width=\linewidth,trim=0 0 0 76bp,clip]{gaussian_reverse_calibration.jpg}
\par\smallskip{\small (c) Gaussian reverse: attained $\delta$ vs. target $h$}
\end{minipage}\hfill
\begin{minipage}[t]{0.485\textwidth}
\centering
\includegraphics[width=\linewidth,trim=0 0 0 76bp,clip]{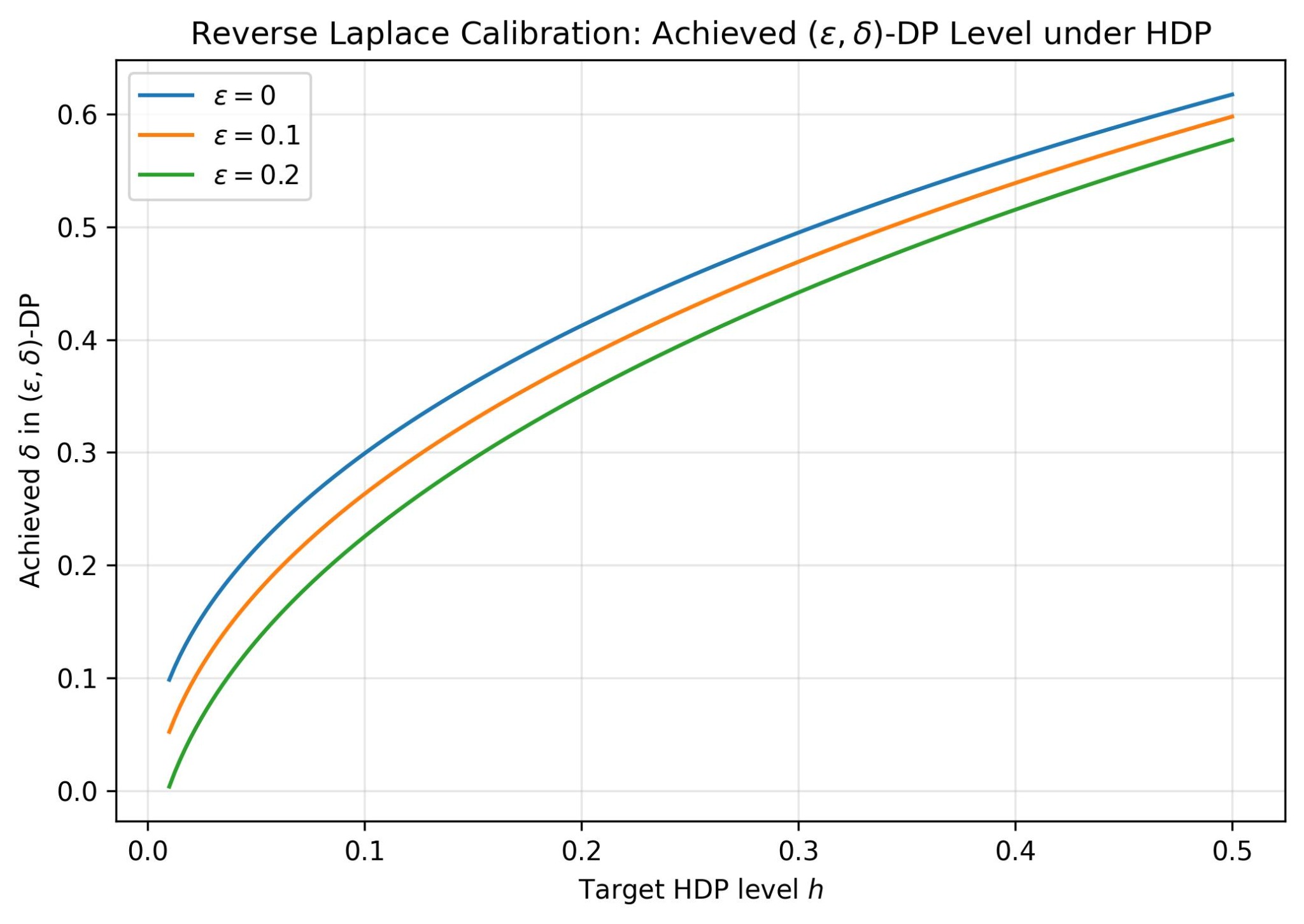}
\par\smallskip{\small (d) Laplace reverse: attained $\delta$ vs. target $h$}
\end{minipage}
\caption{Exact forward and reverse privacy calibrations for the additive Gaussian and scalar Laplace families, using the noise level calibrated to the target guarantee. Forward panels show the achieved Hellinger level at target $(\epsilon,\delta)$-DP; reverse panels show the achieved $\delta$ for target $h$-HDP at a fixed $\epsilon$. The four supplied numerical curves are unchanged.}
\label{fig:calibration-full}
\end{figure}

\clearpage

\section{Uniform rare-event bounds and large-deviation consequences}
\label{app:hdp-rare-events}

We derive the consequences of the Hellinger constraint for sequences
of neighboring output laws. These statements concern the actual output
probabilities. They are not applications of an empirical-divergence
large-deviation theorem. The partition viewpoint is related to
\citet{Bei01} and \citet{BiauGyrofi2005}, whose sampling experiments
involve empirical distributions on partitions.

For each $n$, let $M_n$ satisfy $h_n$-HDP, with $0\leq h_n\leq2$,
and let $\mathcal I_n$ be a nonempty collection of triples
$\iota=(D,D',A)$ consisting of an ordered neighboring pair and a
measurable output event. Write
\[
p_{n,\iota}=P_{n,D}(A),\qquad
q_{n,\iota}=P_{n,D'}(A).
\]
The event classes and neighboring pairs may vary with $n$.
No restriction on the cardinality of $\mathcal I_n$ is needed.

\begin{proposition}[Uniform probability comparison]
\label{prop:hdp-app-uniform}
For every $n$ and $\iota\in\mathcal I_n$,
\begin{equation}
(\sqrt{q_{n,\iota}}-\sqrt{h_n})_+^2
\leq p_{n,\iota}
\leq\min\{1,(\sqrt{q_{n,\iota}}+\sqrt{h_n})^2\}.
\label{eq:hdp-app-square-root}
\end{equation}
If $\inf_{\iota\in\mathcal I_n}q_{n,\iota}\geq a_n>0$, then
\begin{equation}
\sup_{\iota\in\mathcal I_n}
\left|\frac{p_{n,\iota}}{q_{n,\iota}}-1\right|
\leq2\sqrt{\frac{h_n}{a_n}}+\frac{h_n}{a_n}.
\label{eq:hdp-app-uniform-relative}
\end{equation}
In particular, $h_n/a_n\to0$ implies uniform convergence of these
probability ratios to one.
\end{proposition}

\begin{proof}
Apply the binary-partition inequality
\eqref{eq:hdp-op-binary} to each triple and drop the complementary-event
term. This yields
$|\sqrt{p_{n,\iota}}-\sqrt{q_{n,\iota}}|\leq\sqrt{h_n}$
and proves \eqref{eq:hdp-app-square-root}.
Multiplying by
$\sqrt{p_{n,\iota}}+\sqrt{q_{n,\iota}}$ gives
\[
|p_{n,\iota}-q_{n,\iota}|
\leq2\sqrt{h_nq_{n,\iota}}+h_n.
\]
Division by $q_{n,\iota}$ and the assumed lower bound prove
\eqref{eq:hdp-app-uniform-relative}.
These are simultaneous inequalities between the output laws; no union
bound over the event class is involved.
\end{proof}

\begin{corollary}[Exponential bounds]
\label{cor:hdp-app-exponential}
Let $v_n\to\infty$ and $a,b>0$. If
$\sup_{\iota\in\mathcal I_n}q_{n,\iota}\leq e^{-v_na}$
and $h_n\leq e^{-v_nb}$, then
\begin{equation}
\sup_{\iota\in\mathcal I_n}p_{n,\iota}
\leq4e^{-v_n\min\{a,b\}}.
\label{eq:hdp-app-exponential}
\end{equation}
More generally, if the corresponding logarithmic upper limits are
at most $-a$ and $-b$, then
\[
\limsup_{n\to\infty}\frac1{v_n}
\log\sup_{\iota\in\mathcal I_n}p_{n,\iota}
\leq-\min\{a,b\},
\]
where $\log0=-\infty$.
If instead
$\inf_{\iota\in\mathcal I_n}q_{n,\iota}\geq e^{-v_na}$
and $h_n\leq e^{-v_nb}$ with $b>a$, then
\begin{equation}
\sup_{\iota\in\mathcal I_n}
\left|\frac{p_{n,\iota}}{q_{n,\iota}}-1\right|
\leq2e^{-v_n(b-a)/2}+e^{-v_n(b-a)}.
\label{eq:hdp-app-exponential-relative}
\end{equation}
\end{corollary}

\begin{proof}
By \eqref{eq:hdp-app-square-root}, the first supremum is bounded by
$(e^{-v_na/2}+e^{-v_nb/2})^2$, which is at most the right side of
\eqref{eq:hdp-app-exponential}.
For logarithmic upper limits, apply this argument with
$a-\eta,b-\eta$ for $0<\eta<\min(a,b)$ and then let $\eta\downarrow0$.
The last assertion follows by substituting $a_n=e^{-v_na}$ into
\eqref{eq:hdp-app-uniform-relative}.
\end{proof}

\begin{proposition}[Transfer of relative approximations]
\label{prop:hdp-app-sharp-transfer}
Suppose positive reference quantities $r_{n,\iota}$ satisfy
\[
\eta_n:=
\sup_{\iota\in\mathcal I_n}
\left|\frac{q_{n,\iota}}{r_{n,\iota}}-1\right|
\longrightarrow0,
\qquad
r_{n,*}:=
\inf_{\iota\in\mathcal I_n}r_{n,\iota}>0.
\]
Then
\begin{equation}
\sup_{\iota\in\mathcal I_n}
\left|\frac{p_{n,\iota}}{r_{n,\iota}}-1\right|
\leq\eta_n+
2\sqrt{\frac{(1+\eta_n)h_n}{r_{n,*}}}
+\frac{h_n}{r_{n,*}}.
\label{eq:hdp-app-sharp-transfer}
\end{equation}
Hence $h_n/r_{n,*}\to0$ transfers the reference approximation uniformly,
including any prefactor present in $r_{n,\iota}$.
\end{proposition}

\begin{proof}
Use
$|p_{n,\iota}-q_{n,\iota}|
\leq2\sqrt{h_nq_{n,\iota}}+h_n$,
divide by $r_{n,\iota}$, and apply
$q_{n,\iota}/r_{n,\iota}\leq1+\eta_n$.
The triangle inequality proves \eqref{eq:hdp-app-sharp-transfer}.
\end{proof}

For example, suppose
$r_{n,\iota}=C_{n,\iota}v_n^{-\kappa_{n,\iota}}
e^{-v_n I_{n,\iota}}$, where
$C_{n,\iota}\geq c>0$,
$0\leq\kappa_{n,\iota}\leq K$, and
$I_{n,\iota}\leq I_*<b$.
If $h_n\leq e^{-v_nb}$, then, for $v_n\geq1$,
\[
\frac{h_n}{r_{n,*}}
\leq c^{-1}v_n^K e^{-v_n(b-I_*)}\longrightarrow0.
\]
Thus a separately established sharp reference approximation transfers
to the neighboring output laws. This conclusion does not assume that
the reference probabilities arise from an empirical partition statistic.

\begin{remark}[Sharpness of the exponential scale]
\label{rem:hdp-app-rate-sharpness}
Define
$\mathcal R(h,q)=
\sup\{P_1(A):D_{HD}(P_1,P_2)\leq h,\ P_2(A)\leq q\}$,
where the supremum is over probability laws and events.
Since $U_h$ is nondecreasing,
Proposition 16 of the main paper and its binary attaining examples
give $\mathcal R(h,q)=U_h(q)$. In particular,
\[
\max\{q,h-h^2/4\}
\leq\mathcal R(h,q)
\leq(\sqrt q+\sqrt h)^2.
\]
The lower bounds follow by taking $P_1=P_2$ and by the
distinguishing-event example, respectively.
If $q_n,h_n\to0$ and
$-v_n^{-1}\log q_n\to a>0$,
$-v_n^{-1}\log h_n\to b>0$, then
\[
-\frac1{v_n}\log\mathcal R(h_n,q_n)
\longrightarrow\min\{a,b\}.
\]
Indeed, for $0<h\leq2$ the lower bound is at least
$\max\{q,h/2\}$, while the upper bound is at most
$4\max\{q,h\}$.
For fixed $h>0$, instead,
$\lim_{q\downarrow0}\mathcal R(h,q)=h-h^2/4$.
Thus relative preservation of arbitrarily rare events requires a
corresponding restriction on the budget scale.
\end{remark}

\section{An example of crossing privacy profiles}
\label{app:hdp-crossing}

The following two neighboring-output experiments have the same
attained Hellinger distance but different orderings at different
privacy-profile thresholds. Consider
\[
P_1^R=(3/4,1/4,0),\qquad
P_2^R=(1/4,3/4,0).
\]
Their likelihood ratio takes the values $3$ and $1/3$ on the common
support, so this mechanism is pure $\log3$-DP.
Direct calculation gives
$D_{HD}(P_1^R,P_2^R)=2-\sqrt3$.
Now set $r=1-\sqrt3/2$ and take
\[
P_1^H=(1-r,r,0),\qquad
P_2^H=(1-r,0,r).
\]
Here $D_{HD}(P_1^H,P_2^H)=2r=2-\sqrt3$ as well.

For $t\geq0$, summing the positive cell contributions in
\eqref{eq:hdp-op-profiles} gives the two directional profiles, which
coincide within each experiment:
\begin{equation}
\delta_R(t)=\frac{(3-e^t)_+}{4},
\qquad
\delta_H(t)=r.
\label{eq:hdp-app-crossing}
\end{equation}
Indeed, in the first experiment only the cell with likelihood ratio
$3$ can have a positive contribution, and it contributes
$(3-e^t)_+/4$.
In the second, the common cell contributes zero, while the
distinguishing cell contributes $r$.
Reversing the laws gives the same quantities.

Thus $\delta_H(0)=r<1/2=\delta_R(0)$, whereas
$\delta_R(t)=0<r=\delta_H(t)$ for $t\geq\log3$.
The crossing occurs at $t=\log(2\sqrt3-1)$.
Consequently, equality of attained Hellinger distances does not impose
a pointwise ordering of privacy profiles.
The uniform likelihood-ratio cutoff of the first experiment describes
one aspect of its protection, not an ordering at every disclosure
threshold.

Both experiments can be realized as mechanisms on a binary database
domain. They also admit versions with Lebesgue densities by replacing
the three points by uniform laws on three disjoint unit intervals,
using the same weights.
The example compares these specified mechanisms and does not assert
a universal ordering of privacy classes.

\section{Influence functions and breakdown of PMHDE}
\label{app:robustness-proofs}

Write $T(g)=\bta_g$, let $\hat\bta_n$ denote MHDE, and let
$\widetilde\bta_{A,n}=\hat\bta_n^{(K_{A,n})}$ be its private version,
where $A\in\{\mathrm{PGD},\mathrm{PNR}\}$. The estimators are computed
from the same sample and use the same density estimate, bandwidth,
truncation, and criterion. Set
$R_{A,n}=\widetilde\bta_{A,n}-\hat\bta_n$.
The probability $\mathbf P_g$ includes sampling from $g$ and the
algorithmic randomization.

The inputs are the asymptotic-equivalence and consistency conclusions
of Theorem 42, and the approximation bounds in
Theorems 36 and 39. Each is applied
under its hypotheses at the stated sampling law and with the stipulated
initialization and tuning. In particular, whenever Theorem 42
applies under $g$,
\begin{equation}
\sqrt n R_{A,n}\longrightarrow0,
\qquad
\widetilde\bta_{A,n}\longrightarrow T(g)
\quad\text{in probability}.
\label{eq:rob-efficiency-input}
\end{equation}
The results below are consequences of these inputs. Uniform statements
and full-parameter-space breakdown specify their additional domain
conditions separately.

\subsection[Influence and the finite-alpha response]{Influence and the finite-$\alpha$ response}
\label{app:rob-influence}

\begin{proposition}[Asymptotic linear influence of PMHDE]
\label{prop:rob-influence}
Assume the hypotheses of Theorem 42 under $g$, and suppose
the ordinary MHDE has the representation
\[
\sqrt n\{\hat\bta_n-T(g)\}
=\frac1{\sqrt n}\sum_{i=1}^n\psi_g(X_i)+r_{n,g},
\qquad r_{n,g}\longrightarrow0\quad\text{in probability},
\]
where $\mathbf E_g\psi_g(X)=0$ and
$\mathbf E_g\|\psi_g(X)\|_2^2<\infty$.
Then PGD and PNR PMHDE have the same asymptotic linear influence
function $\psi_g$.
\end{proposition}

\begin{proof}
By Theorem 42(i), $\sqrt nR_{A,n}\to0$ in probability.
Hence
\[
\sqrt n\{\widetilde\bta_{A,n}-T(g)\}
=\frac1{\sqrt n}\sum_{i=1}^n\psi_g(X_i)
+r_{n,g}+\sqrt nR_{A,n}.
\]
The last two terms tend to zero, proving the assertion.
\end{proof}

Let $g_0=f_{\bta_0}$ and
$g_{\alpha,z}=(1-\alpha)g_0+\alpha q_z$, where $q_z$ is a density.
Following the finite-contamination formulation of \citet{Beran1977},
write
\begin{equation}
\operatorname{IF}_\alpha(z;T,g_0)
=\frac{T(g_{\alpha,z})-T(g_0)}{\alpha},
\qquad 0<\alpha<1.
\label{eq:rob-alpha-influence}
\end{equation}
For each sampling law on which the private estimator has a deterministic
probability limit, denote this limit by $T_A(g)$.
Define $\operatorname{IF}_\alpha(z;T_A,g_0)$ by the same expression
with $T_A$ in place of $T$.

For a sample-level version, let superscripts $(\alpha,z)$ and $0$
denote samples from $g_{\alpha,z}$ and $g_0$, respectively. Use any
coupling with these prescribed marginal sampling and randomization
laws, and put
\begin{align}
\widehat{\operatorname{IF}}_{\alpha,A,n}(z)
&=\frac{\widetilde\bta_{A,n}^{\alpha,z}
             -\widetilde\bta_{A,n}^{0}}{\alpha},\notag\\
\widehat{\operatorname{IF}}_{\alpha,H,n}(z)
&=\frac{\hat\bta_n^{\alpha,z}-\hat\bta_n^{0}}{\alpha}.
\label{eq:rob-alpha-empirical}
\end{align}
Here $H$ labels the ordinary MHDE, not another algorithm.

\begin{proposition}[$\alpha$-influence function of PMHDE]
\label{prop:rob-alpha-influence}
Fix $0<\alpha<1$ and $z$. Assume the hypotheses of
Theorem 42 at both $g_0$ and $g_{\alpha,z}$, with the same
Hellinger minimizer selection used for the private and non-private
procedures. Then, for $A\in\{\mathrm{PGD},\mathrm{PNR}\}$,
\begin{equation}
\operatorname{IF}_\alpha(z;T_A,g_0)
=\operatorname{IF}_\alpha(z;T,g_0),
\qquad
\widehat{\operatorname{IF}}_{\alpha,A,n}(z)
\longrightarrow\operatorname{IF}_\alpha(z;T,g_0)
\quad\text{in probability}.
\label{eq:rob-alpha-identity}
\end{equation}
Moreover,
\begin{equation}
\sqrt n\,
\bigl\|\widehat{\operatorname{IF}}_{\alpha,A,n}(z)
       -\widehat{\operatorname{IF}}_{\alpha,H,n}(z)\bigr\|_2
\longrightarrow0\quad\text{in probability}.
\label{eq:rob-alpha-root-n}
\end{equation}
For any set of locations on which the hypotheses hold pointwise,
the two limiting $\alpha$-influence curves therefore have identical
suprema of their norms, allowing the value $+\infty$.
\end{proposition}

\begin{proof}
Apply Theorem 42(ii) at the two sampling laws. It gives
$T_A(g_0)=T(g_0)$ and $T_A(g_{\alpha,z})=T(g_{\alpha,z})$,
which prove the functional identity. The sample response differs from
that common limit in norm by at most
\[
\frac{
\|\widetilde\bta_{A,n}^{\alpha,z}-T(g_{\alpha,z})\|_2
+\|\widetilde\bta_{A,n}^{0}-T(g_0)\|_2}{\alpha}.
\]
Both terms in the numerator tend to zero in probability, proving the
sample convergence for fixed $\alpha,z$. Finally, exactly,
\begin{equation}
\bigl\|\widehat{\operatorname{IF}}_{\alpha,A,n}(z)
       -\widehat{\operatorname{IF}}_{\alpha,H,n}(z)\bigr\|_2
\leq\frac{\|R_{A,n}^{\alpha,z}\|_2+\|R_{A,n}^{0}\|_2}{\alpha}.
\label{eq:rob-alpha-sample}
\end{equation}
Multiply by $\sqrt n$ and apply Theorem 42(i) at both
laws. A union bound suffices; independence is not required.
The statement about suprema follows from equality of the limiting
functions at every location in the specified set.
\end{proof}

The convergence statements above refer to the random sample response,
not its expectation. The functional identity requires only the
consistency conclusion of the efficiency theorem; the root-$n$
comparison uses its asymptotic-equivalence conclusion. If the clean and
contaminated samples and the two algorithmic randomizations are
independent, Theorem 42(iii) also gives, for fixed $\alpha,z$,
\begin{equation}
\sqrt n\left\{
\widehat{\operatorname{IF}}_{\alpha,A,n}(z)
-\operatorname{IF}_\alpha(z;T,g_0)\right\}
\ \Longrightarrow\
N\left(0,\frac{\Sigma_{g_{\alpha,z}}+\Sigma_{g_0}}{\alpha^2}\right).
\label{eq:rob-alpha-clt}
\end{equation}
Indeed, apply the two central limit statements jointly using
independence and take their difference divided by $\alpha$.
Equation~\eqref{eq:rob-alpha-root-n} transfers the same limit to the
ordinary sample response. Equality of the limiting curves does not
assert uniform convergence of the random curves over all $z$.

\paragraph{The local functional derivative.}
For the manuscript's criterion $2HD^2(f_{\bta},g)$, the population
estimating function is
$\Psi(\bta,g)=-2\int\sqrt{g f_{\bta}}\,\mathbf u_{\bta}$.
Consider a dominated density path $g_t=g_0+t a$ with $\int a=0$,
and a differentiable local minimizer $T(g_t)$.
Suppose the required derivatives may be taken under the integral,
$I(\bta_0)$ is nonsingular, and
$\int\|\mathbf u_{\bta_0}\||a|<\infty$.
The gradient and Hessian formulas in the manuscript give
\[
\partial_{\bta}\Psi(\bta_0,g_0)=I(\bta_0),
\qquad
\partial_g\Psi(\bta_0,g_0)[a]
=-\int\mathbf u_{\bta_0}(x)a(x)\,dx.
\]
Differentiation of $\Psi\{T(g_t),g_t\}=0$ yields
\begin{equation}
\left.\frac{d}{dt}T(g_t)\right|_{t=0}
=I(\bta_0)^{-1}\int\mathbf u_{\bta_0}(x)a(x)\,dx.
\label{eq:rob-functional-derivative}
\end{equation}
Thus the influence at the model is
$\psi_{g_0}=I(\bta_0)^{-1}\mathbf u_{\bta_0}$ on these regular
paths. For $a=q_z-g_0$, the derivative is
$\int\psi_{g_0}(x)q_z(x)\,dx$.
The asymptotic linear representation also includes the density-estimation
remainder, as in Proposition~\ref{prop:rob-influence}.
The finite-$\alpha$ identity in
Proposition~\ref{prop:rob-alpha-influence} requires no interchange
between $n\to\infty$ and $\alpha\downarrow0$.

\subsection{The limiting-functional breakdown point}
\label{app:rob-functional-breakdown}

Let $\Theta_*\subseteq\mathbb R^m$ be the full open parameter space,
and let $\mathcal Q$ be a nonempty class of contaminating densities.
For $0<\alpha<1$, set
$\mathcal G_\alpha=\{(1-\alpha)g_0+\alpha q:q\in\mathcal Q\}$.
For a functional $S$ on these laws, define
\begin{equation}
b^*(S;\mathcal Q,g_0)
=\inf\left\{\alpha\in(0,1):
\{S(g):g\in\mathcal G_\alpha\}
\text{ is not contained in a compact subset of }\Theta_*\right\},
\label{eq:rob-functional-breakdown-definition}
\end{equation}
with value one if the set in braces is empty. Leaving compact subsets
includes location escape, scale explosion, and scale collapse for
$\Theta_*=\mathbb R\times(0,\infty)$.

For the private sequence, the limiting-functional breakdown point is
$b^*(T_A;\mathcal Q,g_0)$, where $T_A(g)$ is its probability limit
at each fixed sampling law. Thus this definition takes $n\to\infty$
first. For these breakdown statements the estimators and targets are
considered on $\Theta_*$. The efficiency conclusion must apply at each
specified target through its local assumptions; a fixed compact
constraint is not used to prevent escape.

\begin{proposition}[Limiting-functional breakdown of PMHDE]
\label{prop:rob-functional-breakdown}
Suppose MHDE is consistent for $T(g)$ and the consistency conclusion of
Theorem 42(ii) applies to PGD or PNR PMHDE at every fixed
$g\in\{g_0\}\cup\bigcup_{0<\alpha<1}\mathcal G_\alpha$,
with the same global minimizer selection on $\Theta_*$.
Then
\begin{equation}
b^*(T_A;\mathcal Q,g_0)=b^*(T;\mathcal Q,g_0).
\label{eq:rob-functional-breakdown-identity}
\end{equation}
In particular, whenever the ordinary MHDE functional has breakdown
point $1/2$ on this class, both private procedures have the same
limiting-functional breakdown point $1/2$.
\end{proposition}

\begin{proof}
Theorem 42(ii) gives $T_A(g)=T(g)$ at every stated
sampling law. Fix $\alpha$. If the ordinary targets lie in a compact
set $K\subset\Theta_*$, all the private probability limits lie in
that same $K$. Conversely, if no compact subset contains the ordinary
targets, then for every compact $K$ there is
$g\in\mathcal G_\alpha$ with $T(g)\notin K$; its private probability
limit $T_A(g)$ is also outside $K$.
The sets of contamination fractions in
\eqref{eq:rob-functional-breakdown-definition} are therefore identical,
which proves \eqref{eq:rob-functional-breakdown-identity}.
\end{proof}

\begin{lemma}[A classical Hellinger benchmark]
\label{lem:rob-half-benchmark}
Let $g_0=f_{\bta_0}$ and $\rho(f,g)=\int\sqrt{fg}$.
Assume global maximizers
$T(g)\in\arg\max_{\bta\in\Theta_*}\rho(f_{\bta},g)$ exist
for the mixtures in the specified class, with a fixed selection.
Suppose that, for each $c>0$,
$\{\bta:\rho(f_{\bta},g_0)\geq c\}$ is contained in a compact
subset of $\Theta_*$, and that $\mathcal Q$ contains a sequence
$q_j=f_{\bta_j}$ satisfying
\begin{equation}
\sup_{\bta\in K}\rho(f_{\bta},q_j)\longrightarrow0
\quad\text{for every compact }K\subset\Theta_*.
\label{eq:rob-classical-separation}
\end{equation}
Then $b^*(T;\mathcal Q,g_0)=1/2$.
\end{lemma}

\begin{proof}
For $g_\alpha=(1-\alpha)g_0+\alpha q$, comparison with
$f_{\bta_0}$ and $\sqrt{a+b}\leq\sqrt a+\sqrt b$ give
\[
\sqrt{1-\alpha}
\leq\rho(f_{T(g_\alpha)},g_\alpha)
\leq\sqrt{1-\alpha}\rho(f_{T(g_\alpha)},g_0)+\sqrt\alpha.
\]
For $\alpha<1/2$, it follows that
\begin{equation}
\rho(f_{T(g_\alpha)},g_0)
\geq1-\sqrt{\frac{\alpha}{1-\alpha}}>0.
\label{eq:rob-classical-affinity}
\end{equation}
All the targets therefore lie in a compact set independent of $q$.
For $\alpha>1/2$, take $q=q_j$. Fitting $q_j$ itself gives affinity
at least $\sqrt\alpha$. On any fixed compact $K$,
\[
\sup_{\bta\in K}\rho(f_{\bta},g_{\alpha,j})
\leq\sqrt{1-\alpha}
+\sqrt\alpha\sup_{\bta\in K}\rho(f_{\bta},q_j).
\]
For sufficiently large $j$, the right side is strictly less than
$\sqrt\alpha$. Every global maximizer must then lie outside $K$.
Thus all fractions below $1/2$ preserve compact containment and every
fraction above $1/2$ admits an escape sequence. Taking the infimum in
\eqref{eq:rob-functional-breakdown-definition} proves the result.
\end{proof}

This is a population affinity argument of the classical minimum-distance
type \citep{Beran1977,Basu2011}; its sufficient conditions are stated
here to identify the benchmark used in the private transfer.

\begin{corollary}[Half breakdown for the limiting PMHDE functional]
\label{cor:rob-half-PMHDE}
Suppose the conditions of Lemma~\ref{lem:rob-half-benchmark} hold,
and the consistency conditions of
Proposition~\ref{prop:rob-functional-breakdown} hold for both PGD
and PNR. Then
\[
b^*(T_{\mathrm{PGD}};\mathcal Q,g_0)
=b^*(T_{\mathrm{PNR}};\mathcal Q,g_0)
=b^*(T;\mathcal Q,g_0)=\tfrac12.
\]
For every fixed $0<\alpha<1/2$, the limiting $\alpha$-influence
functions of the two private procedures are uniformly bounded over
$q_z\in\mathcal Q$ and equal the MHDE curve.
\end{corollary}

\begin{proof}
The breakdown identity follows by applying
Proposition~\ref{prop:rob-functional-breakdown} to
Lemma~\ref{lem:rob-half-benchmark}.
For fixed $\alpha<1/2$, \eqref{eq:rob-classical-affinity} puts all
$T(g_{\alpha,z})$ in one compact set. Their Euclidean distances from
$T(g_0)$ are uniformly bounded, as are those distances divided by the
fixed positive $\alpha$.
The identity $T_A(g)=T(g)$ from
Proposition~\ref{prop:rob-functional-breakdown} identifies the same
limiting curves for PMHDE.
\end{proof}

\subsection{Probability bounds and simultaneous breakdown transfer}
\label{app:rob-breakdown}

For a fixed privacy budget, the approximation bounds of
Theorems 36 and 39 have the form
\begin{equation}
\mathbf P_g\{\|R_{A,n}\|_2>e_{A,n}(g;\xi)\}\leq\xi,
\qquad
 e_{A,n}(g;\xi)=C_A(g)n^{-1/p}
\{K_{A,n}\log(K_{A,n}/\xi)\}^{1/2},
\label{eq:rob-source-bound}
\end{equation}
for $0<\xi<1$ on the domains and sample sizes covered by those
results. The endpoint $p=2$ refers to the consistency-scale use of
the weak sensitivity bound. Local constants, initialization, and
sample-size cutoffs are retained when the sampling laws vary.

\begin{lemma}[Finite-sample contamination comparison]
\label{lem:rob-probability-comparison}
Put $B_n=\|\hat\bta_n-\bta_0\|_2$ and
$\widetilde B_{A,n}=\|\widetilde\bta_{A,n}-\bta_0\|_2$.
Whenever \eqref{eq:rob-source-bound} holds, let $e=e_{A,n}(g;\xi)$.
For $t>e$,
\begin{equation}
\mathbf P_g\{B_n>t+e\}-\xi
\leq\mathbf P_g\{\widetilde B_{A,n}>t\}
\leq\mathbf P_g\{B_n>t-e\}+\xi.
\label{eq:rob-tail-comparison}
\end{equation}
If the bound is applied at both $g_{\alpha,z}$ and $g_0$, then with
probability at least $1-2\xi$,
\begin{equation}
\bigl\|\widehat{\operatorname{IF}}_{\alpha,A,n}(z)
       -\widehat{\operatorname{IF}}_{\alpha,H,n}(z)\bigr\|_2
\leq\frac{e_{A,n}(g_{\alpha,z};\xi)+e_{A,n}(g_0;\xi)}{\alpha}.
\label{eq:rob-alpha-probability}
\end{equation}
\end{lemma}

\begin{proof}
The reverse triangle inequality gives
$|\widetilde B_{A,n}-B_n|\leq\|R_{A,n}\|_2$.
On $\{\|R_{A,n}\|_2\leq e\}$,
$\{B_n>t+e\}\subseteq\{\widetilde B_{A,n}>t\}
\subseteq\{B_n>t-e\}$.
Adding the probability of the complementary event proves
\eqref{eq:rob-tail-comparison}.
Equation~\eqref{eq:rob-alpha-sample} and a union bound prove
\eqref{eq:rob-alpha-probability}.
\end{proof}

If displacement is measured from a clean estimate instead of
$\bta_0$, put
$V_{H,n}=\hat\bta_n^{\alpha,z}-\hat\bta_n^0$ and
$V_{A,n}=\widetilde\bta_{A,n}^{\alpha,z}-\widetilde\bta_{A,n}^0$.
Then
$|\|V_{A,n}\|_2-\|V_{H,n}\|_2|
\leq\|R_{A,n}^{\alpha,z}\|_2+\|R_{A,n}^0\|_2$.
Thus \eqref{eq:rob-tail-comparison} holds with radius
$e_{A,n}(g_{\alpha,z};\xi)+e_{A,n}(g_0;\xi)$ and failure probability
$2\xi$. This comparison does not require independent clean and
contaminated samples.

For a simultaneous sample-size and contamination limit, use the
compact exhaustion of the full open parameter space
\[
\mathcal K_R=\left\{\bta\in\Theta_*:
\|\bta\|_2\leq R,\quad
\operatorname{dist}(\bta,\mathbb R^m\setminus\Theta_*)\geq R^{-1}
\right\},\qquad R\geq1.
\]
The distance to an empty complement is infinite. These compact sets
increase to $\Theta_*$ and contain every compact subset for large
$R$. An estimator sequence $S_n$ is asymptotically tight uniformly
over $\mathcal G_\alpha$ when
\begin{equation}
\lim_{R\to\infty}\limsup_{n\to\infty}
\sup_{g\in\mathcal G_\alpha}
\mathbf P_g\{S_n\notin\mathcal K_R\}=0.
\label{eq:rob-tightness-definition}
\end{equation}
Its simultaneous asymptotic breakdown threshold is the infimum of
those $\alpha\in(0,1)$ for which this limit is positive, with value
one if there are none. An escape witness is a sequence
$g_{\alpha,n}\in\mathcal G_\alpha$ for which
$\mathbf P_{g_{\alpha,n}}\{S_n\in K\}\to0$ for every compact
$K\subset\Theta_*$.
A failed or inadmissible output is counted as outside every compact
set; the approximation error is assigned value $+\infty$ on failures.

\begin{proposition}[Simultaneous asymptotic breakdown transfer]
\label{prop:rob-breakdown-transfer}
Fix a threshold $b\in(0,1)$. For every $\alpha<b$, suppose MHDE
is asymptotically tight uniformly over $\mathcal G_\alpha$ and there
are $r_{A,n,\alpha}\to0$, $\xi_{A,n,\alpha}\to0$ such that
\begin{equation}
\sup_{g\in\mathcal G_\alpha}
\mathbf P_g\{\|R_{A,n}\|_2>r_{A,n,\alpha}\}
\leq\xi_{A,n,\alpha}.
\label{eq:rob-uniform-error}
\end{equation}
For every $\alpha>b$, suppose there is an MHDE escape witness
$g_{\alpha,n}\in\mathcal G_\alpha$ along which
$R_{A,n}\to0$ in probability. Then PMHDE has simultaneous asymptotic
breakdown threshold $b$. In particular, the threshold is $1/2$ for
PGD and PNR whenever these hypotheses hold with $b=1/2$.
\end{proposition}

\begin{proof}
Fix $R\geq1$. If $\bta\in\mathcal K_R$ and
$\|\bta'-\bta\|_2\leq(2R)^{-1}$, then
$\bta'\in\mathcal K_{2R+1}$: its norm is at most
$R+(2R)^{-1}$ and its distance from the parameter-space complement
is at least $(2R)^{-1}$.
For $\alpha<b$, the error radius is eventually at most $(2R)^{-1}$.
Hence
\begin{align*}
\sup_{g\in\mathcal G_\alpha}
\mathbf P_g\{\widetilde\bta_{A,n}\notin\mathcal K_{2R+1}\}
&\leq\sup_{g\in\mathcal G_\alpha}
\mathbf P_g\{\hat\bta_n\notin\mathcal K_R\}
+\xi_{A,n,\alpha}.
\end{align*}
Take $\limsup_n$, then $R\to\infty$, to prove tightness below $b$.
For the escape witness above $b$, reverse the roles of the estimates:
\[
\mathbf P_{g_{\alpha,n}}
\{\widetilde\bta_{A,n}\in\mathcal K_R\}
\leq\mathbf P_{g_{\alpha,n}}
\{\hat\bta_n\in\mathcal K_{2R+1}\}
+\mathbf P_{g_{\alpha,n}}
\{\|R_{A,n}\|_2>(2R)^{-1}\}.
\]
Both terms tend to zero. Every $\alpha>b$ therefore admits escape
and the threshold equals $b$. No claim at $\alpha=b$ is needed.
\end{proof}

The same compact-set argument also transfers uniform tightness from
bounded population targets whenever an ordinary estimation error
bound holds uniformly over the contamination class. Along a population
escape sequence it transfers escape if both the ordinary estimation
error and the private approximation error vanish. This is the extra
sampling step needed to pass from
Corollary~\ref{cor:rob-half-PMHDE} to the simultaneous statement.

For location--scale parameters, a lower bound $\hat\sigma_n\geq c>0$
and private error at most $c/2$ give
$\widetilde\sigma_{A,n}\geq c/2$. The compact-set proof therefore
includes scale collapse. For an exploding Euclidean displacement, the
weaker requirement that the private error be negligible relative to
that displacement suffices by \eqref{eq:rob-tail-comparison}.

\subsection{PGD and PNR error radii}
\label{app:rob-algorithm-radii}

The following local calculations display the constants in the
probability comparisons. They use the manuscript's Hessian bounds and
perturbation structure without changing either algorithm.
They are stated on a good-data event $\mathcal E_{A,n}(g)$ on which
an interior MHDE root exists, the starting point is in the stated
neighborhood, and the local derivative bounds hold. Write
$\ell_{A,n}(g)=\mathbf P_g\{\mathcal E_{A,n}(g)^c\}$.
If the original approximation theorem is used as an unconditional
input, its failure probability already accounts for the required
localization; no additional $\ell_{A,n}$ is then added.

Use a closed ball $\overline B_r(\hat\bta_n)\subset\Theta_*$ on which
$2\tau_1 I\preceq H_n(\bta)\preceq2\tau_2 I$, with
$0<\tau_1\leq\tau_2<\infty$.
Constants and radii below may depend on $g$, but are deterministic
bounds on the specified good-data event. Let $e_{0,n}$ be a deterministic
upper bound for $\|\hat\bta_n^{(0)}-\hat\bta_n\|_2$ on this event. Let
$c_u=\{-8\log(1-u/2)\}^{-1/2}$, $0<u<2$, as in the manuscript.
For total privacy budget $0<\epsilon<2$ and $K$ iterations, the existing
allocation uses $\epsilon'=\epsilon/K$.

\paragraph{PGD.}
The Gaussian vector radius used in the PGD proof is
\begin{equation}
u_{G,n}(\xi)=\Delta_n c_{\epsilon'}
\{4\sqrt m+2\sqrt{2\log(K/\xi)}\}.
\label{eq:rob-pgd-noise-radius}
\end{equation}
The corresponding Gaussian norm estimate and a union bound give
$\max_{k<K}\|N_{n,k}\|_2\leq u_{G,n}(\xi)$ with probability at least
$1-\xi$. For $0<\eta\leq(2\tau_2)^{-1}$, put
$q_G=1-2\eta\tau_1\in[0,1)$.
The integral Hessian formula along the line from $\hat\bta_n$ to an
iterate gives
\[
\|\hat\bta_n^{(k+1)}-\hat\bta_n\|_2
\leq q_G\|\hat\bta_n^{(k)}-\hat\bta_n\|_2+\eta u_{G,n}(\xi).
\]
If $e_{0,n}\leq r$ and
$u_{G,n}(\xi)/(2\tau_1)\leq r$, induction keeps the iterates in
this ball and yields
\begin{equation}
e_{G,n}^{\mathrm{loc}}(\xi)
=q_G^K e_{0,n}
 +\frac{u_{G,n}(\xi)}{2\tau_1}(1-q_G^K).
\label{eq:rob-pgd-local-error}
\end{equation}
The probability that the actual PGD error exceeds this radius is at
most $\ell_{G,n}(g)+\xi$. This expression separates finite iteration
and privacy perturbation.

\paragraph{PNR.}
For the unit-step local Newton phase used in the PNR proof, let
\begin{align}
u_{N,n}(\xi)&=\Delta_n c_{\epsilon'/2}
 \{4\sqrt m+2\sqrt{2\log(2K/\xi)}\},\notag\\
v_{N,n}(\xi)&=\Delta_n^{(H)}c_{\epsilon'/2}
 \{2m\log(4Km/\xi)\}^{1/2}.
\label{eq:rob-pnr-noise-radii}
\end{align}
These are the gradient and symmetric-matrix Gaussian norm radii used
in the supplementary PNR proof; the concentration tools are cited there
as \citet{Rigollet2023,Tropp2015}. The independent entries of the upper
triangle are reflected into the lower triangle. A union bound over the
two noise sequences gives simultaneous bounds
$\|N_{n,k}\|_2\leq u_{N,n}$, $\|W_{n,k}\|_2\leq v_{N,n}$
for all $k<K$, with probability at least $1-\xi$.
Suppose $v_{N,n}\leq\tau_1$ and
$\|\nabla L_n(\bta)\|_2\leq B_1$ on the local ball.
The inverse identity
$(H+W)^{-1}-H^{-1}=-(H+W)^{-1}WH^{-1}$ gives
\[
\|(H+W)^{-1}-H^{-1}\|_2\leq\frac{v_{N,n}}{2\tau_1^2}.
\]
Writing the noisy direction as the exact Newton direction plus an
error, as in the manuscript's PNR decomposition, gives the bound
\begin{equation}
a_{N,n}(\xi)
=\frac{u_{N,n}}{2\tau_1}
 +\frac{B_1v_{N,n}}{2\tau_1^2}
 +\frac{u_{N,n}v_{N,n}}{2\tau_1^2}.
\label{eq:rob-pnr-direction-error}
\end{equation}
In particular, both the gradient and Hessian perturbations and their
interaction remain in the calculation.

If $H_n$ is Lipschitz with constant $L_H$ on the local ball, set
$c_N=L_H/(4\tau_1)$. Integration of the Hessian between an iterate
and the MHDE root proves
\begin{equation}
E_{k+1}\leq c_NE_k^2+a_{N,n}(\xi),
\qquad E_k=\|\hat\bta_n^{(k)}-\hat\bta_n\|_2.
\label{eq:rob-pnr-recursion}
\end{equation}
Indeed, the exact Newton remainder has norm at most
$\|H_n^{-1}\|(L_H/2)E_k^2\leq c_NE_k^2$.
For a concrete sufficient local choice, suppose
$c_Nr\leq1/4$, $e_{0,n}\leq r/2$, and $a_{N,n}\leq r/8$.
If $c_N>0$, induction gives
\begin{equation}
e_{N,n}^{\mathrm{loc}}(\xi)
=c_N^{-1}(c_Ne_{0,n})^{2^K}+2a_{N,n}(\xi).
\label{eq:rob-pnr-local-error}
\end{equation}
To check the induction, write $x_k=c_NE_k$ and $b=c_Na_{N,n}$.
Then $x_0\leq1/8$, $b\leq1/32$, and
$x_{k+1}\leq x_k^2+b$. If $x_k\leq x_0^{2^k}+2b$,
putting $t=x_0^{2^k}\leq1/8$ yields
$(t+2b)^2+b\leq t^2+2b$ because $4t+4b\leq1$.
Also the resulting radius is at most $e_{0,n}+2a_{N,n}\leq3r/4$,
so the iterates remain in the ball. If $c_N=0$, one Newton step leaves
only the direction error and $E_K\leq a_{N,n}$ for $K\geq1$.
The failure probability is at most $\ell_{N,n}(g)+\xi$.
For a fixed damped step $0<\eta<1$, the corresponding verified
recursion is
$E_{k+1}\leq(1-\eta)E_k+\eta c_NE_k^2+\eta a_{N,n}$;
the displayed quadratic-iteration bound pertains to $\eta=1$.

\subsection{Varying contamination and privacy levels}
\label{app:rob-rate-conditions}

\begin{proposition}[Private response at a changing contamination fraction]
\label{prop:rob-varying-alpha}
Let $\alpha_n>0$, allow $z=z_n$, and suppose the two approximation
bounds \eqref{eq:rob-source-bound} hold under $g_{\alpha_n,z_n}$ and
$g_0$ with failure probabilities at most $\xi_n$. Put
\begin{equation}
d_{A,n}=\frac{e_{A,n}(g_{\alpha_n,z_n};\xi_n)
                    +e_{A,n}(g_0;\xi_n)}{\alpha_n}.
\label{eq:rob-varying-alpha}
\end{equation}
Then
\[
\mathbf P\!\left\{
\bigl\|\widehat{\operatorname{IF}}_{\alpha_n,A,n}(z_n)
-\widehat{\operatorname{IF}}_{\alpha_n,H,n}(z_n)\bigr\|_2>d_{A,n}
\right\}\leq2\xi_n.
\]
If $d_{A,n}\to0$ and $\xi_n\to0$, the private and non-private
sample responses differ by a term tending to zero in probability.
\end{proposition}

\begin{proof}
Apply \eqref{eq:rob-alpha-sample} with $(\alpha,z)=(\alpha_n,z_n)$
and bound the two approximation errors by a union bound.
\end{proof}

This result compares the sample responses at the stated contamination
scale. At $\alpha_n=c/\sqrt n$, $c>0$, root-$n$ negligible errors
along the two sampling sequences suffice. For fixed positive
$\alpha$, consistency-scale errors suffice. In all cases the exact
decomposition is
\begin{equation}
\widetilde\bta_{A,n}-\bta_0
=\{T(g_{\alpha_n,z_n})-\bta_0\}
 +\{\hat\bta_n-T(g_{\alpha_n,z_n})\}+R_{A,n}.
\label{eq:rob-three-terms}
\end{equation}

Take $\xi_n=n^{-2}$ and a fixed privacy budget. On a family where the
sensitivity and local constants have the required uniform bounds,
choose $K_{G,n}$ between fixed positive multiples of $\log n$ and
$K_{N,n}$ between fixed positive multiples of $\log\log n$ for large
$n$, with the lower constants large enough that the initialization
terms are no larger than the noise terms. The existing approximation
rates, or the local bounds above, give
\begin{align}
e_{G,n}(g;\xi_n)&\leq C_G(g)n^{-1/p}\log n,\notag\\
e_{N,n}(g;\xi_n)&\leq C_N(g)n^{-1/p}
 \{\log n\,\log\log n\}^{1/2}.
\label{eq:rob-specialized-rates}
\end{align}
Both vanish even at $p=2$ when the stated constants stay bounded and
localization failure tends to zero. Thus consistency-scale sensitivity
can suffice for breakdown; the root-$n$ statements require the stronger
rate. Along moving contamination laws, the displayed constants and
sample-size cutoffs must satisfy the stated bounds along that sequence.

The noise-scale dependence on a varying privacy budget can be retained
directly. Since $-\log(1-x)\geq x$ for $0<x<1$,
\[
c_{\epsilon_n/K}\leq\sqrt{\frac{K}{4\epsilon_n}},
\qquad
c_{\epsilon_n/(2K)}\leq\sqrt{\frac{K}{2\epsilon_n}}.
\]
If both sensitivities are bounded by constants times $n^{-1/p}$ and
$0<\epsilon_n<2$ satisfies
$\epsilon_n\geq c_\epsilon n^{-\gamma}$, with $c_\epsilon>0$ and
$0\leq\gamma<2/p$, the local error bounds give, once the initialization
terms are smaller,
\begin{equation}
e_{A,n}(g;\xi_n)\leq C_A(g)n^{-1/p+\gamma/2}L_A(n),
\quad
L_G(n)=\log n,\quad
L_N(n)=\{\log n\,\log\log n\}^{1/2}.
\label{eq:rob-varying-privacy-rate}
\end{equation}
The PNR gradient--Hessian interaction is a product of vanishing radii
and is absorbed for sufficiently large $n$. With bounded constants,
$\gamma<2/p$ gives vanishing error, and $\gamma<2/p-1$ gives
root-$n$ negligible error. Uniform versions also require uniform
localization and a common prescribed tuning rule.

If $\alpha_n\geq c_\alpha n^{-r}$, where $c_\alpha>0$ and $r\geq0$,
substitution into \eqref{eq:rob-varying-alpha} gives
\[
d_{A,n}\leq
\frac{C_A(g_{\alpha_n,z_n})+C_A(g_0)}{c_\alpha}
 n^{r-1/p+\gamma/2}L_A(n).
\]
For bounded constants and vanishing failure probabilities,
$r<1/p-\gamma/2$ is sufficient for equality of the normalized sample
responses in probability. This bound retains the interaction of
sample size, privacy budget, and contamination fraction.

Finally, suppose for fixed $q_z$ that
$\{T(g_{\alpha,z})-\bta_0\}/\alpha\to b_z$ as $\alpha\downarrow0$,
and that under $g_{\alpha_n,z}$ the non-private centered estimator
$\sqrt n\{\hat\bta_n-T(g_{\alpha_n,z})\}$ converges to $N(0,V_0)$.
If $\sqrt n\alpha_n\to c<\infty$ and
$\sqrt nR_{A,n}\to0$ along those laws, then
\[
\sqrt n(\widetilde\bta_{A,n}-\bta_0)
\ \Longrightarrow\ N(cb_z,V_0).
\]
Multiply \eqref{eq:rob-three-terms} by $\sqrt n$ and use Slutsky's
theorem. Here the non-private triangular-array limit and the private
approximation along that sequence are the stated premises.

\section{GDP calibration for the scalar Laplace mechanism}
\label{app:laplace-gdp}

Consider $M_b(D)=W(D)+Z_b$, where $W$ is scalar,
$Z_b$ has density $(2b)^{-1}e^{-|z|/b}$, and the common scale $b>0$
does not depend on the dataset. Let
$\Delta=\sup_{D\sim D'}|W(D)-W(D')|\in(0,\infty)$ be its exact
sensitivity and put $r=\Delta/(2b)$. The exact Laplace calculation gives
\begin{equation}
\delta_L(t;r)=[1-e^{t/2-r}]_+,
\qquad \eta_L(r)=1-e^{-r},
\qquad h_L(b)=2\{1-(1+r)e^{-r}\},
\label{eq:lap-gdp-source-quantities}
\end{equation}
where $\delta_L$ is the worst-case profile over ordered neighboring
pairs, $\eta_L$ is the worst-case TV distance, and $h_L$ is the
attained HDP level. The profile formula is also given by
\citet[Theorem~3]{BalleBartheGaboardi2020}. The affinity formula
follows by splitting the integral at the two neighboring query values.
All suprema can be interpreted as limits along pairs approaching
$\Delta$; attainment of the sensitivity is not required.

For $\mu>0$, write
\begin{equation}
\delta_G(t;\mu)
=\Phi\left(\frac\mu2-\frac t\mu\right)
 -e^t\Phi\left(-\frac\mu2-\frac t\mu\right),
\qquad t\geq0.
\label{eq:lap-gdp-gaussian-profile}
\end{equation}
The full-profile characterization of \citet{Dong2022} states that
$M_b$ is $\mu$-GDP if and only if
$\delta_L(t;r)\leq\delta_G(t;\mu)$ for every $t\geq0$.
The following argument verifies this inequality over the entire
profile, not only at its TV value $t=0$.

\begin{proposition}[Smallest GDP parameter for scalar Laplace noise]
\label{prop:laplace-gdp-exact}
The scalar Laplace mechanism is $\mu$-GDP if and only if
\begin{equation}
\mu\geq\mu_L(r):=
2\Phi^{-1}\left(1-\frac12e^{-r}\right)
=2\Phi^{-1}\left(\frac{1+\eta_L(r)}2\right).
\label{eq:lap-gdp-minimum}
\end{equation}
Here $r=\Delta/(2b)>0$. Thus $\mu_L(r)$ is the smallest admissible
GDP parameter of the mechanism.
\end{proposition}

\begin{proof}
Necessity follows at $t=0$: any GDP parameter must satisfy
\[
1-e^{-r}=\delta_L(0;r)
\leq\delta_G(0;\mu)=2\Phi(\mu/2)-1,
\]
which is equivalent to \eqref{eq:lap-gdp-minimum}.

For sufficiency, put
$k=\Phi^{-1}(1-e^{-r}/2)>0$ and $q=\Phi(-k)=e^{-r}/2$.
The profiles $\delta_G(t;2k)$ and $\delta_L(t;r)$ agree at $t=0$.
For $0<t<2r$, differentiation gives
\begin{equation}
\frac{d}{dt}\{\delta_G(t;2k)-\delta_L(t;r)\}
=e^{t/2}\left\{
q-e^{t/2}\Phi\left(-k-\frac{t}{2k}\right)
\right\}.
\label{eq:lap-gdp-profile-derivative}
\end{equation}
Indeed, the derivative of the Gaussian profile is
$-e^t\Phi(-k-t/(2k))$; the normal-density terms cancel since
$\phi(k-t/(2k))=e^t\phi(k+t/(2k))$.

To determine the sign in \eqref{eq:lap-gdp-profile-derivative}, define
$v(t)=e^{t/2}\Phi(-k-t/(2k))$. With $x=k+t/(2k)>0$,
\[
\frac{v'(t)}{v(t)}
=\frac12-\frac{\phi(x)}{2k\Phi(-x)}
<\frac12-\frac{x}{2k}\leq0.
\]
The strict inequality follows from
\[
\phi(x)=\int_x^\infty u\phi(u)\,du
>x\int_x^\infty\phi(u)\,du=x\Phi(-x).
\]
Thus $v(t)\leq v(0)=q$, and
\eqref{eq:lap-gdp-profile-derivative} is nonnegative. It follows that
$\delta_L(t;r)\leq\delta_G(t;2k)$ on $[0,2r]$.
For $t\geq2r$ the Laplace profile is zero, so the inequality continues
to hold. This proves $2k$-GDP by the full-profile characterization.
Since $\partial_\mu\delta_G(t;\mu)=\phi(\mu/2-t/\mu)>0$, every
$\mu\geq2k$ is also admissible. Necessity proves minimality.
\end{proof}

Although the profiles agree at $t=0$ for this smallest GDP parameter,
the proof gives
$\delta_L(t;r)<\delta_G(t;\mu_L(r))$ for every finite $t>0$.
Thus the result identifies the tightest Gaussian trade-off lower
bound for this Laplace family; it does not identify its actual
trade-off function with a Gaussian one.

\begin{corollary}[HDP correspondence and inverse GDP calibration]
\label{cor:laplace-gdp-inverse}
For $0<h<2$, let $r_h>0$ be the unique solution of
$(1+r_h)e^{-r_h}=1-h/2$. The minimum scalar Laplace scale attaining
$h$-HDP is $b_H(h)=\Delta/(2r_h)$. At that scale, the smallest
GDP parameter is
\begin{equation}
\mu_h^L=2\Phi^{-1}\left(1-\frac12e^{-r_h}\right).
\label{eq:lap-gdp-from-h}
\end{equation}
Every scalar Laplace mechanism satisfying $h$-HDP therefore also
satisfies $\mu_h^L$-GDP.

Conversely, for a prescribed $\mu>0$, put
\begin{equation}
q_\mu=2\Phi(-\mu/2),
\qquad r_\mu=-\log q_\mu.
\label{eq:lap-gdp-inverse-r}
\end{equation}
Then $M_b$ is $\mu$-GDP if and only if
\begin{equation}
b\geq b_{\mathrm{GDP}}(\mu):=\frac{\Delta}{2r_\mu}.
\label{eq:lap-gdp-inverse-b}
\end{equation}
The noise variance and attained Hellinger level at this minimum
scale are, respectively,
\begin{equation}
\operatorname{Var}(Z_{b_{\mathrm{GDP}}})
=\frac{\Delta^2}{2r_\mu^2},
\qquad
h_L^{\mathrm{GDP}}(\mu)
=2\{1-(1+r_\mu)e^{-r_\mu}\}
=2\{1-q_\mu(1-\log q_\mu)\}.
\label{eq:lap-gdp-inverse-h}
\end{equation}
These inverse calibrations recover one another:
$b_{\mathrm{GDP}}(\mu_h^L)=b_H(h)$ and
$h_L^{\mathrm{GDP}}(\mu_h^L)=h$.
\end{corollary}

\begin{proof}
The function $r\mapsto2\{1-(1+r)e^{-r}\}$ increases strictly on
$(0,\infty)$ from $0$ to $2$, and $r\mapsto\mu_L(r)$ is strictly
increasing. Thus the condition $h_L(b)\leq h$ is equivalent to
$r\leq r_h$, and the first assertion follows from
Proposition~\ref{prop:laplace-gdp-exact}.

For the inverse, \eqref{eq:lap-gdp-minimum} is equivalent to
$e^{-r}\geq2\Phi(-\mu/2)$, hence to $r\leq r_\mu$ and
$b\geq\Delta/(2r_\mu)$. Substitution into the variance $2b^2$ and
\eqref{eq:lap-gdp-source-quantities} proves
\eqref{eq:lap-gdp-inverse-h}. Finally,
$2\Phi(-\mu_h^L/2)=e^{-r_h}$, which proves the inverse identities.
\end{proof}

\begin{remark}[Scope of the inverse-TV expression]
\label{rem:lap-gdp-inverse-tv-scope}
For an additive Gaussian mechanism with a common dataset-independent
covariance $\sigma^2I$, the smallest GDP parameter is
$\Delta/\sigma$, and its exact TV distance is
$2\Phi\{\Delta/(2\sigma)\}-1$.
Together with Proposition~\ref{prop:laplace-gdp-exact}, this shows that
in either this Gaussian family or the common-scale scalar Laplace
family, the smallest GDP parameter equals
$2\Phi^{-1}\{(1+\eta)/2\}$, with $\eta$ the exact worst-case TV distance.

Consequently, within either specified family, $h$-HDP implies GDP
with parameter
\[
2\Phi^{-1}\left(
\frac{1+\sqrt{h-h^2/4}}2
\right),\qquad 0<h<2.
\]
For $0<h<1$, replacing $\sqrt{h-h^2/4}$ by $\sqrt h$ gives a further
valid relaxation. These conclusions use the verified family-specific
profile comparisons; they are not consequences of TV control alone.
For the scalar Laplace family, \eqref{eq:lap-gdp-from-h} retains the
exact Hellinger information and can be used instead of this relaxation.
\end{remark}

\clearpage

\end{document}